\documentclass[paper=a4,fontsize=11pt,DIV=14,headinclude,headings=normal]{scrbook}
\usepackage{latexsym, amssymb, amsthm, amsmath, makeidx, ngerman}
\usepackage[ansinew]{inputenc}

\setkomafont{sectioning}{\normalfont\normalcolor\bfseries}

\allowdisplaybreaks

\newtheorem{defin}{}[chapter]
\newtheorem{saetze}[defin]{}
\newtheorem{lemmas}[defin]{}
\newtheorem{folger}[defin]{}
\newtheorem{verein}[defin]{}
\newtheorem{bemerk}[defin]{}
\newtheorem{defini}[defin]{}
\newtheorem{exampl}[defin]{}

\newenvironment{satz}      {\begin{saetze}\rmfamily {\bfseries Satz:}}{\end{saetze}}
\newenvironment{lemma}     {\begin{lemmas}\rmfamily {\bfseries  Lemma:}}{\end{lemmas}}
\newenvironment{folg}      {\begin{folger}\rmfamily {\bfseries  Folgerung:}}{\end{folger}}
\newenvironment{defi}      {\begin{defini}\rmfamily {\bfseries  Definition:}}{\end{defini}}
\newenvironment{ver}       {\begin{verein}\rmfamily {\bfseries  Vereinbarung:}}{\end{verein}}
\newenvironment{beme}      {\begin{bemerk}\rmfamily {\bfseries  Bemerkung:}}{\end{bemerk}}

\newenvironment{beweis}    {\noindent{\scshape Beweis}:}    {{\hfill $\Box$ \medskip}}

\newcommand{\ra}{\rightarrow}

\newcommand{\mt}{\mapsto}
\newcommand{\ol}{\overline}

\newcommand{\gap}{{\scshape Gap}}

\newcommand{\Z}{\mathbb Z}
\newcommand{\N}{\mathbb N}

\newcommand{\F}{\mathbb F}
\newcommand{\f}{\F_p}
\newcommand{\mal}{\cdot \ldots \cdot}
\newcommand{\dd}{, \ldots, }

\newcommand{\1}{^{-1}}

\newcommand{\vi}{{\varphi}}
\newcommand{\al}{{\alpha}}
\newcommand{\skp}[2]{\langle #1 \,|\, #2  \rangle}
\newcommand{\erz}[1]{\langle #1 \rangle}
\newcommand{\sw}{\sqrt{w}}
\newcommand{\fsw}{\F_p(\sqrt{w})}

\begin{document}



\pagenumbering{Roman}

\title{Hausarbeit im Rahmen\\ der Ersten Staatsprüfung\\ für das Lehramt an Gymnasien}

\author{vorgelegt von: Boris Girnat\\ Thema: \emph{Die Klassifikation der Gruppen bis zur Ordnung
$p^5$} \\ 1. Gutachter: Prof. Dr. Bettina Eick}

\thispagestyle{empty}

\date{}

\maketitle

\tableofcontents

\newpage

\pagenumbering{arabic}

\chapter{Themenstellung und historische Bemerkungen}

\section{Thema der Arbeit und Übersicht über die Ergebnisse}

Diese Arbeit beschäftigt sich mit der Klassifikation von Gruppen der Ordnung $p^n$ bis auf
Isomorphie, wobei $p$ eine Primzahl größer als 3 und $n$ eine natürliche Zahl mit $1 \leq n
\leq 5$ ist. Durch ein konstruktives Beweisverfahren wird für jeden Isomorphietyp explizit
eine Gruppe konstruiert und durch eine endliche Präsentation dargestellt. Die Liste dieser
Präsentationen findet man in der Zusammenfassung. Der folgende Satz gibt einen Überblick über
die Anzahl der Isomorphieklassen (er wird später bewiesen).

\begin{satz}
Es sei $p$ eine Primzahl und $p>3$. Dann gilt:
\begin{enumerate}
  \item Die Gruppen der Ordnung $p$ bilden eine Isomorphieklasse.
  \item Die Gruppen der Ordnung $p^2$ bilden zwei Isomorphieklassen.
  \item Die Gruppen der Ordnung $p^3$ bilden fünf Isomorphieklassen.
  \item Die Gruppen der Ordnung $p^4$ bilden fünfzehn Isomorphieklassen.
  \item Die Gruppen der Ordnung $p^5$ bilden $61 + 2 \cdot p + ggT(4, p-1) + 2 \cdot ggT(3, p-1)$ Isomorphieklassen.
\end{enumerate}
\end{satz}

Der theoretische Hintergrund für das konstruktive Beweisverfahren stammt aus der
algorithmischen Gruppentheorie. Dort hat man einen Algorithmus entwickelt, mit dem man von
$p$"=Gruppen niedrigerer Ordnung aus solche höherer Ordnung konstruieren kann. Beginnt man mit
der Konstruktion bei elementarabelschen $p$"=Gruppen, dann kann man zu den Gruppen einer
beliebig vorgegebenen Ordnung $p^n$ ein Vertretersystem der Isomorphieklassen konstruieren. So
hat man beispielsweise die Gruppen der Ordnung $2^5$ und $3^5$ klassifiziert (vgl.
\cite{bes02}).

Anders als in der algorithmischen Gruppentheorie wird dieser Ansatz hier nicht verwendet, um
$p$"=Gruppen für eine \emph{vorab fest gewählte} Primzahl zu klassifizieren, sondern um die
$p$"=Gruppen bis zur Ordnung $p^5$ für \emph{variables} $p$ zu klassifizieren. Durch den
algorithmischen Ansatz ist die Beweisführung systematischer als für die Fälle $p, p^2, p^3$
und $p^4$, zu denen Anzahl und Vertretersysteme der Isomorphieklassen seit längerem bekannt
sind. Der Fall $p^5$ tritt an dieser Stelle neu hinzu.

Als Nebenresultat wird für die Ordnungen $p$, $p^2$, $p^3$ und $p^4$ die Automorphismengruppe
zu jedem Vertreter der Isomorphieklassen angegeben. Für einige Gruppen der Ordnung $p^5$ lässt
sich die Automorphismengruppen aus den Beweisschritten unmittelbar ablesen.

\section{Historische Bemerkungen}

Bis zur Mitte des 19. Jahrhunderts war der Begriff der Gruppe, so wie er in der heutigen
Mathematik verwendet wird, allenfalls in Ansätzen vorhanden. Die Beschäftigung mit endlichen
Gruppen beschränkte sich im Wesentlichen auf Gruppen von Permutationen, die im Anschluss an
die Arbeiten Galois' recht früh eine Bedeutung in der Algebra gewonnen hatten. Die Aufsätze
\cite{cay54} und \cite{cay59} von Arthur Cayley kann man als erste Schritte dazu ansehen, den
modernen Gruppenbegriff auszubilden und die Anfänge und Ziele der Gruppentheorie als einer
neuen mathematischen Disziplin zu prägen.

\subsection{Die Klassifikation von Gruppen bis auf Isomorphie}

Cayley stellt mit diesen beiden Veröffentlichungen ein Konzept vor, das großen Einfluss
erworben hat: Statt sich mit konkreten mathematischen Objekten wie etwa Permutationen zu
beschäftigen, fragt Cayley auf einer abstrakteren Ebenen nach der "`Struktur"' einer Gruppe,
die durch unterschiedliche Objekte wie etwa Zahlen, Matrizen oder Permutationen realisiert
werden kann. Die Frage nach der "`Struktur"' wird über den Begriff der Isomorphie exakt
gefasst und bietet damit die Möglichkeit, eine Äquivalenzrelation für Gruppen zu definieren
und Gruppen in Isomorphieklassen einzuteilen. Diese Art der Klassifikation ist die
feinkörnigste, die aus Sicht der abstrakten Gruppentheorie sinnvoll und wünschenswert ist, und
stellt vor allem im Bereich der endlichen Gruppen neben vielen anderen Zielen eines der
Hauptanliegen der Gruppentheorie dar.

\subsection{Das Aufkommen der algorithmischen Gruppentheorie}

In der zweiten Hälfte des 20. Jahrhunderts hat sich ein neuer Zweig der Gruppentheorie
entwickelt, der mit algorithmischen Methoden arbeitet und beispielsweise die Klassifikation
endlicher Gruppen durch rechnergestützte Verfahren vorantreibt (vgl. \cite{bes02}).

Sofern möglich, werden in der algorithmischen Gruppentheorie neben den klassischen Objekten
der frühen Gruppentheorie wie Matrix- und Permutationsgruppen häufig so genannte endliche
Präsentationen von Gruppen betrachtet, mit denen die Struktur einer Gruppe auf eine Menge von
Wörtern abgebildet und durch eine endliche Teilmenge dieser Wörter repräsentiert werden kann.
Damit führt man nicht nur Cayleys "`Tendenz zur abstrakten Sichtweise"' fort, sondern kommt
dem algorithmischen Konzept entgegen, Rechenoperationen als Manipulationen von Zeichenfolgen
aufzufassen. Außerdem lassen sich manche unendlichen Gruppen endlich präsentieren und dadurch
einer algorithmischen Betrachtung zugänglich machen.

In dieser Arbeit werden ausschließlich endliche Gruppen eine Rolle spielen. Endliche
Präsentationen dieser Gruppen werden vor allem aus zwei Gründen verwendet: Mit ihnen lässt
sich einerseits eine Gruppe kompakt darstellen und andererseits wird hier für die
Klassifikation ein Verfahren verwendet, das häufig über Manipulationen endlicher
Präsentationen voranschreitet. Im Kapitel über die Grundlagen aus der Gruppentheorie werden
endliche Präsentationen und einige ihrer Manipulationsmöglichkeiten in aller Kürze
vorgestellt.

\subsection{Abelsche $p$"=Gruppen}

Für abelsche Gruppen lassen sich die Isomorphieklassen endlicher $p$"=Gruppen unmittelbar aus
dem Hauptsatz über endlich erzeugte abelsche Gruppen entnehmen (vgl. etwa \cite{mey80}, S.
78): Ist $G$ eine abelsche $p$"=Gruppe der Ordnung $p^n$, dann gibt es eine Partition $(n_1
\dd n_l)$ von $n$ mit $n_1 \geq \ldots \geq n_l$, sodass $G$ isomorph zum direkten Produkt
$C_{p^{n_1}} \times \ldots \times C_{p^{n_l}}$ der zyklischen Gruppen $C_{p^{n_1}} \dd
C_{p^{n_l}}$ ist.

Die Gesamtheit der Isomorphieklassen abelscher $p$"=Gruppen der Ordnung $p^n$ ist also durch
die Menge der absteigend (oder aufsteigend) geordneten Partitionen von $n$ unmittelbar
gegeben. Für die Gruppenordnungen, die im Rahmen dieser Arbeit von Interesse sind, kann man
daher ohne weitere Überlegungen ein Vertretersystem der Isomorphieklassen abelscher Gruppen
aufschreiben:
\begin{table} [h] \centering
\begin{tabular}{l|l}
  Ordnung & Vertretersystem der Isomorphieklassen \\ \hline
  $p$   & $C_p$ \\
  $p^2$ & $C_{p^2}, C_p^2$ \\
  $p^3$ & $C_{p^3}, C_{p^2} \times C_p, C_p^3$ \\
  $p^4$ & $C_{p^4}, C_{p^3} \times C_p, C_{p^2}^2, C_{p^2} \times C_p^2, C_p^4$ \\
  $p^5$ & $C_{p^5}, C_{p^4} \times C_p, C_{p^3} \times C_{p^2}, C_{p^3} \times C_p^2,
  C_{p^2}^2 \times C_p, C_{p^2} \times C_p^3, C_p^5$ \\
\end{tabular}
\end{table}

\subsection{Nichtabelsche $p$"=Gruppen}

Anders als im Fall der abelschen Gruppen verfügt man für nichtabelsche Gruppen keineswegs über
eine vollständige theoretische Übersicht, die dem Hauptsatz über endlich erzeugte abelsche
Gruppen entsprechen könnte. Man ist weit davon entfernt, für jede beliebige Gruppenordnung ein
Vertretersystem der Isomorphieklassen nichtabelscher Gruppen ohne weiteres Nachdenken angeben
zu können.

Im Fall endlicher $p$"=Gruppen haben Higman \cite{hig60} und Sims \cite{sim65} eine
asymptotische Näherung für die Anzahl der Isomorphieklassen angegeben. Für Gruppen der Ordnung
$p^n$ wird die Zahl der Isomorphietypen durch $$p^{2n^3/27+o\left(n^{8/3}\right)}$$
abgeschätzt, wobei $\lim\limits_{n \ra \infty} o(n^{8/3})=0$ ist.

Man hat auf unterschiedlichen Wegen versucht, $p$"=Gruppen zu klassifizieren. An dieser Stelle
kann kein Überblick über die umfangreiche Literatur gegeben werden. Statt dessen wird die
Entwicklung in groben Zügen skizziert. Nur einige zentrale Arbeiten werden genannt.

Die Isomorphieklassen der $p$"=Gruppen bis zur Ordnung $p^3$ waren bereits im 19. Jahrhundert
bekannt ~-- vor allem in Gestalt von Symmetriegruppen geometrischer Objekte. Diese Liste wurde
kurz vor der Wende zum 20. Jahrhundert durch Hölder mit dem Artikel \cite{hoe93} um die
Isomorphieklassen zur Gruppenordnung $p^4$ ergänzt. Seit diesem Zeitpunkt liegt die folgende
Übersicht über die Isomorphieklassen der $p$"=Gruppen bis zur Ordnung $p^4$ vor:
\begin{table} [h] \centering
\begin{tabular}{l|r|l}
  Ordnung & Anzahl & Vertreter \\ \hline
  $p$   & 1  & $C_p$ \\
  $p^2$ & 2  & $C_{p^2}, C_p^2$ \\
  $p^3$ & 5  & $C_{p^3}, C_{p^2} \times C_p, C_p^3, D_p, Q_p$ \\
  $p^4$ & 15 & keine Standardbezeichnung üblich \\
\end{tabular}
\end{table}

Über die Gruppen der Ordnung $p^5$ sind einige Arbeiten mit verschiedenen Zielen erschienen.
Hall legt in \cite{hal40} ein Klassifikationsverfahren für $p$"=Gruppen vor, das diese Gruppen
nicht bis auf Isomorphie, sondern bis auf eine dort eingeführte Äquivalenzrelation einteilt,
die Hall Isoklinismus nennt. Für die Gruppen der Ordnung $p^5$ gibt es nach Hall zehn
Isoklinismenklassen. Bagnera veröffentlicht mit \cite{bag98} eine Arbeit über die Gruppen der
Ordnung $p^5$. Dieser in Rio de Janeiro erschienene Artikel ist schwer zugänglich und konnte
zum Vergleich mit der vorliegenden Arbeit leider nicht beschafft werden.

\section{Algorithmische Konstruktion von $p$"=Gruppen für festes $p$}

Anders als die frühen Aufsätze über $p$"=Gruppen stützt sich diese Arbeit nicht auf
Ad"=hoc"=Argumente oder auf Anleihen aus der Darstellungstheorie, sondern benutzt den
theoretischen Hintergrund, der in der algorithmischen Gruppentheorie in den siebziger und
achtziger Jahren des zwanzigsten Jahrhunderts von Havas, Newman und O'Brien entwickelt worden
ist, um Probleme im Umfeld der Burnside"=Fragen algorithmisch zu lösen. Die ursprünglichen
Ideen gingen von Higman aus, wurden aber von ihm selbst nicht weiterentwickelt.

Als Nebenprodukt dieser Überlegungen ist ein Algorithmus entstanden, mit dem man aus endlichen
Präsentationen von $p$"=Gruppen niedrigerer Ordnung solche höherer Ordnung konstruieren kann.
Wenn man eine Primzahl $p$ fest auswählt und wenn man mit der Konstruktion bei den
elementarabelschen Gruppen anfängt, erhält man durch dieses Verfahren eine Liste endlicher
Präsentationen von $p$"=Gruppen, in der jede Isomorphieklasse von $p$ an über jede Potenz von
$p$ bis zu einer beliebig vorgegebenen Schranke $p^n$ genau einmal auftritt. Damit ist es
möglich, für eine vorgegebene Primzahl $p$ die Gruppen der Ordnung $p, p^2 \dd p^n$ durch
algorithmische Methoden bis auf Isomorphie zu klassifizieren. Die Grenzen dieses Verfahrens
sind allein durch den gegenwärtigen Stand der Rechnertechnik gesetzt.

Für einen ersten Überblick über das algorithmische Verfahren und über die Art und Weise, wie
sich die spätere Argumentation entwickeln wird, werden hier die Grundideen des Ansatzes
informell und ohne technische Details vorgestellt. Einer ausführlichen Darstellung ist ein
späteres Kapitel gewidmet.

\subsection{Die grundlegenden Begriffe}

Der Algorithmus baut auf zwei Merkmalen einer endlichen $p$"=Gruppe $G$ auf, die unter jedem
Automorphismus von $G$ erhalten bleiben: Einerseits lässt sich nach dem Basissatz von Burnside
der Gruppe $G$ eine natürliche Zahl $d$ zuordnen, sodass jedes \emph{minimale
Erzeugendensystem} von $G$ genau $d$ Elemente enthält. Andererseits besitzt jede endliche
$p$"=Gruppe $G$ eine charakteristische Folge von Untergruppen $$G=\gamma_0(G) \geq \gamma_1(G)
\geq \ldots,$$ die in dieser Arbeit in Anlehnung an den englischen Sprachgebrauch als
\emph{absteigende $p$"=Zentralreihe} bezeichnet wird. Da die absteigende $p$"=Zentralreihe
einer endlichen $p$"=Gruppe nach endlichen vielen Schritten auf der trivialen Untergruppe
stationär wird, lässt sich $G$ eindeutig der Wert $c$ zuordnen, unter dem sie zum ersten Mal
die triviale Untergruppe erreicht. Dieser Wert $c$ wird die \emph{$p$"=Klasse} von $G$
genannt.

Mit Hilfe dieser beiden charakteristischen Daten über $G$ lässt sich der Schlüsselbegriff
definieren, der dem gesamten Algorithmus zugrunde liegt, nämlich der Begriff des
\emph{unmittelbaren Nachfolgers} einer $p$"=Gruppe: Hat $G$ die $p$"=Klasse $c$ und $d$ als
minimale Anzahl von Erzeugern, so ist eine $p$"=Gruppe $H$ ein \emph{Nachfolger} von $G$, wenn
$H$ ebenfalls $d$ als minimale Anzahl von Erzeugern hat und wenn die Faktorgruppe
$H/\gamma_c(H)$ isomorph zu $G$ ist. Eine Gruppe $H$ ist ein \emph{unmittelbarer Nachfolger}
von $G$, wenn $H$ ein Nachfolger von $G$ ist und außerdem die $p$"=Klasse $c+1$ hat.

Über den Begriff des Nachfolgers lässt sich ein Verfahren beschreiben, mit dem man ausgehend
von der elementarabelschen Gruppe $C_p^d$ jede endliche $p$"=Gruppe ermitteln kann, die $d$
als minimale Anzahl von Erzeugern hat. Ist $G=\gamma_0(G) > \gamma_1(G) > \dots >
\gamma_{c-1}(G) > \gamma_c(G) = \{1\}$ die absteigende $p$"=Zentralreihe von $G$ bis zum
ersten Auftreten der trivialen Untergruppe und hat $G$ den Wert $d$ als minimale Anzahl von
Erzeugern, so gilt:
\begin{enumerate}
  \item Nach dem Basissatz von Burnside ist die Faktorgruppe
  $G/\gamma_1(G)$ isomorph zu $C_p^d$. Daher ist $G$ ein Nachfolger der
  elementarabelschen Gruppe $C_p^d$.
  \item Für jedes $i$ mit $1 \leq i \leq c-1$ gilt: Die Faktorgruppe $G/\gamma_{i+1}(G)$
  ist ein unmittelbarer Nachfolger von $G/\gamma_i(G)$.
\end{enumerate}

\subsection{Die Idee des Algorithmus}

Mit den bisherigen Überlegungen zeichnet sich der Weg ab, über den man endliche $p$"=Gruppen
berechnen kann: Jede endliche $p$"=Gruppe mit $d$ als minimaler Anzahl von Erzeugern ist ein
Nachfolger --~wenn auch nicht unbedingt ein unmittelbarer~-- der elementarabelschen Gruppe
$C_p^d$. Außerdem gibt es für $G$ eine endliche Folge von $p$"=Gruppen, sodass aufsteigend von
$C_p^d$ bis zu $G$ jedes Glied dieser Folge ein unmittelbarer Nachfolger des vorangegangenen
ist und dass mit jedem Schritt in dieser Folge die $p$"=Klasse der Gruppe um 1 größer ist als
die der vorangegangenen. Daher steigt auch mit jedem Schritt die Ordnung der Gruppen
wenigstens um den Faktor $p$.

Man kann also auf folgendem Weg alle Isomorphieklassen von $p$"=Gruppen bis zu einer
vorgegebenen $p$"=Klasse $c$ und $d$ als minimaler Anzahl von Erzeugern erreichen:
\begin{enumerate}
  \item Man beginne mit der elementarabelschen Gruppe $C_p^d$, die die
  $p$"=Klasse $1$ hat, d.\,h. man setze für eine Menge $L_1$ als
  Anfangswert des Verfahrens $L_1=\{C_p^d\}$.
  \item Man bestimme alle unmittelbaren Nachfolger der Gruppen aus
  $L_i$ und fasse sie in einer Menge $L_{i+1}$ zusammen.
  Alle Elemente von $L_{i+1}$ haben die $p$"=Klasse $i+1$.
  \item Man wiederhole Schritt $2$, bis $i$ den Wert $c$
  erreicht.
\end{enumerate}
Wir werden sehen, dass es einfacher ist, unmittelbare Nachfolger zu ermitteln als beliebige.
Daher ist der Algorithmus so kleinschrittig über die unmittelbaren Nachfolger aufgebaut.

\subsection{Die einzelnen Schritte: Unmittelbare Nachfolger}

Bisher ist die Rede von einer Ermittlung der unmittelbaren Nachfolger gewesen und damit offen
geblieben, ob mit einer Ermittlung die Berechnung im Sinne Turings oder eine Spielart des
Kaffeesatzlesens gemeint ist. In der Tat ist das Verfahren eine algorithmische Methode, die
gleich etwas näher beschrieben wird und bei deren Beschreibung besonders auf die beiden
folgenden Fragen eingegangen werden soll:
\begin{enumerate}
  \item Wie lassen sich unmittelbare Nachfolger berechnen?
  \item Wie lässt sich Redundanz vermeiden? Da in der abstrakten Gruppentheorie nicht die Realisierung der Gruppen von Interesse ist,
sondern nur ihre Isomorphieklassen, ist es wünschenswert, eine  Liste unmittelbarer
Nachfolgern zu erhalten, in der keine Gruppe zu einer anderen isomorph ist.
\end{enumerate}

Beide Fragen werden im Artikel \cite{obr90} zusammenfassend beantwortet, indem O'Brien etliche
Erkenntnisse benutzt, die im Zusammenhang mit den Burnside"=Fragen insbesondere durch Havas
und Newman entstanden sind.

Zu jeder endlichen $p$"=Gruppe $G$ gibt es eine endliche $p$"=Gruppe $P(G)$, welche die
folgende Eigenschaft hat: Jede Gruppe $H$, die eine elementarabelsche Untergruppe $Z$ enthält,
sodass $H/Z$ isomorph zu $G$ ist, ist selbst isomorph zu einer Faktorgruppe von $P(G)$. Da die
letzte nichttriviale Stufe einer absteigenden $p$"=Zentralreihe eine elementarabelsche
Untergruppe ist, ist insbesondere jeder unmittelbare Nachfolger von $G$ zu einer Faktorgruppe
von $P(G)$ isomorph. Die Gruppe $P(G)$ wird im Englischen das \emph{$p$"=Cover} von $G$
genannt. Dieser Ausdruck wird hier unübersetzt ebenfalls benutzt.

Die Verwendung des bestimmten Artikels legt es nahe, dass die Gruppe $P(G)$ eindeutig bestimmt
sei. Sie ist in der Tat (bis auf Isomorphie) eindeutig bestimmt und kann nach \cite{hav77} aus
einer endlichen Präsentation von $G$ berechnet werden, indem man eine endliche Präsentation
von $G$ nach dem algorithmischen Verfahren aus \cite{hav77} erweitert und --~sofern es nötig
ist~-- nach dem Algorithmus von Knuth"=Bendix in eine konsistente endliche Präsentation
überführt. Das Knuth"=Bendix"=Verfahren steht allgemein für endliche Präsentationen zur
Verfügung und wird beispielsweise in \cite{sim94} auf den Seiten 43 bis 95 ausführlich
beschrieben.

Für die weitere Schilderung des Verfahrens sind zwei Untergruppen des $p$"=Covers $P(G)$ von
$G$ von Bedeutung, nämlich der \emph{Multiplikator} $M(G)$ und der \emph{Nukleus} $N(G)$.

Es lässt sich zeigen, dass auch $P(G)$ eine elementarabelsche Untergruppe $M(G)$ enthält,
sodass die Faktorgruppe $P(G)/M(G)$ isomorph zu $G$ ist. Die Untergruppe $M(G)$ ist bis auf
die Wahl eines Epimorphismus von $P(G)$ auf $G$ eindeutig bestimmt. Das Konstruktionsverfahren
der Präsentation von $P(G)$ aus der von $G$ zeichnet einen dieser Epimorphismen aus. Man kann
daher $M(G)$ als eindeutig bestimmt annehmen (die Isomorphieklassen der unmittelbaren
Nachfolger von $G$, die über $P(G)$ und $M(G)$ berechnet werden, sind von der Wahl von $M(G)$
unabhängig).

Die zweite wichtige Untergruppe von $P(G)$ ist der Nukleus $N(G)$. Die Untergruppe $N(G)$ ist
als $c$"=te Stufe der absteigenden $p$"=Zentralreihe von $P(G)$ definiert. Da $P(G)$
einerseits mindestens die $p$"=Klasse $c$ hat (andernfalls gäbe es keinen Epimorphismus von
$P(G)$ auf $G$) und da $P(G)$ andererseits höchstens die $p$"=Klasse $c+1$ hat (andernfalls
wäre $M(G)$ nicht elementarabelsch), ist $N(G)$ trivial oder eine Untergruppe von $M(G)$ (denn
sonst hätte $G$ eine größere $p$"=Klasse als $c$).

Jede Untergruppe $U$ von $M(G)$, die ein Supplement zu $N(G)$ in $M(G)$ ist, liefert in
Gestalt der Faktorgruppe $P(G)/U$ einen unmittelbaren Nachfolger von $G$, sofern $P(G)$ die
$p$"=Klasse $c+1$ hat. Diese Aussage ist ein zentrales Ergebnis der Arbeit \cite{obr90}.
Untergruppen von $M(G)$ mit dieser Eigenschaft werden als \emph{zulässige Untergruppen} von
$P(G)$ bezeichnet. Falls also $G$ überhaupt unmittelbare Nachfolger besitzt, erhält man
dadurch eine vollständige Liste ihrer Isomorphieklassen, dass man $P(G)$ nach jeder zulässigen
Untergruppe von $P(G)$ faktorisiert.

\subsection{Vermeidung von Redundanz: Bahnen zulässiger Untergruppen}

Nach dem bisher geschilderten Verfahren erhält man zwar eine \emph{vollständige} Liste
unmittelbarer Nachfolger von $G$, aber andererseits auch eine Liste, die möglicherweise
mehrere unmittelbare Nachfolger derselben Isomorphieklasse enthält. Das $p$"=Cover zur Gruppe
$C_p^3$ enthält beispielsweise $(p^6-1)/(p-1)$ zulässige Untergruppen einer solchen Ordnung,
dass die Faktorgruppen von $P(C_p^3)$ nach diesen Gruppen unmittelbare Nachfolger der Ordnung
$p^4$ sind. Aber die unmittelbaren Nachfolger der Ordnung $p^4$ bilden nur vier verschiedene
Isomorphieklassen. Da es im allgemeinen schwer festzustellen ist, ob zwei Gruppen, die durch
endliche Präsentationen gegeben sind, isomorph zueinander sind, stellt sich die Frage, ob man
vor dem Faktorisieren eine Auswahl unter den zulässigen Untergruppen treffen kann, sodass sich
eine Redundanz unter den unmittelbaren Nachfolgern nicht ergibt.

Als weitere zentrale Erkenntnis stellt der Artikel \cite{obr90} eine Methode zur Verfügung,
mit der man die zulässigen Untergruppen in Äquivalenzklassen einteilen kann, die den
Isomorphieklassen der unmittelbaren Nachfolger entsprechen. Diese Äquivalenzklassen sind die
Bahnen der zulässigen Untergruppen unter einer bestimmten Operation der Automorphismengruppe
von $G$ auf $M(G)$: Jedem Automorphismus $\al$ von $G$ lässt sich ein so genannter
\emph{Erweiterungsautomorphismus} $\al^* \in Aut(P(G))$ zuordnen, der $M(G)$ invariant lässt.
Daher operiert die Automorphismengruppe $Aut(G)$ von $G$ über Erweiterungsautomorphismen auf
$M(G)$. Die Bahnen der zulässigen Untergruppen unter dieser Operation entsprechen den
Äquivalenzklassen der unmittelbaren Nachfolger, die sich durch Faktorisierung von $P(G)$ nach
zulässigen Untergruppen ergeben.

Da der Nukleus als Teil der absteigenden $p$"=Zentralreihe eine charakteristische Untergruppe
von $P(G)$ ist, braucht man die Operation von $Aut(G)$ nicht auf die Menge der zulässigen
Untergruppen einzuschränken, sondern kann die Bahnen auf der Menge aller Untergruppen des
Multiplikators betrachten, denn jede dieser Bahnen enthält keine oder ausschließlich zulässige
Untergruppen.

Zur Berechnung solcher Bahnen kann man Methoden aus der linearen Algebra zur verwenden: Der
Multiplikator $M(G)$ ist eine elementarabelsche Gruppe und der endlich Körper $\f$ aus $p$
Elementen operiert durch $\f \times M(G) \ra M(G): (f,m) \mt m^f$ im Sinn einer
Skalarmultiplikation auf $M(G)$. Daher kann man $M(G)$ als Vektorraum über $\f$ auffassen. Hat
$M(G)$ die Ordnung $p^k$, so ist $M(G)$ isomorph zur additiven Gruppe $\f^k$. Daher gibt es
einen Operationshomomorphismus $\vi: Aut(G) \ra GL(k,p)$, sodass die Operation von $Aut(G)$
über Erweiterungsautomorphismen auf $M(G)$ der Operation von $Aut(G)$ über $\vi$ auf $\f^k$
entspricht. Aus diesem Grunde kann man Matrixgruppen verwenden, um die Bahnen der zulässigen
Untergruppen zu ermitteln.

\subsection{Zusammenfassende Übersicht}

Insgesamt setzt sich die Berechnung der unmittelbaren Nachfolger einer $p$"=Gruppe $G$ aus den
folgenden Schritten zusammen:
\begin{enumerate} \label{algunmittnachf}
  \item Man berechne aus einer Präsentation von $G$ nach \cite{hav77}
  und dem Knuth"=Bendix"=Algorithmus eine Präsentation von $P(G)$.
  Man lese $M(G)$ unmittelbar aus der Präsentation von $P(G)$ ab
  und bestimme $N(G)$ und damit die $p$"=Klasse von $P(G)$. Unter
  günstigen Umständen lässt sich $N(G)$ ebenfalls unmittelbar aus
  der Präsentation von $P(G)$ ablesen, andernfalls stehen
  algorithmische Verfahren zur Verfügung.
  \item Falls $P(G)$ genauso wie $G$ die $p$"=Klasse $c$ hat, so
  hat $G$ keine unmittelbaren Nachfolger. In diesem Fall ist das
  Verfahren abgeschlossen. Andernfalls fahre man fort.
  \item Man berechne, wie $Aut(G)$ über Erweiterungsautomorphismen auf
  $M(G) \cong \f^k$ operiert, bzw. man berechne den entsprechenden
  Operationshomomorphismus $\vi: Aut(G) \ra GL(k,p)$. Für diesen
  Schritt sind zu jedem Element aus $Aut(G)$ oder --~falls
  bekannt~-- auch nur zu jedem Erzeuger von $Aut(G)$
  endlich viele Wortoperationen in $P(G)$ durchzuführen.
  \item Man ermittle die Bahnen der Untergruppen in $M(G)$ oder
  der entsprechenden Unterräume in $\f^k$ unter der Operation von
  $Aut(G)$ und bilde aus den Bahnen, die zulässige Untergruppen
  enthalten, ein Vertretersystem $V$. Für diesen Schritt stehen die
  Bahn"=Stabilisator"=Algorithmen zur Verfügung. Falls nur unmittelbare
  Nachfolger einer bestimmten Ordnung von Interesse sind, kann
  dieser Schritt auf Untergruppen von $M(G)$ einer entsprechenden
  Ordnung oder Unterräume von $\f^k$ einer entsprechenden Dimension
  eingeschränkt werden.
  \item Man bilde für jedes $U$ aus $V$ die Faktorgruppe $P(G)/U$.
  So erhält man eine vollständige Liste unmittelbarer
  Nachfolger von $G$, in der jede Isomorphieklasse genau einmal
  vertreten ist.
\end{enumerate}

Mit dieser Befehlsfolge ist der Algorithmus skizziert, der zu einer vorgegebenen $p$"=Gruppe
bzw. einer ihrer endlichen Präsentationen alle Isomorphieklassen unmittelbarer Nachfolger
redundanzfrei berechnet. Zu diesem Algorithmus liegen Implementationen für einige
Computer"=Algebra"=Systeme vor. So hat Newman beispielsweise einige Ergänzungspakete für
\gap{} \cite{gap} initiiert, die diesen Algorithmus umsetzen.

Dass jeder Schritt der oben aufgeführten Befehlsfolge algorithmisch durchführbar ist, lässt
sich mit einer Ausnahme unmittelbar erkennen: Es wird stillschweigend vorausgesetzt, dass die
Automorphismengruppe von $G$ bekannt ist. Diese Bedingung scheint für die meisten $p$"=Gruppen
nicht erfüllt zu sein. Gerade der oben dargestellte Algorithmus erlaubt es aber, die
Automorphismengruppen der unmittelbaren Nachfolger einer $p$"=Gruppe aus der
Automorphismengruppe von $G$ zu berechnen. Dies hat O'Brien bereits im Artikel \cite{obr90}
dargelegt. Spätere Veröffentlichungen von O'Brien \cite{obr94} und Eick u.\,a. \cite{eic02}
beschreiben dieses Verfahren genauer und schlagen einige Optimierungsmöglichkeiten vor. Die
Lösung des Problems ist überraschend einfach: Ist $U$ eine zulässige Untergruppe von $P(G)$,
so lässt sich die Automorphismengruppe von $P(G)/U$ unmittelbar aus dem Stabilisator von $U$
unter der Operation von $Aut(G)$ über Erweiterungsautomorphismen auf $M(G)$ ablesen. Beginnt
man damit, $p$"=Gruppen ausgehend von einer elementarabelschen Gruppe $C_p^d$ zu berechnen, so
macht der Ausgangspunkt des Verfahrens keine Schwierigkeiten, denn die Automorphismengruppe
von $C_p^d$ ist in der Gestalt von $GL(d,p)$ hinreichend bekannt.

\chapter{Grundlagen aus der Gruppentheorie}

\section{Endliche Gruppen von Primzahlpotenzordnung}

\begin{defi}
Eine Gruppe $G$ ist eine \emph{endliche $p$"=Gruppe} zur Primzahl $p$, wenn es eine natürliche
Zahl $n$ gibt mit $|G|=p^n$.
\end{defi}

\begin{beme}
Da in dieser Arbeit nur endliche $p$"=Gruppen eine Rolle spielen, wird statt "`endliche
$p$"=Gruppe"' nur "`$p$"=Gruppe"' geschrieben.
\end{beme}

\begin{defi}
Die \emph{absteigende $p$"=Zentralreihe} $$\gamma_0(G) \geq \gamma_1(G) \geq \ldots \geq
\gamma_i(G) \geq \gamma_{i+1}(G) \geq \ldots$$ einer $p$"=Gruppe $G$ ist folgendermaßen
rekursiv definiert: Es ist $\gamma_0(G)=G$ und
$$\gamma_{i+1}(G)=[\gamma_i(G),G]\gamma_i(G)^p=\langle [h,g]k^p \mid h,k \in \gamma_i(G)
\text{ und } g \in G \rangle,$$ wobei $[h,g]=h\1 g\1 h g$ der \emph{Kommutator} von $h$ und
$g$ ist. Wenn $\gamma_c(G)=\{1\}$ ist und $c$ die kleinste natürliche Zahl mit dieser
Eigenschaft ist, so ist $c$ die \emph{$p$"=Klasse} von $G$.
\end{defi}

%

\begin{satz}
(Basissatz von Burnside) Ist $G$ eine $p$"=Gruppe mit $|G/\gamma_1(G)|=p^d$, so enthält jedes
minimale Erzeugendensystem von $G$ genau $d$ Elemente.
\end{satz}

\begin{beweis}
Siehe \cite{hup67}, Seite 273.
\end{beweis}

\begin{satz} \label{kommzentreihe} \label{frat} \label{dlem}
Es sei $G=\erz{a_1 \dd a_d}$ eine endliche $p$"=Gruppe und $G/\gamma_1(G)$ mit
$|G/\gamma_1(G)|=p^d$ werde von den Bildern von $a_1 \dd a_d$ unter dem natürlichen
Homomorphismus $\nu: G \ra G/\gamma_1(G)$ erzeugt.
\begin{enumerate}
  \item Ist $\theta$ ein Homomorphismus von $G$, so ist
  $\gamma_i(G)^\theta=\gamma_i(G^\theta)$ für alle $i \in \N_0$. Insbesondere ist jedes Glied
  der absteigenden $p$"=Zentralreihe ein Normalteiler von $G$.
  \item Ist $N$ ein echter Normalteiler von $G$ und hat $G/N$ die $p$"=Klasse
  $c$, so ist $\gamma_c(G)$ eine Untergruppe von $N$.
  \item Für alle $i, j \in \N_0$ ist die Gruppe $[\gamma_i(G),\gamma_j(G)]$ eine
  Untergruppe von $\gamma_{i+j}(G)$.
  \item Die Faktorgruppe $\gamma_1(G)/\gamma_2(G)$ wird von $\{a_i^p \mid 1 \leq i
  \leq d\} \cup \{[a_j,a_i] \mid 1 \leq i < j \leq d\}$ erzeugt.
  \item Es sei $s \in \N$ und $U_s$ eine Teilmenge von $\gamma_s(G)$, sodass
  die Faktorgruppe $\gamma_s(G)/\gamma_{s+1}(G)$ unter dem natürlichen Homomorphismus $\nu_s:
  \gamma_s(G) \ra \gamma_s(G)/\gamma_{s+1}(G)$  von
  den Bildern der Elemente aus $U_s$ erzeugt wird. Dann wird
  $\gamma_{s+1}(G)/\gamma_{s+2}(G)$ von der Menge $\{(u^p)^{\nu_{s+1}} \mid u \in U\} \cup
  \{([u,a_i])^{\nu_{s+1}} \mid u \in U \text{ und } 1 \leq i \leq d\}$ erzeugt, wobei
  $\nu_{s+1}: \gamma_s(G) \ra \gamma_s(G)/\gamma_{s+1}(G)$ der natürliche Homomorphismus ist.
  \item Für jedes $i$, sodass $\gamma_i(G)$ nichttrivial ist, sind die Gruppen
  $\gamma_i(G)^p$ und $[\gamma_i(G),G]$ sowie $\gamma_{i+1}(G)$ echte Untergruppen von $\gamma_i(G)$.
  Insbesondere gibt es ein $c \in \N$, sodass $\gamma_c(G)=\{1\}$ ist.
\end{enumerate}
\end{satz}

\begin{beweis}
Siehe \cite{sim94}, Seite 446.
\end{beweis}

\begin{satz} \label{nzyk}
Ist $G$ eine $p$"=Gruppe und $i \in \N_0$, so ist $\gamma_{i+1}(G)$ der kleinste Normalteiler
$N$ von $G$, der in $\gamma_i(G)$ enthalten ist, sodass $\gamma_i(G)/N$ elementarabelsch und
im Zentrum von $G/N$ enthalten ist.
\end{satz}

\begin{beweis}
Siehe \cite{sim94}, Seite 446.
\end{beweis}

\begin{lemma} \label{autpzent}
Ist $G=\erz{a_1 \dd a_d}$ eine endliche $p$"=Gruppe der Ordnung $p^n$ mit $d$ als minimaler
Anzahl von Erzeugern und ist $N$ eine charakteristische und elementarabelsche Untergruppe des
Zentrums von $G$ und hat $N$ die Ordnung $p^k$, so ist die Untergruppe $V$ von $Aut(G)$, deren
Elemente auf $G/N$ die Identität liefern, ein Normalteiler von $Aut(G)$ und isomorph zur
elementarabelschen Gruppe $C_p^{d(n-k)}$.
\end{lemma}

\begin{beweis}
Es sei $H=G/N$, $\psi: G \ra N$ der natürliche Homomorphismus und $T$ eine Transversale zu $N$
in $G$. Da $N$ eine charakteristische Untergruppe von $G$ ist, induziert jeder Automorphismus
$\theta$ von $G$ einen Automorphismus $\theta^\vi$ von $H$. Die Abbildung $\vi: Aut(G) \ra
Aut(H)$ ist offensichtlich ein Homomorphismus und der Kern von $\vi$ ist gerade $V$. Also ist
$V$ nach den üblichen Homomorphiesätzen ein Normalteiler von $Aut(G)$. Da $\{a_1 \dd a_d\}$
ein Erzeugendensystem von $G$ ist, ist $\{a_1^\psi \dd a_d^\psi\}$ ein Erzeugendensystem von
$H$. Es sei $$X=\left\{(x_1 \dd x_d) \in G^d \mid x_i^\psi=a_i^\psi \text{ für alle } i \text{
mit } 1 \leq i \leq d\right\}.$$ Da es für jedes $g \in G$ genau ein $t \in T$ und genau ein
$m \in N$ mit $g=tm$ sowie $|N|=p^k$ und $|T|=p^{n-k}$ als auch $g^\psi=(tm)^\psi=t^\psi$ ist,
ist $|X|=p^{d(n-k)}$. Weil $V$ der Kern von $V$ ist und damit für alle $\nu \in V$ die
Bedingung $(a_i^\psi)^{\nu^\vi}=(a_i^\psi)$ für alle $i$ mit $1 \leq i \leq d$ erfüllt ist,
ist $(x_1^{\nu} \dd x_d^{\nu}) \in X$ für alle $\nu \in V$. Also operiert $V$ in dieser Weise
auf $X$. Es sei $B$ die Bahn eines Elementes $x=(x_1 \dd x_d)$ aus $X$ unter dieser Operation
von $V$ in $X$. Nach den Bahnen"=Stabilisator"=Sätzen ist $|B|$ gleich dem Index
$[V:Stab_V(x)]$. Es sei $\nu \in Stab_V(x)$. Dann ist $x_i^\nu=x_i$ für alle $i$ mit $1 \leq i
\leq d$. Da $G$ nach dem Basissatz von Burnside von $\{x_1 \dd x_d\}$ erzeugt und damit $\nu$
jedes Element aus $G$ invariant lässt, ist $\nu$ die identische Abbildung. Also ist
$|Stab_V(x)|=1$ und $|B|=|V|$. Da $x$ beliebig aus $X$ gewählt ist, hat jede Bahn in $X$ unter
$V$ die Länge $|V|$, d.\,h. es gibt ein $k \in \N$ mit $|X|= k |V|$. Also ist $|V|$ ein Teiler
von $p^{d(n-k)}$ und $V$ damit eine $p$"=Gruppe, deren Ordnung $p^{d(n-k)}$ teilt. Nun wird
eine Untergruppe $U$ von $V$ angegeben, die die Ordnung $p^{d(n-k)}$ hat. Es sei $\nu(n_1 \dd
n_d): G \ra G$ die Abbildung, die für alle $i$ mit $1 \leq i \leq d$ das Element $a_i$ auf
$a_i n_i$ mit $n_i \in N$ abbildet. Da $\{a_1 \dd a_d\}$ ein Erzeugendensystem von $G$ ist,
ist mit dieser Definition $\nu(n_1 \dd n_d)$ auf ganz $G$ erklärt, und da $N$ zentral in $G$
ist, ist $\nu(n_1 \dd n_d)$ ein Homomorphismus. Außerdem ist $(a_1^{\nu(n_1 \dd n_d)} \dd
a_d^{\nu(n_1 \dd n_d)})$ ein Element von $X$ und damit ein Erzeugendensystem von $G$. Daher
ist $\nu(n_1 \dd n_d)$ ein Automorphismus von $G$. Da $\nu(n_1 \dd n_d)^\vi=id_H$ ist, ist
$\nu(n_1 \dd n_d)$ ein Element von $V$. Es sei $U=\{\nu(n_1 \dd n_d) \mid n_1 \dd n_d \in
N\}$. Da $id_G \in U$ ist und für $\nu(n_1 \dd n_d)$ und $\nu(m_1 \dd m_d) \in U$ die
Bedingung $\nu(n_1 \dd n_d) \circ \nu(m_1 \dd m_d)\1=\nu(n_1 \dd n_d) \circ \nu(m_1\1 \dd
m_d\1)=\nu(n_1 m_1\1 \dd n_d m_d\1) \in U$ erfüllt ist, ist $U$ eine Untergruppe von $V$ und
elementarabelsch, da $N$ elementarabelsch ist. Die Ordnung von $U$ ist
$|N|^d=(p^{n-k})^d=p^{d(n-k)}$. Also ist $U=V$. (Ähnlich \cite{suz82}, Seite 242)
\end{beweis}

\section{Polyzyklische Gruppen und ihre Präsentationen}

Da in den späteren Teilen der Arbeit nur endliche $p$"=Gruppe betrachtet werden, wird auch
dieser Abschnitt auf den endlichen Fall eingeschränkt.

\begin{defi}
Es sei $G$ eine endliche Gruppe. Eine \emph{polyzyklische Folge} der Länge $n$ ist eine Folge
von Untergruppen $G=G_1 \supseteq G_2 \supseteq \ldots \supseteq G_n \supseteq G_{n+1}=\{1\}$,
sodass für alle $i$ mit $1 \leq i \leq n$ die Untergruppe $G_{i+1}$ normal in $G_i$ und die
Faktorgruppe $G_i/G_{i+1}$ zyklisch ist. Eine Gruppe $G$ ist \emph{polyzyklisch}, wenn $G$
eine polyzyklische Folge hat. Eine Menge $A=\{a_1 \dd a_n\} \subseteq G$ ist ein
\emph{polyzyklisches Erzeugendensystem} bzw. die Folge $(a_1 \dd a_n)$ eine
\emph{polyzyklische Erzeugerfolge} von $G$ der Länge $n$, wenn $G$ von der Menge $A$ und die
Faktorgruppe $G_i/G_{i+1}$ für jedes $i$ mit $1 \leq i \leq n$ von $a_i G_{i+1}$ erzeugt wird.
Hat die Faktorgruppe $G_i/G_{i+1}$ die Ordnung $m_i$, so ist $m_i$ die \emph{relative Ordnung}
von $a_i$.
\end{defi}

\begin{lemma} \label{poly}
\begin{enumerate}
  \item Jede polyzyklische Gruppe hat ein polyzyklisches Erzeugendensystem.
  \item Hat ein polyzyklisches Erzeugendensystem einer Gruppe $G$ die Länge
  $n$, so hat jedes polyzyklisches Erzeugendensystem von $G$ die Länge $n$.
  \item Faktorgruppen polyzyklischer Gruppen sind polyzyklisch.
  \item Untergruppen polyzyklischer Gruppen sind polyzyklisch.
  \item Ist $N$ ein Normalteiler einer Gruppe $G$ und sind sowohl $N$ als auch $G/N$
  polyzyklisch, so ist auch $G$ polyzyklisch.
  \item Hat $G$ ein polyzyklisches Erzeugendensystem der Länge $n$, so wird
  jede Untergruppe von $G$ von höchstens $n$ Elementen erzeugt.
  \item Ist $A=\{a_1 \dd a_n\}$ ein polyzyklisches Erzeugendensystem der Gruppe
  $G$ und $g$ ein Element von $G$, so gibt es ganze Zahlen $e_1 \dd e_n$,
  sodass $g=a_1^{e_1} \mal a_n^{e_n}$ ist. Die Exponenten $e_1 \dd e_n$ sind eindeutig bestimmt,
  wenn für alle $i$ die Zahl $e_i$ so gewählt wird, dass $0 \leq e_i < m_i=|\erz{a_i G_{i+1}}|$ ist.
\end{enumerate}
\end{lemma}

\begin{beweis}
Siehe \cite{sim94}, Seite 390 bis 392 sowie 394 und 395.
\end{beweis}

\begin{defi}
Ist $A=\{a_1 \dd a_n\}$ ein polyzyklisches Erzeugendensystem der Gruppe $G$ und ist
$g=a_1^{e_1} \mal a_n^{e_n}$ mit $0 \leq e_i < m_i=|\erz{a_i G_{i+1}}|$, so ist $a_1^{e_1}
\mal a_n^{e_n}$ die nach \ref{poly} eindeutig bestimmte \emph{Normalform} von $g$ bezüglich
$A$. In diesem Fall heißt $(e_1 \dd e_n) \in \Z^n$ der \emph{Exponenten"=} oder
\emph{Koeffizientenvektor} von $g$ bezüglich $A$. Der erste Eintrag des Exponentenvektors von
$g$ ungleich Null ist der \emph{führende Exponent} oder \emph{führende Koeffizient} von $g$.
\end{defi}

\begin{defi} \label{praes}
Es sei $A=\{a_1 \dd a_n\}$ eine endliche Menge, $A^*$ die Menge der Wörter über $A$ und
$\{(u_1,v_1) \dd (u_m,v_m)\}$ eine endliche Teilmenge von $A^* \times A^*$. Ist $G$ eine
Gruppe und ist $G \cong \erz{a_1 \dd a_n \mid u_1=v_1 \dd u_m=v_m}$, so ist $\erz{a_1 \dd a_n
\mid u_1=v_1 \dd u_m=v_m}$ eine \emph{endliche Präsentation} von $G$. Die Elemente der Menge
$\{(u_1,v_1) \dd (u_m,v_m)\}$ bzw. die ihnen entsprechenden Gleichungen in $\erz{a_1 \dd a_n
\mid u_1=v_1 \dd u_m=v_m}$ heißen dann \emph{Relationen}. Ist $u=v$ eine Relation, so ist $u$
die linke Seite von $u=v$ und $v$ die rechte Seite von $u=v$, und $u=v$ ist \emph{trivial},
wenn $v=1$ ist.
\end{defi}

\begin{lemma} \label{stdpraes}
Es sei $A=\{a_1 \dd a_n\}$ ein polyzyklisches Erzeugendensystem einer endlichen Gruppe $G$,
wobei für jedes $i$ mit $1 \leq i \leq n$ das Element $a_i$ die relative Ordnung $m_i$ hat.
Dann gibt es genau eine endliche Präsentation $\erz{a_1 \dd a_n \mid S}$ von $G$, sodass die
Elemente von $S$ für alle $i,j \in \{1 \dd n\}$ folgendermaßen gegeben sind
\begin{eqnarray*}
   a_i^{m_i}   &=&       a_{i+1}^{e_1(i,i+1)} \mal a_n^{e_1(i,n)} \\
   a_j a_i     &=& a_i   a_{i+1}^{e_2(i,j,i+1)} \mal a_n^{e_2(i,j,n)} \text{ für } j>i
\end{eqnarray*}
Dabei sind die Werte der Funktionen $e_1$ und $e_2$ so bestimmt, dass jede rechte Seite einer
Relation von $S$ in Normalform ist.
\end{lemma}

\begin{beweis}
Siehe \cite{sim94}, Seite 395.
\end{beweis}

\begin{defi}
Ist $G$ eine polyzyklische Gruppe und $A=\{a_1 \dd a_n\}$ ein polyzyklisches Erzeugendensystem
von $G$, so heißt die im Sinne von \ref{stdpraes} eindeutig bestimmte Präsentation $\erz{a_1
\dd a_n \mid S}$ von $G$ die \emph{polyzyklische Standardpräsentation} von $G$ bezüglich $A$.
Die Relationen des Typs 1 heißen \emph{Potenzrelationen} und die Relationen des Typs 2 heißen
\emph{Kommutatorrelationen}.
\end{defi}

\begin{lemma} \label{ppoly}
Jede endlich $p$"=Gruppe $G$ ist polyzyklisch.
\end{lemma}

\begin{beweis}
Nach \ref{frat} gibt es ein kleinstes $c \in \N$, sodass $\gamma_c(G)$ trivial ist. Nun wird
durch Induktion über die Glieder der absteigenden $p$"=Zentralreihe gezeigt, dass $G$
polyzyklisch ist. Die triviale Gruppe $\gamma_c(G)$ ist polyzyklisch. Nun sei bereits gezeigt,
dass $\gamma_i(G)$ polyzyklisch ist. Da $H=\gamma_{i-1}(G)/\gamma_i(G)$ nach \ref{nzyk} eine
elementarabelsche Gruppe und von endlichem Index ist, ist $H$ isomorph zu $C_p^s$ mit
$p^s=|H|$ und es gibt daher $h_1 \ldots h_s \in H$, sodass $H_0, H_1 \dd H_s$ mit
$H_i=\erz{h_1 \dd h_i}$ für $1 \leq i \leq s$ und $H_0=\{1\}$ eine polyzyklische Folge $H=H_s
> H_{s-1} > \ldots > H_1 > H_0=\{1\}$ von $H$ ist. Damit ist $H$ polyzyklisch. Da sowohl
$\gamma_i(G)$ als auch die Faktorgruppe $H=\gamma_{i-1}(G)/\gamma_i(G)$ polyzyklisch sind, ist
nach \ref{poly} auch $\gamma_{i-1}(G)$ polyzyklisch.
\end{beweis}

\begin{defi}
Eine endliche Präsentation $\erz{a_1 \dd a_n \mid S}$ ist eine
\emph{Potenz"=Kommutator"=Präsentation zur Primzahl $p$}, wenn sämtliche Elemente von $S$
folgendermaßen gegeben sind
\begin{eqnarray*}
   a_i^p              &=& a_{i+1}^{e_1(i,i+1)} \mal a_n^{e_1(i,n)} \\
   \text{} [a_j, a_i] &=& a_{j+1}^{e_2(i,j,j+1)} \mal a_n^{e_2(i,j,n)} \text{ für } j>i
\end{eqnarray*}
Dabei sind die Werte der Funktionen $e_1$ und $e_2$ so bestimmt, dass jede rechte Seite einer
Relation von $S$ in Normalform ist. Die Relationen vom Typ 1 heißen \emph{Potenzrelationen}
und die Relationen vom Typ 2 heißen \emph{Kommutatorrelationen}.
\end{defi}

\begin{defi}
Eine Potenz"=Kommutator"=Präsentation $\erz{a_1 \dd a_n \mid S}$ zur Primzahl $p$ ist
\emph{konsistent}, wenn $|\erz{a_1 \dd a_n \mid S}|=p^n$ ist.

\end{defi}

\begin{satz}
Jede endliche $p$"=Gruppe $G$ besitzt eine Potenz"=Kommutator"=Präsentation zur
Primzahl $p$.
\end{satz}

\begin{beweis}
Siehe \cite{syl72}.
\end{beweis}

\begin{beme} \label{kurz}
Im Weiteren werden in Potenz"=Kommutator"=Präsentationen alle trivialen Relationen nicht
notiert. Statt beispielsweise $$\erz{a_1, a_2, a_3 \mid [a_2,a_1]=a_3, [a_3, a_1]=1, [a_3,
a_2]=1, a_1^p=1, a_2^p=1, a_3^p=1}$$ oder $$\erz{a_1, a_2, a_3 \mid [a_2,a_1]=a_3, [a_3, a_1],
[a_3, a_2], a_1^p, a_2^p, a_3^p}$$ wird nur $$\erz{a_1, a_2, a_3 \mid [a_2,a_1]=a_3}$$
geschrieben.
\end{beme}


\begin{defi} \label{gew}
Eine Potenz"=Kommutator"=Präsentation $\erz{A \mid S}$ zur Primzahl $p$ ist mit $\omega$ und
$D$ \emph{gewichtet}, wenn $D$ eine Teilmenge von $S$ und $\omega: A \ra \N$ eine Abbildung
ist, sodass die folgenden Bedingungen erfüllt sind:
\begin{enumerate}
  \item Es ist $\omega(a_1)=1$ und $\omega(a_i) \leq \omega(a_{i+1})$ für alle
  $i$ mit $1 \leq i \leq n-1$.
  \item Zu jedem $i$ mit $1 \leq i \leq n$ gibt es eine Relation $a_i^p=w_i$
  und in $w_i$ kommen nur solche Elemente $a$ von $A$ vor, sodass
  $\omega(a) \geq \omega(a_i)+1$ ist.
  \item Zu jedem $i$ und $j$, sodass es eine Relation $[a_j, a_i]=u_{ij}$ gibt,
  kommen in $u_{ij}$ nur solche Elemente $a$ von $A$ vor, sodass
  $\omega(a) \geq \omega(a_i)+ \omega(a_j)$ ist.
  \item Ist $\omega(a_k) > 1$, dann ist eine der beiden folgenden Bedingungen
  erfüllt:
\begin{enumerate}
  \item Es gibt ein $i$ und ein $j$ mit $1 \leq i < j \leq n$, sodass
  $\omega(a_i)=1$ und $\omega(a_j)=\omega(a_k)-1$ und $u_{ij}=a_k$ ist.
  \item Es gibt ein $i$ mit $1 \leq i \leq n$, sodass
  $\omega(a_i)=\omega(a_k)-1$ und $w_i=a_k$ ist.
\end{enumerate}
Genau eine der Relationen $[a_j, a_i]=a_k$ oder $a_i^p=a_k$ wird als \emph{Definition} von
$a_k$ ausgezeichnet. Die Menge $D$ besteht genau aus den Definitionen von $\erz{A \mid S}$.

\end{enumerate}
Die Abbildung $\omega$ wird \emph{Gewichtungsfunktion} genannt. Ist $\erz{A \mid S}$ durch
$\omega$ und $D$ gewichtet, so ist der Maximalwert von $\omega$ die \emph{Gewichtsklasse} von
$\erz{A \mid S}$.
\end{defi}

\begin{lemma} \label{gewpzentral}
Es sei $G$ eine $p$"=Gruppe. Weiterhin sei $\erz{a_1 \dd a_n \mid S}$ eine mit $\omega$
gewichtete Potenz"=Kommutator"=Präsentation von $G$, sodass $\omega(a_1) \leq \ldots \leq
\omega(a_n)$ ist, und $A_h$ sei definiert durch $A_h=\{a_i \mid a_i \in A \text{ und }
\omega(a_i) \geq h+1\}$ für alle $h \in \{1 \dd \omega(n)\}$. Dann ist $\gamma_h(G)=\langle
A_h \rangle$ für alle $h \in \N_0$.
\end{lemma}

\begin{beweis}
Die Aussage wird über vollständige Induktion bewiesen. Wenn $h=0$ ist, so ist nach
$\gamma_0(G)=G=\langle A \rangle = \langle A_0 \rangle$ die Behauptung erfüllt. Es sei $h>0$,
und $\gamma_h(G)$ werde nach Induktionsvoraussetzung von $A_h$ erzeugt. Dann ist
$\gamma_{h+1}(G)=[\gamma_h(G),G]\gamma_h(G)^p$, d.\,h. $\gamma_{h+1}(G)$ wird von allen
Elementen $a_l^p$ und $[a_l,a_i]$ erzeugt, sodass $a_l \in A_h$ und $a_i \in A$ ist. Da
$\omega(a_l^p)=\omega(a_l)+1 \geq h+1$ und $\omega([a_l,a_i])= \omega(a_l)+\omega(a_i) \geq
h+1$ ist und damit $a_l^p$ sowie $[a_l,a_i]$ in $A_{h+1}$ liegen, ist $\gamma_{h+1}(G)$ eine
Untergruppe von $\langle A_{h+1} \rangle$.

Umgekehrt ist aber auch $\erz{A_{h+1}}$ eine Untergruppe von $\gamma_{h+1}(G)$: Es sei $a_k
\in A_{h+1}$ und damit $\omega(a_k) \geq h+1 > 1$. Dann ist entweder $a_k=a_i^p$ mit
$\omega(a_k)=\omega(a_i)+1$ oder $\omega(a_k)=[a_l,a_i]$ mit
$\omega(a_k)=\omega(a_i)+\omega(a_j)$ für geeignete Elemente $a_i, a_l \in A_h$ und $a_i \in
A$. Nach Induktionsvoraussetzung und \ref{kommzentreihe} ist dann $a_k$ ein Element aus
$\gamma_{h+1}(G)$. (Beweis nach \cite{hav77}, Seite 217)
\end{beweis}

\begin{folg}
Ist $c$ die $p$"=Klasse von $G$ und hat $G$ eine gewichtete Potenz"=Kommutator"=Präsentation,
so ist $c$ ebenfalls die Gewichtsklasse von $G$.
\end{folg}

\begin{folg}
Ist $\erz{a_1 \dd a_n \mid S}$ eine konsistente und gewichtete
Potenz"=Kommutator"=Präsentation einer endlichen $p$"=Gruppe $G$ und $\{a_1 \dd a_d\}$ die
Teilmenge von $\{a_1 \dd a_n\}$ mit der Gewichtung 1, so enthält jedes nicht zu verkleinernde
Erzeugendensystem von $G$ genau $d$ Elemente.
\end{folg}

%
%
%

\section{Faktorgruppen polyzyklischer Gruppen und ihre Präsentationen}

\begin{defi}
Es sei $A=\{a_1 \dd a_n\}$ ein polyzyklisches Erzeugendensystem einer endliche Gruppe $G$ und
$U=(g_1 \dd g_s)$ eine endliche Folge von Elementen aus $G$, sodass $a_1^{e_{i1}} \mal
a_n^{e_{in}}$ für alle $i$ mit $1 \leq i \leq s$ die Normalform von $g_i$ sei.
\begin{enumerate}
  \item Die Matrix $(e_{ij}) \in Z^{s \times n}$ ist die
  \emph{Exponentenmatrix} von $U$ (bezüglich $A$).
  \item Die Folge $U$ ist in \emph{Standardform}, wenn die Exponentenmatrix $M$
  von $U$ die folgenden Bedingungen erfüllt:
  \begin{enumerate}
    \item Keine Zeile von $M$ ist der Nullvektor.
    \item Die Matrix $M$ ist zeilenreduziert über $\Z$.
    \item Ist $e_{ij}$ ein Eckeintrag von $M$, dann ist $e_{ij}$ ein Teiler der relativen
    Ordnung von $a_j$.
  \end{enumerate}
  \item Ist $U$ ist Standardform und $M=(e_{ij})$ die Exponentenmatrix von $U$,
  so ist $(f_1 \dd f_s) \in \Z^s$ ein \emph{zulässiger Exponentenvektor} von
  $U$, wenn $0 \leq f_k < m_k/e_{ik}$ ist, sofern $a_k$ die relative Ordnung $m_k$
  hat und $e_{ik}$ ein Eckeintrag von $M$ ist.
  \item Ist $U$ in Standardform, so ist $E(U)$ die Menge der zulässigen Exponentenvektoren von $U$ und
  $S(U)=\{g_1^{f_1} \mal g_s^{f_s} \mid (f_1 \dd f_s) \in E(U)\}$, sofern $U$
  in Standardform ist.
  \item Die Folge $U$ ist \emph{voll}, wenn $U$ in Standardform ist und die
  folgenden beiden Bedingungen erfüllt sind:
  \begin{enumerate}
    \item Für alle $i$ und $j$ mit $1 \leq i < j \leq s$ enthält $S(U)$ das
    Element $g_i\1 g_j g_i$.
    \item Ist $e_{ij}$ ein Eckeintrag der Exponentenmatrix von $U$, dann enthält $S(U)$
    das Element $g_i^q$, wobei $q=m_j/e_{ij}$ und $m_j$ die relative Ordnung von $a_j$ ist.
  \end{enumerate}
  \item Die folgenden Operationen sind \emph{elementare Operationen} auf $U$ (mit $i,j
  \in \{1 \dd s\}$):
  \begin{enumerate}
    \item Vertauschen von $g_i$ und $g_j$ für $i \neq j$.
    \item Ersetzen von $g_i$ durch $g_1\1$.
    \item Ersetzen von $g_i$ durch $g_i g_j^e$ für $i \neq j$ und $e \in \Z$.
    \item Hinzufügen eines beliebigen Elementes aus $\erz{g_1 \dd g_s}$ als
    $g_{s+1}$.
    \item Entfernen von $g_s$, wenn $g_s=1$ ist.
  \end{enumerate}
  \item Eine Folge $V=(h_1 \dd h_t)$ ist \emph{äquivalent} zu $U$, wenn $V$ durch
  endlich viele elementare Operationen in $U$ umgeformt werden kann.
\end{enumerate}
\end{defi}

\begin{lemma}
Ist $G$ eine polyzyklische Gruppe und $U$ eine eine endliche Folge von
Elementen aus $G$ in Standardform, so ist $S(U)$ genau dann eine Gruppe, wenn
$U$ voll ist. In diesem Fall ist $U$ eine polyzyklische Erzeugerfolge von
$S(U)$.
\end{lemma}

\begin{beweis}
Siehe \cite{sim94}, Seite 409.
\end{beweis}

\begin{lemma} \label{vollegruppe}
Es sei $A=\{a_1 \dd a_n\}$ ein polyzyklisches Erzeugendensystem
einer endlichen $p$"=Gruppe $G$ und $U=(g_1 \dd g_s)$ eine
endliche Folge aus Elementen von $G$ in Standardform.
\begin{enumerate}
  \item Alle Eckeinträge der Exponentenmatrix von $U$ sind 1.
  \item Ist $H$ eine elementarabelsche Untergruppe von $G$ und sind $g_1 \dd
  g_s$ Elemente von $H$, so ist $S(U)$ eine Gruppe.
\end{enumerate}
\end{lemma}

\begin{beweis}
Da $G$ eine $p$"=Gruppe ist, hat jedes Element aus $A$ die relative Ordnung
$p$, und da $p$ eine Primzahl ist, ist 1 der einzige Teiler von $p$, der als
Eckeintrag in der Exponentenmatrix von $U$ auftritt (denn im Fall von $p$ als
Eckeintrag wäre der zu diesem Eckeintrag gehörende Exponentenvektor kein
Exponentenvektor eines Wortes in Normalform).

Da $H$ abelsch ist, enthält $S(U)$ für jedes $i$ und $j$ mit $1 \leq i < j \leq
s$ das Element $g_i\1 g_j g_i=g_j$. Da jedes Element von $A$ die relative
Ordnung $p$ hat und nach der ersten Teilaussage alle Eckeinträge $e_{ij}$ der
zu $U$ gehörenden Exponentenmatrix 1 sind, ist $p$ für alle $j$ mit $1 \leq j
\leq n$ die relative Ordnung von $a_j$, geteilt durch $e_{ij}$. Daher, und weil
$H$ elementarabelsch ist, ist $g_i^p=1$ für alle $i$ mit $1 \leq i \leq s$ ein
Element von $H$. Also ist die Exponentenmatrix von $U$ voll und $S(U)$ nach
\ref{vollegruppe} eine Gruppe.
\end{beweis}

\begin{satz} \label{faktor}
Es sei $A=(a_1 \dd a_n)$ eine polyzyklische Erzeugerfolge einer endlichen Gruppe $G$ und
$\erz{A \mid S}$ eine polyzyklische Präsentation von $G$ und $U=(g_1 \dd g_s)$ eine Folge von
Elementen aus $G$, sodass $U$ voll und $N=S(U)$ ein Normalteiler von $G$ ist. Weiterhin sei
für alle $i$ mit $1 \leq i \leq n$ die Bedingung $b_i = a_i N$ erfüllt und es sei $B=\{b_1 \dd
b_n\}$.
\begin{enumerate}
  \item Die Folge $(b_1 \dd b_n)$ ist eine polyzyklische Erzeugerfolge von
  $G/N$.
  \item Es sei $T=S$. Man forme $T$ in folgender Weise um:
  \begin{enumerate}
    \item Ist $u=v$ eine Kommutatorrelation von $T$, so ersetze in $u$ und in
    $v$ für alle $i$ mit $1 \leq i \leq n$ jedes Vorkommnis von $a_i$ durch
    $b_i$.
    \item Für alle $i$ mit $1 \leq i \leq n$ gilt: Gibt es ein $j$ mit $1 \leq
    j \leq s$, sodass $g_j=a_i^{e_i} \mal a_n^{e_n}$ mit $e_i > 0$ ist, so hat
    $b_i$ die relative Ordnung $e_i$ und es gilt $b_i^{e_i}=b_n^{-e_n} \mal
    b_{i+1}^{-e_{i+1}}$. Falls es eine Potenzrelation $r$ zu $a_i$ gibt, so
    ersetze $r$ in $T$ durch $b_i^{e_i}=[b_n^{-e_n} \mal b_{i+1}^{-e_{i+1}}]$,
    wobei $[b_n^{-e_n} \mal b_{i+1}^{-e_{i+1}}]$ das Wort in Normalform von
    $b_n^{-e_n} \mal b_{i+1}^{-e_{i+1}}$ ist. Gibt es hingegen kein $j$ mit $1 \leq
    j \leq s$, sodass $g_j=a_i^{e_i} \mal a_n^{e_n}$ mit $e_i > 0$ ist, und
    gibt es eine Potenzrelation $r$ zu $a_i$, so ersetze in $r$ jedes Vorkommnis von $a_i$ durch
    $b_i$.
  \end{enumerate}
\end{enumerate}
Dann ist $(B,T)$ eine polyzyklische Präsentation von $G/N$.
\end{satz}

\begin{beweis}
Siehe \cite{sim94}, Seite 414.
\end{beweis}


%

\chapter{Algorithmisches Erzeugen von $p$"=Gruppen für fest gewähltes $p$}

In diesem Kapitel wird der theoretische Hintergrund des Algorithmus dargestellt, mit dem man
von endlichen $p$"=Gruppen niedrigerer Ordnung aus solche höherer Ordnung konstruieren kann
und dadurch eine Liste von Gruppen erhält, in der jede Isomorphieklasse zur Gruppenordnung $p,
p^2 \dd p^n$ genau einmal auftritt, wobei $p^n$ eine beliebig vorgegebenen Schranke ist.

Das Verfahren geht auf eine Idee von Higman zurück und ist hauptsächlich von Havas, New\-man
und O'Brien seit Mitte der siebziger Jahre ausgearbeitet worden (siehe \cite{hav77},
\cite{obr90} und \cite{obr94}). Am Anfang stand nicht die Erzeugung von $p$"=Gruppen selbst im
Vordergrund, sondern die Lösung einiger Probleme der Burnside"=Fragen. Erst mit den Aufsätzen
\cite{obr90} und \cite{obr94} von O'Brien wird die algorithmische Konstruktion endlicher
$p$"=Gruppen zum selbständigen Interessenschwerpunkt und wird seitdem in der algorithmischen
Gruppentheorie zur Klassifikation endlicher $p$"=Gruppen eingesetzt. Einen kleinen Überblick
über diese Entwicklung bieten beispielsweise die Aufsätze \cite{bes02} und \cite{eic99}.

Der Algorithmus ist bereits im ersten Kapitel informell skizziert worden. An dieser Stelle
wird es ausführlich dargestellt. Die Darstellung fasst die Ergebnisse der Aufsätze
\cite{hav77}, \cite{obr90} und \cite{obr94} in einer einheitlichen Notation zusammen und
konzentriert sich dabei auf die Erzeugung von $p$"=Gruppen. Die ursprüngliche Motivation über
die Burnside"=Fragen wird nicht weiter verfolgt. Die Beweise oder zumindest die Beweisideen
der zentralen Sätzen stammen größtenteils aus diesen Arbeiten und sind im Folgenden mit ihren
Quellen gekennzeichnet.

\section{Nachfolger, $p$"=Cover, Nukleus und Multiplikator}

In diesem Abschnitt werden die Schlüsselbegriffe des algorithmischen Verfahrens definiert und
einige der zentralen Sätze zitiert.

\begin{defi}
Es sei $G$ eine endliche $p$"=Gruppe mit der $p$"=Klasse $c$ und mit $d$ als
minimaler Anzahl von Erzeugern.
\begin{enumerate}
  \item Eine endliche $p$"=Gruppe $H$ ist ein \emph{Nachfolger} von $G$, wenn die
  Faktorgruppe $H/\gamma_c(H)$ isomorph zu $G$ ist und $H$ ebenfalls $d$ als minimale Anzahl
  von Erzeugern hat.
  \item Eine endliche $p$"=Gruppe $H$ ist ein \emph{unmittelbarer Nachfolger} von $G$,
  wenn $H$ ein Nachfolger von $G$ ist und die $p$"=Klasse $c+1$ hat.
  \item Die Gruppe $G$ ist \emph{erweiterbar}, wenn $G$ unmittelbare Nachfolger
  hat, andernfalls ist $G$ \emph{abschließend}.
\end{enumerate}
\end{defi}

\begin{lemma}\label{des1}
Jede endliche $p$"=Gruppe mit $d$ als minimaler Anzahl von Erzeugern ist ein
Nachfolger der elementarabelschen Gruppe $C_p^d$.
\end{lemma}

\begin{beweis}
Die elementarabelsche Gruppe $C_p^d$ ist nach \ref{frat} isomorph zu $G/\gamma_1(G)$ und hat
die $p$"=Klasse 1 sowie $d$ als minimale Anzahl von Erzeugern.
\end{beweis}

\begin{lemma}\label{des2}
Für jede endliche $p$"=Gruppe $G$ ist $G/\gamma_{i+1}(G)$ ein unmittelbarer Nachfolger von
$G/\gamma_i(G)$, sofern $i \in \N$ kleiner ist als die $p$"=Klasse von $G$.
\end{lemma}

\begin{beweis}
Die Gruppe $G/\gamma_{i}(G)$ hat die $p$"=Klasse $i$ und die Gruppe $G/\gamma_{i+1}(G)$ die
$p$"=Klasse $i+1$. Es sei $\theta: G \ra \gamma_{i+1}(G)$ der natürliche Homomorphismus. Nach
\ref{dlem} ist daher
$\gamma_i(G/\gamma_{i+1}(G))=\gamma_i(G^\theta)=\gamma_i(G)^\theta=\gamma_i(G)/\gamma_{i+1}(G)$.
Mit den üblichen Homomorphiesätzen erhält man
$$(G/\gamma_{i+1}(G))/\gamma_i(G/\gamma_{i+1}(G))=(G/\gamma_{i+1}(G))/(\gamma_i(G)/\gamma_{i+1}(G))
\cong G/\gamma_{i+1}(G).$$ Damit ist $G/\gamma_{i+1}(G)$ ein unmittelbarer Nachfolger von
$G/\gamma_i(G)$.
\end{beweis}

\begin{satz} \label{vollst}
Es sei $G$ eine $p$"=Gruppe mit der $p$"=Klasse $c$. Dann gibt es eine endliche Folge von
Gruppen $G_1 \dd G_c$, sodass gilt: Es ist erstens $G_c=G$, zweitens ist $G_1$
elementarabelsch und drittens ist $G_{i+1}$ für alle $i \in \{1 \dd c-1\}$ ein unmittelbarer
Nachfolger von $G_i$.
\end{satz}

\begin{beweis}
Dieser Satz ergibt sich unmittelbar aus \ref{des1} und \ref{des2}.
\end{beweis}

\begin{defi}
Es sei $G$ eine $p$"=Gruppe der $p$"=Klasse $c$ mit $d$ als minimaler Anzahl von Erzeugern,
$F$ eine freie Gruppe vom Rang $d$, $R$ der Kern eines Epimorphismus $\theta$ von $F$ auf $G$
und $R^*=[R,F]R^p$.
\begin{enumerate}
  \item Die Gruppe $P(G)=F/R^*=F/[R,F]R^p$ ist das \emph{$p$"=Cover} von $G$.
  \item Die Gruppe $M(G)=R/R^*=R/[R,F]R^p$ ist der \emph{Multiplikator} von $G$.
  \item Die Gruppe $N(G)=\gamma_c(P(G))$ ist der \emph{Nukleus} von $G$.
\end{enumerate}
Diese Gruppen sind nach \ref{pmniso} bis auf Isomorphie eindeutig bestimmt. Da in dieser
Arbeit nur die Isomorphieklassen der Gruppen von Interesse sind, wird die Bezugnahme auf $F$
und $\theta$ im Weiteren vernachlässigt.
\end{defi}

\begin{satz} \label{expcover}
Es sei $G$ eine endliche $p$"=Gruppe mit $d$ als minimaler Anzahl von Erzeugern, $H$ ebenfalls
eine $p$"=Gruppe mit $d$ als minimaler Anzahl von Erzeugern und $Z \leq H$ eine zentrale,
elementarabelsche Untergruppe von $H$, sodass $G$ isomorph zu $H/Z$ ist. Dann ist $H$ ein
homomorphes Bild von $P(G)$.
\end{satz}

\begin{beweis}
Die Gruppe $G$ habe die $p$"=Klasse $c$. Es sei $F$ eine freie Gruppe vom Rang $d$ und $R$ der
Kern eines Epimorphismus $\theta: F \ra G$. Dann ist $P(G)=F/[F,R]R^p$. Da $H$ eine Gruppe mit
$d$ als minimaler Anzahl von Erzeugern ist, gibt es einen Epimorphismus $\psi: F \ra H$; und
da $H$ eine Faktorgruppe enthält, die isomorph zu $G$ ist, gibt es auch einen Epimorphismus
$\varphi: H \ra G$, wobei $Z$ der Kern von $\vi$ ist. Da $F$ eine freie Gruppe ist, kann man
annehmen, dass $\theta=\psi \circ \varphi$ erfüllt ist. Da $R$ der Kern von $\theta$ ist, wird
$R$ nach den üblichen Homomorphiesätzen von $\psi$ auf eine Untergruppe von $Z$ abgebildet.
Nach \ref{dlem} ist $R^*=[R,F]R^p$ seinerseits eine Untergruppe von $R$ und wird damit
ebenfalls auf eine Untergruppe von $Z$ abgebildet. Da aber $Z$ zentral und elementarabelsch
ist, werden sowohl $R^p$ als auch $[R,F]$ auf die triviale Untergruppe von $Z$ abgebildet.
Damit ist insgesamt $(R^*)^\psi=\{1\} \leq Z$. Also ist $H$ ein homomorphes Bild von $F/R^*$.
(Beweis nach \cite{obr90}, Seite 679)
\end{beweis}

\begin{lemma} \label{pmniso}
Ist $G$ eine endliche $p$"=Gruppe mit der $p$"=Klasse $c$ und mit $d$ als minimaler Anzahl von
Erzeugern und $R$ der Kern eines Epimorphismus $\theta$ von der freien Gruppe $F$ vom Rang $d$
auf $G$, so ist die Isomorphieklasse der Gruppe $P(G)=F/([F,R]R^p)$ unabhängig von der Wahl
des Epimorphismus $\theta$.
\end{lemma}

\begin{beweis}
Die Gruppen $G_1$ und $G_2$ seien isomorph, die Abbildungen $\theta_1: F \ra G_1$ und
$\theta_2: F \ra G_2$ Epimorphismen und $R_1=Kern(\theta_1)$ sowie $R_2=Kern(\theta_2)$. Es
seien $P(G_1)=F/[R_1,F]R_1^p$ und $P(G_2)=F/[R_2,F]R_2^p$. Dann sind $P(G_1)$ und $P(G_2)$
isomorph, da beide Gruppen nach \ref{expcover} homomorphe Bilder voneinander sind. (Beweis
nach \cite{obr90}, Seite 679)
\end{beweis}

\begin{lemma} \label{nukinmult}
Der Nukleus $N(G)$ des $p$"=Covers $P(G)$ einer endlichen $p$"=Gruppe $G$ ist eine Untergruppe
des Multiplikators $M(G)$.
\end{lemma}

\begin{beweis}
Die Gruppe $G$ habe die $p$"=Klasse $c$. Analog zum Beweis von \ref{expcover} erhält man, dass
$R/R^*$ ein Normalteiler von $F/R^*$ ist. Nach den Isomorphiesätzen ist $G \cong F/R \cong
(F/R^*)/(R/R^*)$. Also hat $(F/R^*)/(R/R^*)$ die $p$"=Klasse $c$; und nach \ref{dlem} ist
daher $\gamma_c(F/R^*)=\gamma_c(P(G))=N(G)$ eine Untergruppe von $F/R^*=M(G)$.
\end{beweis}

\begin{lemma}
Es ist $G \cong P(G)/M(G)$ und $M(G)$ ist zentral und elementarabelsch.
Außerdem hat $P(G)$ höchstens die $p$"=Klasse $c+1$.
\end{lemma}

\begin{beweis}
Aus dem Beweis zu \ref{nukinmult} ergibt sich, dass $G$ isomorph zu $P(G)/M(G)$ ist, und mit
\ref{expcover} erhält man, dass $M(G)$ zentral und elementarabelsch ist. Da
$N(G)=\gamma_c(P(G))$ nach \ref{nukinmult} eine Untergruppe von $M(G)$ ist, wäre $M(G)$ im
Widerspruch zur ersten Teilaussage dieses Lemmas nicht elementarabelsch, wenn $P(G)$ eine
höhere $p$"=Klasse als $c+1$ hätte.
\end{beweis}

\begin{folg} \label{pfort}
Die Gruppe $G$ ist genau dann erweiterbar, wenn $P(G)$ die $p$"=Klasse $c+1$
hat.
\end{folg}

\section{Die Berechnung des $p$"=Covers}

\begin{defi} \label{cov}
Es sei $\skp{A}{S}$ mit $A=\{a_1 \dd a_n\}$ eine mit $\omega$ und der Definitionenmenge $D
\subseteq S$ gewichtete konsistente Potenz"=Kommutator"=Präsentation zur Primzahl $p$.
Weiterhin sei $\omega(a_i) > 1$ für $i \in \{d+1 \dd n\}$. Dann ist $\skp{A'}{S'}$ mit
$A'=\{a_1 \dd a_{d+n(n+1)/2}\}$ die \emph{$p$"=Cover"=Präsentation} von $\skp{A}{S}$ (bzw. von
der durch $\skp{A}{S}$ gegebenen Gruppe), sofern $S'$ genau aus den folgenden Relationen
besteht:
\begin{enumerate}
  \item $S'$ enthält die $n-d$ Definitionen zu $\{a_{d+1} \dd
  a_n\}$ aus $D \subseteq S$.
  \item Für $k \in \{n+1 \dd d+n(n+1)/2\}$ enthält $S'$ genau $d+n(n-1)/2$
  Relationen der Art $u_k=v_k a_k$, wobei $u_k=v_k$ eine Relation
  in $S$ ist, die keine Definition von $\skp{A}{S}$ ist. Diese Relationen heißen
  \emph{Kernrelationen}.
  \item $S'$ enthält triviale Relationen, durch die $a_{n+1} \dd
  a_{d+n(n+1)/2}$ zentral und von der Ordnung $p$ sind.
\end{enumerate}
Die Relationen der ersten beiden Arten sind Definitionen von $\skp{A'}{S'}$.
\end{defi}

\begin{satz} \label{newman}
Hat $G$ die gewichtete Potenz"=Kommutator"=Präsentation $\skp{A}{S}$ und ist $\skp{A'}{S'}$
die $p$"=Cover"=Präsentation  von $\skp{A}{S}$, so ist $\skp{A'}{S'}$ eine
Potenz"=Kommutator"=Präsentation von $P(G)$, sofern $\skp{A'}{S'}$ konsistent ist. Ist
hingegen $\skp{A'}{S'}$ nicht konsistent, so lässt sich $\skp{A'}{S'}$ mit dem
Knuth"=Bendix"=Verfahren in eine konsistente Potenz"=Kommutator"=Präsentation $\skp{A''}{S''}$
überführen, sodass $\skp{A''}{S''}$ eine Potenz"=Kommutator"=Präsentation von $P(G)$ ist.
\end{satz}

\begin{beweis}
Siehe \cite{hav77}, Seite 220. Das Knuth"=Bendix"=Verfahren wird beispielsweise in
\cite{sim94}, Seite 43 bis 95, dargestellt.
\end{beweis}

\begin{lemma}
Hat $G$ eine konsistente $p$"=Cover"=Präsentation $\skp{A'}{S'}$ mit der Gewichtsklasse $c$,
so wird $M(G)$ von den Erzeugern von $P(G)$ erzeugt, die über Kernrelationen definiert sind,
und $N(G)$ von den Erzeugern, die mit $c$ gewichtet sind.
\end{lemma}

\begin{beweis}
Siehe \cite{hav77}, Seite 220 und \ref{gewpzentral}.
\end{beweis}

\section{Zulässige Untergruppen und Isomorphieklassen unmittelbarer Nachfolger}

Mit dem vorangegangenen Abschnitt ist deutlich geworden, dass alle unmittelbaren Nachfolger
einer $p$"=Gruppe $G$ zu einer Faktorgruppe von $P(G)$ isomorph sind (falls $G$ überhaupt
unmittelbare Nachfolger hat). Nun werden die folgenden Fragen beantwortet:
\begin{enumerate}
  \item Welche Faktorgruppen von $P(G)$ stellen unmittelbare Nachfolger von $G$ dar?
  \item Wenn $P(G)/U_1$ und $P(G)/U_2$ unmittelbare Nachfolger von $G$ sind, lässt sich dann
  bereits an $U_1$ und $U_2$ erkennen, ob $P(G)/U_1$ und $P(G)/U_2$ isomorph zueinander sind?
\end{enumerate}
Im allgemeinen ist es eine schwierige Aufgabe, zu entscheiden ob zwei Gruppen isomorph zu
einander sind bzw. ob zwei endliche Präsentationen Gruppen derselben Isomorphieklasse
darstellen. Daher ist eines der Ergebnisse für diese Klassifikationsaufgabe besonders
hilfreich: Nur Untergruppen $U$ einer gewisse Teilmenge $Z$ der Untergruppen von $P(G)$
stellen als Faktorgruppen $P(G)/U$ unmittelbare Nachfolger von $G$ dar. Diese Untergruppen
werden zulässige Untergruppen genannt. Auf $Z$ operiert die Automorphismengruppe von $G$,
sodass $Z$ in Bahnen zerfällt. Die Bahnen dieser Operation entsprechen den Isomorphieklassen
unmittelbarer Nachfolger. Damit wird das Problem der Isomorphieprüfung auf das Problem der
Bahnenberechnung zurückgeführt.

\begin{defi}
Eine Untergruppe $U$ von $M(G)$ ist \emph{zulässig}, wenn
$P(G)/U$ ein unmittelbarer Nachfolger von $G$ ist, andernfalls ist
$U$ \emph{unzulässig}.
\end{defi}

\begin{defi} \label{zulgl}
Die zulässigen Untergruppen $U_1$ und $U_2$ sind \emph{äquivalent}, wenn $P(G)/U_1$ isomorph
zu $P(G)/U_2$ ist.
\end{defi}

\begin{satz} \label{suppl}
Eine Untergruppe $U$ ist genau dann zulässig, wenn $U$ eine echte Untergruppe von $M(G)$ und
ein Supplement zu $N(G)$ ist.
\end{satz}

\begin{beweis}
Die Gruppe $U$ sei eine zulässige Untergruppe. Es sei $R$ der Kern eines Epimorphismus
$\theta$ von $P(G)$ auf $G$ und $R^*=[R,F]R^p$. Da $M(G)=R/R^*$ und $U$ eine Untergruppe von
$M(G)$ ist, gibt es eine Untergruppe $M$ von $R$, sodass $U=M/R^*$ ist. Die Gruppe $G \cong
F/R$ habe die $p$"=Klasse $c$. Da $H=P(G)/U=P(G)/(M/R^*)$ ein unmittelbarer Nachfolger von $G$
ist, hat $H$ die $p$"=Klasse $c+1$. Daher ist $M$ eine echte Untergruppe von $R$, bzw.
$U=M/R^*$ eine echte Untergruppe von $M(G)=R/R^*$ (denn im Falle $M=R$ wäre $c=c+1$).

Nun ist noch zu zeigen, dass $U=M/R^*$ ein Supplement zu $N(G)=\gamma_c(P(G))$ ist. Die Gruppe
$R$ ist ein Normalteiler von $F$ und $F/R \cong G$ hat die $p$"=Klasse $c$. Daher ist nach
\ref{dlem} $\gamma_c(F)$ eine Untergruppe von $R$. Da $M$ und $\gamma_c(F)$ Untergruppen von
$R$ und sogar Normalteiler sind, ist $M \gamma_c(F)$ ebenfalls eine Untergruppe von $R$. Da $G
\cong H/\gamma_c(H)$ ist, ist nach demselben Argument wie im Beweis zu \ref{expcover}
$R^\theta$ eine Untergruppe von $\gamma_c(H)$. Andererseits ist nach \ref{dlem} auch
$\gamma_c(F^\theta)$ eine Untergruppe von $R^\theta$. Denn $F/R$ hat die $p$"=Klasse $c$ und
$R$ ist ein Normalteiler von $F$. Also ist nach \ref{dlem} auch $\gamma_c(F)$ eine Untergruppe
von $R$. Aus $\gamma_c(F)^\theta \leq R^\theta$ ergibt sich über
$\gamma_c(F)^\theta=\gamma_c(F^\theta)=\gamma_c(H)$ schließlich $\gamma_c(H) \leq R$.
Insgesamt ist also $R^\theta=\gamma_c(H)$. Außerdem ist $R^\theta=R/M$ und
$\gamma_c(H)=(\gamma_c(F)M)/M$. Denn einerseits ist $G \cong
H/\gamma_c(H)=H/R^\theta=((F/R^*)/(M/R^*))/R^\theta=(F/M)/R^\theta=F/R$ und andererseits ist
$\gamma_c(H)=\gamma_c(F/M)=\gamma_c(F^\theta)=\gamma_c(F)^\theta=(\gamma_c(F)M)/M$. Damit
erhält man insgesamt $R/M=(\gamma_c(F)M)/M$ und darüber $R=\gamma_c(F)M$. Wenn $R^*$
ausfaktorisiert wird, ergibt sich $(M/R^*)\gamma_c(F)R^*/R^*=R/R^*$ und damit, dass $U=M/R^*$
ein Supplement zu $U$ ist.

Nun sei umgekehrt $U=M/R^*$ eine echte Untergruppe von $M(G)$, sodass $U$ ein Supplement zum
Nukleus $N(G)=\gamma_c(F/R^*)$ ist. Da $U$ ein Supplement zu $N(G)$ in $M(G)$ ist, ist
$R/R^*=(\gamma_c(F)M)/R^*$ und daher auch $(\gamma_c(F)M)/M=R/M$. Daraus ergibt sich nach
\ref{dlem} $\gamma_c(F/M)=R/M$, wobei der natürliche Homomorphismus nach $M$ betrachtet wird.
Da $F/M$ eine Faktorgruppe von $F/R^*$ ist, hat $F/M$ ebenfalls $d$ als minimale Anzahl von
Erzeugern und die Faktorgruppe $(F/M)/\gamma_c(F/M)=(F/M)/(R/M)$ ist nach den üblichen
Homomorphiesätzen isomorph zu $F/R \cong G$. Damit ist $F/M$ ein Nachfolger von $G$.

Nun ist noch zu zeigen, dass $F/M$ auch ein unmittelbarer Nachfolger ist. Da $U=M/R^*$ eine
echte Untergruppe von $M(G)=R/R^*$ ist, ist die Gruppe $\gamma_c(F/M)=R/M$ nicht trivial und
$F/M$ hat eine höhere $p$"=Klasse als $c$. Da $F/M$ aber isomorph zu einer Faktorgruppe von
$F/R^*$ ist und $F/R^*$ selbst die $p$"=Klasse $c+1$ hat, hat auch $F/M$ die $p$"=Klasse
$c+1$. Also ist $F/M$ ein unmittelbarer Nachfolger von $G$. (Beweis nach \cite{obr90}, Seite
680)
\end{beweis}

\begin{defi} \label{erwaut}
Es sei $G$ eine $p$"=Gruppe der $p$"=Klasse $c$ mit $d$ als minimaler Anzahl von Erzeugern,
$F$ eine freie Gruppe vom Rang $d$ mit den Erzeugern $a_1 \dd a_d$ und $R$ der Kern eines
Epimorphismus von $F$ auf $G$. Weiterhin sei $\alpha$ ein Automorphismus von $F/R$. Eine
Abbildung $\alpha^*: P(G) \ra P(G)$ heißt \emph{Erweiterungsautomorphismus} von $\alpha$, wenn
$\alpha^*$ in folgender Weise definiert ist: Ist zu jedem $i$ mit $1 \leq i \leq d$ ein Wort
$u_i$ in $F$ ausgewählt, sodass die Bedingung $(a_i R)^\alpha=u_i R$ erfüllt ist, dann ist für
jedes Wort $w(a_1 \dd a_d) \in F$ die Abbildung $\alpha^*$ erklärt durch $(w(a_1 \dd
a_d)R^*)^{\alpha^*}=w(u_1 \dd u_d)R^*$, wobei $w(u_1 \dd u_d)$ das Wort ist, das dadurch
entsteht, dass man jedes Vorkommnis von $a_i$ für $i$ mit $1 \leq i \leq d$ in $w(a_1 \dd
a_d)$ durch $u_i$ ersetzt.
\end{defi}

\begin{lemma} \label{erwautaut}
Jeder Erweiterungsautomorphismus ist ein Automorphismus.
\end{lemma}

\begin{beweis}
Zuerst wird gezeigt, dass $\alpha^*$ eine wohldefinierte Abbildung ist, indem nachgewiesen
wird, dass erstens $R^*$ unter $\alpha^*$ invariant ist und zweitens dass jede Nebenklassen
von $R^*$ in $F$ unabhängig von der Wahl ihres Repräsentanten stets auf dieselbe Nebenklasse
abgebildet wird. Für den ersten Schritt sei $w(a_1 \dd a_d)$ ein Element aus $R$. Die
Nebenklasse $R=w(a_1 \dd a_d)R$ wird durch $\alpha$ auf $R^\alpha=w(u_1 \dd u_d)R$ abgebildet.
Da $\alpha$ ein Automorphismus von $F/R$ ist, ist $R$ unter $\alpha$ invariant und daher
$R^\alpha=R=w(u_1 \dd u_d)R$, d.\,h. es ist $w(u_1 \dd u_d) \in R$. Ist sogar $w(a_1 \dd a_d)
\in R^*$, dann ist auch $w(u_1 \dd u_d) \in R^*$, denn $R^*$ ist nach \ref{dlem} eine
charakteristische Untergruppe von $R$ und damit unter jedem Automorphismus invariant.

Nun wird in einem zweiten Schritt gezeigt, dass die Nebenklassen von $R^*$
unabhängig von der Wahl der Repräsentanten abgebildet werden. Es seien $w_1(a_1
\dd a_d)$ und $w_2(a_1 \dd a_d)$ Worte in $F$, sodass $w_1(a_1 \dd
a_d)R^*=w_2(a_1 \dd a_d)R^*$ erfüllt ist. Dann ist nach dem Identitätskriterium
für Nebenklassen $w_2(a_1 \dd a_d)\1w_1(a_1 \dd a_d) \in R^*$. Analog zur
vorangegangenen Überlegung erhält man, dass $w_2(u_1 \dd u_d)\1w_1(u_1 \dd u_d)
\in R^*$ ist bzw. dass die beiden Nebenklassen $w_1(u_1 \dd u_d)R^*$ und
$w_2(u_1 \dd u_d)R^*$ identisch sind. Damit ist $\alpha^*$ wohldefiniert.

Da $\alpha^*$ über eine Wortersetzung definiert ist, ist $\alpha^*$ ein injektiver
Homomorphismus. Nun ist noch zu zeigen, dass $\alpha^*$ surjektiv ist: Da $\alpha$ ein
Automorphismus von $F/R$ und $(a_iR)^\alpha=u_iR$ für alle $i$ mit $1 \leq i \leq d$ ist, wird
$F/R$ nicht nur von $\{a_1 R \dd a_d R\}$, sondern auch von $\{u_1 R \dd u_d R\}$ erzeugt.
Daher ist $F/R^*=\langle u_1 R^* \dd u_d R^*, R/R^* \rangle$. Da aber $R/R^*$ eine Untergruppe
der Frattinigruppe $\gamma_1(F/R^*)$ von $F/R^*$ ist, ist $R/R^*$ in jedem Erzeugendensystem
überflüssig, d.\,h. $F/R^*$ wird bereits von $\{u_1 R^* \dd u_d R^*\}=\{(a_1 R^*)^{\alpha^*}
\dd (a_d R^*)^{\alpha^*}\}$ erzeugt und $\alpha^*$ ist daher surjektiv. (Beweis nach
\cite{obr90}, Seite 681)
\end{beweis}

\begin{lemma}
Sind $\alpha_1^*$ und $\alpha_2^*$ Erweiterungsautomorphismen von $\alpha \in
Aut(G)$, so sind die Einschränkungen von $\alpha_1^*$ und $\alpha_2^*$ auf
$M(G)$ identisch.
\end{lemma}

\begin{beweis}
Zu den Erzeugern $a_1 R \dd a_d R$ der Gruppe $F/R$ seien die beiden
Erweiterungsautomorphismen $\alpha^*_1$ und $\alpha^*_2$ so definiert, dass $(w(a_1 \dd
a_d)R^*)^{\alpha^*_1}=w(u_1 \dd u_d)R^*$ und $(w(a_1 \dd a_d)R^*)^{\alpha^*_2}=w(v_1 \dd
v_d)R^*$ für jedes Wort $w(a_1 \dd a_d)$ aus $F$ gelte. Man nehme an, dass für alle $i$ mit $1
\leq i \leq d$ die Wörter $u_i$ und $v_i$ verschieden seien, d.\,h. dass es ein nichttriviales
$r_i \in R$ gebe, sodass $v_i = u_i r_i$ erfüllt sei. Da nach dem Beweis zu \ref{erwautaut}
$R$ und $R^*$ unter $\alpha^*_1$ wie auch unter $\alpha^*_2$ invariant sind, ist
$(R/R^*)^{\alpha^*_1}=(R/R^*)^{\alpha^*_2}=R/R^*=M(G)$. Ist also $w(a_1 \dd a_d)$ ein Wort in
$R$, so sind auch $w(u_1 \dd u_d)$ und $w(v_1 \dd v_d)$ Wörter in $R$, d.\,h. die
Einschränkungen von $\alpha^*_1$ und $\alpha^*_2$ sind Automorphismen von $M(G)$. Nun ist noch
zu zeigen, dass in diesem Fall die Nebenklassen $w(u_1 \dd u_d)R^*$ und $w(v_1 \dd v_d)R^*$
identisch sind. Ist $w(a_1 \dd a_d)$ ein Wort in $R$, so besteht $w(a_1 \dd a_d)$ aus
Kommutatoren und $p$"=ten Potenzen über $a_1 \dd a_d$ und ebenso bestehen dann $w(u_1 \dd
u_d)$ und $w(v_1 \dd v_d)$ aus Kommutatoren und $p$"=ten Potenzen über $u_1 \dd u_d$ bzw. $v_1
\dd v_d$. Für alle $i,j$ mit $1 \leq i,j \leq d$ gilt daher $[v_j, v_i] R^*=[u_j r_j, u_i r_i]
R^* = [u_j, u_i] R^*$ und $v_i R^* = (u_i r_i)^p R^* = u_i^p r_i^p R^*= u_i R^*$, da $R^*$ in
$R$ zentral und $R/R^*$ zu $p$ elementarabelsch ist. Durch Induktion über die Wortlänge ergibt
sich damit, dass $w(u_1 \dd u_d) R^* = w(v_1 \dd v_d) R^*$ ist, womit die Identität der
Einschränkungen von $\al_1^*$ und $\al_2^*$ auf $M(G)$ gezeigt ist. (Beweis nach \cite{obr90},
Seite 681)
\end{beweis}

\begin{satz} \label{exterwaut}
Sind $M_1/R^*$ und $M_2/R^*$ Untergruppen von $P(G)$, die in $M(G)$ enthalten
sind, und gibt es einen Isomorphismus $\theta$ von $F/M_1$ auf $F/M_2$, dann
existiert ein von $\theta$ induzierter Automorphismus $\alpha \in Aut(G)$,
sodass der Erweiterungsautomorphismus $\alpha^*$ von $\alpha$ die Gruppe
$M_1/R^*$ auf $M_2/R^*$ abbildet.
\end{satz}

\begin{beweis}
Zuerst wird gezeigt, dass ein $\alpha \in Aut(G)$ von $\theta$ induziert wird: Zu jedem $i$
mit $1 \leq i \leq d$ wähle man ein Wort $b_i$ aus $F$, sodass $(a_i M_1)^\theta= b_i M_2$
erfüllt ist. Nach \ref{dlem} ist dann $(R/M_1)^\theta = \gamma_c(F/M_1)^\theta =
\gamma_c((F/M_1)^\theta)=\gamma_c(F/M_2)=R/M_2$. Nach den Homomorphiesätzen und den
Voraussetzungen über $\theta$ ist $G \cong F/R \cong (F/M_1)/(R/M_1) \cong (F/M_2)/(R/M_2)$
und damit induziert $\theta$ einen Automorphismus $\alpha$ von $G$.

Nun ist noch zu zeigen, dass der Erweiterungsautomorphismus $\alpha^*$ die Gruppe $M_1/R^*$
auf Gruppe $M_2/R^*$ abbildet: Es sei $w(a_1 \dd a_d)$ ein Wort in $M_1$ und $\hat{\alpha}^*$
die Einschränkung von $\alpha^*$ auf $M(G)$. Dann ist $(w(a_1 \dd
a_d)R^*)^{\hat{\alpha}^*}=w(b_1 \dd b_d)R^*$. Für die Behauptung reicht der Nachweis, dass
$w(b_1 \dd b_d)$ ein Wort aus $M_2$ ist: Es ist $$w(b_1 \dd b_d)M_2 = w(b_1 M2 \dd b_d M_2) =
w((a_1 M_1)^\theta \dd (a_d M_1)^\theta)$$ $$=(w(a_1 \dd a_d) M_1)^\theta=M_1^\theta=M_2.$$
Damit ist $w(b_1 \dd b_d)$ ein Wort aus $M_2$ und daher $(M_1/R^*)^{\hat{\alpha}^*}$ eine
Untergruppe von $M_2/R^*$. Da aber beide Gruppen denselben Index in $F/R^*$ haben, ist
$M_1/R^*=M_2/R^*$. (Beweis nach \cite{obr90}, Seite 681)
\end{beweis}

\begin{lemma} \label{indperm}
Jeder Erweiterungsautomorphismus eines Automorphismus von $G$ induziert eine Permutation der
zulässigen Untergruppen in $M(G)$.
\end{lemma}

\begin{beweis}
Automorphismen einer Gruppe induzieren Permutationen der Untergruppen dieser Gruppe. Damit
reicht der Nachweis, dass zulässige Untergruppen wieder auf zulässige abgebildet werden: Die
Untergruppe $M/R^*$ von $M(G)$ sei zulässig. Da der Nukleus $N(G)=\gamma_c(F/R^*)$ als Teil
der absteigenden $p$"=Zentralreihe charakteristisch ist und da $M(G)=R/R^*$ nach dem Beweis
von \ref{erwautaut} von jedem Erweiterungsautomorphismus invariant gelassen wird, gilt für
einen Erweiterungsautomorphismus $\alpha^*$ von $\al \in Aut(G)$, dass $(M/R^*)^{\alpha^*}
\gamma_c(F/R^*)= ((M/R^*) \gamma_c(F/R^*))^{\alpha^*}=R/R^*$ ist. Damit ist gezeigt, dass auch
$(M/R^*)^{\alpha^*}$ ein Supplement des Nukleus und nach \ref{suppl} ebenfalls eine zulässige
Untergruppe von $M(G)$ ist. (Beweis nach \cite{obr90}, Seite 682)
\end{beweis}

\begin{defi}
Für eine endliche $p$"=Gruppe $G$ ist $Aut(G)^*$ die Menge der Erweiterungsautomorphismen und
$Perm(G)^*$ die Menge der Permutationen, die durch $Aut(G)^*$ auf die Menge der zulässigen
Untergruppen von $M(G)$ induziert wird.
\end{defi}

\begin{satz} \label{bahnenzu}
Die Bahnen der zulässigen Untergruppen unter der Operation von $Perm(G)^*$ sind
die Äquivalenzklassen im Sinne von \ref{zulgl}.
\end{satz}

\begin{beweis}
Die Untergruppen $M_1/R^*$ und $M_2/R^*$ von $M(G)$ seien zulässig und äquivalent im Sinne von
\ref{zulgl}, d.\,h. die Gruppen $F/M_1$ und $F/M_2$ sind zueinander isomorph. Nach
\ref{exterwaut} gibt es dann einen Erweiterungsautomorphismus $\alpha^*$, der $M_1/R^*$ auf
$M_2/R^*$ abbildet, und nach \ref{indperm} induziert $\alpha^*$ eine Permutation $\alpha' \in
Perm(G)^*$, sodass $(M_1/R^*)^{\alpha'}=M_2/R^*$ ist. Damit liegen $M_1/R^*$ und $M_2/R^*$ in
derselben Bahn unter der Operation von $Perm(G)^*$.

Umgekehrt seien $M_1/R^*$ und $M_2/R^*$ Untergruppen von $M(G)$, die unter der Operation von
$Perm(G)^*$ in derselben Bahn liegen. Dann gibt es eine Permutation $\alpha' \in Perm(G)^*$,
sodass $(M_1/R^*)^{\alpha'}=M_2/R^*$. Da $\alpha' \in Perm(G)^*$ ist gibt es einen
Erweiterungsautomorphismus $\alpha^*$, sodass $\alpha'$ von $\alpha^*$ induziert wird und
$(M_1/R^*)^{\alpha^*}=M_2/R^*$ ist. Da $\alpha^*$ einem Automorphismus von $F/M_1$ entspricht,
ist $(F/M_1)/(M_1/R^*)$ nach den üblichen Homomorphiesätzen isomorph zu
$((F/M_1)/(M_1/R^*))^{\alpha^*}= (F/M_2)/(M_2/R^*)$. (Beweis nach \cite{obr90}, Seite 682)
\end{beweis}

\begin{folg}
Zerfällt die Menge der Untergruppen von $M(G)$ unter der Operation von $Aut(G)$
über die Erweiterungsautomorphismen von $Aut(G)$ in Bahnen, so enthält jede
dieser Bahnen entweder nur zulässige Untergruppe oder nur unzulässige
Untergruppen.
\end{folg}

\section{Die Berechnung der Automorphismengruppen unmittelbarer Nachfolger}

Mit den beiden vorangegangenen Abschnitten des Kapitels ist geschildert worden, wie sich von
einer elementarabelschen Gruppe ausgehend eine endliche Folge unmittelbarer Nachfolger bis zu
jeder abstrakten $p$"=Gruppen einer gewissen Ordnung ermitteln lässt. Für jeden Schritt dieser
Folge wird die Automorphismengruppen der jeweiligen Gruppe benötigt, von der man unmittelbare
Nachfolger berechnen möchte. Der folgende Satz zeigt, wie man die fraglichen
Automorphismengruppe berechnen kann: Ist $P(G)/U$ ein unmittelbarer Nachfolger von $G$, so
lässt sich die Automorphismengruppe von $P(G)/U$ unmittelbar aus dem Stabilisator $S \leq
Aut(G)$ der zulässigen Untergruppe $U$ ablesen.


\begin{defi}
Ist $G$ eine Gruppe, $\alpha$ ein Automorphismus von $G$ und $N$ ein Normalteiler von $G$,
sodass $N$ unter $\alpha$ invariant ist, und $\theta: G \ra G/N$ der natürliche Homomorphismus
bezüglich $N$, so ist die Abbildung $\bar{\alpha}: G/N \ra G/N$ die \emph{Projektion} von
$\alpha$ auf $G/N$, wenn für alle $g \in G$ die Bedingung
$(g^\theta)^{\bar{\alpha}}=(g^\alpha)^\theta$ erfüllt ist.
\end{defi}

\begin{satz} \label{aut}
Es sei $M/R^*$ eine zulässige Untergruppe von $M(G)$ und $S$ der Stabilisator von $M/R^*$
unter der Operation von $Aut(G)$ auf $M(G)$ über Erweiterungsautomorphismen, d.\,h. es sei
$$S_G=Stab_{Aut(G)}(M/R^*)= \{\alpha \in Aut(G) \mid (M/R^*)^{\alpha^*}=M/R^*\}$$ und
$$S=\left\{\bar{\alpha}^* \mid \alpha \in S_G \right\} \leq Aut(F/M)$$ die Menge der
Projektionen von Erweiterungsautomorphismen von $S_G$ auf $F/M$ und $V$ die Gruppe aller
Automorphismen von $F/M$, deren Projektion auf $F/R \cong G$ die Identität ist. Dann ist
$Aut(F/M)=SV$ und $V$ ist isomorph zu $C_p^k$, wobei $p^k=|M/R^*|$ ist.
\end{satz}

\begin{beweis}
Es sei $\theta$ ein Automorphismus von $F/M$. Für jedes $i$ mit $1 \leq i \leq d$ wähle man
ein Wort $u_i$ in $F$, sodass $(a_i M)^\theta = u_i M$ ist. Da $F/M$ ein unmittelbarer
Nachfolger von $G \cong F/R$ ist, ist $(F/M)/\gamma_c(F/M)$ isomorph zu $F/R$, und da
$\gamma_c(F/M)=R/M$ eine charakteristische Untergruppe von $F/M$ ist, kann $\theta$ auf
$(F/M)/(R/M) \cong F/R \cong G$ eingeschränkt werden und induziert damit einen Automorphismus
$\alpha$ von $G$. Ist $\alpha$ die Identität von $Aut(G)$, so ist $\theta$ ein Element aus
$V$. Ist hingegen $\alpha$ nicht die Identität von $Aut(G)$, so stabilisiert die Erweiterung
$\alpha^*$ die Untergruppe $M/R^*$: Die Erweiterung $\alpha^*$ sei über diese Bedingungen
$(a_i R)^\alpha=u_i R$ für jedes $i$ mit $1 \leq i \leq d$ nach \ref{erwaut} definiert. Es sei
$w(a_1 \dd a_d)$ ein Wort aus $M$. Dann ist $w(u_1 \dd u_d)M= w(u_1 M \dd u_d M)=w((a_1
M)^\theta \dd (a_d M)^\theta)=(w(a_1 \dd a_d)M)^\theta = M^\theta = M$, d.\,h. $w(u_1 \dd
u_d)$ ist ein Wort aus $M$. Daher ist $M/R^*$ sowohl unter $\theta$ als auch unter $\alpha^*$
invariant. Also ist $\bar{\alpha}$, die Projektion von $\alpha^*$ auf $F/M$, ein Element aus
$S$.

Nun ist noch zu zeigen, dass $\theta$ als Produkt eines Elementes aus $S$ und eines Elementes
aus $V$ dargestellt werden kann. Für alle $i$ mit $1 \leq i \leq d$ wähle man $u_i \in F$ und
$r_i \in R$, sodass $(a_i R^*)^{\al^*}=u_i r_i R^*$ erfüllt ist (man kann $r_i$ derart wählen,
da $R$ unter $\al^*$ invariant ist). Dann ist $\bar{\al}$ durch $(a_i M)^{\bar{\al}}=u_i r_i
M$ definiert. Nun sei der Automorphismus $\nu$ von $F/M$ durch $(u_i M)^\nu= u_i r_i M$
definiert. Die Einschränkung von $\nu$ auf $F/R$ ist die Identität. Also ist $\nu$ aus $V$. Es
ist $\bar{\al}=\theta \nu$ und damit $\theta=\bar{\al} \nu\1$ mit $\bar{\al} \in S$ und $\nu\1
\in V$. Da $M/R^*$ elementarabelsch und $\gamma_c(F/M) \cong M/R^*$ zentral in $F/M$ ist, ist
$V$ nach \ref{autpzent} isomorph zu $C_p^k$, wobei $k$ diejenige natürliche Zahl mit
$p^k=|M/R^*|$ ist. (Beweis ähnlich \cite{obr90}, Seite 683)
\end{beweis}

\chapter{Gruppen und Operationen}

\section{Gruppenoperationen}


\begin{defi}
Eine Gruppe $G$ \emph{operiert} (von rechts) auf einer nichtleeren Menge $X$, wenn es eine
Abbildung $\vi: G \times X \ra X: (g,x) \mt x^g$ gibt, sodass für alle $g,h \in G$ und $x \in
X$ sowohl $x^{gh}=(x^g)^h$ als auch $x^1=x$ gilt, wobei $1$ das Neutralelement von $G$ ist.
\end{defi}

\begin{defi}
Operiert eine Gruppe $G$ auf einer Menge $X$ und ist $x \in X$, so ist $x^G=\{x^g \mid g \in
G\}$ die \emph{Bahn} und $Stab_G(x)=\{g \in G \mid x^g=x\}$ der \emph{Stabilisator} von $x$
unter $G$.
\end{defi}

\begin{satz} \label{stab}
Operiert eine Gruppe $G$ auf einer Menge $X$, so ist für alle $x \in X$ der Stabilisator
$Stab_G(x)$ eine Untergruppe von $G$ und es gilt $\left|x^G\right|=[G:Stab_G(x)],$ wobei
$[G:Stab_G(x)]$ der Index von $Stab_G(x)$ in $G$ ist.
\end{satz}

\begin{beweis}
Siehe \cite{mey80}, Seite 67.
\end{beweis}

\begin{lemma} \label{komm1}
Es sei $G$ eine Gruppe und $a,b \in G$. Ist der Kommutator von $a$ und $b$ zentral, so gilt
für alle $n,m,l \in \N$, dass $a^n b^m a^l=a^{n+l} b^m [b,a]^{lm}$ ist.
\end{lemma}

\begin{beweis}
Die Behauptung wird durch vollständige Induktion über $l$ bewiesen. Für $l=1$ gilt $a^n b^m a
= a^n b^{m-1} a b [b,a] = a^{n+1} b^m [b,a]^m$, wie man durch vollständige Induktion über $m$
beweisen kann. Für $l>1$ ist $a^n b^m a^l = a^n b^m a^{l-1} a = a^{n+l-1} b^m [b,a]^{m(l-1)} a
= a^{n+l-1} b^m a [b,a]^{m(l-1)} = a^{n+l} b^m [b,a]^m [b,a]^{m(l-1)}=a^{n+l} b^m [b,a]^{lm}$.
\end{beweis}

\section{Bemerkungen über endliche Körper}

\begin{lemma}\label{kzyklisch}
Ist $K$ ein Körper und $G$ eine endliche Untergruppe der multiplikativen Gruppe von $K$, so
ist $G$ zyklisch.
\end{lemma}

\begin{beweis}
Es sei $n=kgV\{|g| \mid g \in G\}$. Dann gilt $g^n=1$ für alle $g \in G$. Also ist jedes $g
\in G$ eine Nullstelle des Polynoms $f=x^n-1 \in K[x]$. Da $f$ höchstens $n$ Nullstellen hat,
ist $|G| \leq n$. Daher ist nur noch zu zeigen, dass es ein $g \in G$ gibt, sodass $|g|=n$
ist. Denn in diesem Fall ist $G$ zyklisch. Es reicht zu zeigen, dass das Erzeugnis $\langle a,
b \rangle$ zweier Elemente $a,b \in G$ ein Element $c$ enthält, so daß $|c|=kgV(|a|,|b|)$ ist.
Die Behauptung ergibt sich dann über vollständige Induktion. Es sei also $a,b \in G$ und
$|a|=p_1^{a_1} \mal p_r^{a_r} a'$ und $|b|=p_1^{b_1} \mal p_r^{b_r} b'$ mit $ggT(a',b')=1$ und
$a_i \geq b_i$ für alle $i$ mit $1 \leq i \leq s \leq r$ und $a_i < b_i$ andernfalls für
geeignete Primzahlen $p_1 \dd p_r$. Es sei $x=p_{s+1}^{a_{s+1}} \mal p_r^{a_r}$ und
$y=p_1^{b_1} \mal p_s^{b_s}$. Weiterhin seien $c_1=a^x$ und $c_2=b^y$. Dann ist
$|c_1|=p_1^{a_1} \mal p_s^{a_s} a'$ und $|c_2|=p_{s+1}^{b_{s+1}} \mal p_r^{b_r}b'$. Für das
Produkt $c=c_1 c_2$ erhält man also $|c|=|c_1 c_2|=|c_1| \cdot |c_2|=kgV(|a|,|b|)$.
\end{beweis}

\begin{defi}
Ist $w$ ein Erzeuger der multiplikativen Gruppe von $\f$, so ist $$l: \f \ra \Z: x \mt
  \begin{cases}
    \min \{n \in \N_0 \mid w^n=x \} & \text{ für } x \neq 0, \\
    -1 & \text{ für } x=0
  \end{cases}$$
die \emph{Exponentenabbildung} von $\f$ bezüglich $w$.
\end{defi}

\begin{lemma} \label{indexadd}
Ist $m, n \in \Z$, so ist $\erz{\ol{m}}$ ein Ideal in $\Z_n$ vom Index $ggT(m,n)$, wobei
$\ol{m}$ das Bild von $m$ unter dem natürlichen Homomorphismus von $\Z$ auf $\Z_n$ ist.
\end{lemma}

\begin{beweis}
Wie aus der Zahlentheorie bekannt ist, ist die von $\ol{m}$ erzeugte Untergruppe ein Ideal in
$\Z_n$. Man kann annehmen, dass $0 \leq m <n$ ist (andernfalls betrachte man statt $m$ ein $m'
\in \Z$ mit $0 \leq m' < n$ und $m' \equiv m \text{ mod } n$). Dann ergibt sich für die
Ordnung von $\ol{m}$, dass
\begin{eqnarray*}
  |\ol{m}|&=& \min \{k \mid  k \in \N \text{ und } km \equiv n \text{ mod }n\} \\
  &=& \min \{k \mid  k \in \N \text{ und es gibt ein } l \in \N \text{ mit } km=ln\} \\
  &=& \frac{kgV(m,n)}{m}=\frac{\frac{mn}{ggT(m,n)}}{m}=\frac{n}{ggT(m,n)}
\end{eqnarray*}
ist. Da $\erz{\ol{m}}$ zyklisch ist, ist $|\erz{\ol{m}}|=|\ol{m}|$. Nach dem Satz von Lagrange
erhält man daher
$$[\Z_n:\erz{\ol{m}}]=\frac{|\Z_n|}{|\erz{\ol{m}}|}=\frac{n}{\frac{n}{ggT(m,n)}}=ggT(m,n).$$
\end{beweis}

\begin{lemma} \label{indexmult}
Wird die multiplikative Gruppe $\F_p^*$ des endlichen Körpers $\F_p$ von $w$ erzeugt, so ist
$\erz{w^a}$ eine Untergruppe vom Index $ggT(a, p-1)$ in $\F_p^*$. Inbesondere ist die
Abbildung $\theta: \F_p \ra \F_p: x \mt x^a$ genau dann bijektiv, wenn $ggT(a, p-1)=1$ ist,
und $(\F_p^*)^\theta$ ist eine Untergruppe von $\f^*$ vom Index $ggT(a, p-1)$.
\end{lemma}

\begin{beweis}
Die Aussage ergibt sich unmittelbar aus \ref{kzyklisch} und \ref{indexadd}.
\end{beweis}

\begin{defi}
Ist $w$ ein Erzeuger der multiplikativen Gruppe $\f^*$ von $\f$, so ist $\f^+=\{ w^a \mid 0
\leq a < \frac{p-1}{2}\}$ und $\f^-=\{-z \mid z \in \f^+\}$ bezüglich $w$.
\end{defi}

\begin{lemma} \label{pm}
Es ist $\f^-=\{w^b \mid \frac{p-1}{2} \leq b < p-1\}$ und damit $|\f^+|=|\f^-|=\frac{p-1}{2}$
sowie $\f^+ \dot{\cup} \f^- = \f^*$.
\end{lemma}

\begin{beweis}
Es sei $w^a \in \f^+$ und $w^a + w^b=0$. O.\,B.\,d.\,a. sei $a < b$. Dann ist
$w^a(1+w^{b-a})=0$. Da $w^a \neq 0$ ist, ist $1+w^{b-a}=0$ bzw. $w^{b-a}=-1$. Also ist $b-a
\equiv \frac{p-1}{2}$ mod $p-1$, denn $-1$ ist nach \ref{indexmult} das einzige Element $z \in
\f$ mit $z^2=1$ und $z \neq 1$ und diese Bedingung ist nur für $w^{(p-1)/2}$ erfüllt. Da $0
\leq a < \frac{p-1}{2}$ und $a \leq b \leq p-2$ ist, ist $b$ als $b = \frac{p-1}{2} + a$
eindeutig bestimmt.
\end{beweis}

\begin{lemma} \label{aequiv}
Es sei $F=\{(a,b) \in \f^2 \mid a \neq b\}$, sodass $\f^*$ von $w$ erzeugt wird, und $\approx$
die Äquivalenzrelation auf $F^2$, die durch $(a,b) \approx (c,d) \Leftrightarrow((a=c \wedge
b=d) \vee (a=d \wedge b=c))$ definiert ist. Operiert $\f^*$ auf $F$ durch $(\f^* \times F) \ra
F: (z,(a,b)) \mt (za,zb)$, so ist $B=\{(1,w^k) \mid 1 \leq k \leq \frac{p-1}{2}\}$ ein
Vertretersystem der Bahnen unter der Operation von $\f^*$ modulo $\approx$.
\end{lemma}

\begin{beweis}
Es sei $(a,b) \in F$. Da $a$ und $b$ Elemente von $\f^*$ sind, gibt es ein $u$ und ein $v$ aus
$\{0 \dd p-2\}$, sodass $a=w^u$ und $b=w^v$ ist. Es ist dann
$(a,b)^{\f^*}=(w^u,w^v)^{\f^*}=\{(w^{u+q}, w^{v+q}) \mid 0 \leq q \leq p-2\}=\{(w^{q},
w^{v-u+q}) \mid 0 \leq q \leq p-2\}$. Also ist $(1,w^{v-u})$ ein Repräsentant von
$(a,b)^{\f^*}$, d.\,h. $B'=\{(1,w^k) \mid 1 \leq k \leq p-2\}$ ist ein Vertretersystem der
Bahnen in $F$ unter der Operation von $\f^*$. Da $(1, w^k) \approx (w^k,1)$ ist und $(w^k,1)$
in derselben Bahn wie $(1, w^{-k})$ liegt, ergibt sich mit \ref{pm} die Behauptung.
\end{beweis}

\section{Vektorräume und Matrixgruppen über $\f$}

\begin{defi}
Es ist $M(m \times n,p)$ die Menge der $(m \times n)$"=Matrizen über dem Körper $\f$ und
$GL(n,p)=\{m \in M(n \times n, p) \mid det(m) \neq 0\}$ die generelle lineare Gruppe von $M(n
\times n, p)$ sowie $SL(n,p)=\{m \in M(n \times n, p) \mid det(m) = 1 \text{ oder }
det(m)=-1\}$ die spezielle lineare Gruppe von $M(n \times n, p)$. Weiterhin ist $SL^+(n,p)=\{m
\in M(n \times n, p) \mid det(m) = 1\}$.
\end{defi}

\begin{lemma} \label{glslz}
Es ist $$|GL(n,p)|= \prod_{i=0}^{n-1} (p^n-p^i)$$ und $$|SL^+(n,p)|= \frac{|GL(n,p)|}{p-1},$$
und die Menge der Diagonalmatrizen von $GL(n,p)$ ist das Zentrum von $GL(n,p)$.
\end{lemma}

\begin{beweis}
\cite{hup67}, Seite 178.
\end{beweis}

\begin{lemma} \label{einvek}
Operiert eine Untergruppe $A$ von $GL(n,p)$ auf $V=\F_p^n$, ist $Z$ das Zentrum von $GL(n,p)$
und sind $B_1 \dd B_m$ die Bahnen der Elemente von $V$ unter der Operation von $AZ$ sowie $C_1
\dd C_l$ die Bahnen der eindimensionalen Unterräume von $V$ unter Operation von $A$, so ist
$m=l$ und jedem $B_i \in \{B_1 \dd B_m\}$ ist eineindeutig ein $C_j \in \{C_1 \dd C_l\}$
zugeordnet, indem ein Vertreter von $B_i$ einen Vertreter von $C_j$ aufspannt und $B_i$ die
Vereinigung der eindimensionalen Unterräume von $C_j$ abzüglich dem Nullvektor ist.
\end{lemma}

\begin{beweis}
Nach \cite{hup67}, Seite 177, sind die Elemente von $Z$ genau die Elemente von $GL(n,p)$, die
jeden eindimensionalen Unterraum von $V$ auf sich selbst abbildet, und nach \ref{glslz}
operiert $Z$ transitiv auf jedem eindimensionalen Unterraum von $V$. Da $Z$ nach Definition
eine Untergruppe von $AZ$ ist, operiert auch $AZ$ transitiv auf jedem Unterraum von $V$. Daher
gilt für jedes $v$ aus $V$, dass $v^Z=\erz{v}$ und damit dass $\erz{v}^A=(v^Z)^A=v^{AZ}$ ist.
\end{beweis}

\begin{folg} \label{zentin}
Operiert eine Untergruppe $A$ von $GL(n,p)$ auf $V=\F_p^n \backslash\{0\}$ und enthält $A$ das
Zentrum von $GL(n,p)$, so ist die Anzahl der Bahnen der Elemente von $V$ unter $A$ gleich der
Anzahl der Bahnen der eindimensionalen Unterräume von $V$ unter $A$.
\end{folg}


\section{Operationen von $GL(n,p)$}

\begin{lemma} \label{tensor}
Operiert eine Gruppe $S$ auf einer elementarabelschen Gruppe $M$ und ist $U$ eine unter der
Operation von $S$ invariante Untergruppe, dann gibt es eine Bijektion zwischen den Bahnen der
Komplemente zu $U$ in $M$ unter der Operation von $S$ und den Bahnen der Elemente des
Tensorproduktes $M/U \otimes U$ unter der diagonalen Operation von $S$ vermittels $(a \otimes
b)^s = a^{s\1} \otimes b^s$.
\end{lemma}

\begin{beweis}
Die Komplemente der Untergruppe $U$ in $M$ stehen in Bijektion zur Kohomologiegruppe
$Z^1(M/U,U)$. Da $M$ elemetarabelsch ist, ist $Z^1(M/U,U)=Hom(M/U,U)$ und daher kann
$Z^1(M/U,U)$ mit $M/U \otimes U$ identifiziert werden. Den Bahnen der Komplemente von $U$ in
$M$ unter der Operation von $S$ entsprechen den Bahnen der Kozykel, auf denen $S$ durch
$\delta^s(a)=\delta(a^{s\1})^s$ operiert.
\end{beweis}

\begin{lemma} \label{jordan2}
Wird die multiplikative Gruppe des Körpers $\F_p$ von $w$ erzeugt und ist $\F_p(\sw)$ die
Körpererweiterung von $\F_p$ um $\sw$, so zerfällt jedes Polynom aus $\F_p[x]$ von Grad 2 in
$\fsw[x]$ in Linearfaktoren. Ist $l: \F_p \ra \Z$ eine injektive Abbildung, so sind alle
normierten Polynome aus $\F_p[x]$ vom Grad 2 durch die disjunkte Vereinigung der folgenden
Mengen gegeben:
\begin{enumerate}
  \item $R_1=\{ (x-a)^2 \mid a \in \f\}$ mit $|R_1|=p$.
  \item $R_2=\{ (x-a)(x-b) \mid a,b \in \f \text{ und } l(a)<l(b)\}$ mit
  $|R_2|=\frac{p(p-1)}{2}$.
  \item $I=\{ (x-a-b \sw)(x-a+b \sw) \mid a \in \f \text{ und } b \in \f^+\}$
  mit $|I|=\frac{p(p-1)}{2}$.
\end{enumerate}
Dabei ist $R_1 \cup R_2$ genau die Menge der reduziblen und normierten Polynome vom Grad 2 und
$I$ die der irreduziblen und normierten.
\end{lemma}

\begin{beweis}
Dass die Mengen $R_1$ und $R_2$ disjunkt sind und ihre Elemente reduzible Polynome sind, lässt
sich unmittelbar erkennen. Da die Faktorisierung der Polynome nur bis auf Permutationen der
Linearfaktoren eindeutig ist, sind mit $R_2$ bereits alle reduziblen und normierten Polynome
mit zwei verschiedenen Nullstellen gegeben und mit $R_1$ und $R_2$ überhaupt alle reduziblen
und normierten Polynome vom Grad 2. Mit demselben Argument ergibt sich, dass in der Definition
von $I$ die Werte für $b$ auf $\f^+$ eingeschränkt werden können. Jedes Element aus $I$ ist
irreduzibel, da durch $b \neq 0$ beide Nullstellen nicht in $\f$ liegen. Andererseits sind
alle Elemente aus $I$ tatsächlich normierte Polynome aus $\f[x]$, denn für alle $a,b \in \f$
ist $(x-a-b \sw)(x-a+b \sw)= x^2 - (a+b\sw +a -b\sw) + (a+b\sw)(a-b\sw)=x^2 - 2a x + a^2 +
b^2w$ in $\f[x]$. Da die Mächtigkeit von $R_1 \cup R_2 \cup I$ gleich $p + \frac{p(p-1)}{2} +
\frac{p(p-1)}{2}=p^2$ ist, sind mit der disjunkten Vereinigung von $R_1$, $R_2$ und $I$ alle
normierten und irreduziblen Polynome vom Grad 2 in $\f[x]$ gegeben.
\end{beweis}

\begin{satz}
Operiert $GL(2,p)$ durch Konjugation auf der Menge $M(2 \times 2, p)$, so ist durch die
folgende Liste von Matrizen ein Vertretersystem der $p(p-1)+2p$ Konjugiertenklassen angegeben,
wobei $w$ ein Erzeuger der multiplikativen Gruppe von $\f$ und $l: \F_p \ra \Z$ die
Exponentenabbildung zu $w$ ist:
\begin{enumerate}
  \item $\begin{pmatrix}
    a     & 0 \\
    0     & a
  \end{pmatrix}$ mit $a \in \f$ entsprechend $p$ Konjugiertenklassen.
  \item $\begin{pmatrix}
    a     & 1 \\
    0     & a
  \end{pmatrix}$ mit $a \in \f$ entsprechend $p$ Konjugiertenklassen.
  \item $\begin{pmatrix}
    a     & 0 \\
    0     & b
  \end{pmatrix}$ mit $a,b \in \f$ und $l(a) < l(b)$ entsprechend $\frac{p(p-1)}{2}$ Konjugiertenklassen.
  \item $\begin{pmatrix}
    a + b\sw   & 0 \\
    0          & a - b\sw
  \end{pmatrix}$ mit $a \in \f$ und $b \in \f^+$ entsprechend $\frac{p(p-1)}{2}$ Konjugiertenklassen.
\end{enumerate}
\end{satz}

\begin{beweis}
Unter Bezug auf \ref{jordan2} ergibt sich die Behauptung nach \cite{bri85}, Seite 98, 111 und
112, auf denen die Berechnung der Jordanschen Normalform einer Matrix beschrieben ist.
\end{beweis}

\begin{lemma} \label{spiegel}
Operiert $GL(2,p)$ durch Konjugation auf $M(2 \times 2,p)$, so sind für alle $a$ und $b$ aus
$\f$ die Matrizen $$\begin{pmatrix}
 a  & 0 \\
 0  & b
\end{pmatrix} \quad \text{ und } \quad \begin{pmatrix}
 b  & 0 \\
 0  & a
\end{pmatrix}$$
zueinander konjugiert.
\end{lemma}

\begin{beweis}
Die Jordanblockzerlegung einer Matrix ist nur bis auf ihre Reihenfolge der Jordanblöcke
eindeutig bestimmt.
\end{beweis}

\begin{lemma} \label{kro}
Ist $Z$ das Zentrum von $GL(2,p)$ und operiert die Gruppe $G=Z \times GL(2,p)$ vermittels
$\delta: G \times M(2 \times 2,p): ((z,g),m) \mt m^{(z,g)}=z g\1 m g$ auf $M(2 \times 2, p)$,
so ist durch die folgende Liste von Matrizen ein Vertretersystem der insgesamt $p+5$ Bahnen
unter dieser Operation angegeben, wobei $w$ ein Erzeuger der multiplikativen Gruppe von $\f$
ist:
\begin{enumerate}
  \item[1a.] $\begin{pmatrix}
    0     & 0 \\
    0     & 0
  \end{pmatrix}$ entsprechend einer Bahn.
  \item[1b.] $\begin{pmatrix}
    1     & 0 \\
    0     & 1
  \end{pmatrix}$ entsprechend einer Bahn.
  \item[2a.] $\begin{pmatrix}
    0     & 1 \\
    0     & 0
  \end{pmatrix}$ entsprechend einer Bahn.
  \item[2b.] $\begin{pmatrix}
    1     & 1 \\
    0     & 1
  \end{pmatrix}$ entsprechend einer Bahn.
  \item[3a.] $\begin{pmatrix}
    0     & 0 \\
    0     & 1
  \end{pmatrix}$ entsprechend einer Bahn.
  \item[3b.] $\begin{pmatrix}
    1     & 0 \\
    0     & w^k
  \end{pmatrix}$ mit $k \in \{1 \dd \frac{p-1}{2}\}$ entsprechend $\frac{p-1}{2}$ Bahnen.
  \item[4a.] $\begin{pmatrix}
    0     & w \\
    1     & 0
  \end{pmatrix}$ entsprechend einer Bahn.
  \item[4b.] $\begin{pmatrix}
    1     & bw \\
    b     & 1
  \end{pmatrix}$ mit $b \in \f^+$ entsprechend $\frac{p-1}{2}$ Bahnen.
\end{enumerate}
\end{lemma}

\begin{beweis}
Aus \ref{jordan2} ist bekannt, wie $M=M(2 \times 2,p)$ unter der Operation von $GL(2,p)$
vermittels $\theta: GL(2,p) \times M: (g,m) \mt m^g=g\1 m g$ in Bahnen zerfällt. Da jedes $z
\in Z$ mit jedem $g \in GL(2,p)$ kommutiert, ist $m^{(z,g)}=z g\1 m g = g\1 z m g$, d.\,h. die
Bahnen unter der Operation von $G$ auf $M$ sind disjunkte Vereinigungen der Bahnen, die in
$\ref{jordan2}$ durch ein Vertretersystem aufgelistet und in vier Klassen eingeteilt worden
sind. Jede dieser vier Klassen wird nun unter der Operation von $Z$ vermittels $\vi: Z \times
M \ra M: (z,m) \mt mz$ einzeln betrachtet.

Für die Bahnen der ersten Klasse gilt: Ist $E$ die Einheitsmatrix von $M$, so liegen $aE$ und
$bE$ mit $a,b \in \f$ genau dann unter der Operation von $Z$ über $\vi$ in einer Bahn, wenn
entweder $a=b=0$ oder $a,b \in \f^*$ ist. Also gibt es genau zwei Bahnen unter der Operation
von $Z$ und damit auch von $G$ auf $M$, die durch die Matrizen $0 \cdot E$ und $E$
repräsentiert werden.

Für die Bahnen der zweiten Klasse gilt: Ist $$g= \begin{pmatrix}
    l     & b \\
    c     & d
  \end{pmatrix}$$ aus $GL(2,p)$, $z$ aus $Z$ und $m$ eine Matrix der zweiten Klasse, so ist
\begin{eqnarray*}
  z g\1 m g &=& z g\1 \begin{pmatrix}
    a     & 1 \\
    0     & a
  \end{pmatrix} g \\
  &=& z g\1 \left( \begin{pmatrix}
    a     & 0 \\
    0     & a
  \end{pmatrix} + \begin{pmatrix}
    0     & 1 \\
    0     & 0
  \end{pmatrix} \right) g\\
  &=& z g\1 \begin{pmatrix}
    a     & 0 \\
    0     & a
  \end{pmatrix} g + z g\1 \begin{pmatrix}
    0     & 1 \\
    0     & 0
  \end{pmatrix} g \\
  &=& z \begin{pmatrix}
    a     & 0 \\
    0     & a
  \end{pmatrix} g\1 g + z \frac{1}{det(g)} \begin{pmatrix}
    d     & -b \\
    -c    & l
  \end{pmatrix} \begin{pmatrix}
    0     & 1 \\
    0     & 0
  \end{pmatrix} \begin{pmatrix}
    l     & b \\
    c     & d
  \end{pmatrix} \\
  &=& z \begin{pmatrix}
    a     & 0 \\
    0     & a
  \end{pmatrix} + z \frac{1}{det(g)} \begin{pmatrix}
    cd    & d^2 \\
    -c^2  & -cd
  \end{pmatrix} \\
\end{eqnarray*}
Ist $c=0$ und $d=1$, so nimmt $det(g)=l$ unabhängig von $z$ jeden Wert aus $\f^*$ an. Daher
können die beiden angegebenen Repräsentanten analog zur Fallunterscheidung $a=0$ und $a \neq
0$ gewählt werden.

Für die Bahnen der dritten Klasse gilt: Unter der Operation von $Z$ über $\vi$ auf $M$ liegen
die Matrizen $$\begin{pmatrix}
 0  & 0 \\
 0  & b
\end{pmatrix} \quad \text{ und } \quad \begin{pmatrix}
 0  & 0 \\
 0  & 1
\end{pmatrix}$$
für alle $b \in \f^*$ offensichtlich in derselben Bahn. Damit bleiben nur noch die Matrizen
der $$\begin{pmatrix}
 a  & 0 \\
 0  & b
\end{pmatrix}$$
aus der dritten Klasse mit $a,b \in \f^*$ zu betrachten. Nach \ref{spiegel} liegen die
Matrizen $$\begin{pmatrix}
 a  & 0 \\
 0  & b
\end{pmatrix} \quad \text{ und } \quad \begin{pmatrix}
 b  & 0 \\
 0  & a
\end{pmatrix}$$
unter der Operation von $GL(2,p)$ durch Konjugation in derselben Bahn. Nach \ref{aequiv} sind
die Matrizen $$\begin{pmatrix}
 1  & 0 \\
 0  & w^k
\end{pmatrix}$$ mit $k \in \{1 \dd \frac{p-1}{2}\}$ ein Vertretersystem der
Bahnen unter der Operation von $G$ auf $M$ über $\vi$, die nach \ref{jordan2} in der dritten
Klasse in derselben Bahn liegen.

Für die vierte Klasse gilt: Die Matrizen $$\begin{pmatrix}
 b\sw  & 0 \\
 0     & -b \sw
\end{pmatrix} \quad \text{ und } \quad \begin{pmatrix}
 \sw  & 0 \\
 0    & -\sw
\end{pmatrix}\quad \text{ sowie } \quad \begin{pmatrix}
 0  & w \\
 1  & 0
\end{pmatrix}$$
auf der einen Seite und ebenso die Matrizen $$\begin{pmatrix}
 a+b\sw  & 0 \\
 0       & a-b \sw
\end{pmatrix} \quad \text{ und } \quad \begin{pmatrix}
 1+ a\1 b\sw  & 0 \\
 0       & 1- a\1 b \sw
\end{pmatrix},$$ $$\text{sowie } \quad \begin{pmatrix}
 1 + b' \sw & 0 \\
 0          & 1-b' \sw
\end{pmatrix} \quad \text{ und } \quad \begin{pmatrix}
 1   & b'w \\
 b'  & 1
\end{pmatrix}$$
mit $b$ und $b'$ aus $\f^+$ auf der anderen Seite sind unter Konjugation äquivalent.
\end{beweis}

\begin{lemma} \label{det}
Operiert die Gruppe $G=GL(2,p)$ vermittels $\theta: G \times M(2 \times 2,p) \ra M(2 \times
2,p): (g,m) \mt det(g) g\1 m g$ auf $M(2 \times 2, p)$, so ist durch die folgende Liste von
Matrizen ein Vertretersystem der insgesamt $p+7$ Bahnen unter dieser Operation angegeben,
wobei $w$ ein Erzeuger der multiplikativen Gruppe von $\f$ ist:
\begin{enumerate}
  \item[1a.] $\begin{pmatrix}
    0     & 0 \\
    0     & 0
  \end{pmatrix}$ entsprechend einer Bahn.
  \item[1b.] $\begin{pmatrix}
    1     & 0 \\
    0     & 1
  \end{pmatrix}$ entsprechend einer Bahn.
  \item[2a.] $\begin{pmatrix}
    0     & 1 \\
    0     & 0
  \end{pmatrix}$ entsprechend einer Bahn.
  \item[2a'.] $\begin{pmatrix}
    0     & w \\
    0     & 0
  \end{pmatrix}$ entsprechend einer Bahn.
  \item[2b.] $\begin{pmatrix}
    1     & 1 \\
    0     & 1
  \end{pmatrix}$ entsprechend einer Bahn.
  \item[2b'.] $\begin{pmatrix}
    1     & w \\
    0     & 1
  \end{pmatrix}$ entsprechend einer Bahn.
  \item[3a.] $\begin{pmatrix}
    0     & 0 \\
    0     & 1
  \end{pmatrix}$ entsprechend einer Bahn.
  \item[3b.] $\begin{pmatrix}
    1     & 0 \\
    0     & w^k
  \end{pmatrix}$ mit $k \in \{1 \dd \frac{p-1}{2}\}$ entsprechend $\frac{p-1}{2}$ Bahnen.
  \item[4a.] $\begin{pmatrix}
    0     & w \\
    1     & 0
  \end{pmatrix}$ entsprechend einer Bahn.
  \item[4b.] $\begin{pmatrix}
    1     & w^{k+1} \\
    w^k     & 1
  \end{pmatrix}$ mit $k \in \{0 \dd \frac{p-1}{2}-1\}$ entsprechend $\frac{p-1}{2}$ Bahnen.
\end{enumerate}
\end{lemma}

\begin{beweis}
Ist $d \in M(2 \times 2,\f[\sw])$ eine Diagonalmatrix, so ist für alle $g \in GL(p,2)$ die
Bedingung $det(g) g\1 d g= det(g) d g\1 g = det(g) d$ erfüllt. Da dann auch $z g\1 d g = z d
g\1 g = z g$ für alle $z$ aus dem Zentrum $Z$ von $GL(2,p)$ gilt, stimmt die Bahn von $d$
unter der Operation von $G$ über $\theta$ mit der Bahn von $d$ unter der Operation von $Z
\times G$ über $\delta: (Z \times G) \times M(2 \times 2,p) \ra M(2 \times 2,p): ((z,g),m) \mt
m^{(z,g)}=z g\1 m g$ überein. Daher kann das Vertretersystem der Bahnen von Diagonalmatrizen
aus \ref{kro} weit gehend übernommen werden. Es bleibt nur zu zeigen, dass die Matrizen
$$m=\begin{pmatrix}
    k     & 1 \\
    0     & k
  \end{pmatrix}$$ mit $k \in \f$ die vier Bahnen mit den in $2a, 2a', 2b$ und
$2b'$ aufgelisteten Repräsentanten bilden. Es sei $$g=\begin{pmatrix}
    a     & b \\
    c     & d
  \end{pmatrix} \in GL(2,p).$$ Dann ist
\begin{eqnarray*}
  det(g) g\1 m g &=& det(g) g\1 \begin{pmatrix}
    k     & 1 \\
    0     & k
  \end{pmatrix} g \\
  &=& det(g) g\1 \left( \begin{pmatrix}
    k     & 0 \\
    0     & k
  \end{pmatrix} + \begin{pmatrix}
    0     & 1 \\
    0     & 0
  \end{pmatrix} \right) g\\
  &=& det(g) g\1 \begin{pmatrix}
    k     & 0 \\
    0     & k
  \end{pmatrix} g + det(g) g\1 \begin{pmatrix}
    0     & 1 \\
    0     & 0
  \end{pmatrix} g \\
  &=& det(g) \begin{pmatrix}
    k     & 0 \\
    0     & k
  \end{pmatrix} g\1 g + det(g) \frac{1}{det(g)} \begin{pmatrix}
    d     & -b \\
    -c    & a
  \end{pmatrix} \begin{pmatrix}
    0     & 1 \\
    0     & 0
  \end{pmatrix} \begin{pmatrix}
    a     & b \\
    c     & d
  \end{pmatrix} \\
  &=& (ad-bc) \begin{pmatrix}
    k     & 0 \\
    0     & k
  \end{pmatrix} + \begin{pmatrix}
    cd    & d^2 \\
    -c^2  & -cd
  \end{pmatrix} \\
  &=& \begin{pmatrix}
    (ad-bc)k + cd  & d^2 \\
    -c^2           & (ad-bc)k -cd
  \end{pmatrix}
\end{eqnarray*}
Damit ist die Bahn von $m$
\begin{eqnarray*}
  B(k,1) &=&\begin{pmatrix}
    k     & 1 \\
    0     & k
  \end{pmatrix}^G \\
  &=& \left\{ \begin{pmatrix}
    (ad-bc)k + cd  & d^2 \\
    -c^2           & (ad-bc)k -cd
  \end{pmatrix} \mid a,b,c,d \in \f \text{ und } ad-bc \neq 0 \right\} \\
  &=& \left\{ \begin{pmatrix}
    k u + cd  & d^2 \\
    -c^2           & ku -cd
  \end{pmatrix} \mid c,d,u \in \f \text{ und } u \neq 0 \text{ sowie } (c,d) \neq
  (0,0)\right\}
\end{eqnarray*}
Für dieselbe Matrix ergibt sich nach \ref{kro} unter der Operation von $Z \times G$ über
$\delta: (Z \times G) \times M(2 \times 2,p): ((z,g),m) \mt m^{(z,g)}=z g\1 m g$ die Bahn
\begin{eqnarray*}
   C(k) &=& \begin{pmatrix}
    k     & 1 \\
    0     & k
  \end{pmatrix}^{Z \times G} \\
  &=& \left\{ \begin{pmatrix}
    z (k  +  cd)  & z d^2 \\
    - z c^2       & z(k -cd)
  \end{pmatrix} \mid c,d,u \in \f \text{ und } z \neq 0 \text{ sowie } (c,d) \neq
  (0,0)\right\}
\end{eqnarray*}
Es sei weiterhin
\begin{eqnarray*}
  B(k,w) &=& \begin{pmatrix}
    k     & w \\
    0     & k
  \end{pmatrix}^G \\
  &=& \left\{ \begin{pmatrix}
    w(k u + cd)  & w d^2 \\
    -w c^2           & w(ku -cd)
  \end{pmatrix} \mid c,d,u \in \f \text{ und } u \neq 0 \text{ sowie } (c,d) \neq
  (0,0)\right\} \\
  &=& \left\{ \begin{pmatrix}
    k u + cd  & w d^2 \\
    -w c^2           & ku -cd
  \end{pmatrix} \mid c,d,u \in \f \text{ und } u \neq 0 \text{ sowie } (c,d) \neq
  (0,0)\right\}
\end{eqnarray*}
Nach \ref{indexmult} ist für alle $k \in \f$ die Menge $C(k)$ die disjunkte Vereinigung der
Mengen $B(k,1)$ und $B(k,w)$. Damit gibt es neben den Bahnen der Arten $1,3$ und $4$ vier
weitere Bahnen, deren Vertreter in $2a, 2a', 2b$ und $2b'$ aufgelistet sind.
\end{beweis}

\begin{lemma} \label{opdp2}
Der Vektorraum $\f^3$ zerfällt unter der Operation der Gruppe $$A=\left\{ \begin{pmatrix}
  u^2 v  & - u^2 v^2 z & 0 \\
  0      & u v^2       & 0 \\
  0      & 0           & v
\end{pmatrix} \mid u, v \in \f^* \text{ und } z \in \f \right\}$$ in
folgende Bahnen, wobei $w$ ein Erzeuger der multiplikativen Gruppe von $\f$ ist:
\begin{enumerate}
  \item $B_0=(0,0,0)^A$ mit $|B_0|=1$,
  \item $B_1=(0,1,1)^A=(0,1,-1)^A$ mit $|B_1|=(p-1)^2$,
  \item $B_2=(1,0,1)^A=(1,0,-1)^A$ mit $|B_2|=\frac{1}{2}p(p-1)^2$,
  \item $B_3=(1,0,w)^A=(1,0,-w)^A$ mit $|B_3|=\frac{1}{2}p(p-1)^2$,
  \item $B_4=(0,0,1)^A$ mit $|B_4|=p-1$,
  \item $B_5=(0,1,0)^A$ mit $|B_5|=p-1$,
  \item $B_6=(1,0,0)^A$ mit $|B_6|=p(p-1)$.
\end{enumerate}
\end{lemma}

\begin{beweis}
Die Behauptung dürfte bis auf die Bahnen $B_2$ und $B_3$ offensichtlich sein. Der Stabilisator
dieser beiden Bahnen ist $$S=\left\{ \begin{pmatrix}
  u^2    & 0       & 0 \\
  0      & u       & 0 \\
  0      & 0       & 1
\end{pmatrix} \mid u \in \f^* \text{ und } u^2=1 \right\},$$ was sich über die
Lösung des Gleichungssystems aus $1=u^2v$, $0=u^2 v^2 z$ und $v w^a = w^a$ für $a \in \{0,1\}$
ergibt. Aus demselben Gleichungssystem ergibt sich nach \ref{indexmult}, dass $(1,0,w^a)$ und
$(1,0,w^b)$ genau dann in derselben Bahn liegen, wenn $a \equiv b$ mod $2$ ist. Eine Summation
der Bahnenlängen liefert die Vollständigkeit der oben stehenden Liste.
\end{beweis}

\begin{lemma} \label{dreieckvbar}
Es sei $$H=\left\{\begin{pmatrix}
  1   & a \\
  0   & b
\end{pmatrix} \mid a, b \in \f \text{ und } b \neq 0 \right\} \leq GL(2,p)$$ und $G=\erz{ z h\1
\otimes h \mid z \in \f^* \text{ und } h \in H}$. Weiterhin sei $V=\f^4$ und $\{v_1, v_2, v_3,
v_4\}$ die Standardbasis von $V$ sowie $U=\{v_2, v_4\}$ und $\bar{V}=V/U$. Operiert $G$ durch
Multiplikation von rechts auf $V$ und damit in natürlicher Weise auf dem Faktorraum $\bar{V}$,
so zerfällt $\bar{V}$ unter dieser Operation in folgende Bahnen:
\begin{enumerate}
  \item $B_1=(0,0)^G$,
  \item $B_2=(0,1)^G$,
  \item $B_3=(1,0)^G$.
\end{enumerate}
\end{lemma}

\begin{beweis}
Es ist
\begin{eqnarray*}
G &=& \left\langle z \begin{pmatrix}
  1   & -a b\1 \\
  0   & b\1
\end{pmatrix} \otimes \begin{pmatrix}
  1   & a \\
  0   & b
\end{pmatrix} \mid b \in \f \text{ und } z,b \in \f^*  \right\rangle \\
 &=& \left \langle  z \begin{pmatrix}
  1   & a  & -a b\1 & - a^2 b\1 \\
  0   & b  & 0      & - a \\
  0   & 0  & b\1    & a b\1 \\
  0   & 0  & 0      & 1
\end{pmatrix} \mid b \in \f \text{ und } z,b \in \f^*  \right\rangle.
\end{eqnarray*}
Aus dieser Darstellung ist ersichtlich, dass $G$ analog zu $$\bar{G}=\left\{z \begin{pmatrix}
  1   & -a b\1 \\
  0   & b\1
\end{pmatrix} \mid b \in \f \text{ und } z,b \in \f^* \right\} = \left\{\begin{pmatrix}
  u   & x \\
  0   & v
\end{pmatrix} \mid x \in \f \text{ und } u,v \in \f^* \right\}$$ auf $\bar{V}$ operiert,
d.\,h. $\bar{G}$ ist die Untergruppe der oberen Dreiecksmatrizen von $GL(2,p)$, die
bekanntermaßen die Bahnen $B_1$, $B_2$ und $B_3$ bei Multiplikation von rechts auf $\bar{V}$
bildet.
\end{beweis}

\chapter{Die Gruppen bis zur Ordnung $p^5$ für eine beliebige Primzahl $p>3$}

In diesem Kapitel werden die Gruppen der Ordnung $p^5$ klassifiziert. Von den
elementarabelschen Gruppe $C_p$, $C_p^2$, $C_p^3$ und $C_p^4$ aus werden Nachfolger bis zur
Ordnung $p^5$ berechnet und aus jeder Isomorphieklasse der Nachfolger wird ein Vertreter
ausgewählt, d.\,h. für $d \in \{1 \dd 4\}$ wird das folgende Verfahren durchlaufen:
\begin{enumerate}
  \item Man beginne mit der elementarabelschen Gruppe $C_p^d$, welche die
  $p$"=Klasse $1$ hat, d.\,h. man setze für eine Menge $L_1$ als
  Anfangswert des Verfahrens $L_1=\{C_p^d\}$.
  \item Man bestimme alle unmittelbaren Nachfolger der Gruppen aus
  $L_i$ und nehme von jedem Isomorphietyp einen Repräsentanten
  in eine Menge $L_{i+1}$ auf, sofern seine Ordnung maximal $p^5$
  ist.
  Alle Elemente von $L_{i+1}$ haben die $p$"=Klasse $i+1$.
  \item Man wiederhole Schritt $2$, bis $i$ den Wert $5$
  erreicht.
\end{enumerate}
Bei jedem Durchgang des zweiten Schritts erhöht sich die $p$"=Klasse der Gruppen um $1$ und
daher die Ordnung der Gruppen wenigstens um den Faktor $p$. Nach maximal vier Iterationen
liegen daher nur noch Gruppen der Ordnung $p^5$ vor. Die Gruppe $C_p^5$ tritt zu dieser Liste
anschließend hinzu. Die zentralen Sätze \ref{vollst} und \ref{expcover} stellen sicher, dass
mit dieser Methode ein Vertretersystem der Isomorphieklassen konstruiert wird.

Auf Seite \pageref{algunmittnachf} ist bereits beschrieben worden, wie im Schritt $2$ die
unmittelbaren Nachfolger berechnet werden. Dieses Verfahren wird nun für jede Gruppe wie dort
beschrieben durchgeführt. Daher ist jeder Abschnitt dieses Kapitels gleich aufgebaut und folgt
den Schritten des dort angegebenen Verfahrens. Der Satz \ref{bahnenzu} garantiert, dass
dadurch für jede Gruppe ein Vertretersystem der Isomorphieklassen ihrer unmittelbaren
Nachfolger entsteht.

Wie schwierig die Berechnung der unmittelbaren Nachfolger ist, hängt von der jeweiligen Gruppe
ab. Das wesentliche Problem besteht darin, die Bahnen zu finden, in welche die Untergruppen
des Multiplikators $M(G)$ unter der Operation von Erweiterungsautomorphismen zerfallen. In
einigen Fällen ist die Zerlegung offensichtlich oder leicht zu berechnen. In anderen Fällen
braucht man einige weniger triviale Hilfsmittel aus der Geometrie und linearen Algebra. Dies
betrifft vor allem die Gruppen $D_p, C_p^3, C_p^4$ und $(p^4,7)$.

Die Gruppe $D_p$ hat $p+8$ und die Gruppe $C_p^3$ hat $p+14$ Isomorphieklassen unmittelbarer
Nachfolger der Ordnung $p^5$. Der lineare Anstieg erklärt sich folgendermaßen: Die
Automorphismengruppen dieser beiden Gruppen operieren auf ihren Multiplikatoren in ähnlicher
Weise wie $GL(2,p)$ durch Konjugation auf der Menge der $(2 \times 2)$"=Matrizen über $\f$.
Diese Operation hat eine von $p$ abhängige Anzahl von Bahnen. Daher dürfte es nicht wundern,
dass die Bestimmung der Nachfolger von $D_p$ und vor allem von $C_p^3$ Anleihen aus der
Jordanblockzerlegung nimmt und dass das Kapitel über Gruppenoperationen einen längeren Exkurs
über Operationen der Gruppe $GL(2,p)$ enthält, auf den hier zurückgegriffen wird. Die
Automorphismengruppe von $C_p^4$ operiert auf ihrem Multiplikator in derselben Weise wie
$GL(4,p)$ auf der Menge der schiefsymmetrischen Matrizen von $M(4 \times 4,p)$.

Als vierte ungewöhnliche Gruppe ist $(p^4,7)$ dafür verantwortlich, dass in der Anzahlformel
über die Isomorphieklassen der Gruppen der Ordnung $p^5$ die beiden Terme mit größten
gemeinsamen Teilern auftreten. Der Grund dafür ist, dass der Operationshomomorphismus der
Automorphismengruppe von $(p^4,7)$ in die Automorphismengruppe von $M(G)$ aus
Polynomfunktionen höheren Grades zusammengesetzt ist und die Abbildung $\theta: \f^* \ra \f^*:
x \mt x^k$ als Bild eine Untergruppe vom Index $ggT(k,p-1)$ in $\f^*$ hat. Daher hat das Bild
von $Aut(G)$ in $Aut(M(G))$ eine Ordnung, die von $p$ abhängig ist, und liefert je nach Wahl
von $p$ Bahnen unterschiedlicher Länge und verschiedener Anzahl.

Alle anderen Gruppen haben eine Anzahl unmittelbarer Nachfolger, die unabhängig von $p$ ist.
Die Bahnenberechnung lässt sich bei ihnen durch die Lösung einfacher, größtenteils linearer
Gleichungssysteme bewältigen und stellt kein Problem dar.

Um die Sprache nicht zu unübersichtlich werden zu lassen, werden die Vertreter mit ihren
Isomorphieklassen identifiziert. Es wird also beispielsweise statt "`die Isomorphieklasse der
Gruppen der Ordnung $p^5$"' nur "`die Gruppen der Ordnung $p^5$"' oder statt "`die
Isomorphieklasse der Nachfolger von $G$"' nur "`die Nachfolger von $G$"' geschrieben.

\begin{beme}
In diesem Kapitel sei $p$ eine beliebige Primzahl und $p>3$. In allen
Potenz"=Kommutator"=Präsentationen außer bei den elementarabelschen Gruppen werden im Sinne
von \ref{kurz} und \ref{praes} die trivialen Relationen weggelassen.
\end{beme}

\begin{defi}
Es sei $G$ eine Gruppe, $n \in \N$, $p$ eine Primzahl und $\vi: G \ra GL(n,p)$ ein
Operationshomomorphismus. Dann ist $\bar{\vi}: G \ra GL(n,p): g \mapsto (g^\vi)^T$ der zu
$\vi$ \emph{duale Operationshomomorphismus}, wobei $T$ die Transpositionsabbildung von
$GL(n,p)$ ist.
\end{defi}

\section{Nachfolger von $C_p$}

\begin{satz} \label{ncp}
Für jedes $i \in \N$ ist $C_{p^{i+1}}$ der einzige unmittelbaren Nachfolger von $C_{p^i}$.
\end{satz}

\begin{beweis}
Nach dem Hauptsatz über endlich erzeugte abelsche Gruppen (vgl. etwa \cite{mey80}, S. 78) ist
für jedes $i \in \N$ eine Gruppe $G$ der Ordnung $p^i$, sodass $G$ zyklisch ist, isomorph zu
$C_{p^i}$. Daher haben $C_{p^i}$ und $C_{p^{i+1}}$ dieselbe minimale Anzahl von Erzeugern. Da
nach dem Basissatz von Burnside und wiederum nach dem Hauptsatz über endlich erzeugte abelsche
Gruppen die Faktorgruppe $F_j=\gamma_j(C_{p^i})/\gamma_{j+1}(C_{p^i})$ für alle $j \in \N_0$
isomorph zu $C_p$ ist, sofern $F_j$ nicht trivial ist, hat $C_{p^{i+1}}$ nach \ref{dlem} eine
um 1 größere $p$"=Klasse als $C_{p^{i}}$ und aus demselben Grund ist
$C_{p^{i+1}}/\gamma_c(C_{p^{i+1}})$ isomorph zu $C_{p^{i}}$, sofern $c$ die $p$"=Klasse von
$C_{p^{i}}$ ist. Damit ist $C_{p^{i+1}}$ ein unmittelbarer Nachfolger von $C_{p^{i}}$ und nach
dem Hauptsatz über endlich erzeugte abelsche Gruppen Repräsentant der einzigen
Isomorphieklasse unmittelbarer Nachfolger von $C_{p^{i}}$.
\end{beweis}

\section{Nachfolger von $C_p^2$}

\begin{ver}
In diesem Abschnitt sei $G=\erz{a_1, a_2 \mid [a_2, a_1]=1, a_1^p=1, a_2^p=1}$. Damit ist $G
\cong C_p^2$.
\end{ver}

\subsection{Die Automorphismengruppe}

\begin{beme}
Die Gruppe $G$ ist isomorph zu additiven Gruppe des Vektorraums $\F_p^2$. Die
Automorphismengruppe $Aut(G)$ kann daher mit mit der generellen linearen Gruppe $GL(2,p)=:A$
identifiziert werden.
\end{beme}

\subsection{$p$"=Cover, Multiplikator und Nukleus}

\begin{lemma}
Für $G$ gilt:
\begin{itemize}
  \item[] $P(G)=\erz{a_1, a_2, a_3, a_4, a_5 \mid [a_2, a_1]=a_3, a_1^p=a_4,
  a_2^p=a_5}$.
  \item[] $M(G)=\erz{a_3, a_4, a_5}$.
  \item[] $N(G)=\erz{a_3, a_4, a_5}$.
\end{itemize}
Die Gruppe $P(G)$ hat die $p$"=Klasse 2, $G$ ist fortsetzbar und jede
Untergruppe von $M(G)$ ist zulässig.
\end{lemma}

\begin{beweis}
Die Präsentation von $P(G)$ ergibt sich aus \ref{cov} und \ref{newman} und dem
Reduktionsverfahren von Knuth"=Bendix, aus der sich $M(G)$ ablesen lässt. Die Gewichtung von
$P(G)$ und damit die Untergruppe $N(G)$ lassen sich nach \ref{gew} aus der Präsentation von
$P(G)$ ablesen: $1=\omega(a_1)=\omega(a_2)$ und $2=\omega(a_3)=\omega(a_4)=\omega(a_5)$.
\end{beweis}

\subsection{Nachfolger der Ordnung $p^5$}

\begin{folg} \label{ncp2p5}
Die Gruppe $G$ hat genau einen unmittelbaren Nachfolger der Ordnung $p^5$, nämlich $P(G)$.
\end{folg}


\subsection{Operation der Erweiterungsautomorphismen}

\begin{lemma}
Mit der Identifikation von $Aut(G)$ mit $GL(2,p)$ ist die Operation von $A=GL(2,p)$ auf $\f^3
\cong M(G)$ durch $$\vi: GL(2,p) \to GL(3,p):m=\begin{pmatrix}
  m_{1} & m_{2}  \\
  m_{3} & m_{4}  \\
\end{pmatrix} \mt M = \begin{pmatrix}
  |m| & 0 & 0 \\
  0 & m_1 & m_2 \\
  0 & m_3 & m_4
\end{pmatrix}$$ gegeben, wobei $|m|$ die Determinante von $m$ ist.
\end{lemma}

\begin{beweis}
Die Bilder der Erzeuger von $M(G)$ ergeben sich unter $m$ vermittels der
Operation über Erweiterungsautomorphismen folgendermaßen:
\begin{eqnarray*}
  a_3^m &=& [a_2,a_1]^m=[a_2^m, a_1^m] \\
  &=& [a_1^{m_3} a_2^{m_4}, a_1^{m_1} a_2^{m_2}]= [a_1^{m_3}, a_1^{m_1} a_2^{m_2}]^{a_2^{m_4}} [a_2^{m_4}, a_1^{m_1}
  a_2^{m_2}]\\
  &=& ([a_1^{m_3}, a_2^{m_2}][a_1^{m_3},
  a_1^{m_1}]^{a_2^{m_2}})^{a_2^{m_4}} [a_2^{m_4}, a_2^{m_2}] [a_2^{m_4},
  a_1^{m_1}]^{a_2^{m_2}}\\
  &=& [a_1^{m_3}, a_2^{m_2}]^{a_2^{m_4}} [a_2^{m_4},
  a_1^{m_1}]^{a_2^{m_2}} = ([a_2, a_1]^{-m_3 m_2})^{a_2^{m_4}} ([a_2,
  a_1]^{m_4 m_1})^{a_2^{m_2}}\\
  &=& [a_2, a_1]^{m_1 m_4 - m_3 m_2} = a_3^{|m|},
\end{eqnarray*}
da $[a_2, a_1]=a_3$ als Element von $M(G)$ zentral in $P(G)$ ist. Weiterhin ist
\begin{eqnarray*}
  a_4^m &=& (a_1^p)^m=(a_1^m)^p=(a_1^{m_1} a_2^{m_2})^p \\
  &=& (a_1^{m_1})^p (a_2^{m_2})^p =(a_1^p)^{m_1} (a_2^p)^{m_2} \\
  &=& a_4^{m_1} a_5^{m_2},
\end{eqnarray*}
da $a_3$ der Kommutator von $a_1$ und $a_2$ ist und die Ordnung $p$ hat. Daher
gibt es ein $k \in \{0 \dd p-1\}$, sodass $(a_1^{m_1} a_2^{m_2})^p =
(a_1^{m_1})^p (a_3^k)^p (a_2^{m_2})^p = (a_1^p)^{m_1} (a_3^p)^k (a_2^p)^{m_2} =
(a_1^p)^{m_1}$ ist. Analog ergibt sich $a_5^m=a_4^{m_3} a_5^{m_4}$.
\end{beweis}

\subsection{Bahnen zulässiger Untergruppen}

\begin{lemma} \label{bahnen_cp_2}
Über den Operationshomomorphismus $\vi$ zerfällt $\F_p^3 \cong M(G)$ unter der Operation von
$A$ in folgende Bahnen nichttrivialer Vektoren.
\begin{itemize}
  \item[] $B_1=(1,0,0)^A$ mit $|B_1|=p-1$.
  \item[] $B_2=(0,0,1)^A$ mit $|B_2|=p^2-1$.
  \item[] $B_3=(1,0,-1)^A$ mit $|B_3|=(p-1)(p^2-1)$.
\end{itemize}
Die zugehörigen Stabilisatoren sind
\begin{eqnarray*}
 & \bar{S}_1= & SL(2,p) \\
 & \bar{S}_2= & \left\{ \begin{pmatrix}
    m_{1} & m_2  \\
    0     & 1  \\
  \end{pmatrix} \mid m_1, m_2 \in \F_p \text{ und } m_1 \neq 0
  \right\} \\
 & \bar{S}_3= & \left\{ \begin{pmatrix}
    1 & m_2  \\
    0 & 1  \\
  \end{pmatrix} \mid m_2 \in \F_p \right\}
\end{eqnarray*}
\end{lemma}

\begin{beweis}
Die Unterräume $\erz{(1,0,0)}$ und $\erz{(0,1,0),(0,0,1)}$ sind unter der Operation von $A$
irreduzible Teilmoduln. Ihnen entsprechen die Bahnen $B_1$ und $B_2$. Es ist offensichtlich
$\bar{S}_1=SL(2,p)$. Den Stabilisator $\bar{S}_2$ erhält man als Lösungsmenge des
Gleichungssystems $\{m_3=0, m_4=1\}$ und $\bar{S}_3$ als Lösungsmenge von $\{m_1 m_4 - m_3
m_2=1, m_3=0, m_4=1\}$. Aus $|\bar{S}_3|=p$ ergibt sich nach \ref{glslz} und \ref{stab} die
Länge von $B_3$. Addiert man die Bahnenlängen, so ergibt sich $p^3-1$. Damit ist gezeigt, dass
die Liste der Bahnen vollständig ist.
\end{beweis}

\begin{lemma} \label{bahnen_cp_2a}
Über den Operationshomomorphismus $\vi$ zerfällt die Menge der eindimensionalen Unterräume
$\F_p^3 \cong M(G)$ unter der Operation von $A$ in folgende Bahnen.
\begin{itemize}
  \item[] $B_1'=\erz{(1,0,0)}^A$.
  \item[] $B_2'=\erz{(0,0,1)}^A$.
  \item[] $B_3'=\erz{(1,0,-1)}^A$.
\end{itemize}
Die Stabilisatoren $S_1$, $S_2$ und $S_3$ dieser Unterräume lauten:
\begin{eqnarray*}
 & S_1= & GL(2,p)=A \\
 & S_2= & \left\{ \begin{pmatrix}
    m_{1} & m_2  \\
    0     & m_{4}  \\
  \end{pmatrix} \mid m_1, m_2, m_4 \in \F_p \text{ und } m_1 \neq 0, m_4 \neq 0
  \right\} \\
 & S_3= & \left\{ \begin{pmatrix}
    1     & m_2  \\
    0     & m_4  \\
  \end{pmatrix} \mid m_2, m_4 \in \F_p \text{ und } m_4 \neq 0 \right\}
\end{eqnarray*}

\end{lemma}

\begin{beweis}
Die Behauptung ergibt sich unmittelbar aus \ref{bahnen_cp_2} und \ref{einvek}, da $A^\vi$
offensichtlich das Zentrum von $GL(3,p)$ enthält. Den Stabilisator $S_2$ erhält man als
Lösungsmenge des Gleichungssystems $\{m_3=0, m_4=a\}$ und $S_3$ als Lösungsmenge von $\{m_1
m_4 - m_3 m_2=a, m_3=0, m_4=a\}$, wobei $a \in \f^*$ ist.
\end{beweis}

\subsection{Nachfolger der Ordnung $p^4$}

\begin{folg} \label{zu1}
Für die Nachfolger von $G$ der Ordnung $p^4$ ergibt sich das folgende
Vertretersystem zulässiger Untergruppen:
\begin{enumerate}
  \item[] $M_1=\erz{a_3}$ entsprechend $\erz{(1,0,0)}$.
  \item[] $M_2=\erz{a_5}$ entsprechend $\erz{(0,0,1)}$.
  \item[] $M_3=\erz{a_3 a_5\1}$ entsprechend $\erz{(1,0,-1)}$.
\end{enumerate}
\end{folg}

\begin{satz} \label{ncp2p4}
In der folgenden Liste sind sämtliche Nachfolger von $G$ der Ordnung $p^4$ angegeben:
\begin{enumerate}
  \item[] $G_1=\erz{g_1, g_2, g_3, g_4 \mid g_1^p=g_3,
  g_2^p=g_4} \cong C_{p^2}^2$ \hfill \gap"=Typ $(p^4,2)$
  \item[] $G_2=\erz{g_1, g_2, g_3, g_4 \mid [g_2, g_1]=g_3, a_1^p=g_4}$ \hfill
  \gap"=Typ $(p^4,3)$
  \item[] $G_3=\erz{g_1, g_2, g_3, g_4 \mid [g_2, g_1]=g_4, g_1^p=g_3,
  g_2^p=g_4}$ \hfill \gap"=Typ $(p^4,4)$
\end{enumerate}
\end{satz}

\begin{beweis}
Dieser Satz ergibt sich nach \ref{bahnenzu}, indem man $P(G)$ nach jedem Repräsentanten des
Vertretersystems zulässiger Untergruppen faktorisiert, das in \ref{zu1} aufgelistet ist. Die
Präsentation der Faktorgruppe erhält man über das Verfahren aus \ref{faktor}.
\end{beweis}

\subsection{Nachfolger der Ordnung $p^3$}

Mit den Gruppen $G_1$, $G_2$ und $G_3$ sind alle Nachfolger der Ordnung $p^4$ ermittelt. Die
Nachfolger der Ordnung $p^3$ erhält man, indem man zulässige Untergruppen aus $P(G)$
ausfaktorisiert, denen zweidimensionale Unterräume in $M(G)$ entsprechen, wenn man $M(G)$ als
Vektorraum über $\F_p$ auffasst. Da $M(G)$ als solcher die Dimension 3 hat, sind die
zweidimensionalen Unterräume von $M(G)$ Komplemente der eindimensionalen. Nach dem
Dualitätsprinzip entsprechen die Bahnen der zweidimensionalen Unterräume denen ihrer
eindimensionalen orthogonalen Komplementen. Daher erhält man über die Bahnen $B_1$, $B_2$ und
$B_3$ unmittelbar die Bahnen $C_1$, $C_2$ und $C_3$ der zweidimensionalen Unterräume. Um
ansprechendere Präsentationen zu bekommen, werden hier andere Vertreter von $B_1$, $B_2$ und
$B_3$ gewählt als für die Nachfolger der Ordnung $p^4$.

\begin{lemma} \label{stab_cp2_p3}
Die Menge der zweidimensionalen Unterräume von $\F_p^3 \cong M(G)$ zerfällt
unter der Operation von $A$ in folgende Bahnen:
\begin{enumerate}
  \item[] $C_1'=\erz{(0,1,0),(0,0,1)}^A$ , orthogonales Komplement zu $\erz{(1,0,0)}$.
  \item[] $C_2'=\erz{(1,0,0),(0,0,1)}^A$, orthogonales Komplement zu $\erz{(0,1,0)}$.
  \item[] $C_3'=\erz{(1,-1,0),(0,0,1)}^A$, orthogonales Komplement zu $\erz{(1,1,0)}$.
\end{enumerate}
Die Stabilisator $T_1$, $T_2$ und $T_3$ zu $C_1'$, $C_2'$ und $C_3'$ sind:
\begin{eqnarray*}
 & T_1= & GL(2,p)=A \\
 & T_2= & \left\{ \begin{pmatrix}
    m_{1} & m_2  \\
    0     & m_{4}  \\
  \end{pmatrix} \mid m_1, m_2, m_4 \in \F_p \text{ und } m_1 \neq 0, m_4 \neq 0
  \right\} \\
 & T_3= & \left\{ \begin{pmatrix}
    m_1     & m_2  \\
    0       & 1  \\
  \end{pmatrix} \mid m_1, m_2 \in \F_p \text{ und } m_1 \neq 0 \right\}
\end{eqnarray*}
\end{lemma}

\begin{beweis}
Diese Behauptung folgt unmittelbar aus \ref{bahnen_cp_2} und dem Dualitätsprinzip. Um für
$C_2'$ und $C_3'$ Stabilisatoren in oberer Dreiecksform zu erhalten, sind hier andere
Vertreter von $B_2'$ und $B_3'$ gewählt worden als in \ref{bahnen_cp_2}.
\end{beweis}

\begin{folg} \label{zu2}
Für Nachfolger der Ordnung $p^3$ ergeben sich die folgenden Repräsentanten der
Bahnen zulässiger Untergruppen:
\begin{enumerate}
  \item[] $N_1=\erz{a_4, a_5}$ entsprechend $C_1'=\erz{(0,1,0),(0,0,1)}^A$.
  \item[] $N_2=\erz{a_3, a_5}$ entsprechend $C_2'=\erz{(1,0,0),(0,0,1)}^A$.
  \item[] $N_3=\erz{a_3 a_4\1, a_5}$ entsprechend $C_3'=\erz{(1,-1,0),(0,0,1)}^A$.
\end{enumerate}
\end{folg}

\begin{satz} \label{ncp2p3}
In der folgenden Liste sind sämtliche unmittelbaren Nachfolger von $G$ der
Ordnung $p^3$ bis auf Isomorphie angegeben:
\begin{enumerate}
  \item[] $H_1=\erz{h_1, h_2, h_3 \mid [h_2, h_1]=h_3} \cong D_p$.
  \item[] $H_2=\erz{h_1, h_2, h_3 \mid h_1^p=h_3} \cong C_{p^2} \times C_p$.
  \item[] $H_3=\erz{h_1, h_2, h_3 \mid [h_2, h_1]=h_3, h_1^p=h_3} \cong Q_p$.
\end{enumerate}
\end{satz}

\begin{beweis}
Dieser Satz ergibt sich nach \ref{bahnenzu}, indem man $P(G)$ nach jedem Repräsentanten des
Vertretersystems zulässiger Untergruppen faktorisiert, das in \ref{zu2} aufgelistet ist. Die
Präsentation der Faktorgruppe erhält man über das Verfahren aus \ref{faktor}.
\end{beweis}

Damit liegt für $C_p^2$ eine vollständige Liste der unmittelbaren Nachfolger bis zur Ordnung
$p^5$ vor. Als nächstes wird untersucht, ob die sechs unmittelbaren Nachfolger von $C_p^2$,
die nicht die Ordnung $p^5$ haben, ihrerseits unmittelbare Nachfolger bis zur Ordnung $p^5$
besitzen. Diese Untersuchung beginnt mir den Gruppen der Ordnung $p^3$, also $H_1$, $H_2$ und
$H_3$, in der Reihenfolge, in der sie hier aufgetreten sind.

\section{Nachfolger von $D_p$}

\begin{ver}
In diesem Abschnitt sei $G=\erz{a_1, a_2, a_3 \mid [a_2, a_1]=a_3}$. Damit ist $G \cong D_p$.
\end{ver}

\begin{lemma}
Die Gruppe $G$ hat die Gewichtung $\omega(a_1)=\omega(a_2)=1$ und $\omega(a_3)=2$.
\end{lemma}

\begin{beweis}
Die Gewichtung lässt sich nach \ref{gew} ermitteln.
\end{beweis}

\subsection{Die Automorphismengruppe}

\begin{lemma} \label{autdp}
Die Automorphismengruppe von $G$ besteht genau aus folgenden Abbildungen: $$\al(u_1, u_2, u_3,
v_1, v_2, v_3):
  \begin{cases}
    a_1 \mt a_1^{u_1} a_2^{u_2} a_3^{u_3}                & \text{wobei } u_1, u_2, u_3, v_1, v_2, v_3 \\
    a_2 \mt a_1^{v_1} a_2^{v_2} a_3^{v_3}                & \in \{0 \dd p-1\} \text{ und }\\
    a_3 \mt                     a_3^{u_1 v_2 - v_1 u_2}  & u_1 v_2 - v_1 u_2 \neq 0 \text{ ist.}
  \end{cases}$$

\end{lemma}

\begin{beweis}
Nach \ref{aut} lässt sich $Aut(G)$ unmittelbar aus dem Stabilisator $T_1$ (siehe
\ref{stab_cp2_p3}) der zu $G$ gehörenden zulässigen Untergruppe ablesen. Damit ergeben sich
genau die oben angegebenen Bilder von $a_1$ und $a_2$. Das Bild von $a_3$ ist durch die Bilder
von $a_1$ und $a_2$ vollständig festgelegt und lässt sich folgendermaßen berechnen:
\begin{eqnarray*}
  a_3^\al &=& [a_2, a_1]^\al=[a_1^\al, a_2^\al] \\
          &=& [a_1^{u_1} a_2^{u_2} a_3^{u_3}, a_1^{v_1} a_2^{v_2} a_3^{v_3}]\\
          &=& [a_1^{u_1} a_2^{u_2}, a_1^{v_1} a_2^{v_2}] = [a_1, a_2]^{u_2 v_1} [a_2, a_1]^{u_1 v_2}\\
          &=& a_3^{u_1 v_2 - v_1 u_2},
\end{eqnarray*}
da $a_3=[a_2, a_1]$ nach der Gewichtung $\omega$ zentral ist.
\end{beweis}

\subsection{$p$"=Cover, Multiplikator und Nukleus}

\begin{lemma}
Für $G$ ergibt sich das $p$"=Cover, der Nukleus und der Multiplikator in folgender Weise:
\begin{enumerate}
  \item[] $P(G)=\erz{a_1,a_2,a_3,a_4,a_5,a_6,a_7 \mid [a_2,a_1]=a_3,
  [a_3,a_1]=a_4, [a_3,a_2]=a_5, a_1^p=a_6, a_2^p=a_7}$.
  \item[] $M(G)=\erz{a_4, a_5, a_6, a_7}$.
  \item[] $N(G)=\erz{a_4, a_5}$.
\end{enumerate}
\end{lemma}

\begin{beweis}
Die Präsentation von $P(G)$ ergibt sich aus \ref{cov} und \ref{newman} und dem
Reduktionsverfahren von Knuth"=Bendix, aus der sich $M(G)$ ablesen lässt. Die Gewichtung von
$P(G)$ und damit die Untergruppe $N(G)$ lässt sich nach \ref{gew} aus der Präsentation von
$P(G)$ ablesen:
\begin{itemize}
  \item[] $1=\omega(a_1)=\omega(a_2)$.
  \item[] $2=\omega(a_3)=\omega(a_6)=\omega(a_7)$.
  \item[] $3=\omega(a_4)=\omega(a_5)$.
\end{itemize}
\end{beweis}

\subsection{Operation der Erweiterungsautomorphismen}

\begin{lemma}
Nach \ref{autdp} lässt sich jeder Automorphismus $\al$ von $D_p$ in der Form $$\al(u_1, u_2,
u_3, v_1, v_2, v_3):
  \begin{cases}
    a_1 \mt a_1^{u_1} a_2^{u_2} a_3^{u_3} & \text{wobei } u_1, u_2, u_3, v_1, v_2, v_3 \\
    a_2 \mt a_1^{v_1} a_2^{v_2} a_3^{v_3} & \in \{0 \dd p-1\} \text{ und }\\
    a_3 \mt a_3^{u_1 v_2 - v_1 u_2} & u_1 v_2 - v_1 u_2 \neq 0 \text{ ist.}
  \end{cases}$$
darstellen. In dieser Darstellung der Automorphismen von $G$ ist der Operationshomomorphismus
von $A:=Aut(G)$ über Erweiterungsautomorphismen auf $\f^4 \cong M(G)$ durch $$\vi: A \to
GL(4,p): \al(u_1, u_2, u_3, v_1, v_2, v_3) \mt M
=
\begin{pmatrix}
  d u_1 & d u_2 & 0   & 0 \\
  d v_1 & d v_2 & 0   & 0 \\
  0     & 0     & u_1 & u_2 \\
  0     & 0     & v_1 & v_2
\end{pmatrix}$$ gegeben, wobei $d= u_1 v_2 - v_1 u_2$ ist.
\end{lemma}

\begin{beweis}
Die Bilder der Erzeuger von $M(G)$ ergeben sich über die Operation durch
Erweiterungsautomorphismen unter $\al$ nach \ref{erwaut} folgendermaßen:
\begin{eqnarray*}
  a_4^\al &=& [a_3,a_1]^\al=[a_3^{u_1 v_2 - v_1 u_2}, a_1^{u_1} a_2^{u_2} a_3^{u_3}]\\
          &=& [a_3^{u_1 v_2 - v_1 u_2}, a_1^{u_1} a_2^{u_2}]\\
          &=& [a_3, a_1]^{(u_1 v_2 - v_1 u_2) u_1} [a_3, a_2]^{(u_1 v_2 - v_1 u_2)u_2}\\
          &=& a_4^{d u_1} a_5^{d u_2},
\end{eqnarray*} da $a_4$ und $a_5$ zentral sind. Ebenso ergibt sich $a_5^\al=
a_4^{d v_1} a_5^{d v_2}$. Weiterhin ist
\begin{eqnarray*}
  a_6^\al &=& (a_1^p)^\al=(a_1^\al)^p=(a_1^{u_1} a_2^{u_2} a_3^{u_3})^p\\
          &=& (a_1^{u_1} a_2^{u_2})^p (a_3^p)^{u_3} = (a_1^{u_1} a_2^{u_2})^p\\
          &=& (a_1^{u_1})^p (a_2^{u_2})^p=a_6^{u_1} a_7^{u_2},
\end{eqnarray*}
da  $a_3$ der Kommutator von $a_1$ und $a_2$ ist und die Ordnung $p$ hat.
Deshalb  gibt es ein $k \in \{0 \dd p-1\}$, sodass $(a_1^{u_1} a_2^{u_2})^p =
(a_1^{u_1})^p (a_3^k)^p (a_2^{u_2})^p = (a_1^p)^{u_1} (a_3^p)^k (a_2^p)^{u_2} =
(a_1^p)^{u_1} (a_2^p)^{u_2}$ ist. Analog erhält man $a_7^\al=a_6^{v_1}
a_7^{v_2}$.
\end{beweis}

\subsection{Bahnen zulässiger Untergruppen für Nachfolger der Ordnung $p^4$}

\begin{lemma} \label{bahnen_1_dp}
Über den Operationshomomorphismus $\vi$ zerfällt $\F_p^4 \cong M(G)$ unter der Operation von
$A$ in folgende Bahnen nichttrivialer Vektoren.
\begin{enumerate}
  \item[] $B_1=(0,0,1,0)^A$ mit $|B_1|= p^2-1$.
  \item[] $B_2=(1,0,0,0)^A$ mit $|B_2|= p^2-1$.
  \item[] $B_3=(1,0,1,0)^A$ mit $|B_3|= (p^2-1)(p-1)$.
  \item[] $B_4=(1,0,0,1)^A$ mit $|B_4|= \frac{(p^2-1)(p^2-p)}{2}$.
  \item[] $B_5=(1,0,0,w)^A$ mit $|B_5|= \frac{(p^2-1)(p^2-p)}{2}$.
\end{enumerate}
Dabei ist $w$ ein Erzeuger der multiplikativen Gruppe von $\f$. Die Stabilisatoren der
Bahnrepräsentanten sind die folgenden. Dabei ist $\hat{S}_i \leq GL(4,p)$ das Bild des
Stabilisators $\bar{S}_i \leq Aut(G)$ unter $\vi$, wobei $\bar{S}_i$ den Vertreter der Bahn
$B_i$ stabilisiert.
\begin{eqnarray*}
  \hat{S}_1 & = & \left\{ \begin{pmatrix}
    v_2     & 0     & 0   & 0   \\
    v_1 v_2 & v_2^2 & 0   & 0   \\
    0       & 0     & 1   & 0   \\
    0       & 0     & v_1 & v_2 \\
  \end{pmatrix} \mid v_1, v_2 \in \F_p \text{ und } v_2 \neq 0
  \right\} \text{ mit } |\hat{S}_1| = p(p-1)\\
  \hat{S}_2 & = & \left\{ \begin{pmatrix}
    1         & 0        & 0   & 0        \\
    v_1 u_1\1 & u_1^{-3} & 0   & 0        \\
    0         & 0        & u_1 & 0        \\
    0         & 0        & v_1 & u_1^{-2} \\
  \end{pmatrix} \mid u_1, v_1 \in \F_p \text{ und } u_1 \neq 0
  \right\} \text{ mit } |\hat{S}_2| = p(p-1)\\
  \hat{S}_3 & = & \left\{ \begin{pmatrix}
    1         & 0        & 0   & 0        \\
    v_1       & 1        & 0   & 0        \\
    0         & 0        & 1   & 0        \\
    0         & 0        & v_1 & 1        \\
  \end{pmatrix} \mid v_1 \in \F_p \right\} \text{ mit } |\hat{S}_3| = p\\
  \hat{S}_4 = \hat{S}_5 & = & \left\{ \begin{pmatrix}
    1         & 0        & 0     & 0        \\
    0         & \pm 1    & 0     & 0        \\
    0         & 0        & \pm 1 & 0        \\
    0         & 0        & 0     & 1         \\
  \end{pmatrix} \right\} \text{ mit } |\hat{S}_4| = |\hat{S}_5| = 2\\
\end{eqnarray*}
\end{lemma}

\begin{beweis}
Dass $\F_p^4$ derart unter $A$ in Bahnen zerfällt, ist nahezu offensichtlich. Es bleibt zu
bemerken: Die Bahnen $B_1$ und $B_2$ sind verschieden, da die Unterräume
$\erz{(1,0,0,0),(0,1,0,0)}$ und $\erz{(0,0,1,0),(0,0,0,1)}$ irreduzible Teilmuduln unter der
Operation von $A$ sind. Dass $B_4$ und $B_5$ nicht identisch sind, erkennt man daran, dass die
Gleichung $(d u_1, d u_2, v_1, v_2)=(1,0,0,w)$ auf die Gleichung $1=d u_1 =(u_1 v_2 - v_1 u_2)
u_1 = w u_1^2$ führt, die ihrerseits nicht lösbar ist, da $w$ als Erzeuger der multiplikativen
Gruppe von $\F_p$ nach \ref{indexmult} kein Quadrat ist. Damit liegt $(1,0,0,w)$ nicht in
$B_4$ und die beiden Bahnen sind daher verschieden. Addiert man die Länge der Bahnen, so sieht
man, dass $\F_p^4$ mit ihnen vollständig aufgeteilt ist.

Die Stabilisatoren erhält man als Lösungsmengen der folgenden Gleichungssysteme
\begin{enumerate}
  \item $\{u_1=1, u_2=0\}$,
  \item $\{(u_1 v_2 - v_1 u_2)u_1=1, (u_1 v_2 - v_1 u_2)u_2=0\}$,
  \item $\{(u_1 v_2 - v_1 u_2)u_1=1, (u_1 v_2 - v_1 u_2)u_2=0, u_1=1, u_2=0\}$,
  \item $\{(u_1 v_2 - v_1 u_2)u_1=1, (u_1 v_2 - v_1 u_2)u_2=0, v_1=0, v_2=1\}$,
  \item $\{(u_1 v_2 - v_1 u_2)u_1=1, (u_1 v_2 - v_1 u_2)u_2=0, v_1=0, v_2=w\}$.
\end{enumerate}

\end{beweis}

\begin{lemma} \label{bahnen_1_dpa}
Über den Operationshomomorphismus $\vi$ zerfällt die Menge der eindimensionalen Unterräume
$\F_p^4 \cong M(G)$ unter der Operation von $A$ in folgende Bahnen.
\begin{enumerate}
  \item[] $B_1'=\erz{(0,0,1,0)}^A$.
  \item[] $B_2'=\erz{(1,0,0,0)}^A$.
  \item[] $B_3'=\erz{(1,0,1,0)}^A$.
  \item[] $B_4'=\erz{(1,0,0,1)}^A$.
  \item[] $B_5'=\erz{(1,0,0,w)}^A$.
\end{enumerate}
Die zugehörigen Stabilisatoren sind
\begin{eqnarray*}
  \bar{S}_1 & = & A \\
  \bar{S}_2 & = & \left\{ \begin{pmatrix}
    u_1^2 v_2   & 0         & 0   & 0        \\
    u_1 v_1 v_2 & u_1 v_2^2 & 0   & 0        \\
    0           & 0         & u_1 & 0        \\
    0           & 0         & v_1 & v_2 \\
  \end{pmatrix} \mid u_1, v_1, v_2 \in \F_p \text{ und } u_1 v_2 \neq 0
  \right\}\\
  \bar{S}_3& = & \left\{ \begin{pmatrix}
    u_1         & 0         & 0   & 0        \\
    v_1         & u_1\1     & 0   & 0        \\
    0           & 0         & u_1 & 0        \\
    0           & 0         & v_1 & u_1\1 \\
  \end{pmatrix} \mid u_1, v_1 \in \F_p \text{ und } u_1 \neq 0
  \right\}\\
  \bar{S}_4 = \bar{S}_5 & =  & \left\{ \begin{pmatrix}
    v_2         & 0         & 0    & 0        \\
    0           & \pm v_2^2 & 0    & 0        \\
    0           & 0         & \pm 1 & 0        \\
    0           & 0         & 0    & v_2 \\
  \end{pmatrix} \mid v_2 \in \F_p \text{ und } v_2 \neq 0
  \right\}\\
\end{eqnarray*}
\end{lemma}

\begin{beweis}
Die Behauptung ergibt sich unmittelbar aus \ref{bahnen_1_dp} und \ref{einvek}, da $A^\vi$
offensichtlich das Zentrum von $GL(4,p)$ enthält.

Die Stabilisatoren erhält man als Lösungsmengen der folgenden Gleichungssysteme
\begin{enumerate}
  \item $\{u_1=u_1, u_2=0\}$,
  \item $\{(u_1 v_2 - v_1 u_2)u_1=u_1, (u_1 v_2 - v_1 u_2)u_2=0\}$,
  \item $\{(u_1 v_2 - v_1 u_2)u_1=u_1, (u_1 v_2 - v_1 u_2)u_2=0, u_1=u_1, u_2=0\}$,
  \item $\{(u_1 v_2 - v_1 u_2)u_1=v_2, (u_1 v_2 - v_1 u_2)u_2=0, v_1=0, v_2=v_2\}$,
  \item $\{(u_1 v_2 - v_1 u_2)u_1=v_2, (u_1 v_2 - v_1 u_2)u_2=0, v_1=0, v_2=w v_2\}$.
\end{enumerate}
(Auf den rechten Seiten werden $u_1$ und $v_2$ der Einfachheit halber verwendet, anstatt neue
Variablen einzuführen).
\end{beweis}

\begin{lemma} \label{zug_dp_p4}
Für die Nachfolger von $G$ der Ordnung $p^4$ ergibt sich das folgende Vertretersystem
zulässiger Untergruppen:
\begin{enumerate}
  \item $M_2=\erz{a_5,a_6,a_7}$ entsprechend $\erz{(0,1,0,0),(0,0,1,0),(0,0,0,1)}$.
  \item $M_3=\erz{a_5,a_7,a_4 a_6\1}$ entsprechend $\erz{(0,1,0,0),(0,0,0,1),(1,0,-1,0)}$.
  \item $M_4=\erz{a_5,a_6,a_4 a_7\1}$ entsprechend $\erz{(0,1,0,0),(0,0,1,0),(1,0,0,-1)}$.
  \item $M_5=\erz{a_5,a_6,a_4 a_7^{-w}}$ entsprechend $\erz{(0,1,0,0),(0,0,1,0),(1,0,0,-w)}$.
\end{enumerate}
\end{lemma}

\begin{beweis}
Da $P(G)$ die Ordnung $p^7$ hat, werden für die unmittelbaren Nachfolger der Ordnung $p^4$ die
Bahnen der dreidimensionalen Unterräume von $\F_p^4 \cong M(G)$ unter der Operation von $A$
benötigt. Da $\F_p^4$ die Dimension $4$ hat, können statt der dreidimensionalen Unterräume
ihre eindimensionalen orthogonalen Komplemente und die Operationen von $A^T$ auf dieser Menge
betrachtet werden. Wegen der Symmetrie von $A^\vi$ erübrigt sich das Transponieren, d.\,h. es
ist $A^\vi=(A^\vi)^T$. Aus den Repräsentanten der Bahnen zu eindimensionalen Unterräumen in
\ref{bahnen_1_dpa} lassen sich die Vertreter der Bahnen der dreidimensionalen Unterräume
unmittelbar ablesen, indem man orthogonale Komplemente zu ihnen bildet. Mit der Bahn $B_1'$
sind alle Unterräume gegeben, denen keine zulässigen Untergruppen entsprechen, da das
Standardskalarprodukt der Elemente von $B_1'$ mit den Basisvektoren des Unterraumes
$\bar{N}=\erz{(1,0,0,0),(0,1,0,0)}$, der dem Nukleus entspricht, jeweils Null ist. Daher sind
die dreidimensionalen Komplemente der eindimensionalen Unterräume mit Basisvektoren aus $B_1$
nach keine Supplemente von $\bar{N}$ und ihnen entsprechen nach \ref{suppl} keine zulässigen
Untergruppen. Also ist $B_1$ zu vernachlässigen. Alle anderen Bahnen hingegen bestehen nach
demselben Kriterium aus zulässigen Untergruppen.
\end{beweis}

\begin{lemma}
Zu den zulässigen Untergruppen $M_2 \dd M_5$ gehören die folgenden
Stabilisatoren als Untergruppen von $Aut(G)$:
\begin{enumerate}
  \item[] $S_2=\{\al(u_1,u_2,u_3,0,v_2,v_3) \mid u_1, u_2, u_3, v_2, v_3
  \in \F_p \text{ und } u_1 v_2 \neq 0 \}$.
  \item[] $S_3=\{\al(u_1,u_2,u_3,0,u_1\1,v_3) \mid u_1, u_2, u_3, v_3
  \in \F_p \text{ und } u_1 \neq 0 \}$.
  \item[] $S_4=S_5=\{\al(\pm 1,0,u_3,0,v_2,v_3) \mid u_3, v_2, v_3
  \in \F_p \text{ und } v_2 \neq 0 \}$.
\end{enumerate}
\end{lemma}

\begin{beweis}
Die Stabilisatoren $S_2 \dd S_5$ sind aus den Stabilisatoren $\bar{S}_2 \dd
\bar{S}_5$ unter Berücksichtigung des Transponierens unmittelbar ablesbar.
\end{beweis}

\subsection{Nachfolger der Ordnung $p^4$}

\begin{satz} \label{ndpp4}
In der folgenden Liste sind sämtliche unmittelbaren Nachfolger von $G$ der
Ordnung $p^4$ bis auf Isomorphie angegeben:
\begin{enumerate}
  \item[] $G_1=\erz{g_1,g_2,g_3,g_4, \mid [g_2,g_1]=g_3, [g_3,g_1]=g_4}$
  \hfill \gap"=Typ $(p^4,7)$
  \item[] $G_2= \erz{g_1,g_2,g_3,g_4, \mid [g_2,g_1]=g_3,
  [g_3,g_1]=g_4, g_1^p=g_4}$ \hfill  \gap"=Typ $(p^4,8)$
  \item[] $G_3=\erz{g_1,g_2,g_3,g_4 \mid [g_2,g_1]=g_3,
  [g_3,g_1]=g_4, g_2^p=g_4}$ \hfill  \gap"=Typ $(p^4,9)$
  \item[] $G_4=\erz{g_1,g_2,g_3,g_4 \mid [g_2,g_1]=g_3,
  [g_3,g_1]=g_4, g_2^p=g_4^w}$ \hfill  \gap"=Typ $(p^4,10)$
\end{enumerate}
\end{satz}

\begin{beweis}
Dieser Satz ergibt sich nach \ref{bahnenzu}, indem man $P(G)$ nach jedem Repräsentanten des
Vertretersystems zulässiger Untergruppen faktorisiert, das in \ref{zug_dp_p4} aufgelistet ist.
Die Präsentation der Faktorgruppe erhält man über das Verfahren aus \ref{faktor}.
\end{beweis}

\subsection{Bahnen zulässiger Untergruppen für Nachfolger der Ordnung $p^5$}

\begin{lemma} \label{zug_dp_p5}
Für die unmittelbaren Nachfolger von $G$ der Ordnung $p^5$ ist die folgende
Liste ein Vertretersystem der zulässigen Untergruppen, wobei $w$ ein Erzeuger
der multiplikativen Gruppe von $\f$ ist:
\begin{enumerate}
  \item[] $M_1=\erz{a_6,a_7}$ entsprechend $\erz{(0,0,1,0),(0,0,0,1)}$,
  \item[] $M_2=\erz{a_5 a_7\1, a_6}$ entsprechend $\erz{(0,1,0,-1),(0,0,1,0)}$,
  \item[] $M_3=\erz{a_4 a_7\1, a_6}$ entsprechend $\erz{(1,0,0,-1),(0,0,1,0)}$,
  \item[] $M_4=\erz{a_4 a_7^{-w}, a_6}$ entsprechend $\erz{(1,0,0,-w),(0,0,1,0)}$,
  \item[] $M_5=\erz{a_4 a_6\1, a_5 a_7\1}$ entsprechend $\erz{(1,0,-1,0),(0,1,0,-1)}$,
  \item[] $M_6=\erz{a_4 a_6\1 a_7\1, a_5 a_7\1}$ entsprechend $\erz{(1,0,-1,-1),(0,1,0,-1)}$,
  \item[] $M_7=\erz{a_4 a_6\1 a_7^{-w}, a_5 a_7\1}$ entsprechend $\erz{(1,0,-1,-w),(0,1,0,-1)}$,
  \item[] $M_8=\erz{a_4 a_7^{-w}, a_5 a_6\1}$ entsprechend $\erz{(1,0,0,-w),(0,1,-1,0)}$,
  \item[] $M_9^k=\erz{a_4 a_6\1, a_5 a_7^{-w^k}}$ entsprechend $\erz{(1,0,-1,0),(0,1,0,-w^k)}$ mit $k \in
  \{1 \dd \frac{p-1}{2}\}$.
  \item[] $M_{10}^k=\erz{a_4 a_6\1 a_7^{-w^{k+1}}, a_5 a_6^{-w^k} a_7\1}$ entsprechend $\erz{(1,0,-1,-w^{k+1}),(0,1,-w^k,-1)}$ mit $k
  \in \{0 \dd \frac{p-1}{2}-1\}$.
\end{enumerate}
\end{lemma}

\begin{beweis}
Um die unmittelbaren Nachfolger von $G$ der Ordnung $p^5$ bis auf Isomorphie zu
klassifizieren, wird, da der Multiplikator von $G$ die Ordnung $p^7$ hat, ein Vertretersystem
der Bahnen unter Erweiterungsautomorphismen von $G$ von Untergruppen gesucht, die die Ordnung
$p^2$ haben und ein Supplement des Nukleus $N(G)=\erz{a_5, a_6}$ sind. Diese Problemstellung
wird nun in die Sprache der linearen Algebra übersetzt: Es sei $\{f_4, f_5, f_6, f_7\}$ die
Standardbasis von $\f^4$. Dann entspricht dem Nukleus der Unterraum $N=\erz{f_4, f_5}$. Ein
Unterraum $W$ von $\f$ entspricht nur dann einer zulässigen Untergruppe der Ordnung $p^2$,
wenn $W$ die Dimension $2$ hat und $dim(N \cap W)=0$ ist, denn andernfalls wäre $W$ kein
Supplement zu $N$. Entspricht $W$ einer zulässigen Untergruppe der Ordnung $p^2$, so wird $W$
innerhalb dieses Beweises zulässiger Unterraum genannt. Es sei $K=\{f_6, f_7\}$. Dann ist
$\f^4$ die direkte Summe von $N$ und $K$. Ist $W$ ein zulässiger Unterraum, so ist $dim(K \cap
W)=2$, $dim(K \cap W)=1$ oder $dim(K \cap W)=0$. Diese drei Fälle werden nun getrennt
untersucht.

Erster Fall: Es sei $dim(K \cap W)=2$. Dann ist $W=K$ und damit ist $\erz{f_6, f_7}$ ein
zulässiger Unterraum, der in der Liste als $M_1$ geführt wird.

Zweiter Fall: Es sei $dim(K \cap W)=1$. Da $\f^4$ die direkte Summe von $N$ und $K$ ist, gibt
es dann ein $n \in N$ und $k_1, k_2 \in K$ mit $W=\erz{n + k_1, k_2}$. Da $A$ über $\vi$ auf
der Menge der eindimensionalen Unterräume von $K$ transitiv operiert, kann man annehmen, dass
$k_2=f_6$ und dass $k_1 \notin U=\erz{f_6}$ ist, da $W$ andernfalls kein Supplement zu $N$
wäre. Die Gruppe $$S=\left \{\begin{pmatrix}
  u_1^2 v_2   & 0         & 0   & 0  \\
  u_1 v_1 v_2 & u_1 v_2^2 & 0   & 0  \\
  0           & 0         & u_1 & 0  \\
  0           & 0         & v_1 & v_2
\end{pmatrix} \mid u_1, v_2 \in \f^* \text{ und } v_1 \in \f\right\}$$ ist der
Stabilisator von $U$ in $A^\vi$. Nach \ref{tensor} gibt es eine Bijektion zwischen den Bahnen
der Komplemente zu $U$ in $\f^4$ unter der Operation von $S$ und den Bahnen der Elemente des
Tensorproduktes $\bar{M} \otimes U$ unter der Operation von $S$ vermittels $m^{s\1} \otimes
u^s$. Die Operation von $S$ auf dem Faktorraum $\bar{M}=\f^4/U$ vermittels $S \times \bar{M}:
(s,m) \mt m^{s\1}$ ist durch die Gruppe $$\bar{S}=\left \{\begin{pmatrix}
  u_1^2 v_2    & 0           & 0  \\
  -u_1 v_1 v_2 & u_1 v_2^2   & 0  \\
  0            & 0           & v_2
\end{pmatrix} \mid u_1, v_2 \in \f^* \text{ und } v_1 \in \f\right\}$$ gegeben.
Da nicht die Bahnen der dreidimensionalen Komplemente von $U$ selbst von Interesse sind,
sondern deren eindimensionalen orthogonalen Komplemente, wird statt $\bar{S}$ die aus
$\bar{S}$ durch Transponieren hervorgegangene Gruppe $$\hat{S}=\left \{\begin{pmatrix}
  u^2 v    & -u v z & 0  \\
  0        & u v^2  & 0  \\
  0        & 0      & v
\end{pmatrix} \mid u, v \in \f^* \text{ und } z \in \f\right\}$$ betrachtet (der Übersichtlichkeit
wegen werden andere Variablen verwendet). In \ref{opdp2} ist angegeben, wie $\f^3 \cong
\bar{M}$ unter$\hat{S}$ in Bahnen zerfällt. Da $W=\erz{n + k_1, f_6}$ genau dann ein
Komplement zu $N$ ist, wenn die Projektion von $W$ auf $U$ die Dimension 2 hat, und da $W$
genau dann einen eindimensionalen Schnitt mit $U$ hat, wenn die Projektion von $W$ auf $N$ die
Dimemsion 1 hat, ist $W$ genau dann ein zulässiger Unterraum, wenn sowohl $n$ als auch $k_1$
nicht der Nullvektor sind. Da einem in \ref{opdp2} aufgeführten Bahnenvertreter $(a,b,c)$ der
Unterraum $W=\erz{(a,b,0,c),(0,0,1,0)}= \erz{(a,b,0,c),(0,0,-1,0)}$ entspricht, enthalten
genau die Bahnen $B_1, B_2$ und $B_3$ aus \ref{opdp2} zulässige Untergruppen. Ihren Vertreter
entsprechen die Gruppen $M_2, M_3$ und $M_4$ in der oben aufgeführten Liste.

Dritter Fall: Es sei $dim(K \cap W)=0$. Da $\f^4$ die direkte Summe von $N$ und $K$ ist, gibt
es dann $n_1, n_2 \in N$ und $k_1, k_2 \in K$ mit $W=\erz{n_1 + k_1, k_2 + k_2}$. Da $N$ unter
der Operation von $A$ vermittels $\vi$ invariant ist und $A$ über $\vi$ auf $N$ transitiv
operiert, kann man annehmen, dass $n_1=f_4$ und $n_2=f_5$ erfüllt ist. Es sei $U=\erz{f_6,
f_7}$ und $S=A^\vi$. Nach \ref{tensor} entsprechen die Bahnen der zulässigen Unterräume den
Bahnen von $S$ in $N \otimes U$ unter der Operation $\theta: S \times (N \otimes U) \ra N
\otimes U: (s, n \otimes u) \mt n^{s\1} \otimes u^s$. Diese Operation wiederum entspricht
offensichtlich der Operation von $H=GL(2,p)$ auf $M=M(2 \times 2,p)$ vermittels $\delta: H
\times M \ra M: (h,m) \mt det(h)\1 h\1 m h$, wobei statt $\delta$ ebenso die Operation $\nu: H
\times M \ra M: (h,m) \mt det(h)h\1 m h$ verwendet werden kann. Die Bahnen von $M$ unter der
Operation von $H$ über $\nu$ sind in \ref{det} festgehalten. Da einem der dort aufgeführten
Bahnenrepräsentanten $$m=\begin{pmatrix}
  a   & b \\
  c   & d
\end{pmatrix}$$ der Unterraum $W=\erz{(1,0,a,b),(0,1,c,d)}$ entspricht, ist
$W$ genau dann ein zulässiger Unterraum, wenn $det(m) \neq 0$ ist, da
andernfalls $W$ kein Supplement zu $N$ wäre. Die zulässigen Unteräume dieser
Art sind in der oben angeführten Liste als $M_5 \dd M_{10}^k$ zu finden.
\end{beweis}

\subsection{Nachfolger der Ordnung $p^5$}

\begin{satz} \label{ndpp5}
In der folgenden Liste sind sämtliche unmittelbaren Nachfolger
von $G$ der Ordnung $p^5$ bis auf Isomorphie angegeben.
\begin{enumerate}
  \item[] $H_1=\erz{h_1,h_2,h_3,h_4,h_5 \mid [h_2,h_1]=h_3,
  [h_3,h_1]=h_4, [h_3,h_2]=h_5}$,
  \item[] $H_2=\erz{h_1,h_2,h_3,h_4,h_5 \mid [h_2,h_1]=h_3,
  [h_3,h_1]=h_4, [h_3,h_2]=h_5, h_2^p=h_5}$,
  \item[] $H_3=\erz{h_1,h_2,h_3,h_4,h_5 \mid [h_2,h_1]=h_3,
  [h_3,h_1]=h_4, [h_3,h_2]=h_5, h_2^p=h_4}$,
  \item[] $H_4=\erz{h_1,h_2,h_3,h_4,h_5 \mid [h_2,h_1]=h_3,
  [h_3,h_1]=h_4, [h_3,h_2]=h_5, h_2^p=h_4^w}$,
  \item[] $H_5=\erz{h_1,h_2,h_3,h_4,h_5 \mid [h_2,h_1]=h_3,
  [h_3,h_1]=h_4, [h_3,h_2]=h_5, h_1^p=h_4, h_2^p=h_5}$,
  \item[] $H_6=\erz{h_1,h_2,h_3,h_4,h_5 \mid [h_2,h_1]=h_3,
  [h_3,h_1]=h_4 h_5, [h_3,h_2]=h_5, h_1^p=h_4, h_2^p=h_5}$,
  \item[] $H_7=\erz{h_1,h_2,h_3,h_4,h_5 \mid [h_2,h_1]=h_3,
  [h_3,h_1]=h_4 h_5^w, [h_3,h_2]=h_5, h_1^p=h_4, h_2^p=h_5}$,
  \item[] $H_8=\erz{h_1,h_2,h_3,h_4,h_5 \mid [h_2,h_1]=h_3,
  [h_3,h_1]=h_5^w, [h_3,h_2]=h_4, h_1^p=h_4, h_2^p=h_5}$,
  \item[] $H_9^k=\erz{h_1,h_2,h_3,h_4,h_5 \mid [h_2,h_1]=h_3,
  [h_3,h_1]=h_4, [h_3,h_2]=h_5^{w^k}, h_1^p=h_4, h_2^p=h_5}$ mit $k \in
  \{1 \dd \frac{p-1}{2}\}$,
  \item[] $H_{10}^k=\erz{h_1,h_2,h_3,h_4,h_5 \mid [h_2,h_1]=h_3,
  [h_3,h_1]=h_4 h_5^{w^{k+1}}, [h_3,h_2]=h_4^{w^k} h_5, h_1^p=h_4, h_2^p=h_5}$ mit $k
  \in \{0 \dd \frac{p-1}{2}-1\}$.
\end{enumerate}
\end{satz}

\begin{beweis}
Dieser Satz ergibt sich nach \ref{bahnenzu}, indem man $P(G)$ nach jedem Repräsentanten des
Vertretersystems zulässiger Untergruppen faktorisiert, das in \ref{zug_dp_p5} ist. Die
Präsentation der Faktorgruppe erhält man über das Verfahren aus \ref{faktor}.
\end{beweis}

\begin{folg}
Die Gruppe $G$ hat $8+p$ Isomorphieklassen unmittelbarer Nachfolger der Ordnung $p^5$.
\end{folg}

\section{Nachfolger von $(p^4,7)$}

\begin{ver}
In diesem Abschnitt sei $G=\erz{a_1,a_2,a_3,a_4 \mid [a_2,a_1]=a_3, [a_3,a_1]=a_4}$.
\end{ver}

\begin{lemma}
Die Gruppe $G$ ist durch $\omega(a_1)=\omega(a_2)=1$, $\omega(a_3)=2$ und $\omega(a_4)=3$
gewichtet und hat damit die $p$"=Klasse 3.
\end{lemma}

\begin{beweis}
Die Gewichtung lässt sich nach \ref{gew} ermitteln.
\end{beweis}

\subsection{Die Automorphismengruppe}

\begin{lemma}
Jeder Automorphismus $\al$ von $G$ lässt sich eindeutig durch $$\al(u_1, u_2, u_3, u_4, v_2,
v_3, v_4):
  \begin{cases} \label{aut7}
    a_1 \mt a_1^{u_1} a_2^{u_2} a_3^{u_3} a_4^{u_4} & \text{wobei } u_1, u_2, u_3, u_4, v_2, v_4 \\
    a_2 \mt a_2^{v_2} a_3^{v_3} a_4^{v_4} & \in \{0 \dd p-1\} \text{ und }\\
    a_3 \mt a_3^{u_1 v_2} a_4^{u_1 v_3 + (u_1 (u_1-1) v_2)/2} & u_1 v_2 \neq 0 \text{ ist.}\\
    a_4 \mt a_4^{u_1^2 v_2} &
  \end{cases}$$ darstellen.
\end{lemma}

\begin{beweis}
Die Automorphismengruppe von $G$ ist nach \ref{aut} aus dem Stabilisator $\bar{S}_2$ (siehe
\ref{bahnen_1_dpa}) der zu $G$ gehörenden zulässigen Untergruppe unmittelbar ablesbar. Damit
ergeben sich genau die oben angegebenen Bilder von $a_1$ und $a_2$. Die Bilder von $a_3$ und
$a_4$ sind durch die Bilder von $a_1$ und $a_2$ eindeutig bestimmt und lassen sich
folgendermaßen berechnen:
\begin{eqnarray*}
  a_3^\al &=& [a_2,a_1]^\al=[a_2^\al,a_1^\al]=[a_2^{v_2} a_3^{v_3} a_4^{v_4}, a_1^{u_1} a_2^{u_2} a_3^{u_3} a_4^{u_4}]\\
          &=& [a_2^{v_2} a_3^{v_3}, a_1^{u_1} a_2^{u_2} a_3^{u_3}]\\
          &=& (a_2^{v_2} a_3^{v_3})\1 (a_1^{u_1} a_2^{u_2} a_3^{u_3})\1 a_2^{v_2} a_3^{v_3} a_1^{u_1} a_2^{u_2} a_3^{u_3}\\
          &=& a_3^{-v_3} a_2^{-v_2} a_3^{-u_3} a_2^{-u_2} a_1^{-u_1} a_2^{v_2} a_3^{v_3} a_1^{u_1} a_2^{u_2} a_3^{u_3}\\
          &=& a_2^{-v_2 -u_2} a_3^{-v_3 -u_3} a_1^{-u_1} a_2^{v_2} a_3^{v_3} a_1^{u_1} a_2^{u_2} a_3^{u_3}\\
          &=& a_2^{-v_2} a_3^{-v_3} a_1^{-u_1} a_2^{v_2} a_3^{v_3}  a_1^{u_1}\\
          &=& a_3^{u_1 v_2} a_4^{u_1 v_3 + (u_1 (u_1-1) v_2)/2},
\end{eqnarray*}
da $\erz{a_4}$ das Zentrum von $G$ ist. Außerdem ist
\begin{eqnarray*}
  a_4^\al &=& [a_3,a_1]^\al=[a_3^\al, a_1^\al]\\
          &=& [a_3^{u_1 v_2} a_4^{u_1 v_3 + (u_1 (u_1-1) v_2)/2}, a_1^{u_1} a_2^{u_2} a_3^{u_3} a_4^{u_4}]\\
          &=& [a_3^{u_1 v_2}, a_1^{u_1} a_2^{u_2} a_3^{u_3}]\\
          &=& [a_3^{u_1 v_2}, a_1^{u_1}]^{a_3^{u_1 v_2}} [a_3^{u_1 v_2}, a_2^{u_2} a_3^{u_3}]\\
          &=& [a_3^{u_1 v_2}, a_1^{u_1}][a_3^{u_1 v_2}, a_2^{u_2} a_3^{u_3}]\\
          &=& [a_3^{u_1 v_2}, a_1^{u_1}][a_3^{u_1 v_2}, a_3^{u_3}] [a_3^{u_1 v_2}, a_2^{u_2}]^{a_3^{u_3}}\\
          &=& [a_3^{u_1 v_2}, a_1^{u_1}]=[a_3, a_1]^{u_1^2 v_2}=a_4^{u_1^2 v_2},
\end{eqnarray*} da $\erz{a_4}$  das Zentrum von $G$ ist.
\end{beweis}

\subsection{$p$"=Cover, Multiplikator und Nukleus}

\begin{lemma}
Für $G$ ergibt sich das $p$"=Cover, der Nukleus und der Multiplikator in folgender Weise:
\begin{enumerate}
  \item[] $P(G)=\erz{a_1,a_2,a_3,a_4,a_5,a_6,a_7,a_8 \mid [a_2,a_1]=a_3, [a_3,a_1]=a_4,
  [a_4,a_1]=a_5, [a_3,a_2]=a_6, a_1^p=a_7, a_2^p=a_8}$.
  \item[] $M(G)=\erz{a_5, a_6, a_7, a_8}$.
  \item[] $N(G)=\erz{a_5}$.
\end{enumerate}
\end{lemma}

\begin{beweis}
Die Präsentation von $P(G)$ ergibt sich aus \ref{cov} und \ref{newman} sowie dem
Reduktionsverfahren von Knuth"=Bendix, aus der sich die Untergruppe $M(G)$ ablesen lässt. Die
Gewichtung von $P(G)$ und damit die Untergruppe $N(G)$ lassen sich nach \ref{gew} aus der
Präsentation von $P(G)$ ablesen:
\begin{enumerate}
  \item[] $1=\omega(a_1)=\omega(a_2)$.
  \item[] $2=\omega(a_3)=\omega(a_7)=\omega(a_8)$.
  \item[] $3=\omega(a_4)=\omega(a_6)$.
  \item[] $4=\omega(a_5)$.
\end{enumerate}
\end{beweis}

\subsection{Operation der Erweiterungsautomorphismen}

\begin{lemma}
Nach \ref{aut7} lässt sich jeder Automorphismus $\al$ von $G$ in der Form $$\al(u_1, u_2, u_3,
u_4, v_2, v_3, v_4):
  \begin{cases}
    a_1 \mt a_1^{u_1} a_2^{u_2} a_3^{u_3} a_4^{u_4} & \text{wobei } u_1, u_2, u_3, u_4, v_2, v_4 \\
    a_2 \mt a_2^{v_2} a_3^{v_3} a_4^{v_4} & \in \{0 \dd p-1\} \text{ und }\\
    a_3 \mt a_3^{u_1 v_2} a_4^{u_1 v_3 + (u_1 (u_1-1) v_2)/2} & u_1 v_2 \neq 0 \text{ ist.}\\
    a_4 \mt a_4^{u_1^2 v_2} &
  \end{cases}$$ darstellen. In dieser Notation ist der Operationshomomorphismus $\vi$ von $Aut(G)$
über Erweiterungsautomorphismen auf $\f^4 \cong M(G)$ folgendermaßen gegeben $$\vi: Aut(G) \to
GL(4,p): \al(u_1, u_2, u_3, u_4, v_2, v_3, v_4) \mt M
=
 \begin{pmatrix}
  u_1^3 v_2  & 0          & 0          & 0  \\
  0          & u_1 v_2^2  & 0          & 0  \\
  0          & 0          & u_1        & u_2  \\
  0          & 0          & 0          & v_2
\end{pmatrix}$$
\end{lemma}

\begin{beweis}
Die Bilder der Erzeuger von $M(G)$ ergeben sich unter den Erweiterungsautomorphismen
folgendermaßen:
\begin{eqnarray*}
  a_5^\al &=& [a_4,a_1]^\al=[a_4^\al,a_1^\al]\\
          &=& [a_4^{u_1^2 v_2},a_1^{u_1} a_2^{u_2} a_3^{u_3} a_4^{u_4}]\\
          &=& [a_4^{u_1^2 v_2},a_1^{u_1} a_2^{u_2} a_3^{u_3}]\\
          &=& [a_4^{u_1^2 v_2},a_1^{u_1}]=[a_4, a_1]^{u_1^3 v_2}\\
          &=& a_5^{u_1^3 v_2},
\end{eqnarray*}
da $[a_4,a_1]$ zentral und alle anderen auftretenden Kommutatoren trivial sind.
\begin{eqnarray*}
  a_6^\al &=& [a_3,a_2]^\al=[a_3^\al,a_2^\al]\\
          &=& [a_3^{u_1 v_2} a_4^d, a_2^{v_2} a_3^{v_3} a_4^{v_4}]=[a_3^{u_1 v_2}, a_2^{v_2} a_3^{v_3}]\\
          &=& [a_3^{u_1 v_2}, a_2^{v_2}]=[a_3^{u_1 v_2}, a_2^{v_2}]^{u_1 v_2^2}\\
          &=& a_6^{u_1 v_2^2}
\end{eqnarray*}
aus denselben Gründen, wobei $d=u_1 v_3 + (u_1 (u_1-1) v_2)/2$ ist.
\begin{eqnarray*}
  a_7^\al &=& (a_1^p)^\al=(a_1^\al)^p\\
          &=& (a_1^{u_1} a_2^{u_2} a_3^{u_3} a_4^{u_4})^p\\
          &=& (a_1^{u_1})^p (a_2^{u_2})^p = (a_1^p)^{u_1} (a_2^p)^{u_2}\\
          &=& a_7^{u_1} a_8^{u_2}
\end{eqnarray*}
da die Kommutatoren, die aus $a_1, a_2, a_3$ und $a_4$ gebildet sind, trivial
sind oder die Ordnung $p$ haben.
\begin{eqnarray*}
  a_8^\al &=& (a_2^p)^\al=(a_2^\al)^p\\
          &=& (a_2^{v_2} a_3^{v_3} a_4^{v_4})^p\\
          &=& (a_2^{v_2})^p=(a_2^p)^{v_2}\\
          &=& a_8^{v_2}
\end{eqnarray*}
aus demselben Grund.
\end{beweis}

Da für die Nachfolger der Ordnung $p^5$ nicht die dreidimensionalen Unterräume
von $\F_p^4 \cong M(G)$ betrachtet werden, sondern ihre eindimensionalen
Komplemente, wird im weiteren nicht $\vi$, sondern der zu $\vi$ duale
Operationshomomorphismus $\bar{\vi}$ betrachtet, d.\,h. die Operation von
$Aut(G)$ auf dem Dualraum zu $\F_p^4$.

\begin{folg}
Die Gruppe $Aut(G)$ operiert durch $$\bar{\vi}: Aut(G) \to GL(4,p): \al(u_1, u_2, u_3, u_4,
v_2, v_3, v_4) \mt M
=
 \begin{pmatrix}
  u_1^3 v_2  & 0          & 0          & 0  \\
  0          & u_1 v_2^2  & 0          & 0  \\
  0          & 0          & u_1        & 0  \\
  0          & 0          & u_2        & v_2
\end{pmatrix}$$ über Erweiterungsautomorphismen auf dem Dualraum von $\F_p^4$
\end{folg}

\begin{beme}
Die Matrixgruppe $$L=\{\al^{\bar{\vi}} \mid \al \in Aut(G)\}=\left\{
\begin{pmatrix}
  u^3 v      & 0          & 0          & 0  \\
  0          & u v^2      & 0          & 0  \\
  0          & 0          & u          & 0  \\
  0          & 0          & z          & v
\end{pmatrix} \mid u,v,z \in \F_p \text{ und } uv \neq 0
 \right\}$$ operiert nicht derart auf $\F_p^4$, dass die Bahnen der Vektoren von
$\F_p^4 \backslash \{0\}$ in bijektiver Korrespondenz zu den Bahnen der
eindimensionalen Unterräume von $\F_p^4$ stehen.
\end{beme}

\begin{beweis}
Die Gruppe $L$ enthält nicht das volle Zentrum $Z$ von $GL(4,p)$. Denn $Z$ wird von $w E_4$
erzeugt, wobei $w$ ein Erzeuger der multiplikativen Gruppe von $\F_p$ und $E_4$ die
Einheitsmatrix der Dimension $4 \times 4$ ist. Die Gruppe $L$ enthält aber $w E_4$ nicht, da
in diesem Fall $u_1=w$ und $v_2=w$ erfüllt sein müssten. Dann aber wäre auch $w=u_1^3 v_2= w^3
w= w^4$ und $w=u_1 v_2^2=w w^2=w^3$. Nach letzterem wäre $w$ ein Quadrat, was $w$ aber nach
\ref{indexmult} nicht ist (sofern $p \neq 2$ ist --~aber das wird in diesem Kapitel generell
vorausgesetzt). Nach \ref{zentin} lassen sich daher die Bahnen der eindimensionalen Unterräume
unter der Operation von $L$ auf $\F_p^4$ nicht unmittelbar aus den Bahnen der Vektoren
ermitteln.
\end{beweis}

\begin{folg}
Unter der Operation der Gruppe $$A=\left\{
\begin{pmatrix}
  k u^3 v    & 0         & 0        & 0  \\
  0          & k u v^2   & 0        & 0  \\
  0          & 0         & k u      & 0  \\
  0          & 0         & z        & k v
\end{pmatrix} \mid  k,u,v,z \in \F_p \text{ und } kuv \neq
0\right\}$$ auf $\F_p^4$ zerfällt $\F_p^4$ derart in Bahnen, dass sich aus
ihnen die Bahnen der eindimensionalen Unterräume ablesen lassen und dass sie
den Bahnen der eindimensionalen Unterräume im Sinne von \ref{einvek}
entsprechen, die sich unter der Operation von $L$ ergeben.
\end{folg}

\begin{beweis}
Die Behauptung ergibt sich unmittelbar unter Bezug auf \ref{einvek}.
\end{beweis}

\begin{lemma}
Die Gruppe $A$ hat die Ordnung $p(p-1)^3$.
\end{lemma}

\begin{beweis}
Die Ordnung von $A$ wird folgendermaßen bestimmt: Da $A=LZ$ ist, gilt nach den
Homomorphiesätzen $$|LZ|=\frac{|L| |Z|}{|L \cap Z|}.$$ Die Ordnung von $Z$ ist
bekanntlich $p-1$. Da $Z$ von $w E_4$ erzeugt wird und damit $Z=\{w^a E_4 \mid
a \in \{0 \dd p-1\}\}$ ist, kann man den Schnitt von $L$ und $Z$ folgendermaßen
bestimmen: Man erhält, dass für die Elemente des Schnittes $u=w^a$ und $v=w^a$
gilt und daher auch $w^a=u^3 v = w^{3a} w^a = w^{4a}$. Also ist $w^{3a}=1$ bzw.
$3a \equiv 0 \text{ mod } p-1$. Diese Bedingung ist nach \ref{indexadd} bei $3
\nmid p-1$ nur für $a=0$ und andernfalls für jedes $a \in \{0, \frac{1}{3}
(p-1), \frac{2}{3} (p-1)\}$ erfüllt. Andererseits ist auch $w^a=uv^2= w^a
w^{2a}=w^{3a}$, was unabhängig von $p$ nur für $a \in \{0, \frac{1}{2} (p-1)\}$
erfüllt ist. Aus beiden Bedingungen erhält man, dass $a=0$ und damit $|L \cap
Z|=1$ ist.

Als nächster Schritt wird die Ordnung von $L$ bestimmt. Dazu werden die Bahn $B$ und der
Stabilisator $S$ von $(1,0,0,1)$ unter $L$ ermittelt und benutzt, dass $|L|=|B| \cdot |S|$
ist. Die Bahn ist $$B=(1,0,0,1)^L=\{(u^3v,0,z,v) \mid u,v,z \in \f \text{ und } uv \neq 0\}.$$
Da die Abbildung $\theta: \F_p \ra \F_p: x \mt x^3$ nach \ref{indexmult} die multiplikative
Gruppe $\f^*$ auf eine Untergruppe $(\f^*)^\theta \leq \f^*$ vom Index $ggT(3, p-1)$ abbildet,
ist $$|B|=
  \begin{cases}
    p(p-1)^2             & \text{ falls } 3 \nmid p-1, \\
    \frac{1}{3} p(p-1)^2 & \text{ falls } 3 \mid p-1.
  \end{cases}$$ Ist $W_3=\{x \mid x \in \f \text{ und }x^3=1\}$, so ist der
Stabilisator von $(1,0,0,1)$ $$S=\left\{ \begin{pmatrix}
  u^3        & 0          & 0          & 0  \\
  0          & u          & 0          & 0  \\
  0          & 0          & u          & 0  \\
  0          & 0          & 0          & 1
\end{pmatrix} \mid u \in W_3 \right\}.$$ Wiederum nach \ref{indexmult} ergibt
sich, dass $$|S|=
  \begin{cases}
    1      & \text{ falls } 3 \nmid p-1, \\
    3      & \text{ falls } 3 \mid p-1
  \end{cases}$$ ist. Man erhält also insgesamt, dass die
Ordnung von $L$ für jede Primzahl $p > 3$ gleich $p(p-1)^2$ ist. Damit ist $|Z|
\cdot |L| = p(p-1)^3$ die Ordnung von $A$.
\end{beweis}

\subsection{Bahnen zulässiger Untergruppen für Nachfolger der Ordnung $p^4$}

\begin{lemma} \label{bahnen_7}
Der Vektorraum $\f^4$ zerfällt unter der Operation von $A$ in die folgenden Bahnen
nichttrivialer Vektoren mit der jeweils angegebenen Länge.
\begin{itemize}
  \item[] $B_1=(1,0,0,0)^A$ mit $|B_1|=p-1$.
  \item[] $B_2=(1,0,0,1)^A$ mit $|B_2|=p(p-1)^2$, falls $ggT(3, p-1)=1$ ist.
  Wenn hingegen $ggT(3, p-1)=3$ ist, so zerfällt die Bahn $B_2$ in die Bahnen
  $B_{2a}, B_{2b}$ und $B_{2c}$.
  \item[] $B_{2a}=(1,0,0,w^0)^A$ mit $|B_2|=\frac{1}{3} p(p-1)^2$.
  \item[] $B_{2b}=(1,0,0,w^1)^A$ mit $|B_2|=\frac{1}{3} p(p-1)^2$.
  \item[] $B_{2c}=(1,0,0,w^2)^A$ mit $|B_2|=\frac{1}{3} p(p-1)^2$.
  \item[] $B_3=(1,0,1,0)^A$ mit $|B_3|=(p-1)^2$.
  \item[] $B_4=(1,1,0,0)^A$ mit $|B_4|=(p-1)^2$.
  \item[] $B_5=(1,1,0,1)^A$ mit $|B_5|=p(p-1)^3$, falls $ggT(3, p-1)=1$ ist. Für
  $ggT(3, p-1)=3$ zerfällt $B_5$ in die folgenden drei Bahnen.
  \item[] $B_{5a}=(1,1,0,w^0)^A$ mit $|B_{5a}|=\frac{1}{3} p(p-1)^3$
  \item[] $B_{5b}=(1,1,0,w^1)^A$ mit $|B_{5b}|=\frac{1}{3} p(p-1)^3$
  \item[] $B_{5c}=(1,1,0,w^2)^A$ mit $|B_{5c}|=\frac{1}{3} p(p-1)^3$
  \item[] $B_{6a}=(1,1,w^0,0)^A$ mit $|B_{6a}|=\frac{1}{2} (p-1)^3$, wenn
  $ggT(4, p-1)=2$ ist, und $|B_{6a}|=\frac{1}{4} (p-1)^3$, wenn
  $ggT(4, p-1)=4$ ist.
  \item[] $B_{6b}=(1,1,w^1,0)^A$ mit $|B_{6b}|=\frac{1}{2} (p-1)^3$, wenn
  $ggT(4, p-1)=2$ ist, und $|B_{6b}|=\frac{1}{4} (p-1)^3$, wenn
  $ggT(4, p-1)=4$ ist. Wenn $ggT(4, p-1)=4$ ist, gibt es außerdem die Bahnen
  $B_{6c}$ und $B_{6d}$.
  \item[] $B_{6c}=(1,1,w^2,0)^A$ mit $|B_{6c}|=\frac{1}{4} (p-1)^3$.
  \item[] $B_{6d}=(1,1,w^3,0)^A$ mit $|B_{6d}|=\frac{1}{4} (p-1)^3$.
  \item[]  $B_7=(0,0,0,1)^A$ mit $|B_7|=p(p-1)$.
  \item[]  $B_8=(0,0,1,0)^A$ mit $|B_8|=p-1$.
  \item[]  $B_9=(0,1,0,0)^A$ mit $|B_9|=p-1$.
  \item[] $B_{10}=(0,1,0,1)^A$ mit $|B_{10}|=p(p-1)^2$.
  \item[] $B_{11a}=(0,1,1,0)^A$ mit $|B_{11a}|=\frac{1}{2}(p-1)^2$.
  \item[] $B_{11b}=(0,1,w,0)^A$ mit $|B_{11b}|=\frac{1}{2}(p-1)^2$.
\end{itemize}
Es sei $W_2=\{x \in \f \mid x^2=1\}$, $W_3=\{x \in \f \mid x^3=1\}$ und
$W_4=\{x \in \f \mid x^4=1\}$. Die zu den Bahnenvertretern gehörenden
Stabilisatoren sind
\begin{eqnarray*}
  \hat{S}_1 & = & \left\{ \begin{pmatrix}
    1         & 0        & 0           &  0        \\
    0         & u^{-2}   & 0           &  0        \\
    0         & 0        & u^{-2} v\1  &  0        \\
    0         & 0        & z           &  u^{-3}    \\
  \end{pmatrix} \mid u, v, z \in \F_p  \text{ und } u,v \neq 0\right\}
  \text{ mit } |\hat{S}_1|=p(p-1)^2.\\
  \hat{S}_2 & = & \left\{ \begin{pmatrix}
    u^3       & 0        & 0      &  0        \\
    0         & u v      & 0      &  0        \\
    0         & 0        & u v\1  &  0      \\
    0         & 0        & 0      &  1 \\
  \end{pmatrix} \mid u \in W_3 \text{ und } v \in \f^* \right\}
  \text{ mit } |\hat{S}_2|=|W_3|(p-1).\\
  \hat{S}_3 & = & \left\{ \begin{pmatrix}
    1         & 0        & 0      &  0        \\
    0         & u^{-4}   & 0      &  0        \\
    0         & 0        & 1      &  0      \\
    0         & 0        & z      &  u^{-3}        \\
  \end{pmatrix} \mid u, z \in \F_p^* \text{ und } u \neq 0 \right\}
  \text{ mit } |\hat{S}_3|=p(p-1).\\
    \hat{S}_4 & = & \left\{ \begin{pmatrix}
    1         & 0        & 0      &  0        \\
    0         & 1        & 0      &  0        \\
    0         & 0        & u^{-4} &  0        \\
    0         & 0        & z      &  u^{-3}    \\
  \end{pmatrix} \mid u,z \in \F_p \text{ und } u \neq 0 \right\}
  \text{ mit } |\hat{S}_5|=p(p-1).\\
    \hat{S}_5 & = & \left\{ \begin{pmatrix}
    u^3          & 0        & 0        &  0        \\
    0            & u^3      & 0        &  0        \\
    0            & 0        & u\1      &  0        \\
    0            & 0        & 0        &  1   \\
  \end{pmatrix} \mid u \in W_3 \right\}
  \text{ mit } |\hat{S}_5|=|W_3|.\\
    \hat{S}_6 & = & \left\{ \begin{pmatrix}
    u^4       & 0        & 0      &  0        \\
    0         & 1        & 0      &  0        \\
    0         & 0        & 1      &  0        \\
    0         & 0        & z      &  u    \\
  \end{pmatrix} \mid z \in \F_p \text{ und } u \in W_4 \right\}
  \text{ mit } |\hat{S}_6|=p|W_4|.\\
    \hat{S}_7 & = & \left\{ \begin{pmatrix}
    u^3       & 0        & 0      &  0        \\
    0         & u v      & 0      &  0        \\
    0         & 0        & u v\1  &  0        \\
    0         & 0        & 0      &  1    \\
  \end{pmatrix} \mid u,v \in \f^*\right\}
  \text{ mit } |\hat{S}_7|=(p-1)^2.\\
  \hat{S}_{8} & = & \left\{ \begin{pmatrix}
    u^2 v     & 0        & 0      &  0        \\
    0         & v^2      & 0      &  0        \\
    0         & 0        & 1      &  0      \\
    0         & 0        & z      &  u\1 v        \\
  \end{pmatrix} \mid u, v \in \F_p^*  \text{ und } z \in \f^* \right\}
  \text{ mit } |\hat{S}_{8}|=p(p-1)^2.\\
  \hat{S}_{9} & = & \left\{ \begin{pmatrix}
    u^2 v\1   & 0        & 0      &  0        \\
    0         & 1        & 0      &  0        \\
    0         & 0        & v^{-2} &  0      \\
    0         & 0        & z      &  u\1 v\1        \\
  \end{pmatrix} \mid u, v \in \F_p^*  \text{ und } z \in \f^* \right\}
  \text{ mit } |\hat{S}_{9}|=p(p-1)^2.\\
  \hat{S}_{10} & = & \left\{ \begin{pmatrix}
    u^3       & 0        & 0      &  0        \\
    0         & 1        & 0      &  0        \\
    0         & 0        & u^2    &  0      \\
    0         & 0        & 0      &  1        \\
  \end{pmatrix} \mid u \in \F_p^* \right\}
  \text{ mit } |\hat{S}_{10}|=(p-1).\\
  \hat{S}_{11} & = & \left\{ \begin{pmatrix}
    u^2 v     & 0        & 0      &  0        \\
    0         & v^2      & 0      &  0        \\
    0         & 0        & 1      &  0      \\
    0         & 0        & z      &  u\1 v        \\
  \end{pmatrix} \mid u \in \F_p^*, z \in \f \text{ und } v \in W_2 \right\}
  \text{ mit } |\hat{S}_{11}|=|W_2|p(p-1).\\
\end{eqnarray*}
\end{lemma}

\begin{beweis}
Die Längen der Bahnen sind größtenteils unmittelbar ersichtlich oder können aus der Länge der
zugehörigen Stabilisatoren nach den Bahn"=Stabilisator"=Sätzen ermittelt werden. Im Weiteren
werden nur die Fälle behandelt, die nicht allzu offensichtlich sind. Dieser Nachweis betrifft
neben $B_{11}$ die Bahnen, die sich in Abhängigkeit von $p$ "`auffächern"' --~also die Bahnen
$B_2, B_5$ und $B_6$.

Betrachtet man die Bahn $B_2$, die bereits zur Ermittlung der Ordnung von $A$ benutzt worden
ist, so ergibt sich, dass der Vektor $(1,0,0,w^b)$ genau dann in der Bahn des Vektors
$(1,0,0,w^a)$ liegt, wenn die Gleichungen $1=ku^3v$ und $w^b = w^a k v$ erfüllt sind. In
diesem Fall erhält man über $w^{b-a}=kv$ und $1=ku^3 = kvu^3=w^{b-a}$, dass dann $u^3=w^{a-b}$
und damit $a \equiv b$ mod $3$ ist. Nach \ref{indexmult} gibt es genau $ggT(3, p-1)$
Nebenklassen modulo $w^3$ in $\f^*$ und damit liegt im Falle $ggT(3, p-1)=1$ jeder Vektor
$(1,0,0,w^b)$ in der Bahn von $(1,0,0,w^a)$ und im Falle $ggT(3, p-1)=3$ der Vektor
$(1,0,0,w^b)$ genau dann, wenn $a \equiv b$ mod $3$ ist. Also gibt es für $ggT(3, p-1)=1$
genau eine Bahn mit dem Vertreter $(1,0,0,1)$ und für $ggT(3, p-1)=3$ drei Bahnen mit den
Vertretern $(1,0,0,w^0)$, $(1,0,0,w^1)$ und $(1,0,0,w^2)$. Dasselbe Argument trifft auf das
Verhältnis der Bahn $B_5$ zu den Bahnen $B_{5a}, B_{5b}$ und $B_{5c}$ zu.

Im Fall der Bahn $B_6$ ist $(1,1,w^b,0)$ genau dann ein Element der Bahn von $(1,1,w^a,0)$,
wenn die Gleichungen $ku=w^{b-a}$, $1=w^{b-a}v^2$ und $1=w^{b-a}u^2v$ erfüllt sind. In diesem
Fall ist $v=w^{a-b}u^{-2}$ und daher $1=w^{b-a}u^{-4}w^{a-b}=u^{-4}$, d.\,h. $1=u^4$. Also ist
$v^2=u^2v$ und damit $v=u^2$. Man erhält somit $u^4=1=w^{b-a}u^2v=w^{b-a}$. Damit liegen die
Vektoren $(1,1,w^a,0)$ und $(1,1,w^a,0)$ genau dann in derselben Bahn, wenn $a \equiv b$ mod
$4$ ist. Nach \ref{indexmult} gibt es $ggT(4,p-1)$ Äquivalenzklassen in $\f^*$ modulo $w^4$
und damit können als Vertreter dieser Bahnen $(1,1,w^0,0)$ und $(1,1,w^1,0)$ im Falle
$ggT(4,p-1)=2$ und $(1,1,w^0,0) \dd (1,1,w^3,0)$ im Falle $ggT(4, p-1)=4$ gewählt werden. Ein
analoges Argument gilt für die Bahnen $B_{11a}$ und $B_{11b}$ --~allerdings mit dem
Unterschied, dass die dort relevante Bedingung $ggT(2, p-1)=2$ für jede Primzahl $p>3$ erfüllt
ist.
\end{beweis}

\begin{folg} \label{zug_7}
Für die Nachfolger von $G$ der Ordnung $p^5$ ist die folgende Liste ein
Vertretersystem zulässiger Untergruppen. Dabei ist $W_3=\{x \in \f \mid
x^3=1\}$ und $W_4=\{x \in \f \mid x^4=1\}$.
\begin{enumerate}
  \item[] $M_1=\erz{a_6,a_7,a_8}$ entsprechend $\erz{(0,1,0,0),(0,0,1,0),(0,0,0,1)}$.
  \item[] $M_2^a=\erz{a_6,a_7,a_5 a_8^{-w^a}}$ entsprechend $\erz{(0,1,0,0),(0,0,1,0),(1,0,0,-w^{a})}$ mit $a \in W_3$.
  \item[] $M_3=\erz{a_6,a_8,a_5 a_7\1}$ entsprechend $\erz{(0,1,0,0),(0,0,0,1),(1,0,-1,0)}$.
  \item[] $M_4=\erz{a_7,a_8,a_5 a_6\1}$ entsprechend $\erz{(0,0,1,0),(0,0,0,1),(1,-1,0,0)}$.
  \item[] $M_5^a=\erz{a_7,a_5 a_8^{-w^a}, a_6 a_8^{-w^a}}$ entsprechend $\erz{(0,0,1,0),(1,0,0,-w^{a}),(0,1,0,-w^{a})}$ mit $a \in W_3$.
  \item[] $M_6^b=\erz{a_8,a_5 a_7^{-w^b},a_6 a_7^{-w^b}}$ entsprechend $\erz{(0,0,0,1),(1,0,-w^{b},0),(0,1,-w^{b},0)}$ mit $b \in W_4$.
\end{enumerate}
\end{folg}

\begin{beweis}
Das Vertretersystem der zulässigen Untergruppen lässt sich unmittelbar aus \ref{zug_7}
ablesen. Nur die Bahnen $B_1$ bis $B_{6d}$ entsprechen nach \ref{suppl} zulässigen
Untergruppen, da nur sie Unterräume enthalten, die Komplemente zum Unterraum $\erz{(1,0,0,0)}$
sind, der dem Nukleus entspricht.
\end{beweis}

\begin{satz} \label{nt7}
In der folgenden Liste sind sämtliche unmittelbaren Nachfolger von $G$ der
Ordnung $p^5$ bis auf Isomorphie angegeben. Dabei ist $W_3=\{x \in \f \mid
x^3=1\}$ und $W_4=\{x \in \f \mid x^4=1\}$ und $a \in W_3$ sowie $b \in W_4$.
\begin{enumerate}
  \item[] $G_1=\erz{g_1,g_2,g_3,g_4,g_5 \mid [g_2,g_1]=g_3, [g_3,g_1]=g_4,
  [g_4,g_1]=g_5}$.
  \item[] $G_2^a=\erz{g_1,g_2,g_3,g_4,g_5 \mid [g_2,g_1]=g_3, [g_3,g_1]=g_4,
  [g_4,g_1]=g_5^{w^a}, g_2^p=g_5}$.
  \item[] $G_3=\erz{g_1,g_2,g_3,g_4,g_5 \mid [g_2,g_1]=g_3, [g_3,g_1]=g_4,
  [g_4,g_1]=g_5, g_1^p=g_5}$.
  \item[] $G_4=\erz{g_1,g_2,g_3,g_4,g_5 \mid [g_2,g_1]=g_3, [g_3,g_1]=g_4,
  [g_4,g_1]=g_5, [g_3,g_2]=g_5}$.
  \item[] $G_5^a=\erz{g_1,g_2,g_3,g_4,g_5 \mid [g_2,g_1]=g_3, [g_3,g_1]=g_4,
  [g_4,g_1]=g_5^{w^a}, [g_3,g_2]=g_5^{w^a}, g_2^p=g_5}$.
  \item[] $G_6^b=\erz{g_1,g_2,g_3,g_4,g_5 \mid [g_2,g_1]=g_3, [g_3,g_1]=g_4,
  [g_4,g_1]=g_5^{w^b}, [g_3,g_2]=g_5^{w^b}, g_1^p=g_5}$.
\end{enumerate}
\end{satz}

\begin{beweis}
Dieser Satz ergibt sich nach \ref{bahnenzu}, indem man $P(G)$ nach jedem Repräsentanten des
Vertretersystems zulässiger Untergruppen faktorisiert, das in \ref{zug_7} aufgelistet ist. Die
Präsentation der Faktorgruppe erhält man über das Verfahren aus \ref{faktor}.
\end{beweis}

\begin{folg}
Die Gruppe $G$ hat $$3 + ggT(4, p-1) + 2 \cdot ggT(3, p-1)$$ unmittelbare Nachfolger der
Ordnung $p^5$.
\end{folg}

\section{Nachfolger von $(p^4,8)$}

\begin{ver}
In diesem Abschnitt sei $G=\erz{a_1,a_2,a_3,a_4 \mid [a_2,a_1]=a_3, [a_3,a_1]=a_4,
a_1^p=a_4}$.
\end{ver}

\begin{lemma}
Die Gruppe $G$ hat die Gewichtung $\omega(a_1)=\omega(a_2)=1$ und $\omega(a_3)=2$ sowie
$\omega(a_4)=3$ und damit die $p$"=Klasse 3.
\end{lemma}

\begin{beweis}
Die Gewichtung von $G$ lässt sich nach \ref{gew} ermitteln.
\end{beweis}

\subsection{Die Automorphismengruppe}

\begin{lemma}
Jeder Automorphismus $\al$ von $G$ lässt sich durch $$\al(u_1, u_2, u_3, u_4, v_3, v_4):
  \begin{cases} \label{aut8}
    a_1 \mt a_1^{u_1} a_2^{u_2} a_3^{u_3} a_4^{u_4} & \text{wobei } u_1, u_2, u_3, u_4, v_3, v_4 \\
    a_2 \mt a_2^{u_1\1} a_3^{v_3} a_4^{v_4} & \in \{0 \dd p-1\} \text{ und }\\
    a_3 \mt a_3 a_4^{u_1 v_3} & u_1 \neq 0 \text{ ist.}\\
    a_4 \mt a_4^{u_1} &
  \end{cases}$$ darstellen.
\end{lemma}

\begin{beweis}
Die Automorphismengruppe von $G$ ist nach \ref{aut} aus dem Stabilisator $\bar{S}_3$ (siehe
\ref{bahnen_1_dpa}) der zu $G$ gehörenden zulässigen Untergruppe unmittelbar ablesbar. Damit
ergeben sich genau die oben angegebenen Bilder von $a_1$ und $a_2$. Die Bilder von $a_3$ und
$a_4$ sind durch die Bilder von $a_1$ und $a_2$ eindeutig bestimmt und lassen sich
folgendermaßen berechnen:
\begin{eqnarray*}
  a_3^\al &=& [a_2, a_1]^\al=[a_2^\al, a_1^\al] \\
          &=& [a_2^{u_1\1} a_3^{v_3} a_4^{v_4}, a_1^{u_1} a_2^{u_2} a_3^{u_3} a_4^{u_4}]\\
          &=& [a_2^{u_1\1} a_3^{v_3},  a_1^{u_1} a_2^{u_2} a_3^{u_3}]\\
          &=& [a_2^{u_1\1}, a_1^{u_1} a_2^{u_2} a_3^{u_3}]^{a_3^{v_3}} [a_3^{v_3}, a_1^{u_1} a_2^{u_2} a_3^{u_3}]\\
          &=& ([a_2^{u_1\1}, a_3^{u_3}] [a_2^{u_1\1}, a_1^{u_1} a_2^{u_2}]^{a_3^{u_3}})^{a_3^{v_3}} [a_3^{v_3}, a_3^{u_3}] [a_3^{v_3},  a_1^{u_1} a_2^{u_2}]^{a_3^{u_3}}\\
          &=& ([a_2^{u_1\1}, a_1^{u_1}a_2^{u_2}]^{a_3^{u_3}})^{a_3^{v_3}} [a_3^{v_3}, a_1^{u_1} a_2^{u_2}]^{a_3^{u_3}}\\
          &=& (([a_3^{u_1\1}, a_2^{u_2}] [a_2^{u_1\1}, a_1^{u_1}]^{a_2^{u_2}})^{a_3^{u_3}})^{a_3^{v_3}} ([a_3^{v_3}, a_2^{u_2}]  [a_3^{v_3}, a_1^{u_1}]^{a_2^{u_2}})^{a_3^{u_3}}\\
          &=& (([a_2^{u_1\1}, a_1^{u_1}]^{a_2^{u_2}})^{a_3^{u_3}})^{a_3^{v_3}} ([a_3^{v_3}, a_1^{u_1}]^{a_2^{u_2}})^{a_3^{u_3}}\\
          &=& [a_2^{u_1\1}, a_1^{u_1}] [a_3^{v_3}, a_1^{u_1}]\\
          &=& [a_2, a_1] [a_3, a_1]^{u_1 v_3} = a_3 a_4^{u_1 v_3},
\end{eqnarray*} da $a_3$ und $a_4$ zentral sind. Weiterhin ist
\begin{eqnarray*}
  a_4^\al &=& [a_3,a_1]^\al=[a_3^\al, a_1^\al]\\
          &=& [a_3 a_4^{u_1 v_3}, a_1^{u_1} a_2^{u_2} a_3^{u_3} a_4^{u_4}] \\
          &=& [a_3, a_1^{u_1} a_2^{u_2} a_3^{u_3}] = [a_3, a_1^{u_1} a_2^{u_2}]\\
          &=& [a_3, a_1^{u_1}] [a_3, a_2^{u_2}] = [a_3, a_1^{u_1}] \\
          &=& [a_3, a_1]^{u_1}=a_4^{u_1}
\end{eqnarray*}
\end{beweis}

\subsection{$p$"=Cover, Multiplikator und Nukleus}

\begin{lemma}
Für $G$ ergibt sich das $p$"=Cover, der Nukleus und der Multiplikator in folgender Weise mit
der Gewichtung $\omega$:
\begin{enumerate}
  \item[] $P(G)=\erz{a_1 \dd a_7 \mid [a_2,a_1]=a_3,
  [a_3,a_1]=a_4, [a_3,a_2]=a_5, a_1^p = a_4 a_6, a_2^p=a_7}$.
  \item[] $M(G)=\erz{a_5, a_6, a_7}$.
  \item[] $N(G)=\erz{1}$.
\end{enumerate}
Die Gruppe $P(G)$ hat die $p$"=Klasse 3.
\end{lemma}

\begin{beweis}
Die Präsentation von $P(G)$ ergibt sich aus \ref{cov} und \ref{newman} sowie dem
Reduktionsverfahren von Knuth"=Bendix, aus der sich die Untergruppe $M(G)$ ablesen lässt. Die
Gewichtung von $P(G)$
\begin{enumerate}
  \item[] $1=\omega(a_1)=\omega(a_2)$.
  \item[] $2=\omega(a_3)=\omega(a_6)=\omega(a_7)$.
  \item[] $3=\omega(a_4)=\omega(a_5)$.
\end{enumerate}
und damit $N(G)$ lassen sich nach \ref{gew} aus der Präsentation von $P(G)$ ablesen.
\end{beweis}

\begin{satz}
Die Gruppe $G$ ist abschließend.
\end{satz}

\begin{beweis}
Da $P(G)$ dieselbe $p$"=Klasse hat wie $G$, hat $G$ nach \ref{pfort} keine unmittelbaren
Nachfolger, d.\,h. $G$ ist abschließend.
\end{beweis}

\section{Nachfolger von $(p^4,9)$}

\begin{ver}
In diesem Abschnitt sei $G=\erz{a_1,a_2,a_3,a_4 \mid [a_2,a_1]=a_3, [a_3,a_1]=g_4,
a_2^p=g_4}$.
\end{ver}

\begin{lemma}
Die Gruppe $G$ hat die Gewichtung $\omega(a_1)=\omega(a_2)=1$ und $\omega(a_3)=2$ sowie
$\omega(a_4)=3$ und damit die $p$"=Klasse 3.
\end{lemma}

\begin{beweis}
Die Gewichtung von $G$ lässt sich nach \ref{gew} ermitteln.
\end{beweis}

\subsection{Die Automorphismengruppe}

\begin{lemma}
Die Automorphismen $\al$ von $G$ lassen sich in der folgenden Weise darstellen: $$\al(u_3,
u_4, v_2, v_3, v_4):
  \begin{cases} \label{aut9}
    a_1 \mt a_1^{\pm 1} a_3^{u_3} a_4^{u_4} & \text{wobei } u_3, u_4, v_2, v_3, v_4 \\
    a_2 \mt a_2^{v_2} a_3^{v_3} a_4^{v_4} & \in \{0 \dd p-1\} \text{ und }\\
    a_3 \mt a_3^{\pm v_2} a_4^{\pm v_3} & u_1 \neq 0 \text{ ist.}\\
    a_4 \mt a_4^{\pm v_2} &
  \end{cases}$$
\end{lemma}

\begin{beweis}
Die Automorphismengruppe von $G$ ist nach \ref{aut} aus dem Stabilisator $\bar{S}_4$ (siehe
\ref{bahnen_1_dpa}) der zu $G$ gehörenden zulässigen Untergruppe unmittelbar ablesbar. Damit
ergeben sich genau die oben angegebenen Bilder von $a_1$ und $a_2$. Die Bilder von $a_3$ und
$a_4$ sind durch die Bilder von $a_1$ und $a_2$ eindeutig bestimmt und lassen sich
folgendermaßen berechnen:
\begin{eqnarray*}
  a_3^\al &=& [a_2, a_1]^\al=[a_2^\al, a_1^\al]\\
          &=& [a_2^{v_2} a_3^{v_3} a_4^{v_4}, a_1^{\pm 1} a_3^{u_3} a_4^{u_4}]\\
          &=& [a_2^{v_2} a_3^{v_3}, a_1^{\pm 1} a_3^{u_3}]\\
          &=& [a_2^{v_2}, a_1^{\pm 1} a_3^{u_3}]^{a_3^{v_3}} [a_3^{v_3}, a_1^{\pm 1} a_3^{u_3}]\\
          &=& ([a_2^{v_2}, a_3^{u_3}] [a_2^{v_2}, a_1^{\pm 1}]^{a_3^{u_3}})^{a_3^{v_3}} [a_3^{v_3}, a_1^{\pm 1}]\\
          &=& [a_2^{v_2}, a_1^{\pm 1}] [a_3, a_1]^{\pm v_3}\\
          &=& a_3^{\pm v_2} a_4^{\pm v_3}
\end{eqnarray*} Weiterhin ist
\begin{eqnarray*}
  a_4^\al &=& [a_3,a_1]^\al=[a_3^\al, a_1^\al]\\
          &=& [a_3^{\pm v_2} a_4^{\pm v_3}, a_1^{\pm 1} a_3^{u_3} a_4^{u_4}]\\
          &=& [a_3^{\pm v_2}, a_1^{\pm 1} a_3^{u_3}] = [a_3^{\pm v_2}, a_1^{\pm 1}] = a_4^{\pm v_2}
\end{eqnarray*}
\end{beweis}

\subsection{$p$"=Cover, Multiplikator und Nukleus}

\begin{lemma}
Für $G$ ergibt sich das $p$"=Cover, der Nukleus und der Multiplikator in folgender Weise mit
der Gewichtung $\omega$:
\begin{enumerate}
  \item[] $P(G)=\erz{a_1 \dd a_7 \mid [a_2,a_1]=a_3,
  [a_3,a_1]=a_4, [a_3,a_2]=a_5, a_1^p=a_6, a_2^p= a_4 a_7}$.
  \item[] $M(G)=\erz{a_5, a_6, a_7}$.
  \item[] $N(G)=\erz{1}$.
\end{enumerate}
Die Gruppe $P(G)$ hat die $p$"=Klasse 3.
\end{lemma}

\begin{beweis}
Die Präsentation von $P(G)$ ergibt sich aus \ref{cov} und \ref{newman} sowie dem
Reduktionsverfahren von Knuth"=Bendix, aus der sich die Untergruppe $M(G)$ ablesen lässt. Die
Gewichtung von $P(G)$
\begin{enumerate}
  \item[] $1=\omega(a_1)=\omega(a_2)$.
  \item[] $2=\omega(a_3)=\omega(a_6)=\omega(a_7)$.
  \item[] $3=\omega(a_4)=\omega(a_5)$.
\end{enumerate}
und damit $N(G)$ lassen sich nach \ref{gew} aus der Präsentation von $P(G)$ ablesen.
\end{beweis}

\begin{satz}
Die Gruppe $G$ ist abschließend.
\end{satz}

\begin{beweis}
Da $P(G)$ dieselbe $p$"=Klasse hat wie $G$, hat $G$ nach \ref{pfort} keine unmittelbaren
Nachfolger, d.\,h. $G$ ist abschließend.
\end{beweis}

\section{Nachfolger von $(p^4,10)$}

\begin{ver}
In diesem Abschnitt sei $G=\erz{a_1,a_2,a_3,a_4 \mid [a_2,a_1]=a_3, [a_3,a_1]=a_4,
a_2^p=a_4^w}$, wobei $w$ ein Erzeuger der multiplikativen Gruppe von $\f$ ist.
\end{ver}

\begin{lemma}
Die Gruppe $G$ hat die Gewichtung $\omega(a_1)=\omega(a_2)=1$ und $\omega(a_3)=2$ sowie
$\omega(a_4)=3$ und damit die $p$"=Klasse 3
\end{lemma}

\begin{beweis}
Die Gewichtung von $G$ lässt sich nach \ref{gew} ermitteln.
\end{beweis}

\begin{lemma}
Die Automorphismen $\al$ von $G$ lassen sich in der folgenden Weise darstellen: $$\al(u_3,
u_4, v_2, v_3, v_4):
  \begin{cases} \label{aut10}
    a_1 \mt a_1^{\pm 1} a_3^{u_3}     a_4^{u_4}     & \text{wobei } u_3, u_4, v_2, v_3, v_4 \\
    a_2 \mt a_2^{v_2}   a_3^{v_3}     a_4^{v_4}     & \in \{0 \dd p-1\} \text{ und }\\
    a_3 \mt             a_3^{\pm v_2} a_4^{\pm v_3} & u_1 \neq 0 \text{ ist.}\\
    a_4 \mt                           a_4^{\pm v_2} &
  \end{cases}$$
\end{lemma}

\begin{beweis}
Die Automorphismengruppe von $G$ ist nach \ref{aut} aus dem Stabilisator $\bar{S}_4$ (siehe
\ref{bahnen_1_dpa}) der zu $G$ gehörenden zulässigen Untergruppe unmittelbar ablesbar. Damit
ergeben sich genau die oben angegebenen Bilder von $a_1$ und $a_2$. Die Bilder von $a_3$ und
$a_4$ sind durch die Bilder von $a_1$ und $a_2$ eindeutig bestimmt und lassen sich
folgendermaßen berechnen:
\begin{eqnarray*}
  a_3^\al &=& [a_2, a_1]^\al=[a_2^\al, a_1^\al]\\
          &=& [a_2^{v_2} a_3^{v_3} a_4^{v_4}, a_1^{\pm 1} a_3^{u_3} a_4^{u_4}]\\
          &=& [a_2^{v_2} a_3^{v_3}, a_1^{\pm 1} a_3^{u_3}]\\
          &=& [a_2^{v_2}, a_1^{\pm 1} a_3^{u_3}]^{a_3^{v_3}} [a_3^{v_3}, a_1^{\pm 1} a_3^{u_3}]\\
          &=& ([a_2^{v_2}, a_3^{u_3}] [a_2^{v_2}, a_1^{\pm 1}]^{a_3^{u_3}})^{a_3^{v_3}} [a_3^{v_3}, a_1^{\pm 1}]\\
          &=& [a_2^{v_2}, a_1^{\pm 1}] [a_3, a_1]^{\pm v_3}\\
          &=& a_3^{\pm v_2} a_4^{\pm v_3}
\end{eqnarray*} Weiterhin ist
\begin{eqnarray*}
  a_4^\al &=& [a_3,a_1]^\al=[a_3^\al, a_1^\al]\\
          &=& [a_3^{\pm v_2} a_4^{\pm v_3}, a_1^{\pm 1} a_3^{u_3} a_4^{u_4}]\\
          &=& [a_3^{\pm v_2}, a_1^{\pm 1} a_3^{u_3}] = [a_3^{\pm v_2}, a_1^{\pm 1}] = a_4^{\pm v_2}
\end{eqnarray*}
\end{beweis}

\subsection{$p$"=Cover, Multiplikator und Nukleus}

\begin{lemma}
Für $G$ ergibt sich das $p$"=Cover, der Nukleus und der Multiplikator in folgender, wobei $w$
ein Erzeuger der multiplikativen Gruppe von $\f$ ist:
\begin{enumerate}
  \item[] $P(G)=\erz{a_1 \dd a_7 \mid [a_2,a_1]=a_3, [a_3,a_1]=a_4,
  [a_3,a_2]=a_5, a_1^p=a_6, a_2^p=a_4^w a_7}$.
  \item[] $M(G)=\erz{a_5, a_6, a_7}$.
  \item[] $N(G)=\erz{1}$.
\end{enumerate}
Die Gruppe $P(G)$ hat die $p$"=Klasse 3.
\end{lemma}

\begin{beweis}
Die Präsentation von $P(G)$ ergibt sich aus \ref{cov} und \ref{newman} sowie dem
Reduktionsverfahren von Knuth"=Bendix, aus der sich die Untergruppe $M(G)$ ablesen lässt. Die
Gewichtung von $P(G)$
\begin{enumerate}
  \item[] $1=\omega(a_1)=\omega(a_2)$.
  \item[] $2=\omega(a_3)=\omega(a_6)=\omega(a_7)$.
  \item[] $3=\omega(a_4)=\omega(a_5)$.
\end{enumerate}
und damit $N(G)$ lassen sich nach \ref{gew} aus der Präsentation von $P(G)$ ablesen.
\end{beweis}

\begin{satz}
Die Gruppe $G$ ist abschließend.
\end{satz}

\begin{beweis}
Da $P(G)$ dieselbe $p$"=Klasse hat wie $G$, hat $G$ nach \ref{pfort} keine unmittelbaren
Nachfolger, d.\,h. $G$ ist abschließend.
\end{beweis}

\section{Nachfolger von $C_{p^2} \times C_p$}

\begin{ver}
In diesem Abschnitt sei $$G=\erz{a_1, a_2, a_3 \mid a_1^p=a_3}.$$ Damit ist $G \cong C_{p^2}
\times C_p$.
\end{ver}

\begin{lemma}
Die Gruppe $G$ hat die Gewichtung $\omega(a_1)=\omega(a_2)=1$ und $\omega(a_3)=2$.
\end{lemma}

\begin{beweis}
Die Gewichtung von $G$ lässt sich nach \ref{gew} ermitteln.
\end{beweis}

\subsection{Die Automorphismengruppe}

\begin{lemma} \label{autcp2cp}
Für $G$ erhält man die Automorphismen $\al$ in folgender Weise: $$\al(u_1, u_2, u_3, v_2,
v_3):
  \begin{cases}
    a_1 \mt a_1^{u_1} a_2^{u_2} a_3^{u_3}  & \text{wobei } u_1, u_2, u_3, v_2, v_3 \\
    a_2 \mt           a_2^{v_2} a_3^{v_3}  & \in \{0 \dd p-1\} \text{ und }\\
    a_3 \mt                     a_3^{u_1}  & u_1 \neq 0 \text{ sowie } v_2 \neq 0 \text{ ist.}
\end{cases}$$
\end{lemma}

\begin{beweis}
Nach \ref{aut} lässt sich $Aut(G)$ unmittelbar aus dem Stabilisator $T_2$ (siehe
\ref{stab_cp2_p3}) der zu $G$ gehörenden zulässigen Untergruppe ablesen. Damit ergeben sich
genau die oben angegebenen Bilder von $a_1$ und $a_2$. Das Bild von $a_3$ ist durch die Bilder
von $a_1$ und $a_2$ vollständig festgelegt und lässt sich folgendermaßen berechnen:
\begin{eqnarray*}
  a_3^\al &=& (a_1^p)^\al=(a_1^{u_1} a_2^{u_2} a_3^{u_3})^p \\
          &=& (a_1^{u_1})^p (a_2^{u_2})^p (a_3^{u_3})^p = (a_1^p)^{u_1} (a_2^p)^{u_2} (a_3^p)^{u_3}\\
          &=& (a_1^p)^{u_1} = a_3^{u_1},
\end{eqnarray*}
da $G$ kommutativ ist.
\end{beweis}

\subsection{$p$"=Cover, Multiplikator und Nukleus}

\begin{lemma}
Für $G$ ergibt sich das $p$"=Cover, der Nukleus und der Multiplikator in
folgender Weise mit der Gewichtung $\omega$:
\begin{enumerate}
  \item[] $P(G)=\erz{a_1,a_2,a_3,a_4,a_5,a_6 \mid a_1^p=a_3, a_3^p=a_4,
  [a_2,a_1]=a_5, a_2^p=a_6}$.
  \item[] $M(G)=\erz{a_4, a_5, a_6}$.
  \item[] $N(G)=\erz{a_4}$.
\end{enumerate}
\end{lemma}

\begin{beweis}
Die Präsentation von $P(G)$ ergibt sich aus \ref{cov} und \ref{newman} sowie dem
Reduktionsverfahren von Knuth"=Bendix, aus der sich die Untergruppe $M(G)$ ablesen lässt. Die
Gewichtung von $P(G)$
\begin{enumerate}
  \item[] $1=\omega(a_1)=\omega(a_2)$.
  \item[] $2=\omega(a_3)=\omega(a_5)=\omega(a_6)$.
  \item[] $3=\omega(a_4)$.
\end{enumerate}
und damit $N(G)$ lassen sich nach \ref{gew} aus der Präsentation von $P(G)$ ablesen.
\end{beweis}

\subsection{Operation der Erweiterungsautomorphismen}

\begin{lemma}
Werden die Automorphismen $\al$ von $G$ gemäß \ref{autcp2cp} in der Weise $$\al(u_1, u_2, u_3,
v_2, v_3):
  \begin{cases}
    a_1 \mt a_1^{u_1} a_2^{u_2} a_3^{u_3} & \text{wobei } u_1, u_2, u_3, v_2, v_3 \\
    a_2 \mt a_2^{v_2} a_3^{v_3} & \in \{0 \dd p-1\} \text{ und }\\
    a_3 \mt a_3^{u_1}  & u_1 \neq 0 \text{ sowie } v_2 \neq 0 \text{ ist.}
\end{cases}$$ dargestellt, so ist
$$\vi: Aut(G) \to GL(3,p): \al(u_1, u_2, u_3, v_2, v_3) \mt M =
 \begin{pmatrix}
  u_1 & 0       & 0 \\
  0   & u_1 v_2 & 0 \\
  v_3 & 0       & v_2
\end{pmatrix}$$ der Operationshomomorphismus von $A=Aut(G)$ auf $\f^3 \cong M(G)$ über
Erweiterungsautomorphismen.
\end{lemma}

\begin{beweis}
Die Bilder der Erzeuger von $M(G)$ ergeben sich unter $\al$ über $\vi$ folgendermaßen:
\begin{eqnarray*}
  a_4^\al &=& (a_3^p)^\al=(a_3^\al)^p=(a_3^{u_1})^p\\
          &=& (a_3^p)^{u_1}=a_4^{u_1}
\end{eqnarray*} Weiterhin ist
\begin{eqnarray*}
  a_5^\al &=& [a_2,a_1]^\al=[a_2^\al,a_1^\al]= [a_2^{v_2} a_3^{v_3}, a_1^{u_1} a_2^{u_2} a_3^{u_3}]\\
          &=& [a_2^{v_2} a_3^{v_3}, a_1^{u_1} a_2^{u_2} a_3^{u_3}]\\
          &=& [a_2^{v_2}, a_1^{u_1} a_2^{u_2} a_3^{u_3}]^{a_3^{v_3}} [a_3^{v_3}, a_1^{u_1} a_2^{u_2} a_3^{u_3}] \\
          &=& ([a_2^{v_2}, a_3^{u_3}] [a_2^{v_2}, a_1^{u_1} a_2^{u_2}]^{a_3^{u_3}})^{a_3^{v_3}} [a_3^{v_3}, a_3^{u_3}] [a_3^{v_3}, a_1^{u_1} a_2^{u_2}]^{a_3^{u_3}}\\
          &=& ([a_2^{v_2}, a_1^{u_1} a_2^{u_2}]^{a_3^{u_3}})^{a_3^{v_3}} [a_3^{v_3}, a_1^{u_1} a_2^{u_2}]^{a_3^{u_3}}\\
          &=& (([a_2^{v_2}, a_2^{u_2}][a_2^{v_2}, a_1^{u_1}]^{a_2^{u_2}})^{a_3^{u_3}})^{a_3^{v_3}} ([a_3^{v_3}, a_2^{u_2}][a_3^{v_3}, a_1^{u_1}]^{a_2^{u_2}})^{a_3^{u_3}}\\
          &=& (([a_2^{v_2}, a_1^{u_1}]^{a_2^{u_2}})^{a_3^{u_3}})^{a_3^{v_3}}\\
          &=& [a_2, a_1]^{u_1 v_2}=a_5^{u_1 v_2},
\end{eqnarray*} da $a_5$ zentral und die anderen auftretenden Kommutatoren
trivial sind. Außerdem ist
\begin{eqnarray*}
  a_6^\al &=& (a_2^p)^\al = (a_2^\al)^p = (a_2^{v_2} a_3^{v_3})^p\\
          &=& (a_2^{v_2})^p (a_3^{v_3})^p = (a_2^p)^{v_2} (a_3^p)^{v_3} = a_4^{v_3}
          a_6^{v_2},
\end{eqnarray*} da der Kommutator von $a_4$ und $a_6$ trivial ist.
\end{beweis}

Für die Isomorphieklassen der Nachfolger der Ordnung $p^4$ werden die Bahnen zweidimensionaler
Unterräume von $\F_p^3 \cong M(G)$ unter der Operation von $A=Aut(G)$ benötigt. Diese Bahnen
kann man einfacher ermitteln, indem man eindimensionale Komplemente des Dualraumes betrachtet
und den zu $\vi$ dualen Operationshomomorphismus.

\begin{folg}
Der zu $\vi$ duale Operationshomomorphismus ist $$\bar{\vi}: A \to GL(3,p): \al(u_1, u_2, u_3,
v_2, v_3) \mt M =
 \begin{pmatrix}
  u_1 & 0       & v_3 \\
  0   & u_1 v_2 & 0 \\
  0   & 0       & v_2
\end{pmatrix}$$
\end{folg}

\subsection{Bahnen zulässiger Untergruppen für Nachfolger der Ordnung $p^4$}

\begin{lemma} \label{bahnen_cp2cp}
Der Vektorraum $\f^3 \cong M(G)$ zerfällt unter der Operation von $A=Aut(G)$ vermittels
$\bar{\vi}$ in die folgenden Bahnen nichttrivialer Vektoren zerfällt:
\begin{enumerate}
  \item[] $B_1=(0,0,1)^A$ mit $|B_1|=p-1$.
  \item[] $B_2=(0,1,0)^A$ mit $|B_2|=p-1$.
  \item[] $B_3=(0,1,1)^A$ mit $|B_3|=(p-1)^2$.
  \item[] $B_4=(1,0,0)^A$ mit $|B_4|=p(p-1)$.
  \item[] $B_5=(1,1,0)^A$ mit $|B_5|=p(p-1)^2$.
\end{enumerate}
Die Stabilisatoren der Vertreter von $B_4$ und $B_5$ sind
\begin{eqnarray*}
  \hat{S}_4 & = & \left\{ \begin{pmatrix}
    1         & 0        & 0          \\
    0         & v_2      & 0          \\
    0         & 0        & v_2        \\
  \end{pmatrix} \mid v_2 \in \F_p  \text{ und } v_2 \neq 0\right\}\\
  \hat{S}_5 & = & \left\{ \begin{pmatrix}
    1         & 0        & 0           \\
    0         & 1        & 0           \\
    0         & 0        & 1         \\
  \end{pmatrix}  \right\}\\
\end{eqnarray*}
und die Stabilisatoren der eindimensionalen Unterräume, die von den Vertretern
von $B_4$ und $B_5$ aufgespannt werden, sind
\begin{eqnarray*}
  \bar{S}_4 & = & \left\{ \begin{pmatrix}
    u_1       & 0            & 0          \\
    0         & u_1 v_2      & 0          \\
    0         & 0            & v_2        \\
  \end{pmatrix} \mid u_1, v_2 \in \F_p  \text{ und } u_1 v_2 \neq 0\right\}\\
  \bar{S}_5 & = & \left\{ \begin{pmatrix}
    u_1       & 0        & 0           \\
    0         & u_1      & 0           \\
    0         & 0        & 1         \\
  \end{pmatrix} \mid u_1 \in \F_p  \text{ und } u_1 \neq 0 \right\}\\
\end{eqnarray*}
\end{lemma}

\begin{beweis}
Die Behauptung ist offensichtlich.
\end{beweis}

\begin{folg} \label{zug_cp2cp_p4}
Für die Nachfolger von $G$ der Ordnung $p^4$ ergibt sich das folgende
Vertretersystem zulässiger Untergruppen:
\begin{enumerate}
  \item $M_4=\erz{a_5,a_6}$ entsprechend $\erz{(0,1,0),(0,0,1)}$.
  \item $M_5=\erz{a_6,a_4 a_5\1}$ entsprechend $\erz{(0,0,1),(1,-1,0)}$.
\end{enumerate}
Die Stabilisatoren von $M_4$ und $M_5$ unter der Operation von $Aut(G)$ auf
$M(G)$ über Erweiterungsautomorphismen sind
\begin{enumerate}
  \item $S_4=\{\al(u_1, u_2, u_3, v_2, 0) \mid u_1, u_2, u_3, v_2
  \in \{0 \dd p-1\} \text{ und } u_1 \neq 0 \}$.
  \item $S_5=\{\al(u_1, u_2, u_3, 1, 0) \mid u_1, u_2, u_3
  \in \{0 \dd p-1\} \text{ und } u_1 \neq 0 \}$.
\end{enumerate}
\end{folg}

\begin{beweis}
Da nur die Vertreter der Bahnen $B_4$ und $B_5$ mit $(1,0,0)$ ein Skalarprodukt
ungleich Null liefern, entsprechen nur die zweidimensionalen Unterräume den
zulässigen Untergruppen von $M(G)$, zu denen ein Element aus $B_4$ oder $B_5$
der Basisvektor des orthogonalen, eindimensionalen Komplementes ist. Die
Stabilisatoren von $M_4$ und $M_5$ lassen sich nach Transposition aus
$\bar{S}_4$ und $\bar{S}_5$ ablesen.
\end{beweis}

\subsection{Nachfolger der Ordnung $p^4$}

\begin{satz} \label{ncpp2p4}
In der folgenden Liste sind sämtliche unmittelbaren Nachfolger von $G$ der
Ordnung $p^4$ bis auf Isomorphie angegeben:
\begin{enumerate}
  \item $G_1=\erz{g_1,g_2,g_3,g_4 \mid g_1^p=g_3, g_3^p=g_4} \cong C_{p^3} \times C_p$ \hfill
  \gap"=Typ $(p^4,5)$
  \item $G_2=\erz{g_1,g_2,g_3,g_4 \mid g_1^p=g_3, g_3^p=g_4,
  [g_2,g_1]=g_4}$ \hfill \gap"=Typ $(p^4,6)$
\end{enumerate}
\end{satz}

\begin{beweis}
Dieser Satz ergibt sich nach \ref{bahnenzu}, indem man $P(G)$ nach jedem Repräsentanten des
Vertretersystems zulässiger Untergruppen faktorisiert, das in \ref{zug_cp2cp_p4} aufgelistet
ist. Die Präsentation der Faktorgruppe erhält man über das Verfahren aus \ref{faktor}.
\end{beweis}

\begin{satz}
Die Gruppe $G$ hat nur unmittelbare Nachfolger der Ordnung $p^4$.
\end{satz}

\begin{beweis}
Mit den Gruppen $G_1$ und $G_2$ sind alle (Isomorphieklassen der) unmittelbaren
Nachfolger von $G$ angegeben. Unmittelbare Nachfolger der Ordnung $p^5$ gibt es
nicht, da $P(G)$ die Ordnung $p^6$ und sich Nachfolger der Ordnung $p^5$ nur
dann ergäben, wenn man zulässige Untergruppen der Ordnung $p$ ausfaktorisieren
könnte. Da $M(G)$ die Ordnung $p^3$ und $N(G)$ die Ordnung $p$ hat, ist keine
Untergruppe der Ordnung $p$ ein Supplement in $M(G)$ zu $N(G)$. Also gibt es
nach \ref{suppl} keine zulässigen Untergruppen der Ordnung $p$.
\end{beweis}

\section{Nachfolger von $(p^4,5)$ bzw. $C_{p^3} \times C_p$}

\begin{ver}
In diesem Abschnitt sei $$G=\erz{a_1,a_2,a_3,a_4 \mid a_1^p=a_3, a_3^p=a_4}$$ Damit ist $G
\cong C_{p^3} \times C_p$.
\end{ver}

\begin{lemma}
Die Gruppe $G$ hat die Gewichtung $\omega(a_1)=\omega(a_2)=1$, $\omega(a_3)=2$ und
$\omega(a_4)=3$.
\end{lemma}

\begin{beweis}
Die Gewichtung von $G$ lässt sich nach \ref{gew} ermitteln.
\end{beweis}

\subsection{Die Automorphismengruppe}

\begin{lemma}
Jeder Automorphismus $\al$ von $G$ lässt sich in folgender Weise darstellen: $$\al(u_1, u_2,
u_3, u_4, v_2, v_3, v_4): \label{aut5}
  \begin{cases}
    a_1 \mt a_1^{u_1} a_2^{u_2} a_3^{u_3} a_4^{u_4} & \text{wobei } u_1, u_2, u_3, u_4, v_2, v_4 \\
    a_2 \mt a_2^{v_2} a_4^{v_4} & \in \{0 \dd p-1\} \text{ und }\\
    a_3 \mt a_3^{u_1} a_4^{u_3} & u_1 \neq 0 \text{ sowie } v_2 \neq 0 \text{ ist.}\\
    a_4 \mt a_4^{u_1} &
  \end{cases}.$$
\end{lemma}

\begin{beweis}
Die Automorphismengruppe von $G$ ist nach \ref{aut} aus dem Stabilisator $\hat{S}_4$ (siehe
\ref{bahnen_cp2cp}) der zu $G$ gehörenden zulässigen Untergruppe unmittelbar ablesbar. Aus
$\hat{S}_4$ lassen sich die Bilder von $a_1$ und $a_2$ unmittelbar ablesen. Die Bilder von
$a_3$ und $a_4$ sind durch die Bilder von $a_1$ und $a_2$ eindeutig bestimmt und lassen sich
folgendermaßen berechnen:
\begin{eqnarray*}
  a_3^\al &=& (a_1^p)^\al=(a_1^\al)^p=(a_1^{u_1} a_2^{u_2} a_3^{u_3} a_4^{u_4})^p\\
          &=& (a_1^{u_1})^p (a_2^{u_2})^p (a_3^{u_3})^p (a_4^{u_4})^p = (a_1^p)^{u_1} (a_2^p)^{u_2} (a_3^p)^{u_3} (a_4^p)^{u_4}\\
          &=& a_3^{u_1} a_4^{u_3}
\end{eqnarray*} und
\begin{eqnarray*}
  a_4^\al &=& (a_3^p)^\al=(a_3^\al)^p=(a_3^{u_1} a_4^{u_3})^p = (a_3^{u_1})^p (a_4^{u_3})^p\\
          &=& (a_3^p)^{u_1} (a_4^p)^{u_3},
\end{eqnarray*} da $G$ abelsch ist.
\end{beweis}

\subsection{$p$"=Cover, Multiplikator und Nukleus}

\begin{lemma}
Für $G$ ergibt sich das $p$"=Cover, der Nukleus und der Multiplikator in
folgender Weise mit der Gewichtung $\omega$:
\begin{enumerate}
  \item[] $P(G)=\erz{a_1,a_2,a_3,a_4,a_5,a_6,a_7 \mid
  a_1^p=a_3, a_3^p=a_4, a_4^p=a_5, [a_2,a_1]=a_6, a_2^p=a_7}$.
  \item[] $M(G)=\erz{a_5, a_6, a_7}$.
  \item[] $N(G)=\erz{a_5}$.
\end{enumerate}
\end{lemma}

\begin{beweis}
Die Präsentation von $P(G)$ ergibt sich aus \ref{cov} und \ref{newman} und dem
Reduktionsverfahren von Knuth"=Bendix, aus der sich $M(G)$ ablesen lässt. Die Gewichtung von
$P(G)$
\begin{enumerate}
  \item[] $1=\omega(a_1)=\omega(a_2)$.
  \item[] $2=\omega(a_3)=\omega(a_6)=\omega(a_7)$.
  \item[] $3=\omega(a_4)$.
  \item[] $4=\omega(a_5)$.
\end{enumerate}
und damit $N(G)$ lassen sich nach \ref{gew} aus der Präsentation von $P(G)$ ablesen.
\end{beweis}

\subsection{Operation der Erweiterungsautomorphismen}

\begin{lemma}
Werden die Automorphismen $\al$ von $G$ gemäß \ref{aut5} in der Weise $$\al(u_1, u_2, u_3,
u_4, v_2, v_3, v_4):
  \begin{cases}
    a_1 \mt a_1^{u_1} a_2^{u_2} a_3^{u_3} a_4^{u_4} & \text{wobei } u_1, u_2, u_3, u_4, v_2, v_4 \\
    a_2 \mt a_2^{v_2} a_4^{v_4} & \in \{0 \dd p-1\} \text{ und }\\
    a_3 \mt a_3^{u_1} a_4^{u_3} & u_1 \neq 0 \text{ sowie } v_2 \neq 0 \text{ ist.}\\
    a_4 \mt a_4^{u_1} &
  \end{cases}$$
dargestellt, so operiert die Automorphismengruppe von $G$ über den Operationshomomorphismus
$$\vi: Aut(G) \to GL(3,p): \al(u_1, u_2, u_3, u_4, v_2, v_3, v_4) \mt M =
\begin{pmatrix}
  u_1  & 0       & 0 \\
  0    & u_1 v_2 & 0 \\
  v_4  & 0       & v_2
\end{pmatrix}$$ auf $\F_p^3 \cong M(G)$.
\end{lemma}

\begin{beweis}
Nach \ref{aut5} lässt sich jeder Automorphismus $\al$ von $G$ in der oben angegebenen Weise
darstellen. Die Bilder der Erzeuger von $M(G)$ sind unter der Operation von $Aut(G)$ über
Erweiterungsautomorphismen dann folgendermaßen gegeben:
\begin{eqnarray*}
  a_5^\al &=& (a_4^p)^\al=(a_4^\al)^p=(a_4^{u_1})^p=(a_4^p)^{u_1}\\
          &=& a_5^{u_1}
\end{eqnarray*}
sowie
\begin{eqnarray*}
  a_6^\al &=& [a_2, a_1]^\al=[a_2^\al, a_1^\al] \\
          &=& [a_2^{v_2} a_4^{v_4}, a_1^{u_1} a_2^{u_2} a_3^{u_3} a_4^{u_4}]\\
          &=& [a_2, a_1]^{u_1 v_2}=a_6^{u_1 v_2},
\end{eqnarray*}
da $[a_2, a_1]$ zentral und die Kommutatoren $[a_3,a_1] \dd [a_4, a_3]$ trivial
sind. Weiterhin ist
\begin{eqnarray*}
  a_7^\al &=& (a_2^p)^\al=(a_2^\al)^p\\
          &=& (a_2^{v_2} a_4^{v_4})^p =(a_2^p)^{v_2} (a_4^p)^{v_4}\\
          &=& a_5^{v_4} a_7^{v_2},
\end{eqnarray*}
da $[a_4, a_2]=1$ ist.
\end{beweis}

Da für die unmittelbaren Nachfolger der Ordnung $p^5$ zulässige Untergruppen der Ordnung $p^2$
benötigt werden, werden die Bahnen der zweidimensionalen Unterräume von $\F_p^3$ unter der
Operation von $Aut(G)$ vermittels $\vi$ gesucht. Wie üblich werden statt der zweidimensionalen
Unterräume ihre eindimensionalen orthogonalen Komplemente betrachtet und damit statt $\vi$ der
zu $\vi$ duale Operationshomomorphismus verwendet.

\begin{beme}
Der zu $\vi$ duale Operationshomomorphismus ist $$\bar{\vi}: Aut(G) \to
Aut(M(G)): \al(u_1, u_2, u_3, u_4, v_2, v_3, v_4) \mt M = \begin{pmatrix}
  u_1  & 0       &  v_4 \\
  0    & u_1 v_2 & 0 \\
  0    & 0       & v_2
\end{pmatrix}.$$
\end{beme}

\subsection{Bahnen zulässiger Untergruppen für Nachfolger der Ordnung $p^5$}

\begin{lemma}
Der Vektorraum $\F_p^3 \cong M(G)$ zerfällt unter $A=Aut(G)^{\bar{\vi}}$ in folgende Bahnen
nichttrivialer Vektoren:
\begin{enumerate}
  \item[] $B_1=(1,0,0)^A$ mit $|B_1|=p(p-1)$,
  \item[] $B_2=(1,1,0)^A$ mit $|B_2|=p(p-1)^2$,
  \item[] $B_3=(0,1,0)^A$ mit $|B_3|=p-1$,
  \item[] $B_4=(0,0,1)^A$ mit $|B_4|=p-1$,
  \item[] $B_5=(0,1,1)^A$ mit $|B_5|=p(p-1)^2$.
\end{enumerate}
Die Repräsentanten der Bahnen $B_1$ und $B_2$ haben die folgenden
Stabilisatoren:
\begin{eqnarray*}
  \bar{S}_1 & = & \left\{ \begin{pmatrix}
    u_1       & 0        & 0          \\
    0         & u_1 v_2  & 0          \\
    0         & 0        & v_2    \\
  \end{pmatrix} \mid u_1, v_2 \in \F_p \text{ und } u_1 v_2 \neq 0 \right\}\\
  \bar{S}_2 & = & \left\{ \begin{pmatrix}
    u_1       & 0        & 0          \\
    0         & u_1      & 0          \\
    0         & 0        & 1        \\
  \end{pmatrix} \mid u_1 \in \F_p \text{ und } u_1 \neq 0 \right\}\\
\end{eqnarray*}
\end{lemma}

\begin{beweis}
Die Behauptung ist offensichtlich.
\end{beweis}

\begin{folg} \label{zug_5_p5}
Für die Nachfolger von $G$ der Ordnung $p^5$ bilden die folgenden Untergruppen
von $M(G)$ ein Vertretersystem der zulässigen Untergruppen:
\begin{enumerate}
  \item $M_1=\erz{a_6,a_7}$ entsprechend $\erz{(0,1,0),(0,0,1)}$,
  \item $M_2=\erz{a_7,a_5 a_6\1}$ entsprechend $\erz{(0,0,1),(1,-1,0)}$.
\end{enumerate}
\end{folg}

\begin{beweis}
Der Unterraum $\bar{N}=\erz{(1,0,0)}$ entspricht dem Nukleus. Nur die Elemente
von $B_1$ und $B_2$ haben Skalarprodukt ungleich Null mit den Elementen von
$\bar{N}$. Daher entsprechen nur solche zweidimensionale Unterräume von $\f^3$
zulässigen Untergruppen von $M(G)$, die orthogonale Komplemente zu
eindimensionalen Unterräumen sind, deren Basisvektoren in $B_1$ oder $B_2$
enthalten sind.
\end{beweis}

\subsection{Nachfolger der Ordnung $p^5$}

\begin{satz} \label{nt5}
In der folgenden Liste sind sämtliche unmittelbaren Nachfolger von $G$ der
Ordnung $p^5$ bis auf Isomorphie angegeben:
\begin{enumerate}
  \item $G_1=\erz{g_1,g_2,g_3,g_4,g_5 \mid
  g_1^p=g_3, g_3^p=g_4, g_4^p=g_5} \cong C_{p^4} \times C_p$,
  \item $G_2=\erz{g_1,g_2,g_3,g_4,g_5 \mid
  g_1^p=g_3, g_3^p=g_4, g_4^p=g_5, [g_2,g_1]=g_5}$.
\end{enumerate}
\end{satz}

\begin{beweis}
Dieser Satz ergibt sich nach \ref{bahnenzu}, indem man $P(G)$ nach jedem Repräsentanten des
Vertretersystems zulässiger Untergruppen faktorisiert, das in \ref{zug_5_p5} aufgelistet ist.
Die Präsentation der Faktorgruppe erhält man über das Verfahren aus \ref{faktor}.
\end{beweis}

\begin{beme}
Die beiden Gruppen $G_1$ und $G_2$ sind zusammen mit $P(G)$ alle unmittelbaren
Nachfolger von $G$, denn für Nachfolger der Ordnung $p^6$ müsste es Supplemente
der Ordnung $p$ des Nukleus geben, die es offensichtlich nicht gibt.
\end{beme}

\section{Nachfolger von $(p^4,6)$}

\begin{ver}
In diesem Abschnitt sei $$G=\erz{a_1,a_2,a_3,a_4 \mid a_1^p=a_3, a_3^p=a_4, [a_2,a_1]=a_4}.$$
\end{ver}

\begin{lemma}
Die Gruppe $G$ besitzt die Gewichtung $\omega(a_1)=\omega(a_2)=1$, $\omega(a_3)=2$ und
$\omega(a_4)=3$ und damit die $p$"=Klasse 3.
\end{lemma}

\begin{beweis}
Die Gewichtung von $G$ lässt sich nach \ref{gew} ermitteln.
\end{beweis}

\subsection{Die Automorphismengruppe}

\begin{lemma}
Jeder Automorphismus $\al$ von $G$ lässt sich in folgender Weise darstellen: \label{aut6}
$$\al(u_1, u_2, u_3, u_4, v_4):
  \begin{cases}
    a_1 \mt a_1^{u_1} a_2^{u_2} a_3^{u_3} a_4^{u_4} & \text{wobei } u_1, u_2, u_3, u_4, v_4 \\
    a_2 \mt a_2 a_4^{v_4} & \in \{0 \dd p-1\} \text{ und }\\
    a_3 \mt a_3^{u_1} a_4^{u_3}  & v_2 \neq 0 \text{ ist.}\\
    a_4 \mt a_4^{u_1} &
  \end{cases}$$
\end{lemma}

\begin{beweis}
Die Automorphismengruppe von $G$ ist nach \ref{aut} aus dem Stabilisator $\hat{S}_5$ (siehe
\ref{bahnen_cp2cp}) der zu $G$ gehörenden zulässigen Untergruppe unmittelbar ablesbar. Aus
$\hat{S}_5$ lassen sich die Bilder von $a_1$ und $a_2$ unmittelbar ablesen. Die Bilder von
$a_3$ und $a_4$ sind durch die Bilder von $a_1$ und $a_2$ eindeutig bestimmt und lassen sich
folgendermaßen berechnen:
\begin{eqnarray*}
  a_3^\al &=& (a_1^p)^\al=(a_1^\al)^p\\
          &=& (a_1^{u_1} a_2^{u_2} a_3^{u_3} a_4^{u_4})^p= (a_1^{u_1} a_2^{u_2})^p (a_3^{u_3})^p (a_4^{u_4})^p\\
          &=& (a_1^{u_1})^p (a_3^{u_3})^p = a_3^{u_1} a_4^{u_3},
\end{eqnarray*} da $\erz{a_3, a_4}$ das Zentrum von $G$ ist und $a_2$ die Ordnung $p$
hat. Weiterhin ist
\begin{eqnarray*}
  a_4^\al &=& (a_3^p)^\al=(a_3^\al)^p\\
          &=& (a_3^{u_1} a_4^{u_3})^p = (a_3^{u_1})^p (a_4^{u_3})^p\\
          &=& (a_3^p)^{u_1} (a_4^p)^{u_3} =  a_4^{u_1}
\end{eqnarray*} aus demselben Grund.
\end{beweis}

\subsection{$p$"=Cover, Multiplikator und Nukleus}

\begin{lemma}
Für $G$ ergibt sich das $p$"=Cover, der Nukleus und der Multiplikator in folgender Weise:
\begin{enumerate}
  \item $P(G)=\erz{a_1,a_2,a_3,a_4,a_5,a_6 \mid
  a_1^p=a_3, a_3^p=a_4, [a_2,a_1] = a_4 a_5, a_2^p=a_6}$.
  \item $M(G)=\erz{a_5, a_6, a_7}$.
  \item $N(G)=\erz{1}$.
\end{enumerate}
Die Gruppe $P(G)$ hat die $p$"=Klasse 3.
\end{lemma}

\begin{beweis}
Die Präsentation von $P(G)$ ergibt sich aus \ref{cov} und \ref{newman} sowie dem
Reduktionsverfahren von Knuth"=Bendix, aus der sich die Untergruppe $M(G)$ ablesen lässt. Die
Gewichtung von $P(G)$
\begin{enumerate}
  \item $1=\omega(a_1)=\omega(a_2)$.
  \item $2=\omega(a_3)=\omega(a_5)=\omega(a_6)$.
  \item $3=\omega(a_4)$.
\end{enumerate}
und damit $N(G)$ lassen sich nach \ref{gew} aus der Präsentation von $P(G)$ ablesen.
\end{beweis}

\begin{satz}
Die Gruppe $G$ ist abschließend.
\end{satz}

\begin{beweis}
Da $P(G)$ dieselbe $p$"=Klasse hat wie $G$, hat $G$ nach \ref{pfort} keine unmittelbaren
Nachfolger, d.\,h. $G$ ist abschließend.
\end{beweis}

\section{Nachfolger von $Q_p$}

\begin{ver}
In diesem Abschnitt sei $$G=\erz{a_1, a_2, a_3 \mid [a_2, a_1]=a_3, a_1^p=a_3}.$$ Damit ist $G
\cong Q_p$.
\end{ver}

\begin{lemma}
Die Gruppe $G$ hat die Gewichtung $\omega(a_1)=\omega(a_2)=1$ und $\omega(a_3)=2$.
\end{lemma}

\begin{beweis}
Die Gewichtung von $G$ lässt sich nach \ref{gew} ermitteln.
\end{beweis}

\subsection{Die Automorphismengruppe}

\begin{lemma} \label{autqp}
Jeder Automorphismus $\al$ von $G$ lässt sich darstellen als $$\al(u_1, u_2, u_3, v_3):
  \begin{cases}
    a_1 \mt a_1^{u_1} a_2^{u_2} a_3^{u_3} & \text{wobei } u_1, u_2, u_3, v_3 \\
    a_2 \mt a_2 a_3^{v_3} & \in \{0 \dd p-1\} \text{ und }\\
    a_3 \mt a_3^{u_1} & u_1 \neq 0 \text{ ist.}
  \end{cases}$$
\end{lemma}

\begin{beweis}
Nach \ref{aut} lässt sich $Aut(G)$ unmittelbar aus dem Stabilisator $T_3$ (siehe
\ref{stab_cp2_p3}) der zu $G$ gehörenden zulässigen Untergruppe ablesen. Damit ergeben sich
genau die oben angegebenen Bilder von $a_1$ und $a_2$. Das Bild von $a_3$ ist durch die Bilder
von $a_1$ und $a_2$ vollständig festgelegt und lässt sich folgendermaßen berechnen:
\begin{eqnarray*}
  a_3^\al &=& (a_1^p)^\al=(a_1^{u_1} a_2^{u_2} a_3^{u_3})^p \\
          &=& (a_1^{u_1} a_2^{u_2})^p (a_3^{u_3})^p= (a_1^{u_1} a_2^{u_2})^p \\
          &=& (a_1^{u_1})^p (a_2^{u_2})^p =  (a_1^p)^{u_1} = a_3^{u_1},
\end{eqnarray*}
da $a_3$ der Kommutator von $a_1$ und $a_2$ ist und die Ordnung $p$ hat.
\end{beweis}

\subsection{$p$"=Cover, Multiplikator und Nukleus}

\begin{lemma}
Für $G$ ergibt sich das $p$"=Cover, der Nukleus und der Multiplikator in folgender Weise mit
der Gewichtung $\omega$:
\begin{enumerate}
  \item[] $P(G)=\erz{a_1,a_2,a_3,a_4,a_5 \mid [a_2,a_1]=a_3, a_1^p=a_3 a_4, a_2^p=a_5}$.
  \item[] $M(G)=\erz{a_4, a_5}$.
  \item[] $N(G)=\erz{1}$.
\end{enumerate}
Die Gruppe $P(G)$ hat die $p$"=Klasse 2.
\end{lemma}

\begin{beweis}
Die Präsentation von $P(G)$ ergibt sich aus \ref{cov} und \ref{newman} und dem
Reduktionsverfahren von Knuth"=Bendix, aus der sich $M(G)$ ablesen lässt. Die Gewichtung von
$P(G)$
\begin{enumerate}
  \item[] $1=\omega(a_1)=\omega(a_2)$.
  \item[] $2=\omega(a_3)=\omega(a_4)=\omega(a_5)$.
\end{enumerate}
und damit $N(G)$ lassen sich nach \ref{gew} aus der Präsentation von $P(G)$ ablesen.
\end{beweis}

\begin{satz}
Die Gruppe $G$ ist abschließend.
\end{satz}

\begin{beweis}
Da $P(G)$ dieselbe $p$"=Klasse hat wie $G$, hat $G$ nach \ref{pfort} keine unmittelbaren
Nachfolger, d.\,h. $G$ ist abschließend.
\end{beweis}

\section{Nachfolger von $(p^4,2)$ bzw. $C_{p^2} \times C_{p^2}$}

\begin{ver}
In diesem Abschnitt sei $$G=\erz{a_1, a_2, a_3, a_4 \mid a_1^p=a_3, a_2^p=a_4}.$$ Damit ist
$G$ isomorph zu $C_{p^2}^2$, bzw. $G$ ist vom \gap"=Typ $(p^4,2)$.
\end{ver}

\begin{lemma}
Die abelsche Gruppe $G$ hat die Gewichtung $\omega(a_1)=\omega(a_2)=1$ und $\omega(a_3)=
\omega(a_4)=2$ und damit die $p$"=Klasse 2.
\end{lemma}

\begin{beweis}
Die Gewichtung von $G$ lässt sich nach \ref{gew} ermitteln.
\end{beweis}

\subsection{Die Automorphismengruppe}

\begin{lemma} \label{aut2}
Jeder Automorphismus $\al$ von $G$ lässt sich durch $$\al(u_1, u_2, u_3, u_4, v_1, v_2, v_3,
v_4):
  \begin{cases}
    a_1 \mt a_1^{u_1} a_2^{u_2} a_3^{u_3} a_4^{u_4} & \text{wobei } u_1, u_2, u_3, u_4, v_1, v_2, v_3, v_4 \\
    a_2 \mt a_1^{v_1} a_2^{v_2} a_3^{v_3} a_4^{v_4} & \in \{0 \dd p-1\} \text{ und }\\
    a_3 \mt a_3^{u_1} a_4^{u_2} & u_1 v_2 - v_1 u_2 \neq 0 \text{ ist.}\\
    a_4 \mt a_3^{v_1} a_4^{v_2} &
  \end{cases}$$
darstellen.
\end{lemma}

\begin{beweis}
Nach \ref{aut} lässt sich $Aut(G)$ unmittelbar aus dem Stabilisator $S_1$ (siehe
\ref{bahnen_cp_2a}) der zu $G$ gehörenden zulässigen Untergruppe ablesen. Damit ergeben sich
genau die oben angegebenen Bilder von $a_1$ und $a_2$. Die Bilder von $a_3$ und $a_4$ sind
durch die Bilder von $a_1$ und $a_2$ vollständig festgelegt:
\begin{eqnarray*}
  a_3^\al&=& (a_1^p)^\al=(a_1^\al)^p\\
  &=& (a_1^{u_1} a_2^{u_2} a_3^{u_3} a_4^{u_4})^p \\
  &=& (a_1^{u_1})^p (a_2^{u_2})^p (a_3^{u_3})^p (a_4^{u_4})^p\\
  &=& (a_1^{u_1})^p (a_2^{u_2})^p= (a_1^p)^{u_1} (a_2^p)^{u_2}\\
  &=& a_3^{u_1} a_4^{u_2}
\end{eqnarray*}
und $a_4^\al = a_3^{v_1} a_4^{v_2}$ ebenso.
\end{beweis}

\subsection{$p$"=Cover, Multiplikator und Nukleus}

\begin{lemma}
Für $G$ ergibt sich das $p$"=Cover, der Nukleus und der Multiplikator in
folgender Weise mit der Gewichtung $\omega$:
\begin{enumerate}
  \item[] $P(G)=\erz{a_1,a_2,a_3,a_4,a_5,a_6,a_7 \mid a_1^p=a_3,
  a_2^p=a_4, a_3^p=a_5, a_4^p=a_6, [a_2,a_1]=a_7}$.
  \item[] $M(G)=\erz{a_5, a_6, a_7}$.
  \item[] $N(G)=\erz{a_5, a_6}$.
\end{enumerate}
\end{lemma}

\begin{beweis}
Die Präsentation von $P(G)$ ergibt sich aus \ref{cov} und \ref{newman} und dem
Reduktionsverfahren von Knuth"=Bendix, aus der sich $M(G)$ ablesen lässt. Die Gewichtung von
$P(G)$
\begin{enumerate}
  \item[] $1=\omega(a_1)=\omega(a_2)$.
  \item[] $2=\omega(a_3)=\omega(a_4)=\omega(a_7)$.
  \item[] $3=\omega(a_5)=\omega(a_6)$.
\end{enumerate}
und damit $N(G)$ lassen sich nach \ref{gew} aus der Präsentation von $P(G)$ ablesen.
\end{beweis}

\subsection{Operation der Erweiterungsautomorphismen}

\begin{lemma}
Stellt man die Automorphismen $\al$ von $G$ gemäß \ref{aut2}in der Art $$\al(u_1, u_2, u_3,
u_4, v_1, v_2, v_3, v_4):
  \begin{cases}
    a_1 \mt a_1^{u_1} a_2^{u_2} a_3^{u_3} a_4^{u_4} & \text{wobei } u_1, u_2, u_3, u_4, v_1, v_2, v_3, v_4 \\
    a_2 \mt a_1^{v_1} a_2^{v_2} a_3^{v_3} a_4^{v_4} & \in \{0 \dd p-1\} \text{ und }\\
    a_3 \mt a_3^{u_1} a_4^{u_2} & u_1 v_2 - v_1 u_2 \neq 0 \text{ ist.}\\
    a_4 \mt a_3^{v_1} a_4^{v_2} &
  \end{cases}$$
dar, so operiert die Automorphismengruppe über den Operationshomomorphismus $$\vi: Aut(G) \to
GL(3,p): \al(u_1, u_2, u_3, u_4, v_1, v_2, v_3, v_4) \mt M
=
\begin{pmatrix}
  u_1 & u_2 & 0  \\
  v_1 & v_2 & 0  \\
  0 & 0 & u_1 v_2 - v_1 u_2   \\
\end{pmatrix}$$
durch Erweiterungsautomorphismen auf $\f^3 \cong M(G)$.
\end{lemma}

\begin{beweis}
Nach \ref{aut2} lässt sich jeder Automorphismus $\al$ von $G$ in der oben angegebenen Weise
darstellen. Für die Erzeuger von $M(G)$ ergeben sich folgende Bilder:
\begin{eqnarray*}
  a_5^\al &=& (a_3^p)^\al=(a_3^\al)^p=(a_3^{u_1} a_4^{u_2})^p\\
      &=& (a_3^{u_1})^p (a_4^{u_2})^p = (a_3^p)^{u_1} (a_4^p)^{u_2} \\
      &=& a_5^p a_6^{u_2},
\end{eqnarray*}
da der Kommutator von $a_3$ und $a_4$ trivial ist. Aus demselben Grund ist $a_6^\al=a_5^{v_1}
a_6^{v_2}$. Weiterhin gilt
\begin{eqnarray*}
  a_7^\al &=& [a_2, a_1]^\al=[a_1^{v_1} a_2^{v_2} a_3^{v_3} a_4^{v_4}, a_1^{u_1} a_2^{u_2} a_3^{u_3} a_4^{u_4}]\\
      &=& [a_1^{v_1} a_2^{v_2}, a_1^{u_1} a_2^{u_2}]= [a_2, a_1]^{u_1 v_2 - v_1 u_2},
\end{eqnarray*}
da alle anderen neben $[a_2, a_1]$ auftretenden Kommutatoren trivial sind und der Kommutator
$[a_2, a_1]$ zentral ist.
\end{beweis}

Für die Nachfolger der Ordnung $p^5$ werden zulässige Untergruppen der Ordnung $p^2$ gesucht.
Diesen entsprechen zweidimensionale Unterräume in $\F_p^3$. Der Einfachheit halber werden ihre
eindimensionalen Komplemente betrachtet. Wegen der Symmetrie von $Aut(G)^\vi$ erübrigt sich
das Transponieren.

\begin{beme}
Der Operationshomomorphismus $\vi$ ist zu sich selbst dual.
\end{beme}

\begin{lemma}
Der Vektorraum $\f^3 \cong M(G)$ zerfällt unter der Operation von $Aut(G)$ über $\vi$ in
folgende Bahnen nichttrivialer Vektoren:
\begin{enumerate}
  \item[] $B_1=(1,0,0)^A$ mit $|B_1|=p^2-1$.
  \item[] $B_2=(1,0,1)^A$ mit $|B_2|=(p^2-1)(p-1)$.
  \item[] $B_3=(0,0,1)^A$ mit $|B_3|=p-1$.
\end{enumerate}
Die Stabilisatoren zu den Vertretern aus $B_1$ und $B_2$ sind
\begin{eqnarray*}
  \hat{S}_1 & = & \left\{ \begin{pmatrix}
    1         & 0        & 0          \\
    v_1       & v_2      & 0          \\
    0         & 0        & v_2        \\
  \end{pmatrix} \mid v_1, v_2 \in \F_p \text{ und } v_2 \neq 0 \right\}
  \text{ mit } |\hat{S}_1|=p(p-1),\\
  \hat{S}_2 & = & \left\{ \begin{pmatrix}
    1         & 0        & 0          \\
    v_1       & 1        & 0          \\
    0         & 0        & 1        \\
  \end{pmatrix} \mid v_1 \in \F_p \right\} \text{ mit } |\hat{S}_2|=p,\\
\end{eqnarray*}
und die Stabilisatoren der entsprechenden eindimensionalen
Unterräume sind
\begin{eqnarray*}
  \bar{S}_1 & = & \left\{ \begin{pmatrix}
    u_1       & 0        & 0          \\
    v_1       & v_2      & 0          \\
    0         & 0        & u_1 v_2    \\
  \end{pmatrix} \mid u_1, v_1, v_2 \in \F_p \text{ und } u_1 v_2 \neq 0 \right\}\\
  \bar{S}_2 & = & \left\{ \begin{pmatrix}
    u_1       & 0        & 0          \\
    v_1       & 1        & 0          \\
    0         & 0        & u_1        \\
  \end{pmatrix} \mid u_1, v_1 \in \F_p \right\}\\
\end{eqnarray*}
\end{lemma}

\begin{beweis}
Der Beweis erfolgt analog zu dem von \ref{bahnen_cp_2}.
\end{beweis}

\begin{folg} \label{zu3}
Für die Nachfolger von $G$ der Ordnung $p^4$ ergibt sich das folgende
Vertretersystem zulässiger Untergruppen:
\begin{enumerate}
  \item[] $M_1=\erz{a_6, a_7}$ entsprechend $\erz{(0,1,0),(0,0,1)}$.
  \item[] $M_2=\erz{a_6, a_5 a_7\1}$ entsprechend $\erz{(0,1,0),(1,0,-1)}$.
\end{enumerate}
\end{folg}

\begin{beweis}
Die Bahn $B_3$ ist nicht von Interesse, da ihr Vertreter mit $(1,0,0)$ das
Skalarprodukt Null liefert. Also entsprechen den Komplemente zu den
eindimensionalen Unterräumen mit Basisverktoren aus $B_3$ keine zulässigen
Untergruppen.
\end{beweis}

\subsection{Nachfolger der Ordnung $p^5$}

\begin{satz} \label{nt2}
In der folgenden Liste sind sämtliche unmittelbaren Nachfolger von $G$ der
Ordnung $p^5$ bis auf Isomorphie angegeben:
\begin{enumerate}
  \item[] $G_1=\erz{g_1,g_2,g_3,g_4,g_5 \mid g_1^p=g_3,
  g_2^p=g_4, g_3^p=g_5} \cong C_{p^3} \times C_{p^2}$.
  \item[] $G_2=\erz{g_1,g_2,g_3,g_4, g_5 \mid g_1^p=g_3,
  g_2^p=g_4, g_3^p=g_5, [g_2,g_1]=g_5}$.
\end{enumerate}
\end{satz}

\begin{beweis}
Dieser Satz ergibt sich nach \ref{bahnenzu}, indem man $P(G)$ nach jedem Repräsentanten des
Vertretersystems zulässiger Untergruppen faktorisiert, das in \ref{zu3} aufgelistet ist. Die
Präsentation der Faktorgruppe erhält man über das Verfahren aus \ref{faktor}.
\end{beweis}

\section{Nachfolger von $(p^4,3)$}

\begin{ver}
In diesem Abschnitt sei $$G=\erz{a_1, a_2, a_3, a_4 \mid [a_2, a_1]=a_3, a_1^p=a_4}.$$
\end{ver}

\begin{lemma}
Die Gruppe $G$ hat die Gewichtung $\omega(a_1)=\omega(a_2)=1$ und $\omega(a_3)=\omega(a_4)=2$
und damit die $p$"=Klasse 2.
\end{lemma}

\begin{beweis}
Die Gewichtung von $G$ lässt sich nach \ref{gew} ermitteln.
\end{beweis}

\subsection{Die Automorphismengruppe}

\begin{lemma}  \label{aut3}
Jeder Automorphismus $\al$ von $G$ lässt sich durch $$\al(u_1, u_2, u_3, u_4, v_2, v_3, v_4):
  \begin{cases}
    a_1 \mt a_1^{u_1} a_2^{u_2} a_3^{u_3} a_4^{u_4} & \text{wobei } u_1, u_2, u_3, u_4, v_2, v_3, v_4 \\
    a_2 \mt a_2^{v_2} a_3^{v_3} a_4^{v_4} & \in \{0 \dd p-1\} \text{ und }\\
    a_3 \mt a_3^{u_1 v_2}  & u_1 \neq 0 \text{ sowie } v_2 \neq 0 \text{ ist.}\\
    a_4 \mt a_4^{u_1} &
  \end{cases}$$
darstellen.
\end{lemma}

\begin{beweis}
Nach \ref{aut} lässt sich $Aut(G)$ unmittelbar aus dem Stabilisator $S_2$ (siehe
\ref{bahnen_cp_2a}) der zu $G$ gehörenden zulässigen Untergruppe ablesen. Damit ergeben sich
genau die oben angegebenen Bilder von $a_1$ und $a_2$. Die Bilder von $a_3$ und $a_4$ sind
durch die Bilder von $a_1$ und $a_2$ vollständig festgelegt:
\begin{eqnarray*}
  a_3^\al&=&[a_2,a_1]^\al=[a_2^\al,a_1^\al] \\
  &=& [a_2^{v_2} a_3^{v_3} a_4^{v_4}, a_1^{u_1} a_2^{u_2} a_3^{u_3} a_4^{u_4}]\\
  &=& [a_2^{v_2}, a_1^{u_1} a_2^{u_2}] = [a_2^{v_2}, a_1^{u_1}] \\
  &=& [a_2, a_1]^{u_1 v_2} = a_3^{u_1 v_2}
\end{eqnarray*}
da $\langle a_3, a_4 \rangle$ das Zentrum von $G$ ist. Weiterhin ist
\begin{eqnarray*}
  a_4^\al &=& (a_1^p)^\al = (a_1^{u_1} a_2^{u_2} a_3^{u_3} a_4^{u_4})^p\\
  &=& (a_1^{u_1} a_2^{u_2})^p (a_3^{u_3})^p (a_4^{u_4})^p =
  (a_1^{u_1} a_2^{u_2})^p (a_3^p)^{u_3} (a_4^p)^{u_4}\\
  &=& (a_1^{u_1} a_2^{u_2})^p = (a_1^{u_1})^p = a_3^{u_1},
\end{eqnarray*}
da $|a_2|=|a_3|=|a_4|=p$ ist und es ein $k \in \{0 \dd p-1\}$ gibt, sodass $$(a_1^{u_1}
a_2^{u_2})^p = (a_1^{u_1})^p (a_3^k)^p (a_2^{u_2})^p = (a_1^p)^{u_1} (a_3^p)^k (a_2^p)^{u_2} =
(a_1^p)^{u_1}$$ ist, da $a_3$ der Kommutator von $a_1$ und $a_2$ ist.
\end{beweis}

\subsection{$p$"=Cover, Multiplikator und Nukleus}

\begin{lemma}
Für $G$ ergibt sich das $p$"=Cover, der Nukleus und der Multiplikator in folgender Weise:
\begin{enumerate}
  \item[] $P(G)=\langle a_1,a_2,a_3,a_4,a_5,a_6,a_7,a_8
  \mid [a_2,a_1]=a_3, a_1^p=a_4, [a_3,a_1]=a_5, [a_3,a_2]=a_6,
  a_4^p=a_7,a_2^p=a_8 \rangle$.
  \item[] $M(G)=\langle a_5,a_6,a_7,a_8 \rangle$.
  \item[] $N(G)=\langle a_5,a_6,a_7 \rangle$.
\end{enumerate}
\end{lemma}

\begin{beweis}
Die Präsentation von $P(G)$ ergibt sich aus \ref{cov} und \ref{newman} unter Verwendung des
Reduktionsverfahrens von Knuth"=Bendix, aus der sich $M(G)$ ablesen lässt. Die Gewichtung von
$P(G)$
\begin{enumerate}
  \item[] $1=\omega(a_1)=\omega(a_2)$.
  \item[] $2=\omega(a_3)=\omega(a_4)=\omega(a_8)$.
  \item[] $3=\omega(a_5)=\omega(a_6)=\omega(a_7)$.
\end{enumerate}
und damit $N(G)$ lässt sich nach \ref{gew} aus der Präsentation von $P(G)$ ablesen.
\end{beweis}

\subsection{Operation der Erweiterungsautomorphismen}

\begin{lemma}
Stellt man die Automorphismen $\al$ von $G$ gemäß \ref{aut3} in der Form $$\al(u_1, u_2, u_3,
u_4, v_2, v_3, v_4):
  \begin{cases}
    a_1 \mt a_1^{u_1} a_2^{u_2} a_3^{u_3} a_4^{u_4} & \text{wobei } u_1, u_2, u_3, u_4, v_2, v_3, v_4 \\
    a_2 \mt a_2^{v_2} a_3^{v_3} a_4^{v_4} & \in \{0 \dd p-1\} \text{ und }\\
    a_3 \mt a_3^{u_1 v_2}  & u_1 \neq 0 \text{ sowie } v_2 \neq 0 \text{ ist.}\\
    a_4 \mt a_4^{u_1} &
  \end{cases}$$
dar, so operiert die Automorphismengruppe von $G$ über den Operationshomomorphismus $$\vi:
Aut(G) \ra GL(4,p): \al(u_1, u_2, u_3, u_4, v_2, v_3, v_4) \mt
\begin{pmatrix}
  u_1^2 v_2 & u_1 u_2 v_2 & 0   & 0  \\
  0         & u_1 v_2^2   & 0   & 0  \\
  0         & 0           & u_1 & 0  \\
  0         & 0           & v_4 & v_2
\end{pmatrix}$$ durch Erweiterungsautomorphismen auf $\f^4 \cong M(G)$.
\end{lemma}

\begin{beweis}
Nach \ref{aut3} lässt sich jeder Automorphismus $\al$ von $G$ in der oben angegebenen Weise
darstellen. Für die Erzeuger von $M(G)$ ergeben sich folgende Bilder:
\begin{eqnarray*}
  a_5^\al &=& [a_3, a_1]^\al=[a_3^\al, a_1^\al]\\
      &=& [a_3^{u_1 v_2}, a_1^{u_1} a_2^{u_2} a_3^{u_3} a_4^{u_4}]= [a_3^{u_1 v_2}, a_1^{u_1} a_2^{u_2} a_3^{u_3}]^{a_4^{u_4}}\\
      &=& ([a_3^{u_1 v_2}, a_1^{u_1} a_2^{u_2}]^{a_3^{u_3}})^{a_4^{u_4}}\\
      &=& (([a_3^{u_1 v_2}, a_2^{u_2}] [a_3^{u_1 v_2}, a_1^{u_1}]^{a_1^{u_1}})^{a_3^{u_3}})^{a_4^{u_4}}\\
      &=& (([a_3, a_2]^{u_1 u_2 v_2} ([a_3, a_1]^{u_1^2 v_2})^{a_1^{u_1}})^{a_3^{u_3}})^{a_4^{u_4}}\\
      &=& [a_3, a_2]^{u_1 u_2 v_2} [a_3, a_1]^{u_1^2 v_2} = a_6^{u_1 u_2 v_2} a_5^{u_1^2 v_2}\\
      &=& a_5^{u_1^2 v_2} a_6^{u_1 u_2 v_2},
\end{eqnarray*}
da $[a_4,a_3]$ trivial ist und $[a_3,a_1]$ sowie $[a_3,a_2]$ zentral sind.
Ebenso ergibt sich $a_6^\al=a_6^{u_1 v_2^2}$. Weiterhin ist
\begin{eqnarray*}
  a_7^\al &=& (a_4^p)^\al=(a_4^\al)^p=(a_4^{u_1})^p\\
      &=& a_7^{u_1}
\end{eqnarray*} und schließlich ist
\begin{eqnarray*}
  8_3^\al &=& (a_2^p)^\al=(a_2^\al)^p\\
      &=& (a_2^{v_2} a_3^{v_3} a_4^{v_4})^p = (a_2^{v_2})^p (a_3^{v_3})^p (a_4^{v_4})^p\\
      &=& (a_2^p)^{v_2} (a_4^p)^{v_4} = a_7^{v_4} a_8^{v_2},
\end{eqnarray*}
da $|a_3|=p$ und $[a_4,a_2]=[a_4, a_3]=1$ ist und $a_6$, der Kommutator von
$a_2$ und $a_3$, die Ordnung $p$ hat.
\end{beweis}

Da im Weiteren für die unmittelbaren Nachfolger der Ordnung $p^5$ statt der dreidimensionalen
Unterräume ihre eindimensionalen Komplemente betrachtet werden, ist der zu $\vi$ duale
Operationshomomorphismus von Interesse.

\begin{beme}
Der zu $\vi$ duale Operationshomomorphismus ist $$\bar{\vi}: Aut(G) \ra GL(4,p): \al(u_1, u_2,
u_3, u_4, v_2, v_3, v_4) \mt
\begin{pmatrix}
  u_1^2 v_2   & 0           & 0    & 0  \\
  u_1 u_2 v_2 & u_1 v_2^2   & 0    & 0  \\
  0           & 0           & u_1  & v_4  \\
  0           & 0           & 0    & v_2
\end{pmatrix}.$$
\end{beme}

\subsection{Bahnen zulässiger Untergruppen für Nachfolger der Ordnung $p^5$}

\begin{lemma}
Unter der Operation von $A=Aut(G)$ vermittels $\bar{\vi}$ in folgende Bahnen nichttrivialer
Vektoren, wobei $w$ ein Erzeuger der multiplikativen Gruppe von $\F_p$ ist:
\begin{enumerate}
  \item[] $B_1=(0,0,1,0)$ mit $|B_1|=p(p-1)$.
  \item[] $B_2=(1,0,0,0)$ mit $|B_2|=p-1$.
  \item[] $B_3=(1,0,0,1)$ mit $|B_3|=\frac{(p-1)^2}{2}$.
  \item[] $B_4=(w,0,0,1)$ mit $|B_4|=\frac{(p-1)^2}{2}$.
  \item[] $B_5=(1,0,1,0)$ mit $|B_5|=p(p-1)^2$.
  \item[] $B_6=(1,1,0,0)$ mit $|B_6|=p(p-1)$.
  \item[] $B_7=(1,1,0,1)$ mit $|B_7|=p(p-1)^2$.
  \item[] $B_8=(1,1,1,0)$ mit $|B_8|=\frac{p^2(p-1)^2}{2}$.
  \item[] $B_9=(w,w,1,0)$ mit $|B_9|=\frac{p^2(p-1)^2}{2}$.
  \item[] $B_{10}=(0,0,0,1)$ mit $|B_{10}|=p-1$.
\end{enumerate}
Die zu den Bahnenvertretern gehörenden Stabilisatoren sind
\begin{eqnarray*}
  \hat{S}_1 & = & \left\{ \begin{pmatrix}
    v_2       & 0        & 0      &  0        \\
    u_2 v_2   & v_2^2    & 0      &  0        \\
    0         & 0        & 1      &  0        \\
    0         & 0        & 0      &  v_2^2    \\
  \end{pmatrix} \mid u_2, v_2 \in \F_p  \text{ und } v_2 \neq 0\right\}
  \text{ mit } |\hat{S}_1|=p(p-1).\\
    \hat{S}_2 & = & \left\{ \begin{pmatrix}
    1         & 0        & 0      &  0        \\
    u_2 u_1\1 & u_1\1    & 0      &  0        \\
    0         & 0        & u_1    &  v_4      \\
    0         & 0        & 0      &  u_2^{-2} \\
  \end{pmatrix} \mid u_1, u_2, v_4 \in \F_p  \text{ und } u_1 \neq 0\right\}
  \text{ mit } |\hat{S}_2|=p^2(p-1).\\
    \hat{S}_3 = \hat{S}_4 & = & \left\{ \begin{pmatrix}
    1         & 0        & 0      &  0        \\
    \pm u_2   & \pm 1    & 0      &  0        \\
    0         & 0        & \pm 1  &  v_4      \\
    0         & 0        & 0      &  1        \\
  \end{pmatrix} \mid v_4 \in \F_p \right\}
  \text{ mit } |\hat{S}_3|=\hat{S}_4|=2p^2.\\
    \hat{S}_5 & = & \left\{ \begin{pmatrix}
    1         & 0        & 0      &  0        \\
    u_2       & 1        & 0      &  0        \\
    0         & 0        & 1      &  0        \\
    0         & 0        & 0      &  1    \\
  \end{pmatrix} \mid u_2 \in \F_p  \right\}
  \text{ mit } |\hat{S}_5|=p.\\
    \hat{S}_6 & = & \left\{ \begin{pmatrix}
    v_2^{-3}     & 0        & 0        &  0        \\
    1 - v_2^{-3} & 1        & 0        &  0        \\
    0            & 0        & v_2^{-2} &  v_4        \\
    0            & 0        & 0        &  v_2   \\
  \end{pmatrix} \mid v_2, v_4 \in \F_p  \text{ und } v_2 \neq 0\right\}
  \text{ mit } |\hat{S}_6|=p(p-1).\\
    \hat{S}_7 & = & \left\{ \begin{pmatrix}
    1         & 0        & 0      &  0        \\
    0         & 1        & 0      &  0        \\
    0         & 0        & 1      &  v_4        \\
    0         & 0        & 0      &  1    \\
  \end{pmatrix} \mid v_4 \in \F_p \right\}
  \text{ mit } |\hat{S}_7|=p.\\
    \hat{S}_8 = \hat{S}_9 & = & \left\{ \begin{pmatrix}
    \pm 1     & 0        & 0      &  0        \\
    1 \mp 1   & 1        & 0      &  0        \\
    0         & 0        & 1      &  0        \\
    0         & 0        & 0      &  \pm 1    \\
  \end{pmatrix}\right\}
  \text{ mit } |\hat{S}_8|=|\hat{S}_9|=2.\\
  \hat{S}_{10} & = & \left\{ \begin{pmatrix}
    u_1^2     & 0        & 0      &  0        \\
    u_1 u_2   & u_1      & 0      &  0        \\
    0         & 0        & u_1    &  v_4      \\
    0         & 0        & 0      &  1        \\
  \end{pmatrix} \mid u_1, u_2, v_4 \in \F_p  \text{ und } u_1 \neq 0\right\}
  \text{ mit } |\hat{S}_{10}|=p^2(p-1).\\
\end{eqnarray*}
Die Stalisatoren der eindimensionalen Unterräume, die von den Bahnenvertretern
aufgespannt werden sind in Auswahl
\begin{eqnarray*}
  \bar{S}_1 & = & \left\{ \begin{pmatrix}
    u_1^2 v_2   & 0        & 0      &  0        \\
    u_1 u_2 v_2 & v_2^2    & 0      &  0        \\
    0           & 0        & 1      &  0        \\
    0           & 0        & 0      &  v_2^2    \\
  \end{pmatrix} \mid u_1, u_2, v_2 \in \F_p  \text{ und } u_1, v_2 \neq
  0\right\} \\
    \bar{S}_2 & = & A \\
    \bar{S}_3 = \bar{S}_4 & = & \left\{ \begin{pmatrix}
    v_2          & 0         & 0      &  0        \\
    \pm u_2 v_2  & \pm v_2^2 & 0      &  0        \\
    0            & 0         & \pm 1  &  v_4      \\
    0            & 0         & 0      &  v_2        \\
  \end{pmatrix} \mid v_2, v_4 \in \F_p \text{ und } v_2 \neq 0 \right\} \\
    \bar{S}_5 & = & \left\{ \begin{pmatrix}
    u_1       & 0        & 0      &  0        \\
    u_2       & u_1\1    & 0      &  0        \\
    0         & 0        & u_1    &  0        \\
    0         & 0        & 0      &  u_1\1    \\
  \end{pmatrix} \mid u_1, u_2 \in \F_p \text{ und } u_1 \neq 0 \right\} \\
    \bar{S}_6 & = & \left\{ \begin{pmatrix}
    u_1^2 v_2            & 0          & 0        &  0        \\
    u_1 v_2^2 -u_1^2 v_2 & u_1 v_2^2  & 0        &  0        \\
    0                    & 0          & u_1      &  v_4        \\
    0                    & 0          & 0        &  v_2   \\
  \end{pmatrix} \mid u_1, v_2, v_4 \in \F_p  \text{ und } u_1 v_2 \neq
  0\right\}\\
    \bar{S}_7 & = & \left\{ \begin{pmatrix}
    u_1         & 0        & 0      &  0        \\
    u_1\1 -u_1  & u_1\1    & 0      &  0        \\
    0           & 0        & u_1    &  v_4        \\
    0           & 0        & 0      &  u_1\1    \\
  \end{pmatrix} \mid u_1, v_4 \in \F_p \text{ und } u_1 \neq 0 \right\}\\
    \bar{S}_8 = \bar{S}_9 & = & \left\{ \begin{pmatrix}
    \pm u_1^2     & 0        & 0      &  0        \\
    u_1 \mp u_1^2 & u_1      & 0      &  0        \\
    0             & 0        & u_1    &  0        \\
    0             & 0        & 0      &  \pm 1    \\
  \end{pmatrix} \mid u_1 \in \f \text{ und } u_1 \neq 0 \right\}\\
\end{eqnarray*}
\end{lemma}

\begin{beweis}
Die Bahnenlängen erkennt man --~sofern sie nicht unmittelbar ersichtlich sind~-- aus der
Ordnung der Stabilisatoren $\hat{S}_1 \dd \hat{S}_{10}$, die sich aus den Lösungsmengen der
folgenden Gleichungssysteme ergeben:
\begin{enumerate}
  \item $\{u_1=1, v_4=0\}$,
  \item $\{u_1^2 v_2=1\}$,
  \item $\{u_1^2 v_2=1, v_2=1\}$,
  \item $\{u_1^2 v_2 w=w, v_2=1\}$,
  \item $\{u_1^2 v_2=1, u_1=1, v_4=0\}$,
  \item $\{u_1^2 v_2 + u_1 u_2 v_2 = 1, u_1 v_2^2 =1\}$,
  \item $\{u_1^2 v_2 + u_1 u_2 v_2 = 1, u_1 v_2^2 =1, v_2=1\}$,
  \item $\{u_1^2 v_2 + u_1 u_2 v_2 = 1, u_1 v_2^2 =1, u_1=1, v_4=0\}$,
  \item $\{(u_1^2 v_2 + u_1 u_2 v_2) w = w, u_1 v_2^2 w=w, u_1=1, v_4=0\}$,
  \item $\{v_2=1\}$.
\end{enumerate}

Die Stabilisatoren $\bar{S}_1 \dd \bar{S}_9$ erhält man, indem man in den ersten neun
Gleichungssystemen die rechten Seiten mit einer Variable des Wertebereichs $\f^*$
multipliziert und die Lösungsmengen berechnet.

Die Paare von Bahnen gleicher Länge sind nicht identisch, was man daran erkennen kann, dass
der Repräsentant der einen nicht in der Bahn des anderen liegt, denn ansonsten wäre in all
diesen Fällen eine Gleichung ableitbar, nach der $w$ ein Qudarat wäre --~was $w$ als Erzeuger
der multiplikativen Gruppe von $\f$ nach \ref{indexmult} aber nicht ist.
\end{beweis}

\begin{folg} \label{zu4}
Für die unmittelbaren Nachfolger von $G$ der Ordnung $p^5$ ist die folgende
Liste ein Vertretersystem zulässiger Untergruppen:
\begin{enumerate}
  \item[] $M_1=\erz{a_5,a_6,a_8}$ entsprechend $\erz{(1,0,0,0),(0,1,0,0),(0,0,0,1)}$,
  \item[] $M_2=\erz{a_6,a_7,a_8}$ entsprechend $\erz{(0,1,0,0),(0,0,1,0),(0,0,0,1)}$,
  \item[] $M_3=\erz{a_6,a_7,a_5 a_8\1}$ entsprechend $\erz{(0,1,0,0),(0,0,1,0),(1,0,0,-1)}$,
  \item[] $M_4=\erz{a_6,a_7,a_5 a_8^{-w}}$ entsprechend $\erz{(0,1,0,0),(0,0,1,0),(1,0,0,-w)}$,
  \item[] $M_5=\erz{a_6,a_8,a_5 a_7\1}$ entsprechend $\erz{(0,1,0,0),(0,0,0,1),(1,0,-1,0)}$,
  \item[] $M_6=\erz{a_7,a_8,a_5 a_6\1}$ entsprechend $\erz{(0,0,1,0),(0,0,0,1),(1,-1,0,0)}$,
  \item[] $M_7=\erz{a_7,a_5 a_8\1,a_6 a_8\1}$ entsprechend $\erz{(0,0,1,0),(1,0,0,-1),(0,1,0,-1)}$,
  \item[] $M_8=\erz{a_8,a_5 a_7\1,a_6 a_7\1}$ entsprechend $\erz{(0,0,0,1),(1,0,-1,0),(0,1,-1,0)}$,
  \item[] $M_9=\erz{a_8,a_5 a_7^{-w},a_6 a_8^{-w}}$ entsprechend $\erz{(0,0,0,1),(1,0,-w,0),(0,1,-w,0)}$,
\end{enumerate}
\end{folg}

\begin{beweis}
Da der Repräsentant von $B_{10}$ mit allen Basisvektoren das Skalarprodukt Null liefert, die
den Unterraum aufspannen, der dem Nukleus entspricht, stellen die Elemente von $B_{10}$ keine
zulässigen Untergruppen dar. Deshalb wird $B_{10}$ vernachlässigt. Alle anderen Bahnen
bestehen nach demselben Kriterium aus zulässigen Untergruppen.
\end{beweis}

\subsection{Nachfolger der Ordnung $p^5$}

\begin{satz} \label{nt3}
In der folgenden Liste sind sämtliche unmittelbaren Nachfolger von $G$ der
Ordnung $p^5$ bis auf Isomorphie angegeben:
\begin{enumerate}
  \item[] $G_1= \erz{g_1,g_2,g_3,g_4,g_5
  \mid [g_2,g_1]=g_3, g_1^p=g_4, g_4^p=g_5}$,
  \item[] $G_2= \erz{g_1,g_2,g_3,g_4,g_5 \mid [g_2,g_1]=g_3, g_1^p=g_4,
  [g_3,g_1]=g_5}$,
  \item[] $G_3= \erz{g_1,g_2,g_3,g_4, g_5
  \mid [g_2,g_1]=g_3, g_1^p=g_4, [g_3,g_1]=g_5, g_2^p=g_5}$,
  \item[] $G_4= \erz{g_1,g_2,g_3,g_4,g_5
  \mid [g_2,g_1]=g_3, g_1^p=g_4, [g_3,g_1]=g_5^w,g_2^p=g_5}$,
  \item[] $G_5= \erz{g_1,g_2,g_3,g_4,g_5 \mid [g_2,g_1]=g_3, g_1^p=g_4, [g_3,g_1]=g_5,
  g_4^p=g_5}$,
  \item[] $G_6= \erz{g_1,g_2,g_3,g_4,g_5
  \mid [g_2,g_1]=g_3, g_1^p=g_4, [g_3,g_1]=g_5, [g_3,g_2]=g_5}$,
  \item[] $G_7= \erz{g_1,g_2,g_3,g_4,g_5
  \mid [g_2,g_1]=g_3, g_1^p=g_4, [g_3,g_1]=g_5, [g_3,g_2]=g_5,g_2^p=g_5}$,
  \item[] $G_8= \erz{g_1,g_2,g_3,g_4,g_5
  \mid [g_2,g_1]=g_3, g_1^p=g_4, [g_3,g_1]=g_5, [g_3,g_2]=g_5,
  g_4^p=g_5}$,
  \item[] $G_9= \erz{g_1,g_2,g_3,g_4,g_5
  \mid [g_2,g_1]=g_3, g_1^p=g_4, [g_3,g_1]=g_5^w, [g_3,g_2]=g_5^w,
  g_4^p=g_5}$.
\end{enumerate}
\end{satz}

\begin{beweis}
Dieser Satz ergibt sich nach \ref{bahnenzu}, indem man $P(G)$ nach jedem Repräsentanten des
Vertretersystems zulässiger Untergruppen faktorisiert, das in \ref{zu4} aufgelistet ist. Die
Präsentation der Faktorgruppe erhält man über das Verfahren aus \ref{faktor}.
\end{beweis}

\section{Nachfolger von $(p^4,4)$}

\begin{ver}
In diesem Abschnitt sei $$G=\erz{a_1, a_2, a_3, a_4 \mid [a_2, a_1]=a_3, a_2^p=a_3,
a_1^p=a_4}.$$
\end{ver}

\begin{lemma}
Die Gruppe $G$ hat die Gewichtung $\omega(a_1)=\omega(a_2)=1$ und $\omega(a_3)=\omega(a_4)=2$
und daher die $p$"=Klasse 2.
\end{lemma}

\begin{beweis}
Die Gewichtung von $G$ lässt sich nach \ref{gew} ermitteln.
\end{beweis}

\subsection{Die Automorphismengruppe}

\begin{lemma} \label{aut4}
Die Automorphismen $\al$ von $G_3$ lassen sich durch $$\al(u_2, u_3, u_4, v_2, v_3, v_4):
  \begin{cases}
    a_1 \mt a_1 a_2^{u_2} a_3^{u_3} a_4^{u_4} & \text{wobei } u_2, u_3, u_4, v_2, v_3, v_4 \\
    a_2 \mt a_2^{v_2} a_3^{v_3} a_4^{v_4} & \in \{0 \dd p-1\} \text{ und }\\
    a_3 \mt a_3 a_4^{u_2} & v_2 \neq 0 \text{ ist.}\\
    a_4 \mt a_4^{v_2} &
  \end{cases}$$ darstellen.
\end{lemma}

\begin{beweis}
Nach \ref{aut} lässt sich $Aut(G)$ unmittelbar aus dem Stabilisator $S_3$ (siehe
\ref{bahnen_cp_2a}) der zu $G$ gehörenden zulässigen Untergruppe ablesen. Damit ergeben sich
genau die oben angegebenen Bilder von $a_1$ und $a_2$. Die Bilder von $a_3$ und $a_4$ sind
durch die Bilder von $a_1$ und $a_2$ vollständig festgelegt:\begin{eqnarray*}
  a_3^\al&=& (a_1^p)^\al = (a_1^\al)^p = (a_1 a_2^{u_2} a_3^{u_3}
  a_4^{u_4})^p\\
  &=& (a_1 a_2^{u_2})^p = a_1^p (a_2^{u_2})^p = a_3 a_4^{u_2},
\end{eqnarray*}
da $\langle a_3, a_4 \rangle$ das Zentrum von $G_3$ ist und $a_3$, der Kommutator von $a_1$
und $a_2$ die Ordnung $p$ hat. Außerdem ist
\begin{eqnarray*}
  a_4^\al &=& [a_2^\al, a_1^\al]=[a_2^{v_2} a_3^{v_3} a_4^{v_4},
  a_1 a_2^{u_2} a_3^{u_3} a_4^{u_4}] \\
  &=& [a_2^{v_2}, a_1 a_2^{u_2}] =
  [a_2^{v_2}, a_1] = [a_2, a_1]^{v_2} = a_4^{v_2},
\end{eqnarray*}
da $\langle a_3, a_4 \rangle$ das Zentrum von $G_3$ ist.
\end{beweis}

\subsection{$p$"=Cover, Multiplikator und Nukleus}

\begin{lemma}
Für $G$ ergibt sich das $p$"=Cover, der Nukleus und der Multiplikator in folgender Weise:
\begin{enumerate}
  \item[] $P(G)=\erz{a_1,a_2,a_3,a_4,a_5,a_6,a_7 \mid
  [a_2,a_1]=a_3, a_1^p=a_4, [a_3,a_1]=a_5, a_3^p=a_5, [a_4,a_2]=a_5^{p-1}, a_4^p=a_6, a_2^p=a_3 a_7}$.
  \item[] $M(G)=\langle a_5,a_6,a_7 \rangle$.
  \item[] $N(G)=\langle a_5,a_6 \rangle$.
\end{enumerate}
\end{lemma}

\begin{beweis}
Die Präsentation von $P(G)$ ergibt sich aus \ref{cov} und \ref{newman} und dem
Reduktionsverfahren von Knuth"=Bendix, aus der sich $M(G)$ ablesen lässt. Die Gewichtung von
$P(G)$
\begin{enumerate}
  \item[] $1=\omega(a_1)=\omega(a_2)$.
  \item[] $2=\omega(a_3)=\omega(a_4)=\omega(a_7)$.
  \item[] $3=\omega(a_5)=\omega(a_6)$.
\end{enumerate}
und damit $N(G)$ lassen sich nach \ref{gew} aus der Präsentation von $P(G)$ ablesen.
\end{beweis}

\subsection{Operation der Erweiterungsautomorphismen}

\begin{lemma}
Stellt man die Automorphismen $\al$ von $G$ gemäß \ref{aut4} in folgender Weise $$\al(u_2,
u_3, u_4, v_2, v_3, v_4):
  \begin{cases}
    a_1 \mt a_1 a_2^{u_2} a_3^{u_3} a_4^{u_4} & \text{wobei } u_2, u_3, u_4, v_2, v_3, v_4 \\
    a_2 \mt a_2^{v_2} a_3^{v_3} a_4^{v_4} & \in \{0 \dd p-1\} \text{ und }\\
    a_3 \mt a_3^{v_2} & v_2 \neq 0 \text{ ist.}\\
    a_4 \mt a_3^{u_2} a_4 &
  \end{cases}$$
dar, so operiert $Aut(G)$ vermittels Erweiterungsautomorphismen über den
Operationshomomorphismus $$\vi: Aut(G) \to GL(3,p): \al(u_2, u_3, u_4, v_2, v_3, v_4) \mt M =
\begin{pmatrix}
  v_2 & 0   & 0 \\
  u_2 & 1   & 0 \\
  v_3 & v_4 & v_2
\end{pmatrix}$$
auf $\f^3 \cong M(G)$.
\end{lemma}

\begin{beweis}
Gemäß \ref{aut4} lassen sich alle Automorphismen in der oben angegebenen Weise
darstellen. Die Bilder der Erzeuger von $M(G)$ ergeben sich folgendermaßen:
\begin{eqnarray*}
  a_5^\al &=& [a_3, a_1]^\al=[a_3^\al, a_1^\al]\\
      &=& [a_3^{v_2}, a_1 x_2^{u_2} a_3^{u_3} a_4^{u_4}] = [a_3^{v_2}, a_1 a_2^{u_2}]\\
      &=& [a_3^{v_2}, a_2^{u_2}] [a_3^{v_2}, a_1]^{a_2^{u_2}} = [a_3, a_1]^{v_2}\\
      &=& a_5^{v_2},
\end{eqnarray*}
da $[a_3, a_1]$ zentral und $[a_4,a_3]=[a_3,a_2]=1$ ist. Weiterhin ist
\begin{eqnarray*}
  a_6^\al &=& (a_4^p)^\al=(a_4^\al)^p\\
      &=& (a_3^{u_2} a_4)^p=(a_3^{u_2})^p a_4^p = a_5^{u_2} a_6,
\end{eqnarray*}
da $[a_4,a_3]=1$ ist. Außerdem gilt
\begin{eqnarray*}
  a_7^\al &=& (a_3\1 a_2^p)^\al=(a_3\1)^\al (a_2^\al)^p \\
      &=& a_3^{-v_2}(a_2^{v_2} a_3^{v_3} a_4^{v_4})^p = a_3^{-v_2} (a_2^{v_2})^p (a_3^{v_3})^p (a_4^{v_4})^p\\
      &=& (a_3\1 a_2^p)^{v_2} (a_3^p)^{v_3} (a_4^p)^{v_4}\\
      &=& a_5^{v_3} a_6^{v_4} a_7^{v_2},
\end{eqnarray*}
da die Kommutatoren von $a_2, a_3$ und $a_4$ trivial sind oder die Ordnung $p$
haben.
\end{beweis}

Da statt der zweidimensionalen Unterräume von $\F_p^3 \cong M(G)$ ihre eindimensionalen
orthogonalen Komplemente betrachtet werden, ist nicht so sehr $\vi$, sondern der zu $\vi$
duale Operationshomomorphismus von Interesse:

\begin{beme}
Der zu $\vi$ duale Operationshomomorphismus ist $$\bar{\vi}: Aut(G) \to GL(3,p): \al(u_2, u_3,
u_4, v_2, v_3, v_4) \mt M =
\begin{pmatrix}
  v_2 & u_2 & v_3 \\
  0   & 1   & v_4 \\
  0   & 0   & v_2
\end{pmatrix}.$$
\end{beme}

\begin{beme}
Da $Aut(G)^*=\{\alpha^{\bar{\vi}} \mid \alpha \in Aut(G)\}$ nicht das Zentrum von $GL(3,p)$
enthält, entsprechen nach \ref{zentin} die Bahnen der Vektoren nicht den Bahnen der
eindimensionalen Unterräume. Da hier die Bahnen der eindimensionalen Unterräume von Interesse
sind, wird im Einklang mit \ref{einvek} statt $\bar{\vi}$ der Operationshomomorphismus $$\vi':
Aut(G) \to GL(3,p): \al(u_2, u_3, u_4, v_2, v_3, v_4) \mt M =
\begin{pmatrix}
  v_2 & u_2 & v_3 \\
  0   & v_2   & v_4 \\
  0   & 0   & v_2
\end{pmatrix}$$ benutzt und statt $Aut(G)^*$ wird
$A=\{\alpha^{\vi'} \mid \al \in Aut(G)\}$ verwendet.
\end{beme}

\subsection{Bahnen zulässiger Untergruppen für Nachfolger der Ordnung $p^5$}

\begin{lemma}
Unter der Operation von $A$ zerfällt $\f^3 \cong M(G)$ in folgende Bahnen nichttrivialer
Vektoren:
\begin{enumerate}
  \item[] $B_1=(1,0,0)^A$ mit $|B_1|=p^2(p-1)$.
  \item[] $B_2=(0,1,0)^A$ mit $|B_2|=p(p-1)$.
  \item[] $B_3=(0,0,1)^A$ mit $|B_3|=(p-1)$.
\end{enumerate}
Die zugehörigen Stabilisatoren sind
\begin{eqnarray*}
  \hat{S}_1 & = & \left\{ \begin{pmatrix}
    1         & 0        & 0      \\
    0         & 1        & v_4      \\
    0         & 0        & 1      \\
  \end{pmatrix} \mid  v_4 \in \F_p \right\}
  \text{ mit } |\hat{S}_1|=p.\\
  \hat{S}_2 & = & \left\{ \begin{pmatrix}
    1         & u_2      & v_3    \\
    0         & 1        & 0      \\
    0         & 0        & 1      \\
  \end{pmatrix} \mid  u_2, v_3 \in \F_p \right\}
  \text{ mit } |\hat{S}_2|=p^2.\\
  \hat{S}_3 & = & \left\{ \begin{pmatrix}
    1         & u_2      & v_3      \\
    0         & 1        & v_4      \\
    0         & 0        & 1      \\
  \end{pmatrix} \mid  u_2, v_3, v_4 \in \F_p \right\}
  \text{ mit } |\hat{S}_3|=p^3.\\
\end{eqnarray*}
Die folgenden beiden Stabilisatoren lassen die eindimensionalen Unterräume
invariant lassen, die von den Vertretern von $B_1$ bzw. $B_2$ aufgespannt
werden:
\begin{eqnarray*}
  \bar{S}_1 & = & \left\{ \begin{pmatrix}
    v_2       & 0        & 0      \\
    0         & v_2      & v_4    \\
    0         & 0        & v_2    \\
  \end{pmatrix} \mid v_2 v_4 \in \F_p \text{ und } v_2 \neq 0 \right\},\\
  \bar{S}_2 & = & \left\{ \begin{pmatrix}
    v_2       & u_2      & v_3    \\
    0         & v_2      & 0      \\
    0         & 0        & v_2    \\
  \end{pmatrix} \mid u_2, v_2, v_3 \in \F_p \text{ und } v_2 \neq 0\right\}.\\
\end{eqnarray*}
\end{lemma}

\begin{beweis}
Die Behauptung ist unmittelbar ersichtlich.
\end{beweis}

\begin{folg} \label{zu5}
Für die unmittelbaren Nachfolger von $G$ der Ordnung $p^5$ ist die folgende
Liste ein Vertretersystem zulässiger Untergruppen:
\begin{enumerate}
  \item $M_1=\erz{a_6, a_7}$ entsprechend $\erz{(0,1,0),(0,0,1)}$,
  \item $M_2=\erz{a_5, a_7}$ entsprechend $\erz{(1,0,0),(0,0,1)}$,
\end{enumerate}
\end{folg}

\begin{beweis}
Genau die Unterräumen, die von Elementen der Bahnen $B_1$ und $B_2$ aufgespannt
werden, haben orthogonale Komplemente, denen nach \ref{suppl} zulässige
Untergruppen von $M(G)$ entsprechen, da der Vertreter von $B_3$ mit den
Basisvektoren von $\bar{N}=\erz{(1,0,0),(0,1,0)}$ das Skalarprodukt Null
liefert und $\bar{N}$ dem Nukleus entspricht.
\end{beweis}

\subsection{Nachfolger der Ordnung $p^5$}

\begin{satz} \label{nt4}
In der folgenden Liste sind sämtliche unmittelbaren Nachfolger von $G$ der
Ordnung $p^5$ bis auf Isomorphie angegeben:
\begin{enumerate}
  \item $G_1=\erz{g_1,g_2,g_3,g_4,g_5 \mid
  [g_2,g_1]=g_3, g_2^p=g_3, g_1^p=g_4, [g_3,g_1]=g_5, [g_4,g_2]=g_5^{p-1}, g_3^p=g_5}$,
  \item $G_2=\erz{g_1,g_2,g_3,g_4,g_5 \mid
  [g_2,g_1]=g_3, g_2^p=g_3, g_1^p=g_4, g_4^p=g_5}$.
\end{enumerate}
\end{satz}

\begin{beweis}
Dieser Satz ergibt sich nach \ref{bahnenzu}, indem man $P(G)$ nach jedem Repräsentanten des
Vertretersystems zulässiger Untergruppen faktorisiert, das in \ref{zu5} aufgelistet ist. Die
Präsentation der Faktorgruppe erhält man über das Verfahren aus \ref{faktor}.
\end{beweis}

\section{Nachfolger von $C_p^3$}

\begin{ver}
In diesem Abschnitt sei $$G=\erz{a_1, a_2, a_3 \mid a_1^p, a_2^p, a_3^p, [a_2, a_1], [a_3,
a_1], [a_3, a_2]}.$$ Damit ist $G$ isomorph zu $C_p^3$.
\end{ver}

\subsection{Die Automorphismengruppe}

\begin{beme}
Die Gruppe $G$ ist isomorph zu additiven Gruppe des Vektorraumes $\F_p^3$. Die
Automorphismengruppe $Aut(G)$ kann daher mit $GL(3,p)$ identifiziert werden. Im Weiteren sei
$A:=GL(3,p)$.
\end{beme}

\subsection{$p$"=Cover, Multiplikator und Nukleus}

\begin{lemma}
Es gilt:
\begin{enumerate}
  \item[] $P(G)=\erz{a_1 \dd a_9 \mid [a_2,a_1]=a_4, [a_3,a_1]=a_5,
  [a_3,a_2]=a_6, a_1^p=a_7, a_2^p=a_8, a_3^p=a_9}$.
  \item[] $M(G)=\erz{a_4, a_5, a_6, a_7, a_8, a_9}$.
  \item[] $N(G)=\erz{a_4, a_5, a_6, a_7, a_8, a_9}$.
\end{enumerate}
Die Gruppe $P(G)$ hat die $p$"=Klasse 2, $G$ ist fortsetzbar und jede
Untergruppe von $M(G)$ ist zulässig.
\end{lemma}

\begin{beweis}
Die Präsentation von $P(G)$ ergibt sich aus \ref{cov} und \ref{newman} unter Verwendung des
Reduktionsverfahrens von Knuth"=Bendix, aus der sich $M(G)$ ablesen lässt. Die Gewichtung von
$P(G)$
\begin{enumerate}
  \item[] $1=\omega(a_1)=\omega(a_2)=\omega(a_3)$.
  \item[] $2=\omega(a_4)=\omega(a_5)=\omega(a_6)=\omega(a_7)=\omega(a_8)=\omega(a_9)$.
\end{enumerate}
und damit $N(G)$ lassen sich nach \ref{gew} aus der Präsentation von $P(G)$ ablesen.
\end{beweis}

\subsection{Operation der Erweiterungsautomorphismen}

\begin{lemma} \label{autcp3}
Mit der Identifikation von $Aut(G)$ mit $GL(3,p)$ ist die Operation von $A=GL(3,p)$ auf $\f^6
\cong M(G)$ durch den Operationshomomorphismus $\vi$ gegeben: $$\vi: GL(3,p) \to GL(6,p):
m=\begin{pmatrix}
  m_{1} & m_{2} & m_{3} \\
  m_{4} & m_{5} & m_{6} \\
  m_{7} & m_{8} & m_{9}
\end{pmatrix} \mt$$ $$ M = \begin{pmatrix}
  m_1 m_5 - m_4 m_2 & m_1 m_6 - m_4 m_3 & m_2 m_6 - m_5 m_3 & 0 & 0 & 0 \\
  m_1 m_8 - m_7 m_2 & m_1 m_9 - m_7 m_3 & m_2 m_9 - m_8 m_3 & 0 & 0 & 0 \\
  m_4 m_8 - m_7 m_5 & m_4 m_9 - m_7 m_6 & m_5 m_9 - m_8 m_6 & 0 & 0 & 0 \\
  0 & 0 & 0 & m_1 & m_2 & m_3 \\
  0 & 0 & 0 & m_4 & m_5 & m_6 \\
  0 & 0 & 0 & m_7 & m_8 & m_9
\end{pmatrix}$$ $$= \begin{pmatrix}
  |m_{33}| & |m_{32}| & |m_{31}| & 0 & 0 & 0 \\
  |m_{23}| & |m_{22}| & |m_{21}| & 0 & 0 & 0 \\
  |m_{13}| & |m_{12}| & |m_{11}| & 0 & 0 & 0 \\
  0 & 0 & 0 & m_1 & m_2 & m_3 \\
  0 & 0 & 0 & m_4 & m_5 & m_6 \\
  0 & 0 & 0 & m_7 & m_8 & m_9
\end{pmatrix}$$ Dabei ist $m_{ij}$ die Matrix, die aus der Matrix $m$ durch
Streichung der $i$"=ten Zeile und $j$"=ten Spalte entsteht, und es ist $|m_{ij}|=det(m_{ij})$.
\end{lemma}

\begin{beweis}
Für die Erzeuger von $M(G)$ ergeben sich unter der Operation von $A$ auf $M(G)$
über Erweiterungsautomorphismen die folgenden Bilder:
\begin{eqnarray*}
  a_7^m &=& (a_1^p)^m = (a_1^m)^p \\
  &=& (a_1^{m_1} a_2^{m_2} a_3^{m_3})^p\\
  &=& (a_1^{m_1})^p (a_2^{m_2})^p (a_3^{m_3})^p\\
  &=& (a_1^p)^{m_1} (a_2^p)^{m_2} (a_3^p)^{m_3}\\
  &=& a_7^{m_1} a_8^{m_2} a_9^{m_3},
\end{eqnarray*}
nach \ref{komm1} unter der Bedingung, dass die Kommutatoren $[a_2,a_1],
[a_3,a_1]$ und $[a_3,a_2]$ zentral sind und die Ordnung $p$ haben. Ebenso
erhält man $a_8^m=a_7^{m_4} a_8^{m_5} a_9^{m_6}$ und $a_9^m=a_7^{m_7} a_8^{m_8}
a_9^{m_9}$. Weiterhin ist
\begin{eqnarray*}
  a_4^m&=& [a_2^m,a_1^m]= [a_1^{m_4} a_2^{m_5} a_3^{m_6}, a_1^{m_1} a_2^{m_2} a_3^{m_3}]\\
  &=& [a_2^{m_5},a_1^{m_1}] [a_3^{m_6},a_1^{m_1},] [a_1^{m_4},a_2^{m_2}] [a_3^{m_6},a_2^{m_2}] [a_1^{m_4},a_3^{m_3}] [a_2^{m_5},a_3^{m_3}] \\
  &=& [a_2,a_1]^{m_5 m_1} [a_3,a_1]^{m_6 m_1} [a_1,a_2]^{m_4 m_2} [a_3,a_2]^{m_6 m_2} [a_1,a_3]^{m_4 m_3} [a_2,a_3]^{m_5 m_3} \\
  &=& [a_2,a_1]^{m_1 m_5} [a_2,a_1]^{-m_4 m_2} [a_3,a_1]^{m_1 m_6} [a_3,a_1]^{-m_4 {m_3}} [a_3,a_2]^{m_2 m_6} [a_3,a_2]^{-m_5 m_3}\\
  &=& [a_2,a_1]^{m_1 m_5 - m_4 m_2} [a_3,a_1]^{m_1 m_6 - m_4 {m_3}} [a_3,a_2]^{m_2 m_6 - m_5m_3}\\
  &=& a_4^{m_1 m_5 - m_4 m_2} a_5^{m_1 m_6 - m_4 {m_3}} a_6^{m_2 m_6 - m_5m_3},
\end{eqnarray*}
da die Kommutatoren $[a_2,a_1], [a_3,a_1]$ und $[a_3,a_2]$ zentral sind. Ebenso
erhält man
\begin{eqnarray*}
  a_5^m &=& a_4^{m_1 m_8 - m_7 m_2} a_5^{m_1 m_9 - m_7 m_3} a_6^{m_2 m_9 - m_8 m_3}\\
  a_6^m &=& a_4^{m_4 m_8 - m_7 m_5} a_5^{m_4 m_9 - m_7 m_6} a_6^{ m_5 m_9- m_8 m_6}
\end{eqnarray*}
für die übrigen Erzeuger von $M(G)$.
\end{beweis}

\subsection{Bahnen zulässiger Untergruppen für Nachfolger der Ordnung $p^4$}

\begin{lemma} \label{bahnen_cp3_p4}
Der Vektorraum $\f^6 \cong M(G)$ zerfällt unter der Operation von $A=GL(3,p)$ in folgende
Bahnen nichttrivialer Vektoren:
\begin{enumerate}
  \item[] $B_1 = (1,0,0,0,0,0)^A$ mit $|B_1|=p^3-1$.
  \item[] $B_2 = (0,0,0,1,0,0)^A$ mit $|B_2|=p^3-1$.
  \item[] $B_3 = (1,0,0,1,0,0)^A$ mit $|B_3|=(p-1)(p^2-1)(p^2+p+1)$.
  \item[] $B_4 = (1,0,0,0,0,1)^A$mit  $|B_4|=p^2(p^2+p+1)(p-1)^2$.
\end{enumerate}
Zu den Bahnenvertretern gehören die folgenden Stabilisatoren:
\begin{eqnarray*}
  \hat{S}_1 &=& \left\{ \begin{pmatrix}
  1                 & 0       & 0       & 0 & 0 & 0 \\
  m_1 m_8 - m_7 m_2 & m_1 m_9 & m_2 m_9 & 0 & 0 & 0 \\
  m_4 m_8 - m_7 m_5 & m_4 m_9 & m_5 m_9 & 0 & 0 & 0 \\
  0 & 0 & 0 & m_1 & m_2 & 0 \\
  0 & 0 & 0 & m_4 & m_5 & 0 \\
  0 & 0 & 0 & m_7 & m_8 & m_9
\end{pmatrix} \mid \triangle \right\}  \\
  & & \triangle = m_1, m_2, m_4, m_5, m_7, m_8, m_9 \in \f \text{ und } m_9 \neq 0 \text{ sowie }m_1 m_5 - m_4 m_2 = 1\\
  \hat{S}_2 &=& \left\{ \begin{pmatrix}
  m_5               & m_6               & 0                 & 0 & 0 & 0 \\
  m_8               & m_9               & 0                 & 0 & 0 & 0 \\
  m_4 m_8 - m_7 m_5 & m_4 m_9 - m_7 m_6 & m_5 m_9 - m_8 m_6 & 0 & 0 & 0 \\
  0 & 0 & 0 & 1   & 0   & 0   \\
  0 & 0 & 0 & m_4 & m_5 & m_6 \\
  0 & 0 & 0 & m_7 & m_8 & m_9
\end{pmatrix} \mid \triangle \right\}  \\
  & & \triangle = m_4, m_5, m_6, m_7, m_8, m_9 \in \f \text{ und } m_5 m_9 - m_8 m_6 \neq 0 \\
  \hat{S}_3 &=& \left\{ \begin{pmatrix}
  1             & 0        & 0   & 0 & 0 & 0 \\
  m_8           & m_9      & 0   & 0 & 0 & 0 \\
  m_4 m_8 - m_7 & m_4 m_9  & m_9 & 0 & 0 & 0 \\
  0 & 0 & 0 & 1   & 0   & 0 \\
  0 & 0 & 0 & m_4 & 1   & 0 \\
  0 & 0 & 0 & m_7 & m_8 & m_9
  \end{pmatrix} \mid \triangle \right\}  \\
  & & \triangle = m_4, m_7, m_8, m_9 \in \f \text{ und } m_9 \neq 0 \\
  \hat{S}_4 &=& \left\{ \begin{pmatrix}
  1                 & 0                 & 0                 & 0 & 0 & 0 \\
  0                 & m_1               & m_2               & 0 & 0 & 0 \\
  0                 & m_4               & m_5               & 0 & 0 & 0 \\
  0 & 0 & 0 & m_1 & m_2 & 0 \\
  0 & 0 & 0 & m_4 & m_5 & 0 \\
  0 & 0 & 0 & 0   & 0   & 1
\end{pmatrix} \mid \triangle \right\} \\
  & & \triangle = m_1, m_2, m_4, m_5 \in \f \text{ und } m_1 m_5 - m_4 m_2 =1
\end{eqnarray*}
Die folgenden Gruppen sind die Stabilisatoren der eindimensionalen Unterräume, die von den
Bahnenvertretern aufgespannt werden.
\begin{eqnarray*}
  \bar{S}_1 &=& \left\{ \begin{pmatrix}
  m_1 m_5 - m_4 m_2 & 0       & 0       & 0 & 0 & 0 \\
  m_1 m_8 - m_7 m_2 & m_1 m_9 & m_2 m_9 & 0 & 0 & 0 \\
  m_4 m_8 - m_7 m_5 & m_4 m_9 & m_5 m_9 & 0 & 0 & 0 \\
  0 & 0 & 0 & m_1 & m_2 & 0 \\
  0 & 0 & 0 & m_4 & m_5 & 0 \\
  0 & 0 & 0 & m_7 & m_8 & m_9
\end{pmatrix} \mid \triangle \right\}  \\
  & & \triangle = m_1, m_2, m_4, m_5, m_7, m_8, m_9 \in \f \text{ und } m_9 \neq 0 \neq m_1 m_5 - m_4 m_2 \\
  \bar{S}_2 &=& \left\{ \begin{pmatrix}
  m_1 m_5           & m_1 m_6           & 0                 & 0 & 0 & 0 \\
  m_1 m_8           & m_1 m_9           & 0                 & 0 & 0 & 0 \\
  m_4 m_8 - m_7 m_5 & m_4 m_9 - m_7 m_6 & m_5 m_9 - m_8 m_6 & 0 & 0 & 0 \\
  0 & 0 & 0 & m_1   & 0   & 0   \\
  0 & 0 & 0 & m_4 & m_5 & m_6 \\
  0 & 0 & 0 & m_7 & m_8 & m_9
\end{pmatrix} \mid \triangle \right\}  \\
  & & \triangle = m_1, m_4, m_5, m_6, m_7, m_8, m_9 \in \f \text{ und } m_9 \neq 0 \neq m_5 m_9 - m_8 m_6 \\
  \bar{S}_3 &=& \left\{ \begin{pmatrix}
  m_1^2             & 0        & 0        & 0 & 0 & 0 \\
  m_1 m_8           & m_1 m_9  & 0       & 0 & 0 & 0 \\
  m_4 m_8 - m_7 m_1 & m_4 m_9  & m_1 m_9 & 0 & 0 & 0 \\
  0 & 0 & 0 & m_1   & 0   & 0 \\
  0 & 0 & 0 & m_4 & m_1   & 0 \\
  0 & 0 & 0 & m_7 & m_8 & m_9
  \end{pmatrix} \mid \triangle \right\} \\ \text{}
  & & \triangle = m_1, m_4, m_7, m_8, m_9 \in \f \text{ und } m_9 \neq 0 \neq m_1 \\
  \bar{S}_4 &=& \left\{ \begin{pmatrix}
  m_1 m_5 - m_4 m_2 & 0                 & 0                 & 0 & 0 & 0 \\
  0                 & m_1               & m_2               & 0 & 0 & 0 \\
  0                 & m_4               & m_5               & 0 & 0 & 0 \\
  0 & 0 & 0 & m_1 & m_2 & 0 \\
  0 & 0 & 0 & m_4 & m_5 & 0 \\
  0 & 0 & 0 & 0   & 0   & m_1 m_5 - m_4 m_2
\end{pmatrix} \mid \triangle \right\} \\
  & & \triangle = m_1, m_2, m_4, m_5 \in \f \text{ und }  m_1 m_5 - m_4 m_2 \neq 0
\end{eqnarray*}
Die folgenden Gruppen sind die Urbilder von $\bar{S}_1 \dd \bar{S}_4$ unter $\vi$ in
$A=GL(3,p)$:
\begin{eqnarray*}
  S_1 &=& \left\{ \begin{pmatrix}
  m_1 & m_2 & 0 \\
  m_4 & m_5 & 0 \\
  m_7 & m_8 & m_9
\end{pmatrix} \mid m_1, m_2, m_4, m_5, m_7, m_8, m_9 \in \f, m_9 \neq 0 \neq m_1 m_5 - m_4 m_2 \right\}  \\
  S_2 &=& \left\{ \begin{pmatrix}
  m_1   & 0   & 0   \\
  m_4 & m_5 & m_6 \\
  m_7 & m_8 & m_9
\end{pmatrix} m_1, m_4, m_5, m_6, m_7, m_8, m_9 \in \f , m_1 \neq 0 \neq m_5 m_9 - m_8 m_6 \right\}  \\
  S_3 &=& \left\{ \begin{pmatrix}
  m_1   & 0   & 0 \\
  m_4   & m_1   & 0 \\
  m_7   & m_8 & m_9
  \end{pmatrix} \mid m_1, m_4, m_7, m_8, m_9 \in \f \text{ und } m_9 \neq 0 \neq m_1 \right\}  \\
  S_4 &=& \left\{ \begin{pmatrix}
  m_1 & m_2 & 0 \\
  m_4 & m_5 & 0 \\
  0   & 0   & m_1 m_5 - m_4 m_2
\end{pmatrix}  m_1, m_2, m_4, m_5 \in \f, m_1 m_5 - m_4 m_2 \neq 0 \right\}
\end{eqnarray*}
\end{lemma}

\begin{beweis}
Die Stabilisatoren $\hat{S}_1 \dd \hat{S}_4$ ergeben sich aus den Lösungsmengen der folgenden
Gleichungssysteme:
\begin{enumerate}
  \item $\{m_1 m_5 - m_4 m_2 = 1, m_1 m_6 - m_4 m_3 = 0, m_2 m_6 - m_5 m_3 = 0\}$,
  \item $\{m_1 = 1, m_2 = 0, m_3 = 0\}$,
  \item $\{m_1 m_5 - m_4 m_2 = 1, m_1 m_6 - m_4 m_3 = 0, m_2 m_6 - m_5 m_3 = 0, m_1 = 1, m_2 = 0, m_3 = 0\}$,
  \item $\{m_1 m_5 - m_4 m_2 = 1, m_1 m_6 - m_4 m_3 = 0, m_2 m_6 - m_5 m_3 = 0, m_7 = 0, m_8 = 0, m_9 =
  1\}$.
\end{enumerate}

Die Längen der Bahnen $B_1$ und $B_2$ sind unmittelbar ersichtlich, da sie den nichttrivialen
Vektoren der beiden irreduziblen Teilmoduln entsprechen, in die $\f^6$ unter der Operation von
$A$ zerfällt. Die Längen der Bahnen $B_3$ und $B_4$ lassen sich nach dem
Bahn"=Stabilisator"=Satz \ref{stab} ermitteln, indem man die Mächtigkeit von $GL(3,p)$ durch
die Ordnung der Stabilisatoren $\hat{S}_3$ und $\hat{S}_4$ teilt. Dabei ist nach \ref{glslz}
bzw. nach den Lösungsmengen des dritten und vierten Gleichungssystems:
\begin{eqnarray*}
  |\hat{S}_3| &=& p^3(p-1) \\
  |\hat{S}_4| &=& |SL^+(2,p)| = \frac{|GL(2,p)|}{p-1} = \frac{(p^2-1)(p-1)p}{p-1}=(p^2-1)p\\
  |GL(3,p)| &=&  (p^3-1)(p^3-p)(p^3-p^2)=(p^2+p+1)(p+1)(p-1)^3 p^3
\end{eqnarray*}
Addiert man die Länge der vier Bahnen auf, ergibt sich $p^6-1$. Damit ist gezeigt, dass mit
$B_1 \dd B_4$ alle Bahnen erfasst sind. Die Stabilisatoren $\bar{S}_1 \dd \bar{S}_4$ bzw. $S_1
\dd S_4$ ergeben sich als Lösungsmengen von Gleichungssystemen, die aus den oben angegebenen
entstehen, indem man die rechten Seiten mit einer Variable des Wertebereichs $\f^*$
multipliziert.

Es sei $w$ ein Erzeuger von $\f^*$ und $E$ die Einheitsmatrix von $GL(6,p)$. Nun wird gezeigt,
dass die eindimensionalen Unterräume vier Bahnen bilden, deren Vertreter von den
Repräsentanten von $B_1 \dd B_4$ aufgespannt werden (da $A^\vi$ das Zentrum von $GL(6,p)$
nicht enthält, ist das nicht offensichtlich). Nach \ref{bahnenzu} kann zur Ermittlung der
Bahnen eindimensionaler Unterräume statt der Operation von $A$ über $\vi$ die Operation der
Gruppe $G:=\erz{w \cdot E, A^\vi}$ auf $\f^6$ betrachtet werden. Zuerst wird der Index von
$A^\vi$ in $G$ ermittelt. Es ist $\erz{w \cdot E}=\{w^a \cdot E \mid 0 \leq a < p \}$. Da
$m_1=m_5=m_9=w^a$, $m_2=m_3=m_4=m_6=m_7=m_8=0$ und $m_1 m_5 - m_4 m_2 = w^a$ nur für $a=0$
erfüllt ist, ist der Schnitt von $A^\vi$ mit $\erz{w \cdot E}$ trivial. Also ist
$[G:A^\vi]=p-1$ und $|G|=(p-1) \cdot |A|$. Außerdem ist $G=\{w^a M \mid 0 \leq a < p \text{
und } M \in A^\vi\}$. Man betrachte nun die Ordnung der (oben nicht aufgeführten)
Stabilisatoren in $G$ der Repräsentanten von $B_1 \dd B_4$. Sie ergeben sich aus den
Lösungsmengen der folgenden Gleichungssysteme:
\begin{enumerate}
  \item[1'.] $\{w^a(m_1 m_5 - m_4 m_2) = 1, w^a(m_1 m_6 - m_4 m_3) = 0, w^a(m_2 m_6 - m_5 m_3) = 0\}$,
  \item[2'.] $\{w^a m_1 = 1, w^a m_2 = 0, w^a m_3 = 0\}$,
  \item[3'.] $\{w^a(m_1 m_5 - m_4 m_2) = 1, w^a(m_1 m_6 - m_4 m_3) = 0, w^a(m_2 m_6 - m_5 m_3) = 0, w^a m_1 = 1, w^a m_2 = 0, w^a m_3 = 0\}$,
  \item[4'.] $\{w^a(m_1 m_5 - m_4 m_2) = 1, w^a(m_1 m_6 - m_4 m_3) = 0, w^a(m_2 m_6 - m_5 m_3) = 0, w^a m_7 = 0, w^a m_8 = 0, w^a m_9 =
  1\}$.
\end{enumerate}
Teilt man alle vier Gleichungssysteme durch $w^a$, so erkennt man, dass ihre Lösungsmengen den
Stabilisatoren $\bar{S}_1 \dd \bar{S}_4$ entsprechen. Also sind die Stabilisatoren der
Vertreter von $B_1 \dd B_4$ in $G$ jeweils um den Faktor $p-1$ größer als $\hat{S}_1 \dd
\hat{S}_4$. Da die Ordnung von $G$ um demselben Faktor größer ist als die von $A^\vi$, sind
die Bahnen der Vertreter von $B_1 \dd B_4$ unter $G$ genauso lang wie $B_1 \dd B_4$ (man kann
sogar zeigen, dass $G$ dieselben Bahnen liefert wie $A$; aber das ist irrelevant). Daher sind
die Vertreter von $B_1 \dd B_4$ zugleich ein Vertretersystem der Bahnen unter $G$ und die
Unterräume, die von ihnen aufgespannt werden, ein Vertretersystem der Bahnen eindimensionaler
Unterräume unter der Operation von $A$.
\end{beweis}

\begin{folg} \label{zu6}
Für die unmittelbaren Nachfolger von $G$ der Ordnung $p^4$ ist die folgende
Liste ein Vertretersystem zulässiger Untergruppen, wobei $\bot$ die Abbildung
ist, die jedem Unterraum von $\f^6$ sein orthogonales Komplement zuordnet.
\begin{enumerate}
  \item[] $M_1=\erz{a_5,a_6,a_7,a_8,a_9}$ entsprechend $\erz{(1,0,0,0,0,0)}^\bot$,
  \item[] $M_2=\erz{a_4,a_5,a_6,a_8,a_9}$ entsprechend $\erz{(0,0,0,1,0,0)}^\bot$,
  \item[] $M_3=\erz{a_5,a_6,a_8,a_9, a_4 a_7\1}$ entsprechend $\erz{(1,0,0,1,0,0)}^\bot$,
  \item[] $M_4=\erz{a_5,a_6,a_7,a_8, a_4 a_9\1}$ entsprechend $\erz{(1,0,0,0,0,1)}^\bot$.
\end{enumerate}
\end{folg}

\subsection{Nachfolger der Ordnung $p^4$}

\begin{satz} \label{ncp3p4}
In der folgenden Liste sind sämtliche unmittelbaren Nachfolger von $G$ der
Ordnung $p^4$ bis auf Isomorphie angegeben:
\begin{enumerate}
  \item[] $G_1=\erz{g_1, g_2, g_3, g_4 \mid [g_2,g_1]=g_4}$ \hfill \gap"=Typ $(p^4,12)$
  \item[] $G_2=\erz{g_1, g_2, g_3, g_4 \mid g_1^p=g_4}$ \hfill \gap"=Typ $(p^4,11)$
  \item[] $G_3=\erz{g_1, g_2, g_3, g_4 \mid [g_2,g_1]=g_4, g_1^p=g_4}$ \hfill \gap"=Typ $(p^4,13)$
  \item[] $G_4=\erz{g_1, g_2, g_3, g_4 \mid [g_2,g_1]=g_4, g_3^p=g_4}$ \hfill \gap"=Typ $(p^4,14)$
\end{enumerate}
\end{satz}

\begin{beweis}
Dieser Satz ergibt sich nach \ref{bahnenzu}, indem man $P(G)$ nach jedem Repräsentanten des
Vertretersystems zulässiger Untergruppen faktorisiert, das in \ref{zu6} aufgelistet ist. Die
Präsentation der Faktorgruppe erhält man über das Verfahren aus \ref{faktor}.
\end{beweis}

\subsection{Bahnen zulässiger Untergruppen für Nachfolger der Ordnung $p^5$}

\begin{satz} \label{zu7}
Die Menge zulässiger Untergruppen von $M(G)$ der Ordnung $p^4$ zerfällt unter der Operation
von $Aut(G)$ über Erweiterungsautomorphismen in $14+p$ Bahnen. Die folgende Liste ist ein
Vertretersystem dieser Bahnen.
\begin{enumerate}
  \item[] $M_1=\erz{a_6,a_7,a_8,a_9}$,
  \item[] $M_2=\erz{a_5,a_6,a_7,a_8}$,
  \item[] $M_3=\erz{a_5,a_6 a_9\1,a_7,a_8}$,
  \item[] $M_4=\erz{a_4,a_5,a_7,a_8}$,
  \item[] $M_5=\erz{a_4,a_5 a_9\1,a_7,a_8}$,
  \item[] $M_6=\erz{a_4 a_9\1,a_5,a_7,a_8}$,
  \item[] $M_7=\erz{a_4, a_5, a_6 a_9\1, a_7}$,
  \item[] $M_8=\erz{a_4, a_5 a_8\1, a_6 a_9\1, a_7}$,
  \item[] $M_9=\erz{a_4 a_8\1 ,a_5, a_6 a_9\1 ,a_7}$,
  \item[] $M_{10}=\erz{a_4 ,a_5, a_6 ,a_7}$,
  \item[] $M_{11}=\erz{a_4 a_8^{-1}, a_5 a_9^{-1}, a_6 ,a_7}$,
  \item[] $M_{12}=\erz{a_4 a_9^{-1}, a_5, a_6 ,a_7}$,
  \item[] $M_{13}=\erz{a_4 a_8^{-1} a_9^{-1}, a_5 a_9^{-1}, a_6 ,a_7}$,
  \item[] $M_{14}=\erz{a_4, a_5 a_9^{-1}, a_6 ,a_7}$,
  \item[] $M_{15}^k=\erz{a_4 a_8^{-1}, a_5 a_9^{-w^k}, a_6 ,a_7}$ mit $k \in \{1 \dd \frac{p-1}{2}\}$,
  \item[] $M_{16}=\erz{a_4 a_9^{-w}, a_5 a_8^{-1}, a_6 ,a_7}$,
  \item[] $M_{17}^k=\erz{a_4 a_8^{-1} a_9^{-w^k}, a_5 a_8^{-w^{k-1}} a_9^{-1}, a_6 ,a_7}$ mit $k \in \{1 \dd \frac{p-1}{2}\}$.
\end{enumerate}
\end{satz}

\begin{beweis}
Der Multiplikator von $G$ hat die Ordnung $p^6$ und $P(G)$ die Ordnung $p^9$. Da jede
Untergruppe von $M(G)$ zulässig ist, wird für die unmittelbaren Nachfolger der Ordnung $p^5$
ein Vertretersystem der Bahnen gesucht, das die Untergruppen von $M(G)$ der Ordnung $p^4$
unter der Operation von $Aut(G)$ über Erweiterungsautomorphismen erfasst. Diese
Problemstellung wird in die Sprache der linearen Algebra übersetzt: Es sei $\{f_4, f_5, f_6,
f_7, f_8, f_9\}$ die Standardbasis von $V=\f^6$. Entspricht ein Unterraum $W$ von $V$ einer
Untergruppe von $M(G)$ der Ordnung $p^4$, so wird $W$ innerhalb dieses Beweises zulässiger
Unterraum genannt. Ist $W$ ein zulässiger Unterraum, so hat $W$ die Dimension 4.

Es seien $N=\erz{f_4, f_5, f_6}$ und $K=\erz{f_7, f_8, f_9}$. Dann ist $V = N \oplus K$. Es
sei $B=A^\vi$. Die Unterräume $N$ und $K$ sind unter der Operation von $B$ irreduzible
Untermoduln von $V$. Ist $W$ ein zulässiger Unterraum, so ist $dim(W \cap K) \in \{1,2,3\}$.
Diese drei Fälle werden nun getrennt untersucht.

\textbf{Erster Fall:} Es sei $dim(W \cap K)=3$. Dann ist $(W \cap K)^B=W \cap K=K$, da $W \cap
K=K$ und $K$ unter der Operation von $B$ invariant ist. Es seien $\bar{W}=W/K$ und
$\bar{V}=V/K$. Da $\bar{W}$ ein eindimensionaler Unterraum von $\bar{V}$ ist und $B$ transitiv
auf $\bar{V} \backslash\{0\}$ operiert, liegen nach \ref{tensor} alle zulässigen Unterräume
$W$ mit $dim(W \cap K)=3$ in einer Bahn. Als Vertreter dieser Bahn kann $W=\{f_5, f_6, f_7,
f_8\}$ gewählt werden. Dieser Unterraum entspricht der Untergruppe $M_1$.

\textbf{Zweiter Fall:} Es sei $dim(W \cap K)=2$. Im Weiteren wird die folgende Notation
benutzt:
\begin{itemize}
  \item[] $U:=\{f_7, f_8\}$,
  \item[] $\bar{V}:=V/U$, $\bar{K}:=K/U$, $\bar{N}:=N/U$ und $\bar{W}:=W/U$,
  \item[] $\bar{f}_i:=f_i + U$ für $i \in \{4,5,6,9\}$,
  \item[] $\hat{V}:=\bar{V}/\bar{K}$ und $\hat{W}:=\bar{W}/\bar{K}$.
  \item[] $\hat{f}_i:=\bar{f}_i + \bar{K}$ für $i \in \{4,5,6\}$.
\end{itemize}
Die Idee des Beweises lautet folgendermaßen: Der Unterraum $U$ kann als Vertreter der Bahn
gewählt werden, den der Schnitt von $W$ in $K$ hat. Der Stabilisator von $U$ in $B$ wird
berechnet und die Bahn von $W/U$ in $V/U$ ermittelt. Die Bahn von $W$ ergibt sich dann aus der
Bahn von $W/U$ und $U$. Aus technischen Gründen wird der Faktorraum $\hat{V}$ benutzt, um die
Bahnen von $W/U$ zu berechnen. Diese Idee wird nun ausgeführt.

Da $B$ transitiv auf $K \backslash\{0\}$ operiert, hat $W \cap K$ als Bahn unter $A$ die Menge
aller zweidimensionalen Unterräume von $K$. Man kann daher den Unterraum $U=\{f_7, f_8\}$ als
Vertreter von $(W \cap K)^B$ wählen. Nun betrachte man die Operation von $B$ auf $\bar{V}=V/U$
und die Untergruppe $S$ von $B$, die $U$ invariant lässt. Es ist $$S=Stab_B(U)=\left\{
\begin{pmatrix}
  m_1 m_5 - m_4 m_2 & 0                 & 0                 & 0   & 0   & 0 \\
  m_1 m_8 - m_7 m_2 & m_1 m_9           & m_2 m_9           & 0   & 0   & 0 \\
  m_4 m_8 - m_7 m_5 & m_4 m_9           & m_5 m_9           & 0   & 0   & 0 \\
  0                 & 0                 & 0                 & m_1 & m_2 & 0 \\
  0                 & 0                 & 0                 & m_4 & m_5 & 0 \\
  0                 & 0                 & 0                 & m_7 & m_8 & m_9
\end{pmatrix} \mid \triangle\right\}.$$
Dabei ist $\triangle= m_1, m_2, m_4, m_5, m_7, m_8, m_9 \in \f$ und $m_9 \neq 0$ sowie $m_1
m_5 - m_4 m_2 \neq 0$. Es sei $\bar{f}_4=f_4+U, \bar{f}_5=f_5+U, \bar{f}_6=f_6+U$ und
$\bar{f}_9=f_9+U$ sowie $\bar{K}=K/U$, $\bar{W}=W/U$ und $\bar{N}=N/U$. Dann operiert $S$
bezüglich der Basis $\{\bar{f}_4, \bar{f}_5, \bar{f}_6, \bar{f}_9\}$ als $$\bar{S}=\left\{
\begin{pmatrix}
  m_1 m_5 - m_4 m_2 & 0                 & 0          & 0 \\
  m_1 m_8 - m_7 m_2 & m_1 m_9           & m_2 m_9    & 0 \\
  m_4 m_8 - m_7 m_5 & m_4 m_9           & m_5 m_9    & 0 \\
  0                 & 0                 & 0          & m_9
\end{pmatrix} \mid \triangle\right\}$$ auf dem Faktorraum $\bar{V}$. Dabei ist $\triangle=
m_1, m_2, m_4, m_5, m_7, m_8, m_9 \in \f$ und $m_9 \neq 0$ sowie $m_1 m_5 - m_4 m_2 \neq 0$.
Da $m_1 m_8 - m_7 m_2$ und $m_4 m_8 - m_7 m_5$ unabhängig von $m_1, m_2, m_4, m_5$ und $m_9$
jeden Wert in $\f$ annehmen, ist $$\bar{S}=\left\{
\begin{pmatrix}
  m_1 m_5 - m_4 m_2 & 0                 & 0          & 0 \\
  x                 & m_1 m_9           & m_2 m_9    & 0 \\
  y                 & m_4 m_9           & m_5 m_9    & 0 \\
  0                 & 0                 & 0          & m_9
\end{pmatrix} \mid \triangle\right\}.$$ Dabei ist $\triangle= m_1, m_2, m_4, m_5,
m_9, x, y \in \f$ und $m_9 \neq 0$ sowie $m_1 m_5 - m_4 m_2 \neq 0$.

Gesucht werden nun die Bahnen der zweidimensionalen Unterräume $\bar{W}$ von $\bar{V}$ mit
$\bar{W} \cap \bar{K} = \{0\}$, da $dim(W \cap K)$ nach Voraussetzung $2$ ist. Zunächst wird
$\hat{V}=\bar{V}/\bar{K}$ betrachtet und die Bahnen der zweidimensionalen Unterräume von
$\hat{V}$ werden ermittelt, um von diesem Ergebnis aus die Bahnen der zweidimensionalen
Unterräume $\bar{W}$ von $\bar{V}$ zu bestimmen, die die Bedingung $dim(\bar{W} \cap
\bar{K})=0$ erfüllen.

Es sei $\hat{f}_4=\bar{f}_4 + \bar{K}, \hat{f}_5=\bar{f}_5 + \bar{K}$ und $\hat{f}_6 =
\bar{f}_6 + \bar{K}$ sowie $\hat{W}=\bar{W}/\bar{K}$. Dann operiert $\bar{S}$ bezüglich der
Basis $\{\hat{f}_4, \hat{f}_5, \hat{f}_6\}$ als $$\hat{S}=\left\{
\begin{pmatrix}
  m_1 m_5 - m_4 m_2 & 0                 & 0        \\
  x                 & m_1 m_9           & m_2 m_9  \\
  y                 & m_4 m_9           & m_5 m_9
\end{pmatrix} \mid \triangle\right\}$$ auf $\hat{V}$. Dabei ist $\triangle= m_1, m_2, m_4, m_5,
m_9, x, y \in \f$ und $m_9 \neq 0$ sowie $m_1 m_5 - m_4 m_2 \neq 0$. Statt die Bahnen der
zweidimensionalen Unterräume $\hat{W}$ werden die Bahnen der entsprechenden eindimensionalen
orthogonalen Komplemente $\hat{W}^\bot$ betrachtet. Dazu wird die zu $\hat{S}$ transponierte
Matrixgruppe $$\hat{S}^\bot=\left\{
\begin{pmatrix}
  m_1 m_5 - m_4 m_2 & x                 & y        \\
  0                 & m_1 m_9           & m_2 m_9  \\
  0                 & m_4 m_9           & m_5 m_9
\end{pmatrix} \mid \triangle\right\}$$ verwendet (dabei ist $\triangle= m_1, m_2, m_4, m_5,
m_9, x, y \in \f$ und $m_9 \neq 0$ sowie $m_1 m_5 - m_4 m_2 \neq 0$). Da $\hat{S}^\bot$ das
Zentrum von $GL(3,p)$ nicht enthält, lassen sich die Bahnen der eindimensionalen Unterräume
nach \ref{zentin} nicht unbedingt aus den Bahnen der Vektoren ablesen. Daher wird nach
\ref{einvek} statt $\hat{S}^\bot$ die Gruppe $$\hat{S}^*=\left\{
\begin{pmatrix}
  z         & x             & y    \\
  0         & m_1           & m_2  \\
  0         & m_4           & m_5
\end{pmatrix} \mid m_1, m_2, m_4, m_5, x,y,z \in \f \text{ und } z \neq 0 \neq m_1 m_5 - m_4
m_2\right\}$$ verwendet. Unter der Operation von $\hat{S}^*$ bilden die Vektoren mit den
Koordinatenvektoren $(0,0,0)$, $(1,0,0)$ und $(0,0,1)$ offensichtlich ein Vertretersystem der
Bahnen der Vektoren. Durch die Wahl der Basis $\{\hat{f}_4, \hat{f}_5, \hat{f}_6\}$ stellen
nach \ref{einvek} die Unterräume $\erz{\hat{f}_4}$ und $\erz{\hat{f}_6}$ ein Vertretersystem
der Bahnen der eindimensionalen Unterräume unter dieser Operation dar. Also sind nach dem
Dualitätsprinzip die Unterräume $\hat{W}_1=\erz{\hat{f}_5, \hat{f}_6}$ und
$\hat{W}_2=\erz{\hat{f}_4, \hat{f}_5}$ als orthogonale Komplemente von $\erz{\hat{f}_4}$ und
$\erz{\hat{f}_6}$ ein Vertretersystem der zweidimensionalen Unterräume von $\hat{V}$.

Nun wird wieder der Vektorraum $\bar{V}=V/U$ betrachtet. Mit $\hat{W}_1$ und $\hat{W}_2$ kennt
man die Repräsentanten der Bahnen zweidimensionaler Unterräume in $\hat{V}$. Lässt man die
Stabilisatoren von $\hat{W}_1$ und $\hat{W}_2$ auf $\bar{K}$ operieren, dann ergeben sich die
Bahnen der zweidimensionalen Unterräume $\bar{W}$ in $\bar{V}$ aus $\hat{W}_1$ bzw.
$\hat{W}_2$ und den Bahnen ihrer Stabilisatoren in $\bar{K}$.

Die Gruppe $\bar{S}$ operiert auf $\hat{V}$. Daher sind die Stabilisatoren von $\hat{W}_1$ und
$\hat{W}_2$ in $\bar{S}$ definiert. Die Bahnen der zweidimensionalen Unterräume $\bar{W}$ von
$\bar{V}$ erhält man nach \ref{tensor} dadurch, dass man die Operation der Stabilisatoren
$\bar{S}_1=Stab_{\bar{S}}(\hat{W}_1)$ bzw. $\bar{S}_2=Stab_{\bar{S}}(\hat{W}_2)$ auf dem
Komplement $\bar{K}$ zu $\bar{N}$ betrachtet und die Bahnen der Vektoren in $\bar{K}$ unter
der Operation von $\bar{S}_1$ bzw. $\bar{S}_2$ auf $\bar{W}_1/\bar{K} \otimes \bar{K}$ bzw.
$\bar{W}_2/\bar{K} \otimes \bar{K}$ ermittelt, wobei $\bar{W}_1=\erz{\bar{f}_5, \bar{f}_6}$
und $\bar{W}_2=\erz{\bar{f}_4,\bar{f}_5}$ ist. Die Basisvektoren der zweidimensionalen
Unterräume $\bar{W}$ von $\bar{V}$ mit $dim(\bar{W} \cap \bar{K})=0$ ergeben sich dann als
direkte Summen der Basisvektoren von $\bar{W}_1$ bzw. $\bar{W}_2$ und den Vertretern der
jeweiligen Bahnen unter $\bar{S}_1$ bzw. $\bar{S}_2$ in $\bar{K}$.

Es ist $$\bar{S}_1=Stab_{\bar{S}}(\hat{W}_1)=\left\{
\begin{pmatrix}
  m_1 m_5 - m_4 m_2 & 0                 & 0          & 0 \\
  0                 & m_1 m_9           & m_2 m_9    & 0 \\
  0                 & m_4 m_9           & m_5 m_9    & 0 \\
  0                 & 0                 & 0          & m_9
\end{pmatrix} \mid \triangle\right\}.$$ Dabei ist $\triangle= m_1, m_2, m_4, m_5,
m_9 \in \f$ und $m_9 \neq 0$ sowie $m_1 m_5 - m_4 m_2 \neq 0$. Also operiert $\bar{S}_1$ auf
$\bar{W}_1/\bar{K} \otimes \bar{K}$ als
\begin{eqnarray*}
  \tilde{S}_1 &=& \left\langle  \begin{pmatrix}
  m_1 m_9           & m_2 m_9 \\
  m_4 m_9           & m_5 m_9
\end{pmatrix}\1 \otimes (m_9) \mid \triangle  \right\rangle \\
&=& \left\langle m_9\1 \begin{pmatrix}
  m_1            & m_2  \\
  m_4            & m_5
\end{pmatrix}\1 \otimes (m_9) \mid \triangle  \right\rangle \\
&=& \left\langle   \begin{pmatrix}
  m_1            & m_2  \\
  m_4            & m_5
\end{pmatrix}\1 \mid \triangle  \right\rangle \\
&=& \left\{  \begin{pmatrix}
  m_1            & m_2  \\
  m_4            & m_5
\end{pmatrix} \mid \triangle \right\} \\
&=& GL(2,p).
\end{eqnarray*}
Dabei ist $\triangle= m_1, m_2, m_4, m_5, m_9 \in \f$ und $m_9 \neq 0$ sowie $m_1 m_5 - m_4
m_2 \neq 0$. Die Gruppe $GL(2,p)$ hat zwei Bahnen der Vektoren in $\f^2$, deren Vertreter zum
Beispiel $(0,0)$ und $(0,-1)$ sind. Daher können $\bar{W}_{1a}=\erz{\bar{f}_5+0 \cdot
\bar{f}_9, \bar{f}_6 + 0 \cdot \bar{f}_9}$ und $\bar{W}_{1b}=\erz{\bar{f}_5 + 0 \cdot
\bar{f}_9, \bar{f}_6 - 1 \cdot \bar{f}_9}$ als zu $\hat{W}_1$ gehörende Bahnenvertreter der
zweidimensionalen Unterräume von $\bar{V}$ gewählt werden, die einen trivialen Schnitt mit
$\bar{K}$ haben. Ihnen entsprechen die vierdimensionalen Unterräume $W_{1a}=\erz{f_5 + 0 \cdot
f_9, f_6 + 0 \cdot f_9, f_7, f_8}$ und $W_{1b}=\erz{f_5 + 0 \cdot f_9, f_6 - 1 \cdot f_9, f_7,
f_8}$ und damit die zulässigen Untergruppen $M_2$ und $M_3$.

Es ist $$\bar{S}_2=Stab_{\bar{S}}(\hat{W}_2)=\left\{
\begin{pmatrix}
  m_1 m_5 - m_4 m_2 & 0                 & 0          & 0 \\
  x                 & m_1 m_9           & 0          & 0 \\
  y                 & m_4 m_9           & m_5 m_9    & 0 \\
  0                 & 0                 & 0          & m_9
\end{pmatrix} \mid \triangle\right\}.$$ Dabei ist $\triangle= m_1, m_4, m_5,
m_9, x, y \in \f$ und $m_9 \neq 0$ sowie $m_1 m_5 - m_4 m_2 \neq 0$. Also operiert $\bar{S}_2$
auf $\bar{W}_2/\bar{K} \otimes \bar{K}$ als
\begin{eqnarray*}
  \tilde{S}_2 &=& \left\langle  \begin{pmatrix}
  m_1 m_5 - m_4 m_2  & 0       \\
  x                  & m_1 m_9
\end{pmatrix}\1 \otimes (m_9) \mid \triangle \right\rangle \\
&=& \left\langle \begin{pmatrix}
  \frac{1}{m_1 m_5 - m_4 m_2}  & 0  \\
  -\frac{x}{m_1 m_5 - m_4 m_2} & \frac{1}{m_1 m_9}
\end{pmatrix} \otimes (m_9) \mid \triangle  \right\rangle \\
&=& \left\langle   \begin{pmatrix}
  \frac{m_9}{m_1 m_5 - m_4 m_2}    & 0  \\
  -\frac{x m_9}{m_1 m_5 - m_4 m_2} & \frac{1}{m_9}
\end{pmatrix} \mid \triangle  \right\rangle \\
&=& \left\{  \begin{pmatrix}
  a            & 0  \\
  b            & c
\end{pmatrix} \mid b \in \f \text{ und } a,c \in \f^* \right\}
\end{eqnarray*}
Dabei ist $\triangle= m_1, m_2, m_4, m_5, m_9, x \in \f$ und $m_9 \neq 0$ sowie $m_1 m_5 - m_4
m_2 \neq 0$. Also operiert $\tilde{S}_2$ wie die Gruppe der unteren Dreiecksmatrizen von
$GL(2,p)$ auf $\f^2$. Ein Vertretersystem der Bahnen unter dieser Operation ist beispielsweise
$(0,0)$, $(0,-1)$ und $(-1,0)$. Daher können zuzüglich zu $\bar{W}_{1a}$ und $\bar{W}_{1b}$
die Unterräume $\bar{W}_{2a}=\erz{\bar{f}_4+0 \cdot \bar{f}_9, \bar{f}_5 + 0 \cdot
\bar{f}_9}$, $\bar{W}_{2b}=\erz{\bar{f}_4+ 0 \cdot \bar{f}_9, \bar{f}_5 - 1 \cdot \bar{f}_9}$
und $\bar{W}_{2c}=\erz{\bar{f}_4 - 1 \cdot \bar{f}_9, \bar{f}_5 + 0 \cdot \bar{f}_9}$ als
übrige der Bahnenvertreter der zweidimensionalen Unterräume von $\bar{V}$ gewählt werden, die
einen trivialen Schnitt mit $\bar{K}$ haben. Ihnen entsprechen die vierdimensionalen
Unterräume $W_{2a}=\erz{f_4 + 0 \cdot f_9, f_5 + 0 \cdot f_9, f_7, f_8}$, $W_{2b}=\erz{f_4 + 0
\cdot f_9, f_5 - 1 \cdot f_9, f_7, f_8}$ und $W_{2c}=\erz{f_4 - 1 \cdot f_9, f_5 + 0 \cdot
f_9, f_7, f_8}$ und damit die zulässigen Untergruppen $M_4$, $M_5$ und $M_6$.

\textbf{Dritter Fall:} Es sei $dim(W \cap K)=1$. Im Weiteren wird die folgende Notation
benutzt:
\begin{itemize}
  \item[] $U:=\{f_7\}$,
  \item[] $\bar{V}:=V/U$, $\bar{K}:=K/U$, $\bar{N}:=N/U$ und $\bar{W}:=W/U$,
  \item[] $\bar{f}_i:=f_i + U$ für $i \in \{4,5,6,8,9\}$,
\end{itemize}
Die Idee des Beweises ist dieselbe wie im zweiten Fall: Der Unterraum $U$ kann als Vertreter
der Bahn gewählt werden, den der Schnitt von $W$ in $K$ hat. Der Stabilisator von $U$ in $B$
wird berechnet und die Bahn von $W/U$ in $V/U$ ermittelt. Die Bahn von $W$ ergibt sich dann
aus der Bahn von $W/U$ und $U$.

Da $B$ transitiv auf $K \backslash\{0\}$ operiert, hat $W \cap K$ als Bahn unter $B$ die Menge
aller eindimensionalen Unterräume von $K$. Man kann daher als Vertreter von $(W \cap K)^B$ dem
Unterraum $U=\{f_7\}$ wählen. Für $i \in \{4,5,6,8,9\}$ sei $\bar{f}_i=f_i+U$. Nun betrachte
man die Operation von $B$ auf $\bar{V}$ bezüglich der Basis $\{\bar{f}_4, \bar{f}_5,
\bar{f}_6, \bar{f}_8, \bar{f}_9\}$. Der Stabilisator von $U$ in $A$ ist
$$S=\left\{\begin{pmatrix}
  m_1 m_5           & m_1 m_6           & 0                 & 0 & 0 & 0 \\
  m_1 m_8           & m_1 m_9           & 0                 & 0 & 0 & 0 \\
  m_4 m_8 - m_7 m_5 & m_4 m_9 - m_7 m_6 & m_5 m_9 - m_8 m_6 & 0 & 0 & 0 \\
  0                 & 0                 & 0                 & m_1 & 0 & 0 \\
  0                 & 0                 & 0                 & m_4 & m_5 & m_6 \\
  0                 & 0                 & 0                 & m_7 & m_8 & m_9
\end{pmatrix} \mid \triangle \right\},$$ wobei $\triangle = m_1, m_4, m_5, m_6, m_7, m_8, m_9
\in \f \text{ und } m_5 m_9 - m_8 m_6 \neq 0 \neq m_1$ ist. Damit operiert $S$ als $\bar{S}$
folgendermaßen auf $\bar{V}$: $$\bar{S}=\left\{\begin{pmatrix}
  m_1 m_5           & m_1 m_6           & 0                 & 0   & 0 \\
  m_1 m_8           & m_1 m_9           & 0                 & 0   & 0 \\
  m_4 m_8 - m_7 m_5 & m_4 m_9 - m_7 m_6 & m_5 m_9 - m_8 m_6 & 0   & 0 \\
  0                 & 0                 & 0                 & m_5 & m_6 \\
  0                 & 0                 & 0                 & m_8 & m_9
\end{pmatrix} \mid \triangle \right\},$$  wobei $\triangle = m_1, m_4, m_5, m_6, m_7, m_8, m_9
\in \f \text{ und } m_5 m_9 - m_8 m_6 \neq 0 \neq m_1$ ist. Da sowohl $m_4 m_8 - m_7 m_5$ als
auch $m_4 m_9 - m_7 m_6$ unabhängig von $m_5, m_6, m_8$ und $m_9$ jeden beliebigen Wert in
$\f$ annehmen, ist $$\bar{S}=\left\{\begin{pmatrix}
  m_1 m_5  & m_1 m_6  & 0                 & 0   & 0 \\
  m_1 m_8  & m_1 m_9  & 0                 & 0   & 0 \\
  x        & y        & m_5 m_9 - m_8 m_6 & 0   & 0 \\
  0        & 0        & 0                 & m_5 & m_6 \\
  0        & 0        & 0                 & m_8 & m_9
\end{pmatrix} \mid \triangle \right\},$$ wobei $\triangle = m_1, m5, m_6, m_8, m_9, x, y \in \f
\text{ und } m_5 m_9 - m_8 m_6 \neq 0 \neq m_1$ ist.

Der Unterraum $\bar{W}=W/U$ ist ein Komplement zu $\bar{K}=K/U$, denn andernfalls hätte der
Schnitt von $W$ mit $K$ eine höhere Dimension als 1. Daher gibt es $\bar{k}_4, \bar{k}_5$ und
$\bar{k}_6$ in $\bar{K}=K/U$, sodass $\bar{W}=\erz{-\bar{f}_4 + \bar{k}_4, -\bar{f}_5 +
\bar{k}_5, -\bar{f}_6 + \bar{k}_6}=\erz{\bar{f}_4 - \bar{k}_4, \bar{f}_5 - \bar{k}_5,
\bar{f}_6 - \bar{k}_6}$ ist. Da $\bar{K}$ die Dimension 2 hat, gibt es demnach eine Bijektion
zwischen $Hom(\bar{V},\bar{K}) \cong \F_p^6$ und der Menge der Komplemente zu $\bar{K}$ in
$\bar{V}$. Es sei $Kom$ die Menge der Komplemente $\bar{W}$ zu $\bar{K}$. Dann $\beta: \F_6
\ra Kom: (a,b,c,d,e,f) \mt \erz{\bar{f}_4 - a \cdot \bar{f}_8 - b \cdot \bar{f}_9, \bar{f}_5 -
c \cdot \bar{f}_8 - d \cdot \bar{f}_9 , \bar{f}_6 - e \cdot \bar{f}_8 - f \cdot \bar{f}_9}$
eine Bijektion mit der gewünschten Eigenschaft.

Nach \ref{tensor} lassen sich die Bahnen der dreidimensionalen Unterräume $\bar{W}$ mit den
gewünschten Eigenschaften finden, indem man die Operation der Gruppe $G$, die anschließend
definiert wird, durch Multiplikation von rechts auf $\f^6$ betrachtet und den Vektor
$(a,b,c,d,e,f)$ aus $\f^6$ über $\beta$ mit dem Unterraum $\bar{W}=\erz{-\bar{f}_4 + a \cdot
\bar{f}_8 + b \cdot \bar{f}_9, -\bar{f}_5 + c \cdot \bar{f}_8 + d \cdot \bar{f}_9, -\bar{f}_6
+ e \cdot \bar{f}_8 + f \cdot \bar{f}_9}=\erz{\bar{f}_4 - a \cdot \bar{f}_8 - b \cdot
\bar{f}_9, \bar{f}_5 - c \cdot \bar{f}_8 - d \cdot \bar{f}_9 , \bar{f}_6 - e \cdot \bar{f}_8 -
f \cdot \bar{f}_9}$ identifiziert. Die Gruppe $G$ sei im Einklang mit \ref{tensor} definiert
als $$G=\left\langle \begin{pmatrix}
  m_1 m_5  & m_1 m_6  & 0 \\
  m_1 m_8  & m_1 m_9  & 0 \\
  x        & y        & m_5 m_9 - m_8 m_6
\end{pmatrix}\1 \otimes \begin{pmatrix}
  m_5  & m_6  \\
  m_8  & m_9
\end{pmatrix} \mid \triangle  \right\rangle,$$ wobei $\triangle = m_1, m5, m_6, m_8, m_9, x, y \in \f
\text{ und } m_5 m_9 - m_8 m_6 \neq 0 \neq m_1$ ist.

Es sei $\{v_1, v_2, v_3, v_4, v_5, v_6\}$ die Standardbasis von $\f^6$ und $T=\erz{v_1, v_2,
v_3, v_4}$. Der Unterraum $T$ ist invariant unter der Operation von $G$ (invertiert man die
linke Matrix und multipliziert man das Kroneckerprodukt in der Definition von $G$ aus, sieht
man, dass die letzten beiden Spalten in den oberen vier Zeilen nur Nullen enthalten). Die
Operation von $G$ auf $T$ entspricht nach \ref{tensor} der Operation von $Z \times GL(2,p)$
vermittels $\delta: G \times M(2 \times 2,p): ((z,g),m) \mt m^{(z,g)}=z g\1 m g$ auf $M(2
\times 2, p)$, wobei $Z$ das Zentrum von $GL(2,p)$ ist. Daher lassen sich $p+5$ Bahnen der
Operation von $G$ auf $\f^6$ aus \ref{kro} ablesen, indem man die Matrix $$\begin{pmatrix}
  a & b \\
  c & d
\end{pmatrix} \in M(2 \times 2, p)$$ mit dem Vektor $v_{abcd}:=a v_1 + b v_2 + c v_3 + d v_4 \in \f^6$
identifiziert. Dem Koordinatenvektor $(a,b,c,d,0,0)$ von $v_{abcd}$ entspricht nach $\beta$
der Unterraum $$\bar{W}=\erz{\bar{f}_4 - a \bar{f}_8 - b \bar{f}_9, \bar{f}_5 - c \bar{f}_8 -
d \bar{f}_9, \bar{f}_6}$$ von $\bar{V}$ bzw. der Unterraum $$W=\erz{f_4 - a f_8 - b f_9, f_5-
c f_8 - d f_9, f_6, f_7}$$ von $V$ und damit die zulässige Untergruppe $$M_{abcd}=\erz{a_4
a_8^{-a} a_9^{-b}, a_5 a_8^{-c} a_9^{-d}, a_6, a_7}.$$ Liest man die $p+5$ Bahnen aus
\ref{kro} ab und identifiziert einen dort angegebenen Bahnenvertreter $$\begin{pmatrix}
  a & b \\
  c & d
\end{pmatrix} \in M(2 \times 2, p)$$ mit der zulässigen Untergruppe $M_{abcd}=\erz{a_4
a_8^{-a} a_9^{-b}, a_5 a_8^{-c} a_9^{-d}, a_6, a_7}$, so erhält man die Einträge $M_{10}$ bis
$M_{17}^k$ aus dem Vertretersystem zulässiger Untergruppen.

Damit sind die Bahnen aller zulässigen Untergruppen ermittelt worden, denen Vektoren in $T$
entsprechen. Die übrigen Bahnen erhält man folgendermaßen: Es sei $L=\erz{v_5, v_6}$. Alle
Bahnen in $T$ sind bereits bestimmt. Man betrachte nun den Faktorraum $\f^6/T$, auf dem $G$ in
natürlicher Weise operiert. Da $(L/T)^G=L/T$ ist und damit ganz $L/T$ eine einzige Bahn
bildet, kann $v_6 + T$ als Repräsentant dieser Bahn gewählt werden. Daher sind nur noch die
Bahnen in $T$ unter der Operation von $\hat{S}=Stab_G(v_6 + T)$ zu bestimmen. Die Vertreter
dieser Bahnen bilden als direkte Summe mit $v_6$ ein Vertretersystem der noch fehlenden
Bahnen.

Die Bahn $B$ von $v_6$ unter $G$ ist
\begin{eqnarray*}
  B &=& \bigg\{
        \frac{m_9x - m_8y}{m_1 s^2} m_8 v_1 + \frac{m_9x - m_8y}{m_1 s^2} m_9 v_2 + \frac{m_6x -
        m_5y}{m_1 s^2} m_8 v_3 \\
    & & + \frac{m_6x - m_5y}{m_1 s^2} m_9 v_4 + s\1 m_8 v_5 + s\1 m_9 v_6 \mid
        \triangle \bigg\},
\end{eqnarray*}
wobei $s=m_5 m_9 - m_8 m_6$ und $\triangle = m_1, m5, m_6, m_8, m_9, x, y \in \f \text{ und }
m_5 m_9 - m_8 m_6 \neq 0 \neq m_1$ ist. Daher ist der Stabilisator von $v_6 + T$ durch die
Bedingung $m_5=1$ und $m_8=0$ gegeben. Also ist $$\hat{S}=Stab_G(v_6 + T)=\left\langle
\begin{pmatrix}
  m_1      & m_1 m_6  & 0 \\
  0        & m_1 m_9  & 0 \\
  x        & y        & m_9
\end{pmatrix}\1 \otimes \begin{pmatrix}
  1    & m_6  \\
  0    & m_9
\end{pmatrix} \mid \triangle  \right\rangle,$$ wobei $\triangle = m_1, m_6, m_9, x, y \in \f
\text{ und } m_5 m_9 - m_8 m_6 \neq 0 \neq m_1$ ist. Die Bahn $\hat{B}$ von $v_6$ unter
$\hat{S}$ ist
\begin{eqnarray*}
  \hat{B} &=&\bigg\{\frac{m_9x}{m_1} v_2 + \frac{m_6x - y}{m_1 m_9} v_4 + v_6 \mid
             m_6,x,y \in \f \text{ und } m_1, m_9 \in \f^* \bigg\}\\
    &=&\{\bar{x} v_2 + \bar{y} v_4 + v_6 \mid \bar{x}, \bar{y} \in \f\}\\
    &=& v_6 + \erz{v_2, v_4}
\end{eqnarray*}

Da $\hat{B}/L$ den gesamten Unterraum $\erz{v_2, v_4}/L$ enthält, brauchen nicht die Bahnen
von $\hat{S}$ auf $T$ ermittelt zu werden, sondern lediglich die Bahnen, die $\hat{S}$ unter
der natürlichen Operation auf dem Faktorraum $T/\erz{v_2, v_4}$ hat. Dieser Operation von
$\hat{S}$ entspricht die in \ref{dreieckvbar} dargestellte Operation. Zwischen einem dort
angegebenen Bahnenvertreter $(a,b) \in \f^2$ bestehen die folgenden Bijektionen zu den hier
betrachteten Vektoren, Unterräumen und zulässigen Untergruppen:
\begin{eqnarray*}
  (a,b) &\leftrightarrow& a \cdot v_1 + b \cdot v_3 + v_6 \in \f^6 \\
        &\leftrightarrow& \erz{\bar{f}_4 - a \bar{f}_8, \bar{f}_5 - b \bar{f}_8, \bar{f}_6 - \bar{f}_9} \leq \bar{V}\\
        &\leftrightarrow& \erz{f_4 - a f_8, f_5 - b f_8, f_6 - f_9, f_7} \leq V \\
        &\leftrightarrow& \erz{a_4 a_8^{-a}, a_5 a_8^{-b}, a_6 a_9\1, a_7} \leq M(G) \\
\end{eqnarray*}
In \ref{dreieckvbar} sind die Bahnenvertreter $(0,0)$, $(0,1)$ und $(1,0)$ ermittelt worden.
Ihnen entsprechen die zulässigen Untergruppen $M_7$, $M_8$ und $M_9$.
\end{beweis}

\subsection{Nachfolger der Ordnung $p^5$}

\begin{satz} \label{ncp3p5}
Die unmittelbaren Nachfolger von $G$ der Ordnung $p^5$ zerfallen in $14+p$ Isomorphieklassen.
Mit der folgenden Liste von Gruppen ist ein Vertretersystem dieser Isomorphieklassen
angegeben, wobei $w$ ein Erzeuger der multiplikativen Gruppe von $\f$ ist:
\begin{enumerate}
  \item[] $G_1=\erz{g_1, g_2, g_3, g_4, g_5 \mid [g_2,g_1]=g_4, [g_3,g_1]=g_5}$,
  \item[] $G_2=\erz{g_1, g_2, g_3, g_4, g_5 \mid [g_2,g_1]=g_4, g_3^p=g_5}$,
  \item[] $G_3=\erz{g_1, g_2, g_3, g_4, g_5 \mid [g_2,g_1]=g_4,
  [g_3,g_2]=g_5, g_3^p=g_5}$,
  \item[] $G_4=\erz{g_1, g_2, g_3, g_4, g_5 \mid  [g_3,g_2]=g_4, g_3^p=g_5}$,
  \item[] $G_5=\erz{g_1, g_2, g_3, g_4, g_5 \mid  [g_3,g_1]=g_4,
  [g_3,g_2]=g_5, g_3^p=g_5}$,
  \item[] $G_6=\erz{g_1, g_2, g_3, g_4, g_5 \mid [g_2,g_1]=g_4,
  [g_3,g_2]=g_5, g_3^p=g_4}$,
  \item[] $G_7=\erz{g_1, g_2, g_3, g_4, g_5 \mid [g_3,g_2]=g_5, g_2^p=g_4, g_3^p=g_5}$,
  \item[] $G_8=\erz{g_1, g_2, g_3, g_4, g_5 \mid [g_3,g_1]=g_4,
  [g_3,g_2]=g_5, g_2^p=g_4, g_3^p=g_5}$,
  \item[] $G_9=\erz{g_1, g_2, g_3, g_4, g_5 \mid [g_2,g_1]=g_4,
  [g_3,g_2]=g_5, g_2^p=g_4, g_3^p=g_5}$,
  \item[] $G_{10}=\erz{g_1, g_2, g_3, g_4, g_5 \mid g_2^p=g_4, g_3^p=g_5}$,
  \item[] $G_{11}=\erz{g_1, g_2, g_3, g_4, g_5 \mid [g_2,g_1]=g_4, [g_3,g_1]=g_5,
  g_2^p=g_4, g_3^p=g_5}$,
  \item[] $G_{12}=\erz{g_1, g_2, g_3, g_4, g_5 \mid [g_2,g_1]=g_5,
  g_2^p=g_4, g_3^p=g_5}$,
  \item[] $G_{13}=\erz{g_1, g_2, g_3, g_4, g_5 \mid [g_2,g_1]=g_4 g_5, [g_3,g_1]=g_5,
  g_2^p=g_4, g_3^p=g_5}$,
  \item[] $G_{14}=\erz{g_1, g_2, g_3, g_4, g_5 \mid [g_3,g_1]=g_5,
  g_2^p=g_4, g_3^p=g_5}$,
  \item[] $G_{15}^k=\erz{g_1, g_2, g_3, g_4, g_5 \mid [g_2,g_1]=g_4, [g_3,g_1]=g_5^{w^k},
  g_2^p=g_4, g_3^p=g_5}$ mit $k \in \{1 \dd \frac{p-1}{2}\}$,
  \item[] $G_{16}=\erz{g_1, g_2, g_3, g_4, g_5 \mid [g_2,g_1]=g_5^w, [g_3,g_1]=g_4,
  g_2^p=g_4, g_3^p=g_5}$,
  \item[] $G_{17}=\erz{g_1, g_2, g_3, g_4, g_5 \mid [g_2,g_1]=g_4 g_5^{w^k},
  [g_3,g_1]=g_4^{w^{k-1}} g_5, g_2^p=g_4, g_3^p=g_5}$ mit $k \in \{1 \dd \frac{p-1}{2}\}$.
\end{enumerate}
\end{satz}

\begin{beweis}
Dieser Satz ergibt sich nach \ref{bahnenzu}, indem man $P(G)$ nach jedem Repräsentanten des
Vertretersystems zulässiger Untergruppen faktorisiert, das in \ref{zu7} aufgelistet ist. Die
Präsentation der Faktorgruppe erhält man über das Verfahren aus \ref{faktor}.
\end{beweis}

\section{Nachfolger von $(p^4,11)$}

\begin{ver}
In diesem Abschnitt sei $$G=\erz{a_1, a_2, a_3, a_4 \mid a_1^p=a_4}.$$ Damit ist $G \cong
C_{p^2} \times C_p^2$.
\end{ver}

\begin{lemma}
Die Gruppe hat die Gewichtung $\omega(a_1)=\omega(a_2)=\omega(a_3)=1$ und $\omega(a_4)=2$.
\end{lemma}

\begin{beweis}
Die Gewichtung von $G$ lässt sich nach \ref{gew} ermitteln.
\end{beweis}

\subsection{Die Automorphismengruppe}

\begin{lemma} \label{aut11}
Jeder Automorphismus $\al$ von $G_2$ lässt sich durch $$\al(u_1, u_2, u_3, u_4, v_2, v_3, v_4,
w_2, w_3, w_4 ):
  \begin{cases}
    a_1 \mt a_1^{u_1} a_2^{u_2} a_3^{u_3} a_4^{u_4}  & \text{wobei } u_1, u_2, u_3, u_4, v_2, v_3, v_4\\
    a_2 \mt           a_2^{v_2} a_3^{v_3} a_4^{v_4}  & w_2, w_3, w_4 \in \{0 \dd p-1\}\\
    a_3 \mt           a_2^{w_2} a_3^{w_3} a_4^{w_4}  & \text{und } v_2 w_3 - w_2 v_3 \neq 0 \\
    a_4 \mt                               a_4^{u_1}    & \text{wie auch } u_1 \neq 0 \text{ ist.}
  \end{cases}$$ darstellen.
\end{lemma}

\begin{beweis}
Die Automorphismengruppe von $G$ kann nach \ref{aut} aus dem Stabilisator $S_2$ (siehe
\ref{bahnen_cp3_p4}) der zu $G$ gehörenden zulässigen Untergruppe unter Berücksichtigung des
Transponierens unmittelbar abgelesen werden. Damit können auch die Bilder von $a_1$, $a_2$ und
$a_3$ direkt abgelesen werden. Das Bild von $a_4$ ist durch die Bilder von $a_1$, $a_2$ und
$a_3$ eindeutig bestimmt und lässt sich folgendermaßen berechnen:
\begin{eqnarray*}
  a_4^\al &=& (a_1^p)^\al = (a_1^\al)^p \\
  &=& (a_1^{u_1} a_2^{u_2} a_3^{u_3} a_4^{u_4})^p \\
  &=& (a_1^{u_1})^p\\
  &=& (a_1^p)^{u_1}\\
  &=& a_4^{u_1},
\end{eqnarray*}
da $G_2$ abelsch ist und $a_2$, $a_3$ und $a_4$ die Ordnung $p$ haben.
\end{beweis}

\subsection{$p$"=Cover, Multiplikator und Nukleus}

\begin{lemma}
Für $G$ ergibt sich das $p$"=Cover, der Nukleus und der Multiplikator in folgender Weise:
\begin{enumerate}
  \item[] $P(G)=\erz{a_1 \dd a_{10} \mid a_1^p=a_4, a_4^p=a_5, [a_2,a_1]=a_6,
  [a_3,a_1]=a_7, [a_3,a_2]=a_8,a_2^p=a_9, a_3^p=a_{10}}$.
  \item[] $M(G)=\erz{a_5, \dd a_{10}}$.
  \item[] $N(G)=\erz{a_5}$.
\end{enumerate}
\end{lemma}

\begin{beweis}
Die Präsentation von $P(G)$ ergibt sich aus \ref{cov} und \ref{newman} unter Verwendung des
Reduktionsverfahrens von Knuth"=Bendix, aus der sich $M(G)$ ablesen lässt. Die Gewichtung von
$P(G)$
\begin{enumerate}
  \item[] $1=\omega(a_1)=\omega(a_2)=\omega(a_3)$.
  \item[] $2=\omega(a_4)=\omega(a_6)=\omega(a_7)=\omega(a_8)=\omega(a_9)=\omega(a_{10})$.
  \item[] $3=\omega(a_5)$.
\end{enumerate}
und damit $N(G)$ lassen sich nach \ref{gew} aus der Präsentation von $P(G)$ ablesen.
\end{beweis}

\subsection{Operation der Erweiterungsautomorphismen}

\begin{lemma}
Werden die Automorphismen $\al$ von $G$ gemäß \ref{aut11} in der Weise $$\al(u_1, u_2, u_3,
u_4, v_2, v_3, v_4, w_2, w_3, w_4):
  \begin{cases}
    a_1 \mt a_1^{u_1} a_2^{u_2} a_3^{u_3} a_4^{u_4}  & \text{wobei } u_1, u_2, u_3, u_4, v_2, v_3, v_4\\
    a_2 \mt           a_2^{v_2} a_3^{v_3} a_4^{v_4}  & w_2, w_3, w_4 \in \{0 \dd p-1\}\\
    a_3 \mt           a_2^{w_2} a_3^{w_3} a_4^{w_4}  & \text{und } v_2 w_3 - w_2 v_3 \neq 0 \\
    a_4 \mt                               a_4^{u_1}  & \text{wie auch } u_1 \neq 0 \text{ ist.}
  \end{cases}$$ dargestellt, so ist $$\vi:
Aut(G) \to GL(6,p): \al(u_1, u_2, u_3, u_4, v_2, v_3, v_4, w_2, w_3, w_4) \mt $$ $$M =
\begin{pmatrix}
  u_1 & 0       & 0       & 0                 & 0   & 0 \\
  0   & u_1 v_2 & u_1 v_3 & u_2 v_3 - u_3 v_2 & 0   & 0 \\
  0   & u_1 w_2 & u_1 w_3 & u_2 w_3 - u_3 w_2 & 0   & 0 \\
  0   & 0       & 0       & v_2 w_3 - v_3 w_2 & 0   & 0 \\
  v_4 & 0       & 0       & 0                 & v_2 & v_3 \\
  w_4 & 0       & 0       & 0                 & w_2 & w_3
\end{pmatrix}$$ der Operationshomomorphismus von $A=Aut(G)$ auf $\f^6 \cong M(G)$
über Erweiterungsautomorphismen.
\end{lemma}

\begin{beweis}
Nach \ref{aut11} lässt sich jeder Automorphismus $\al$ von $G$ in der oben angegebenen Form
darstellen. Die Bilder der Erzeuger von $M(G)$ ergeben sich unter $\al$ über $\vi$
folgendermaßen: Es ist
\begin{eqnarray*}
  a_5^\al&=& (a_4^p)^\al = (a_4^\al)^p = (a_4^{u_1})^p = (a_4^p)^{u_1} = a_5^{u_1}
\end{eqnarray*}
und
\begin{eqnarray*}
  a_6^\al &=& [a_2,a_1]^\al=[a_2^\al, a_1^\al]\\
  &=& [a_2^{v_2} a_3^{v_3} a_4^{v_4}, a_1^{u_1} a_2^{u_2} a_3^{u_3} a_4^{u_4}]\\
  &=& [a_2, a_1]^{u_1 v_2} [a_3,a_2]^{-u_3 v_2} [a_3, a_1]^{u_1 v_3} [a_3, a_2]^{u_2 v_3} \\
  &=& a_6^{u_1 v_2} a_7^{u_1 v_3} a_8^{u_2 v_3-u_3 v_2},
\end{eqnarray*}
da alle anderen auftretenden Kommutatoren trivial sind. Analog erhält man, dass
auch $a_7^\al = a_6^{u_1 w_2} a_7^{u_1 w_3} a_8^{u_2 w_3-u_3 w_2}$ gilt. Weiterhin ist
\begin{eqnarray*}
  a_8^\al &=& [a_3, a_1]^\al = [a_3^\al, a_2^\al] \\
  &=& [a_2^{w_2} a_3^{w_3} a_4^{w_4}, a_2^{v_2} a_3^{v_3} a_4^{v_4}] \\
  &=& [a_3,a_2]^{v_2 w_3 - v_3 w_2}\\
  &=& a_8^{v_2 w_3 - v_3 w_2},
\end{eqnarray*}
da alle anderen auftretenden Kommutatoren trivial sind. Außerdem ist
\begin{eqnarray*}
  a_9^\al &=& (a_2^p)^\al = (a_2^{v_2} a_3^{v_3} a_4^{v_4})^p \\
  &=& (a_2^{v_2})^p (a_3^{v_3})^p (a_4^{v_4})^p\\
  &=& (a_2^p)^{v_2} (a_3^p)^{v_3} (a_4^p)^{v_4}\\
  &=& a_5^{v_4} a_9^{v_2} a_{10}^{v_3},
\end{eqnarray*}
da die Kommutatoren $[x,y]$ mit $x, y \in \{a_2, a_3, a_4\}$ trivial sind oder
die Ordnung $p$ haben. Ebenso ergibt sich $a_{10}=a_5^{w_4} a_9^{w_2}
a_{10}^{w_3}$.
\end{beweis}

Da im Weiteren für die unmittelbaren Nachfolger der Ordnung $p^5$ statt der
fünfdimensionalen Unterräume ihre eindimensionalen Komplemente betrachtet
werden, ist der zu $\vi$ duale Operationshomomorphismus von Interesse.

\begin{beme}
Der zu $\vi$ duale Operationshomomorphismus ist $$\bar{\vi}: Aut(G) \ra GL(6,p): \al(u_1, u_2,
u_3, u_4, v_2, v_3, v_4, w_2, w_3, w_4) \mt$$ $$
\begin{pmatrix}
  u_1 & 0                 & 0                 & 0                 & v_4 & w_4 \\
  0   & u_1 v_2           & u_1 w_2           & 0                 & 0   & 0 \\
  0   & u_1 v_3           & u_1 w_3           & 0                 & 0   & 0 \\
  0   & u_2 v_3 - u_3 v_2 & u_2 w_3 - u_3 w_2 & v_2 w_3 - v_3 w_2 & 0   & 0 \\
  0   & 0                 & 0                 & 0                 & v_2 & w_2 \\
  0   & 0                 & 0                 & 0                 & v_3 & w_3
\end{pmatrix}.$$
\end{beme}

\subsection{Bahnen zulässiger Untergruppen für Nachfolger der Ordnung $p^5$}

\begin{lemma}
Unter der Operation von $A=Aut(G)$ vermittels $\bar{\vi}$ zerfällt $\f^6$ in folgende Bahnen
nichttrivialer Vektoren:
\begin{enumerate}
  \item[] $B_1=(1,0,0,0,0,0)^A$ mit $|B_1|=p^2(p-1)$,
  \item[] $B_2=(1,0,0,1,0,0)^A$ mit $|B_2|=p^4(p-1)^2$,
  \item[] $B_3=(1,0,1,0,0,0)^A$ mit $|B_3|=p^2(p+1)(p-1)^2$,
  \item[] $B_4=(0,0,0,0,0,1)^A$ mit $|B_4|=p^2-1$,
  \item[] $B_5=(0,0,0,1,0,0)^A$ mit $|B_5|=p^2(p-1)$,
  \item[] $B_6=(0,0,0,1,0,1)^A$ mit $|B_6|=p^2(p-1)(p^2-1)$,
  \item[] $B_7=(0,0,1,0,0,0)^A$ mit $|B_7|=p^2-1$,
  \item[] $B_8=(0,0,1,0,0,1)^A$ mit $|B_8|=(p-1)(p^2-1)$,
  \item[] $B_9=(0,0,1,0,1,0)^A$ mit $|B_9|=p(p-1)(p^2-1)$.
\end{enumerate}
Zu den Vertretern der Bahnen $B_2, B_3$ und $B_9$ gehören folgende
Stabilisatoren:
\begin{eqnarray*}
  \bar{S}_2 &=& \left\{ \begin{pmatrix}
  1   & 0                 & 0                 & 0                 & 0   & 0  \\
  0   & v_2               & w_2               & 0                 & 0   & 0 \\
  0   & v_3               & w_3               & 0                 & 0   & 0 \\
  0   & 0                 & 0                 & v_2 w_3 - v_3 w_2 & 0   & 0 \\
  0   & 0                 & 0                 & 0                 & v_2 & w_2 \\
  0   & 0                 & 0                 & 0                 & v_3 & w_3
  \end{pmatrix} \mid  v_2, v_3, w_2, w_3 \in \f \text{, } v_2 w_3 - v_3 w_2 = 1 \right\}  \\
  \bar{S}_3 &=& \left\{ \begin{pmatrix}
  1   & 0                 & 0                 & 0                 & 0   & 0 \\
  0   & v_2               & w_2               & 0                 & 0   & 0 \\
  0   & 0                 & 1                 & 0                 & 0   & 0 \\
  0   &         - u_3 v_2 & u_2     - u_3 w_2 & v_2               & 0   & 0 \\
  0   & 0                 & 0                 & 0                 & v_2 & w_2 \\
  0   & 0                 & 0                 & 0                 & 0   & 1
  \end{pmatrix} \mid u_2, u_3, v_2, w_2 \in \f \text{ und } v_2 \neq 0 \right\}  \\
\bar{S}_9 &=& \left\{ \begin{pmatrix}
  u_1 & 0                 & 0                 & 0                 & v_4 & w_4 \\
  0   & u_1               & 0                 & 0                 & 0   & 0 \\
  0   & 0                 & 1                 & 0                 & 0   & 0 \\
  0   & - u_3 v_2         & u_2 u_1\1         & v_2 u_1\1         & 0   & 0 \\
  0   & 0                 & 0                 & 0                 & 1   & w_2 \\
  0   & 0                 & 0                 & 0                 & 0   & w_3
  \end{pmatrix} \mid u_1, u_2, u_3, v_4, w_4 \in \f \text{ und } u_1 \neq 0 \right\} \\
\end{eqnarray*}
\end{lemma}

\begin{beweis}
Die meisten Bahnenlängen lassen sich unmittelbar erkennen. Für die weniger offensichtlichen
Fälle sind die Stabilisatoren angegeben, deren Ordnung sich durch Abzählen ermitteln lässt.
Die Stabilisatoren $\bar{S}_2$, $\bar{S}_3$ und $\bar{S}_9$ ergeben sich aus den Lösungsmengen
der folgenden Gleichungssysteme:
\begin{enumerate}
  \item[2.] $\{u_1=1, u_2 v_3 - u_3 v_2 =0, u_2 w_3 - u_3 w_2=0, v_2 w_3 - v_3 w_2=1 \}$,
  \item[3.] $\{u_1=1, u_1 v_3=0, u_1 w_3=1 \}$,
  \item[9.] $\{u_1 v_3=0, u_1 w_3=1, v_2=1, w_2=0 \}$.
\end{enumerate}
Über die Bahnensätze erhält man dann die Länge der Bahnen, wobei $|A^{\bar{\vi}}|=|GL(n,p)|
\cdot p^3 (p-1)=p^5(p+1)(p-1)^2$ ist. Die Vollständigkeit der Liste erkennt man dadurch, dass
die Summation der Bahnenlänge die Mächtigkeit von $\f^6 \backslash \{0\}$ ergibt. Die Bahnen
verschiedener Länge sind verschieden und die beiden Bahnen gleicher Länge, nämlich $B_4$ und
$B_7$, offensichtlich auch, da ihre Vertreter in zwei verschiedenen irreduziblen Teilmoduln
von $\f^6$ unter der Operation von $\bar{\vi}$ liegen.
\end{beweis}

\begin{folg} \label{zu8}
Für die Nachfolger von $G$ der Ordnung $p^5$ ist die folgende Liste ein
Vertretersystem zulässiger Untergruppen.
\begin{enumerate}
  \item[] $M_1=\erz{a_6,a_7,a_8,a_9,a_{10}}$ entsprechend $\erz{(1,0,0,0,0,0)}^\bot$.
  \item[] $M_2=\erz{a_6,a_7,a_5 a_8\1,a_9,a_{10}}$ entsprechend $\erz{(1,0,0,1,0,0)}^\bot$.
  \item[] $M_3=\erz{a_6,a_5 a_7\1, a_8, a_9, a_{10}}$ entsprechend $\erz{(1,0,1,0,0,0)}^\bot$.
\end{enumerate}
\end{folg}

\begin{beweis}
Nach \ref{suppl} sind nur Supplemente des Nukleus $N=\erz{a_5}$, entsprechend
dem Unterraum $\bar{N}=\erz{(1,0,0,0,0,0)}$, zulässige Untergruppen und
wiederum sind genau fünfdimensionalen Unterräume Supplemente zu $\bar{N}$,
deren eindimensionale Komplemente einen Basisvektor haben, der mit dem
Basisvektor von $\bar{N}$ ein Skalarprodukt ungleich Null liefert. Damit
beschränkt sich die Auswahl zulässiger Untergruppen auf $B_1, B_2$ und $B_3$.
\end{beweis}

\subsection{Nachfolger der Ordnung $p^5$}

\begin{satz} \label{nt11}
In der folgenden Liste sind sämtliche Nachfolger von $G$ der Ordnung $p^5$ bis
auf Isomorphie angegeben:
\begin{enumerate}
  \item[] $G_1=\erz{g_1, g_2, g_3, g_4, g_5 \mid g_1^p=g_4, g_4^p=g_5}$,
  \item[] $G_2=\erz{g_1, g_2, g_3, g_4, g_5 \mid g_1^p=g_4, g_4^p=g_5,
  [g_3,g_2]=g_5}$,
  \item[] $G_3=\erz{g_1, g_2, g_3, g_4, g_5 \mid g_1^p=g_4, g_4^p=g_5, [g_3,g_1]=g_5}$.
\end{enumerate}
\end{satz}

\begin{beweis}
Dieser Satz ergibt sich nach \ref{bahnenzu}, indem man $P(G)$ nach jedem Repräsentanten des
Vertretersystems zulässiger Untergruppen faktorisiert, das in \ref{zu8} aufgelistet ist. Die
Präsentation der Faktorgruppe erhält man über das Verfahren aus \ref{faktor}.
\end{beweis}

\section{Nachfolger von $(p^4,12)$}

\begin{ver}
In diesem Abschnitt sei $$G=\erz{a_1, a_2, a_3, a_4 \mid [a_2,a_1]=g_4, a_1^p=a_4}.$$
\end{ver}

\begin{lemma}
Die Gruppe $G$ hat die Gewichtung $\omega(a_1)=\omega(a_2)=\omega(a_3)=1$ und $\omega(a_4)=2$.
\end{lemma}

\begin{beweis}
Die Gewichtung von $G$ lässt sich nach \ref{gew} ermitteln.
\end{beweis}

\subsection{Die Automorphismengruppe}

\begin{lemma} \label{aut12}
Jeder Automorphismus $\al$ von $G$ lässt sich durch $$\al(u_1, u_2, u_3, u_4, v_1, v_2, v_3,
v_4, w_3, w_4 ):
  \begin{cases}
    a_1 \mt a_1^{u_1} a_2^{u_2} a_3^{u_3} a_4^{u_4}  & \text{wobei } u_1, u_2, u_3, u_4, v_1, v_2, v_3,\\
    a_2 \mt a_1^{v_1} a_2^{v_2} a_3^{v_3} a_4^{v_4}  & v_4, w_3, w_4 \in \{0 \dd p-1\}\\
    a_3 \mt                     a_3^{w_3} a_4^{w_4}  & \text{und } u_1 v_2 - v_1 u_2 \neq 0 \\
    a_4 \mt                               a_4^{v_2 u_1 - v_1 u_2}    & \text{wie auch } w_3 \neq 0 \text{ ist.}
  \end{cases}$$ darstellen.
\end{lemma}

\begin{beweis}
Die Automorphismengruppe von $G$ kann nach \ref{aut} aus dem Stabilisator $S_1$ (siehe
\ref{bahnen_cp3_p4}) der zu $G$ gehörenden zulässigen Untergruppe unter Berücksichtigung des
Transponierens unmittelbar abgelesen werden. Damit können auch die Bilder von $a_1$, $a_2$ und
$a_3$ direkt abgelesen werden. Das Bild von $a_4$ ist durch die Bilder von $a_1$, $a_2$ und
$a_3$ eindeutig bestimmt und lässt sich folgendermaßen berechnen:
\begin{eqnarray*}
  a_4^\al &=& [a_2, a_1]^\al = [a_2^\al, a_1^\al] \\
  &=& [a_1^{v_1} a_2^{v_2} a_3^{v_3} a_4^{v_4}, a_1^{u_1} a_2^{u_2} a_3^{u_3} a_4^{u_4}] \\
  &=& [a_1^{v_1} a_2^{v_2} a_3^{v_3}, a_1^{u_1} a_2^{u_2} a_3^{u_3}] \\
  &=& [a_1^{v_1} a_2^{v_2}, a_1^{u_1} a_2^{u_2} a_3^{u_3}]^{a_3^{v_3}} [a_3^{v_3}, a_1^{u_1} a_2^{u_2} a_3^{u_3}]\\
  &=& [a_1^{v_1} a_2^{v_2}, a_1^{u_1} a_2^{u_2} a_3^{u_3}]^{a_3^{v_3}} \\
  &=& ([a_1^{v_1}, a_1^{u_1} a_2^{u_2} a_3^{u_3}]^{a_2^{v_2}} [a_2^{v_2}, a_1^{u_1} a_2^{u_2} a_3^{u_3}])^{a_3^{v_3}} \\
  &=& [a_1,a_2]^{v_1 u_2} [a_2, a_1]^{v_2 u_1}\\
  &=& [a_2, a_1]^{v_2 u_1 - v_1 u_2}\\
  &=& a_4^{v_2 u_1 - v_1 u_2},
\end{eqnarray*}
da alle anderen auftretenden Kommutatoren trivial sind und $[a_2, a_1]$ zentral ist.
\end{beweis}

\subsection{$p$"=Cover, Multiplikator und Nukleus}

\begin{lemma}
Für $G$ ergibt sich das $p$"=Cover, der Nukleus und der Multiplikator in folgender Weise:
\begin{enumerate}
  \item[] $P(G)=\erz{a_1 \dd a_{11} \mid [a_2,a_1]=a_4, [a_4,a_1]=a_5, [a_4,a_2]=a_6,
  [a_3,a_1]=a_7, [a_3,a_2]=a_8, a_1^p=a_9, a_2^p=a_{10}, a_3^p=a_{11}}$.
  \item[] $M(G)=\erz{a_5, \dd a_{11}}$.
  \item[] $N(G)=\erz{a_5,a_6}$.
\end{enumerate}
\end{lemma}

\begin{beweis}
Die Präsentation von $P(G)$ ergibt sich aus \ref{cov} und \ref{newman} und dem
Reduktionsverfahren von Knuth"=Bendix, aus der sich $M(G)$ ablesen lässt. Die Gewichtung von
$P(G)$
\begin{enumerate}
  \item[] $1=\omega(a_1)=\omega(a_2)=\omega(a_3)$.
  \item[] $2=\omega(a_4)=\omega(a_7)=\omega(a_8)=\omega(a_9)=\omega(a_{10})=\omega(a_{11})$.
  \item[] $3=\omega(a_5)=\omega(a_6)$.
\end{enumerate}
und damit $N(G)$ lässt sich nach \ref{gew} aus der Präsentation von $P(G)$ ablesen.
\end{beweis}

\subsection{Operation der Erweiterungsautomorphismen}

\begin{lemma}
Werden die Automorphismen $\al$ von $G$ gemäß \ref{aut12} in der Weise $$\al(u_1, u_2, u_3,
u_4, v_1, v_2, v_3, v_4, w_3, w_4 ):
  \begin{cases}
    a_1 \mt a_1^{u_1} a_2^{u_2} a_3^{u_3} a_4^{u_4}  & \text{wobei } u_1, u_2, u_3, u_4, v_1, v_2, v_3,\\
    a_2 \mt a_1^{v_1} a_2^{v_2} a_3^{v_3} a_4^{v_4}  & v_4, w_3, w_4 \in \{0 \dd p-1\}\\
    a_3 \mt                     a_3^{w_3} a_4^{w_4}  & \text{und } u_1 v_2 - v_1 u_2 \neq 0 \\
    a_4 \mt                               a_4^{v_2 u_1 - v_1 u_2}    & \text{wie auch } w_3 \neq 0 \text{ ist.}
  \end{cases}$$ dargestellt, so ist $$\vi:
Aut(G) \to GL(7,p): \al(u_1, u_2, u_3, u_4, v_1, v_2, v_3, v_4, w_3, w_4) \mt $$ $$M
=
\begin{pmatrix}
  d u_1   & d u_2   & 0       & 0       & 0   & 0   & 0 \\
  d v_1   & d v_2   & 0       & 0       & 0   & 0   & 0 \\
  u_1 w_4 & u_2 w_4 & u_1 w_3 & u_2 w_3 & 0   & 0   & 0 \\
  v_1 w_4 & v_2 w_4 & v_1 w_3 & v_2 w_3 & 0   & 0   & 0 \\
  0       & 0       & 0       & 0       & u_1 & u_2 & u_3 \\
  0       & 0       & 0       & 0       & v_1 & v_2 & v_3 \\
  0       & 0       & 0       & 0       & 0   & 0   & w_3
\end{pmatrix}$$ der Operationshomomorphismus von $A=Aut(G)$ auf $\f^7 \cong M(G)$
über Erweiterungsautomorphismen, wobei $d=u_1 v_2 - v_1 u_2$ ist.
\end{lemma}

\begin{beweis}
Nach \ref{aut12} lässt sich jeder Automorphismus $\al$ von $G$ in der oben angegebenen Form
darstellen. Die Bilder der Erzeuger von $M(G)$ ergeben sich unter $\al$ über $\vi$
folgendermaßen:
\begin{eqnarray*}
  a_5^\al&=& [a_4,a_1]^\al=[a_4^\al, a_1^\al]\\
  &=& [a_4^{v_2 u_1 - v_1 u_2}, a_1^{u_1} a_2^{u_2} a_3^{u_3} a_4^{u_4}]\\
  &=& [a_4,a_1]^{(v_2 u_1 - v_1 u_2)u_1} [a_4,a_2]^{(v_2 u_1 - v_1 u_2)u_2}\\
  &=& a_5^{(v_2 u_1 - v_1 u_2)u_1} a_6^{(v_2 u_1 - v_1 u_2)u_2},
\end{eqnarray*}
da alle anderen auftretenden Kommutatoren trivial sind und $M(G)$ zentral ist.
Analog erhält man, dass auch $a_6^\al = a_5^{(v_2 u_1 - v_1 u_2)v_1} a_6^{(v_2
u_1 - v_1 u_2)v_2}$ gilt. Weiterhin ist
\begin{eqnarray*}
  a_7^\al&=& [a_3, a_1]^\al = [a_3^\al, a_1^\al]\\
  &=& [a_3^{w_3} a_4^{w_4}, a_1^{u_1} a_2^{u_2} a_3^{u_3} a_4^{u_4}]\\
  &=& [a_3, a_1]^{u_1 w_3}[a_3, a_2]^{u_2 w_3}[a_4, a_1]^{u_1 w_4}[a_4, a_2]^{u_2 w_4}\\
  &=& a_5^{u_1 w_4} a_6^{u_2 w_4} a_7^{u_1 w_3} a_8^{u_2 w_3}
\end{eqnarray*}
da alle anderen auftretenden Kommutatoren trivial sind und $M(G)$ zentral ist.
Analog erhält man, dass auch $a_8^\al = a_5^{v_1 w_4} a_6^{v_2 w_4} a_7^{v_1
w_3} a_8^{v_2 w_3}$ ist. Außerdem ist
\begin{eqnarray*}
  a_9^\al &=& (a_1^p)^\al=(a_1^\al)^p\\
  &=& (a_1^{u_1} a_2^{u_2} a_3^{u_3} a_4^{u_4})^p \\
  &=& (a_1^p)^{u_1} (a_2^p)^{u_2} (a_3^p)^{u_3} \\
  &=& a_9^{u_1} a_{10}^{u_2} a_{11}^{u_3},
\end{eqnarray*}
da $a_4$ und der Kommutator von $a_9$ und $a_{10}$ die Ordnung $p$ haben und
alle anderen relevanten Kommutatoren trivial sind. Ebenso erhält man
$a_{10}^\al=a_9^{v_1} a_{10}^{v_2} a_{11}^{v_3}$ und schließlich ist
$a_{11}^\al=(a_3^{w_3} a_4^{w_4})^p=(a_3^p)^{w_3}=a_{11}^{w_3}$, da sowohl
$a_4$ als auch der Kommutator von $a_{10}$ und $a_4$ die Ordnung $p$ haben.
\end{beweis}

Da im Weiteren für die unmittelbaren Nachfolger der Ordnung $p^5$ statt der
sechsdimensionalen Unterräume ihre eindimensionalen Komplemente betrachtet
werden, ist der zu $\vi$ duale Operationshomomorphismus von Interesse.

\begin{beme}
Der zu $\vi$ duale Operationshomomorphismus ist $$\bar{\vi}: Aut(G) \ra GL(7,p) :\al(u_1, u_2,
u_3, u_4, v_1, v_2, v_3, v_4, w_3, w_4) \mt $$ $$M
=
\begin{pmatrix}
  d u_1   & d v_1   & u_1 w_4 & v_1 w_4 & 0   & 0   & 0 \\
  d u_2   & d v_2   & u_2 w_4 & v_2 w_4 & 0   & 0   & 0 \\
  0       & 0       & u_1 w_3 & v_1 w_3 & 0   & 0   & 0 \\
  0       & 0       & u_2 w_3 & v_2 w_3 & 0   & 0   & 0 \\
  0       & 0       & 0       & 0       & u_1 & v_1 & 0 \\
  0       & 0       & 0       & 0       & u_2 & v_2 & 0 \\
  0       & 0       & 0       & 0       & u_3 & v_3 & w_3
\end{pmatrix},$$ wobei $d=u_1 v_2 - v_1 u_2$ ist.
\end{beme}

\begin{beme}
Die Gruppe $L=Aut(G)^{\bar{\vi}}$ enthält nicht das volle Zentrum $Z$ von $GL(7,p)$. Denn $Z$
wird von $w E_7$ erzeugt, wobei $w$ ein Erzeuger der multiplikativen Gruppe von $\F_p$ und
$E_7$ die Einheitsmatrix über $\f$ der Dimension $7 \times 7$ ist. Die Gruppe $L$ enthält aber
$w E_7$ nicht, da in diesem Fall $w_3=w$ und $v_2=w$ erfüllt wäre und zugleich auch $w=v_2 w_3
= w^2$ im Widerspruch zu $w \neq 1$ sein müsste (damit ist der Schnitt von $L$ mit $Z$
trivial). Nach \ref{zentin} lassen sich daher die Bahnen der eindimensionalen Unterräume unter
der Operation von $L$ auf $\F_p^7$ nicht unbedingt aus den Bahnen der Vektoren ablesen.
\end{beme}

\begin{folg}
Es sei $$A=\left\{ k
\begin{pmatrix}
  d u_1   & d v_1   & u_1 w_4 & v_1 w_4 & 0   & 0   & 0 \\
  d u_2   & d v_2   & u_2 w_4 & v_2 w_4 & 0   & 0   & 0 \\
  0       & 0       & u_1 w_3 & v_1 w_3 & 0   & 0   & 0 \\
  0       & 0       & u_2 w_3 & v_2 w_3 & 0   & 0   & 0 \\
  0       & 0       & 0       & 0       & u_1 & v_1 & 0 \\
  0       & 0       & 0       & 0       & u_2 & v_2 & 0 \\
  0       & 0       & 0       & 0       & u_3 & v_3 & w_3
\end{pmatrix} \mid  \triangle \right\}$$ mit $\triangle=u_1, u_2, u_3,
v_1, v_2, v_3, w_3, w_4 \in \f$ und $d=u_1 v_2 - u_2 v_1 \neq 0$ sowie $k \neq 0$ und $w_3
\neq 0$. Unter der Operation von $A$ auf $\F_p^7$ zerfällt $\F_p^7$ derart in Bahnen, dass
sich aus ihnen die Bahnen der eindimensionalen Unterräume ablesen lassen und dass sie den
Bahnen der eindimensionalen Unterräume im Sinne von \ref{einvek} entsprechen, die sich unter
der Operation von $Aut(G)^{\bar{\vi}}$ ergeben. Die Mächtigkeit von $A$ ist
$(p-1)|Aut(G)^{\bar{\vi}}|=p^4(p-1)^3(p^2-1)$.
\end{folg}

\subsection{Bahnen zulässiger Untergruppen für Nachfolger der Ordnung $p^5$}

\begin{lemma}
Unter der Operation von $A$ zerfällt $\f^7$ in folgende Bahnen nichttrivialer Vektoren, wobei
$w$ ein Erzeuger der multiplikativen Gruppe von $\f$ ist:
\begin{enumerate}
  \item[] $B_1   =(0,1,0,0,0,0,0)^A$ mit $|B_1|=p(p^2-1)$,
  \item[] $B_2   =(0,1,0,0,0,0,1)^A$ mit $|B_2|=p^3(p-1)(p^2-1)$,
  \item[] $B_3   =(0,1,0,0,0,1,0)^A$ mit $|B_3|=p(p-1)(p^2-1)$,
  \item[] $B_4   =(0,1,0,0,1,0,0)^A$ mit $|B_4|=\frac{1}{2}p^2(p-1)(p^2-1)$,
  \item[] $B_5   =(0,1,0,0,w,0,0)^A$ mit $|B_5|=\frac{1}{2}p^2(p-1)(p^2-1)$,
  \item[] $B_6   =(0,1,1,0,0,0,0)^A$ mit $|B_6|=p(p-1)(p^2-1)$,
  \item[] $B_7   =(0,1,1,0,0,0,1)^A$ mit $|B_7|=p^3(p-1)^2(p^2-1)$,
  \item[] $B_8   =(0,1,1,0,0,1,0)^A$ mit $|B_8|=p(p-1)^2(p^2-1)$,
  \item[] $B_9   =(0,1,1,0,1,0,0)^A$ mit $|B_9|=\frac{1}{2}p^2(p-1)^2(p^2-1)$,
  \item[] $B_{10}=(0,1,1,0,w,0,0)^A$ mit $|B_{10}|=\frac{1}{2}p^2(p-1)^2(p^2-1)$,
  \item[] $B_{11}=(0,0,0,0,0,0,1)^A$ mit $|B_{11}|=p^2(p-1)$,
  \item[] $B_{12}=(0,0,0,0,0,1,0)^A$ mit $|B_{12}|=p^2-1$,
  \item[] $B_{13}=(0,0,0,1,0,0,0)^A$ mit $|B_{13}|=p^2-1$,
  \item[] $B_{14}=(0,0,0,1,0,0,1)^A$ mit $|B_{14}|=p^2(p-1)(p^2-1)$,
  \item[] $B_{15}=(0,0,0,1,0,1,0)^A$ mit $|B_{15}|=(p-1)(p^2-1)$,
  \item[] $B_{16}=(0,0,0,1,1,0,0)^A$ mit $|B_{16}|=p(p-1)(p^2-1)$.
\end{enumerate}
Dabei ist $w$ ein Erzeuger der multiplikativen Gruppe von $\f$. Zu den Vertretern ausgewählter
Bahnen gehören folgende Stabilisatoren mit ihrer jeweils anschließend angegebenen Ordnung:
\begin{eqnarray*}
  \bar{S}_1 &=& \left\{
  \begin{pmatrix}
    u_1 v_2\1 & d v_1 v_2\1 & 0       & 0       & 0   & 0   & 0 \\
    0         & 1           & 0       & 0       & 0   & 0   & 0 \\
    0         & 0           & u_1 w_3 & v_1 w_3 & 0   & 0   & 0 \\
    0         & 0           & 0       & v_2 w_3 & 0   & 0   & 0 \\
    0         & 0           & 0       & 0       & u_1 & v_1 & 0 \\
    0         & 0           & 0       & 0       & 0   & v_2 & 0 \\
    0         & 0           & 0       & 0       & u_3 & v_3 & w_3
  \end{pmatrix} \mid \triangle \right\}  \\
    & & \triangle = u_3, v_1, v_3 \in \f \text{ und } u_1, v_2, w_3 \in \f^*.
    \text{ }|\bar{S}_1|=p^3(p-1)^3.\\
  \bar{S}_2 &=& \left\{
  \begin{pmatrix}
    u_1 v_2\1 & d v_1 v_2\1 & 0       & 0       & 0   & 0   & 0 \\
    0         & 1           & 0       & 0       & 0   & 0   & 0 \\
    0         & 0           & u_1     & v_1     & 0   & 0   & 0 \\
    0         & 0           & 0       & v_2     & 0   & 0   & 0 \\
    0         & 0           & 0       & 0       & u_1 & v_1 & 0 \\
    0         & 0           & 0       & 0       & 0   & v_2 & 0 \\
    0         & 0           & 0       & 0       & 0   & 0   & 1
  \end{pmatrix} \mid \triangle \right\}  \\
    & & \triangle = v_1 \in \f \text{ und } u_1, v_2 \in \f^*.
    \text{ }|\bar{S}_2|=p(p-1)^2.\\
  \bar{S}_3 &=& \left\{
  \begin{pmatrix}
    v_2\1   & v_1 v_2\1 & 0            & 0             & 0         & 0         & 0 \\
    0       & 1         & 0            & 0             & 0         & 0         & 0 \\
    0       & 0         & v_2^{-2} w_3 & v_1 v_2\1 w_3 & 0         & 0         & 0 \\
    0       & 0         & 0            & w_3           & 0         & 0         & 0 \\
    0       & 0         & 0            & 0             & v_2^{-2}  & v_1 v_2\1 & 0 \\
    0       & 0         & 0            & 0             & 0         & 1         & 0 \\
    0       & 0         & 0            & 0             & v_2\1 u_3 & v_2\1 v_3 & v_2\1 w_3
  \end{pmatrix} \mid  \triangle \right\}  \\
    & & \triangle = u_3, v_1, v_3 \in \f \text{ und } v_2, w_3 \in \f^*.
    \text{ }|\bar{S}_3|=p^3(p-1)^2.\\
  \bar{S}_4 = \bar{S}_5 &=& \left\{
  \begin{pmatrix}
    u_1 v_2 & 0       & 0       & 0             & 0   & 0         & 0 \\
    0       & v_2^2   & 0       & 0             & 0   & 0         & 0 \\
    0       & 0       & w_3     & 0             & 0   & 0         & 0 \\
    0       & 0       & 0       & u_1\1 v_2 w_3 & 0   & 0         & 0 \\
    0       & 0       & 0       & 0             & 1   & 0         & 0 \\
    0       & 0       & 0       & 0             & 0   & u_1\1 v_2 & 0 \\
    0       & 0       & 0       & 0             & u_3 & v_3       & w_3
  \end{pmatrix} \mid  \triangle \right\}  \\
    & & \triangle = u_3, v_3 \in \f \text{ und } u_1, v_2, w_3 \in \f^* \text{ sowie } v^2=1.
    \text{ }|\bar{S}_4|=|\bar{S}_5|=\frac{p^2(p-1)^2}{2}.\\
  \bar{S}_6 &=& \left\{
  \begin{pmatrix}
    u_1 v_2\1 & v_1 v_2\1   & -v_1 v_2\1 & -u_1\1 v_1^2 v_2\1 & 0   & 0   & 0 \\
    0         & 1           & 0          & -u_1\1 v_1         & 0   & 0   & 0 \\
    0         & 0           & 1          & u_1\1 v_1          & 0   & 0   & 0 \\
    0         & 0           & 0          & u_1\1 v_2          & 0   & 0   & 0 \\
    0         & 0           & 0          & 0                  & u_1 & v_1 & 0 \\
    0         & 0           & 0          & 0                  & 0   & v_2 & 0 \\
    0         & 0           & 0          & 0                  & u_3 & v_3 & v_2^2
  \end{pmatrix} \mid  \triangle \right\}  \\
    & & \triangle = u_3, v_1, v_3 \in \f \text{ und } u_1, v_2 \in \f^*.
    \text{ }|\bar{S}_6|=p^3(p-1)^2.\\
  \bar{S}_7 &=& \left\{
  \begin{pmatrix}
    v_2\1   & -w_4 v_2^{-2} & v_2^{-2} w_4 & -v_2^{-3} w_4^{-2} & 0     & 0             & 0 \\
    0       & 1             & 0            & v_2\1 w_4          & 0     & 0             & 0 \\
    0       & 0             & 1            & -v_2\1 w_4         & 0     & 0             & 0 \\
    0       & 0             & 0            & v_2                & 0     & 0             & 0 \\
    0       & 0             & 0            & 0                  & v_2\1 & -v_2^{-3} w_4 & 0 \\
    0       & 0             & 0            & 0                  & 0     & v_2\1         & 0 \\
    0       & 0             & 0            & 0                  & 0     & 0             & 1
  \end{pmatrix} \mid  \triangle \right\}  \\
    & & \triangle = w_4 \in \f \text{ und } v_2 \in \f^*.
    \text{ }|\bar{S}_7|=p(p-1).\\
  \bar{S}_8 &=& \left\{
  \begin{pmatrix}
    v_2^{-3} & -w_4    & v_2^{-2} w_4 & -w_4    & 0         & 0         & 0 \\
    0        & 1       & 0            & w_4     & 0         & 0         & 0 \\
    0        & 0       & v_2^{-2}     & -w_4    & 0         & 0         & 0 \\
    0        & 0       & 0            & 1       & 0         & 0         & 0 \\
    0        & 0       & 0            & 0       & v_2^{-2}  & -w_4      & 0 \\
    0        & 0       & 0            & 0       & 0         & 1         & 0 \\
    0        & 0       & 0            & 0       & u_3 v_2\1 & v_3 v_2\1 & 1
  \end{pmatrix} \mid  \triangle \right\}  \\
    & & \triangle = u_3, v_3, w_4 \in \f \text{ und } v_2 \in \f^*.
    \text{ }|\bar{S}_8|=p^3(p-1).\\
  \bar{S}_9=\bar{S}_{10} &=& \left\{
  \begin{pmatrix}
    v_2     & 0       & 0       & 0       & 0   & 0       & 0 \\
    0       & v_2^2   & 0       & v_2 w_4 & 0   & 0       & 0 \\
    0       & 0       & 1       & 0       & 0   & 0       & 0 \\
    0       & 0       & 0       & u_1 v_2 & 0   & 0       & 0 \\
    0       & 0       & 0       & 0       & 1   & 0       & 0 \\
    0       & 0       & 0       & 0       & 0   & u_1 v_2 & 0 \\
    0       & 0       & 0       & 0       & u_3 & v_3     & 1
  \end{pmatrix} \mid  \triangle \right\}  \\
    & & \triangle = u_3, v_3 \in \f \text{ und } u_1, v_2 \in \f^* \text{ sowie } v^2=1.
    \text{ }|\bar{S}_9|=|\bar{S}_{10}|=p^3(p-1)^3.\\
  \bar{S}_{15} &=& \left\{
  \begin{pmatrix}
    u_1^2   & u_1 v_1 & u_1 v_2\1 w_4 & v_1 v_2\1 w_4 & 0         & 0         & 0 \\
    0       & u_1 v_2 & 0             & w_4           & 0         & 0         & 0 \\
    0       & 0       & u_1 v_2\1     & v_1 v_2\1     & 0         & 0         & 0 \\
    0       & 0       & 0             & 1             & 0         & 0         & 0 \\
    0       & 0       & 0             & 0             & u_1 v_2\1 & v_1 v_1\1 & 0 \\
    0       & 0       & 0             & 0             & 0         & 1         & 0 \\
    0       & 0       & 0             & 0             & u_3 v_2\1 & v_3 v_2\1 & v_2\1
  \end{pmatrix} \mid  \triangle \right\}  \\
    & & \triangle = u_3, v_3, w_4 \in \f \text{ und } u_1, v_2 \in \f^*.
    \text{ }|\bar{S}_{15}|=p^3(p-1)^2.\\
  \bar{S}_{16} &=& \left\{
  \begin{pmatrix}
    u_1 v_2 & 0       & w_4       & 0            & 0         & 0         & 0 \\
    0       & v_2^2   & 0         & u_1\1v_2 w_4 & 0         & 0         & 0 \\
    0       & 0       & u_1 v_2\1 & 0            & 0         & 0         & 0 \\
    0       & 0       & 0         & 1            & 0         & 0         & 0 \\
    0       & 0       & 0         & 0            & 1         & 0         & 0 \\
    0       & 0       & 0         & 0            & 0         & u_1\1 v_2 & 0 \\
    0       & 0       & 0         & 0            & u_1\1 u_3 & u_1\1 v_3 & v_2\1
  \end{pmatrix} \mid  \triangle \right\}  \\
    & & \triangle = u_3, v_3, w_4 \in \f \text{ und } u_1, v_2 \in \f^*.
    \text{ }|\bar{S}_{16}|=p(p-1)^2.\\
\end{eqnarray*}
\end{lemma}

\begin{beweis}
Die Bahnen $B_{11}$ bis $B_{14}$ lassen sich unmittelbar erkennen. Für die übrigen Fälle sind
die Stabilisatoren angegeben, die man aus einfachen Gleichungssystemen berechnen kann und
deren Ordnung sich durch Abzählen ermitteln lässt. Über die Bahnensätze erhält man dann die
Länge der Bahnen. Die Vollständigkeit der Liste erkennt man dadurch, dass die Summation der
Bähnenlängen die Mächtigkeit von $\f^7$ ergibt. Die Bahnen gleicher Länge sind offensichtlich
verschieden, da sie Teilmengen unterschiedlicher irreduzibler Teilmoduln sind.
\end{beweis}

\begin{folg} \label{zu9}
Für die Nachfolger von $G$ der Ordnung $p^5$ ist die folgende Liste ein Vertretersystem
zulässiger Untergruppen.
\begin{enumerate}
  \item[] $M_1   =\erz{a_5, a_7, a_8, a_9, a_{10}, a_{11}}$ entsprechend $\erz{(0,1,0,0,0,0,0)}^\bot$,
  \item[] $M_2   =\erz{a_5, a_7, a_8, a_9, a_{10}, a_6 a_{11}\1}$ entsprechend $\erz{(0,1,0,0,0,0,1)}^\bot$,
  \item[] $M_3   =\erz{a_5, a_7, a_8, a_9, a_6 a_{10}\1, a_{11}}$ entsprechend $\erz{(0,1,0,0,0,1,0)}^\bot$,
  \item[] $M_4   =\erz{a_5, a_7, a_8, a_6 a_9\1, a_{10}, a_{11}}$ entsprechend $\erz{(0,1,0,0,1,0,0)}^\bot$,
  \item[] $M_5   =\erz{a_5, a_7, a_8, a_6 a_9^{-w}, a_{10}, a_{11}}$ entsprechend $\erz{(0,1,0,0,w,0,0)}^\bot$,
  \item[] $M_6   =\erz{a_5, a_6 a_7\1, a_8, a_9, a_{10}, a_{11}}$ entsprechend $\erz{(0,1,1,0,0,0,0)}^\bot$,
  \item[] $M_7   =\erz{a_5, a_6 a_7\1, a_8, a_9, a_{10}, a_6 a_{11}\1}$ entsprechend $\erz{(0,1,1,0,0,0,1)}^\bot$,
  \item[] $M_8   =\erz{a_5, a_6 a_7\1, a_8, a_9, a_6 a_{10}\1, a_{11}}$ entsprechend $\erz{(0,1,1,0,0,1,0)}^\bot$,
  \item[] $M_9   =\erz{a_5, a_6 a_7\1, a_8, a_6 a_9\1, a_{10}, a_{11}}$ entsprechend $\erz{(0,1,1,0,1,0,0)}^\bot$,
  \item[] $M_{10}=\erz{a_5, a_6 a_7\1, a_8, a_6 a_9^{-w}, a_{10}, a_{11}}$ entsprechend $\erz{(0,1,1,0,w,0,0)}^\bot$.
\end{enumerate}
\end{folg}

\begin{beweis}
Nach \ref{suppl} sind nur Supplemente des Nukleus $N=\erz{a_5, a_6}$, entsprechend dem
Unterraum $\bar{N}=\erz{(1,0,0,0,0,0,0), (0,1,0,0,0,0,0)}$, zulässige Untergruppen und
wiederum sind genau sechsdimensionale Unterräume Supplemente zu $\bar{N}$, deren
eindimensionale Komplemente einen Basisvektor haben, der mit den Basisvektoren von $\bar{N}$
ein Skalarprodukt ungleich Null liefert. Damit schränkt sich die Auswahl zulässiger
Untergruppen auf die Bahnen $B_1$ bis $B_{10}$ ein.
\end{beweis}

\subsection{Nachfolger der Ordnung $p^5$}

\begin{satz} \label{nt12}
In der folgenden Liste sind sämtliche Nachfolger von $G$ der Ordnung $p^5$ bis auf Isomorphie
angegeben:
\begin{enumerate}
  \item[] $G_1=\erz{g_1, g_2, g_3, g_4, g_5 \mid [g_2,g_1]=g_4, [g_4,g_2]=g_5}$,
  \item[] $G_2=\erz{g_1, g_2, g_3, g_4, g_5 \mid [g_2,g_1]=g_4, [g_4,g_2]=g_5, g_3^p=g_{5}}$,
  \item[] $G_3=\erz{g_1, g_2, g_3, g_4, g_5 \mid [g_2,g_1]=g_4, [g_4,g_2]=g_5, g_2^p=g_{5}}$,
  \item[] $G_4=\erz{g_1, g_2, g_3, g_4, g_5 \mid [g_2,g_1]=g_4, [g_4,g_2]=g_5, g_1^p=g_5}$,
  \item[] $G_5=\erz{g_1, g_2, g_3, g_4, g_5 \mid [g_2,g_1]=g_4, [g_4,g_2]=g_5^w, g_1^p=g_5}$,
  \item[] $G_6=\erz{g_1, g_2, g_3, g_4, g_5 \mid [g_2,g_1]=g_4, [g_4,g_2]=g_5, [g_3,g_1]=g_5}$,
  \item[] $G_7=\erz{g_1, g_2, g_3, g_4, g_5 \mid [g_2,g_1]=g_4, [g_4,g_2]=g_5, [g_3,g_1]=g_5, g_3^p=g_{5}}$,
  \item[] $G_8=\erz{g_1, g_2, g_3, g_4, g_5 \mid [g_2,g_1]=g_4, [g_4,g_2]=g_5, [g_3,g_1]=g_5, g_2^p=g_{5}}$,
  \item[] $G_9=\erz{g_1, g_2, g_3, g_4, g_5 \mid [g_2,g_1]=g_4, [g_4,g_2]=g_5, [g_3,g_1]=g_5, g_1^p=g_5}$,
  \item[] $G_{10}=\erz{g_1, g_2, g_3, g_4, g_5 \mid [g_2,g_1]=g_4, [g_4,g_2]=g_5^w, [g_3,g_1]=g_5^w, g_1^p=g_5}$,
\end{enumerate}
\end{satz}

\begin{beweis}
Dieser Satz ergibt sich nach \ref{bahnenzu}, indem man $P(G)$ nach jedem Repräsentanten des
Vertretersystems zulässiger Untergruppen faktorisiert, das in \ref{zu9} aufgelistet ist. Die
Präsentation der Faktorgruppe erhält man über das Verfahren aus \ref{faktor}.
\end{beweis}

\section{Nachfolger von $(p^4,13)$}

\begin{ver}
In diesem Abschnitt sei $$G=\erz{a_1, a_2, a_3, a_4 \mid [a_2,a_1]=a_4, a_1^p=a_4}.$$
\end{ver}

\begin{lemma}
Die Gruppe $G$ hat die Gewichtung $\omega(a_1)=\omega(a_2)=\omega(a_3)=1$ und $\omega(a_4)=2$
und damit die $p$"=Klasse 2.
\end{lemma}

\begin{beweis}
Die Gewichtung von $G$ lässt sich nach \ref{gew} ermitteln.
\end{beweis}

\subsection{Die Automorphismengruppe}

\begin{lemma} \label{aut13}
Jeder Automorphismus $\al$ von $G_3$ lässt sich durch $$\al(u_1, u_2, u_3, u_4, v_2, v_3, v_4,
w_2, w_3, w_4 ):
  \begin{cases}
    a_1 \mt a_1^{u_1} a_2^{u_2} a_3^{u_3} a_4^{u_4}  & \text{wobei } u_1, u_2, u_3, u_4, v_2, v_3, v_4\\
    a_2 \mt           a_2^{v_2} a_3^{v_3} a_4^{v_4}  & w_3, w_4 \in \{0 \dd p-1\}\\
    a_3 \mt                     a_3^{w_3} a_4^{w_4}  & \text{und } u_1 \neq 0 \\
    a_4 \mt                               a_4^{u_1 v_2}    & \text{wie auch } v_2 \neq 0 \text{ ist.}
  \end{cases}$$ darstellen.
\end{lemma}

\begin{beweis}
Die Automorphismengruppe von $G$ kann nach \ref{aut} aus dem Stabilisator $S_3$ (siehe
\ref{bahnen_cp3_p4}) der zu $G$ gehörenden zulässigen Untergruppe unter Berücksichtigung des
Transponierens unmittelbar abgelesen werden. Damit können auch die Bilder von $a_1$, $a_2$ und
$a_3$ direkt abgelesen werden. Das Bild von $a_4$ ist durch die Bilder von $a_1$, $a_2$ und
$a_3$ eindeutig bestimmt und lässt sich folgendermaßen berechnen:
\begin{eqnarray*}
  a_4^\al &=& [a_2, a_1]^\al = [a_2^\al, a_1^\al] \\
  &=& [a_2^{v_2} a_3^{v_3} a_4^{v_4}, a_1^{u_1} a_2^{u_2} a_3^{u_3} a_4^{u_4}] \\
  &=& [a_2^{v_2} a_3^{v_3}, a_1^{u_1} a_2^{u_2} a_3^{u_3}] \\
  &=& [a_2^{v_2}, a_1^{u_1} a_2^{u_2} a_3^{u_3}]^{a_3^{v_3}} \\
  &=& ([a_2,a_1]^{u_1 v_2})^{a_3^{v_3}}\\
  &=& [a_2,a_1]^{u_1 v_2} \\
  &=& a_4^{u_1 v_2},
\end{eqnarray*}
da alle anderen auftretenden Kommutatoren trivial sind und $[a_2, a_1]$ zentral ist.
\end{beweis}

\subsection{$p$"=Cover, Multiplikator und Nukleus}

\begin{lemma}
Für $G$ ergibt sich das $p$"=Cover, der Nukleus und der Multiplikator in
folgender Weise mit der Gewichtung $\omega$:
\begin{enumerate}
  \item[] $P(G)=\erz{a_1 \dd a_9 \mid [a_2, a_1] = a_4,
  [a_3, a_1]=a_5, [a_3, a_2]=a_6, a_1^p=a_4 a_7, a_2^p=a_8, a_3^p=a_9}$.
  \item[] $M(G)=\erz{a_5 \dd a_9}$.
  \item[] $N(G)=\erz{1}$.
\end{enumerate}
Die Gruppe $P(G)$ hat die $p$"=Klasse 2.
\end{lemma}

\begin{beweis}
Die Präsentation von $P(G)$ ergibt sich aus \ref{cov} und \ref{newman} unter Verwendung des
Reduktionsverfahrens von Knuth"=Bendix, aus der sich $M(G)$ ablesen lässt. Die Gewichtung von
$P(G)$
\begin{enumerate}
  \item[] $1=\omega(a_1)=\omega(a_2)=\omega(a_3)$.
  \item[] $2=\omega(a_3)= \ldots =\omega(a_9)$.
\end{enumerate}
und damit $N(G)$ lassen sich nach \ref{gew} aus der Präsentation von $P(G)$ ablesen.
\end{beweis}

\begin{satz}
Die Gruppe $G$ ist abschließend.
\end{satz}

\begin{beweis}
Da $P(G)$ dieselbe $p$"=Klasse hat wie $G$, hat $G$ nach \ref{pfort} keine unmittelbaren
Nachfolger, d.\,h. $G$ ist abschließend.
\end{beweis}

\section{Nachfolger von $(p^4,14)$}

\begin{ver}
In diesem Abschnitt sei $$G=\erz{a_1, a_2, a_3, a_4 \mid [a_2,a_1]=a_3, a_3^p=a_4}.$$
\end{ver}

\begin{lemma}
Die Gruppe $G$ hat die Gewichtung $\omega(a_1)=\omega(a_2)=1$, $\omega(a_3)=2$ und
$\omega(a_4)=3$ und damit die $p$"=Klasse 3.
\end{lemma}

\begin{beweis}
Die Gewichtung von $G$ lässt sich nach \ref{gew} ermitteln.
\end{beweis}

\subsection{Die Automorphismengruppe}

\begin{lemma} \label{aut14}
Jeder Automorphismus $\al$ von $G_4$ lässt sich durch $$\al(u_1, u_2, u_3, u_4, v_2, v_3, v_4,
w_2, w_3, w_4 ):
  \begin{cases}
    a_1 \mt a_1^{u_1} a_2^{u_2}                         a_4^{u_4}  & \text{wobei } u_1, u_2, u_4, v_2, v_4\\
    a_2 \mt a_1^{v_1} a_2^{v_2}                         a_4^{v_4}  & \in \{0 \dd p-1\} \text{ und }\\
    a_3 \mt                     a_3^{u_1 v_2 - v_1 u_2} a_4^{w_4}  & u_1 v_2 - v_1 u_2 \neq 0 \text{ ist.}\\
    a_4 \mt                               a_4^{u_1 v_2 - v_2 u_1}  &
  \end{cases}$$ darstellen.
\end{lemma}

\begin{beweis}
Die Automorphismengruppe von $G$ kann nach \ref{aut} aus dem Stabilisator $S_1$ (siehe
\ref{bahnen_cp3_p4}) der zu $G$ gehörenden zulässigen Untergruppe unter Berücksichtigung des
Transponierens unmittelbar abgelesen werden. Damit können auch die Bilder von $a_1$, $a_2$ und
$a_3$ direkt abgelesen werden. Das Bild von $a_4$ ist durch die Bilder von $a_1$, $a_2$ und
$a_3$ eindeutig bestimmt und lässt sich folgendermaßen berechnen:
\begin{eqnarray*}
  a_4^\al &=& [a_2, a_1]^\al=[a_2^\al, a_1^\al]\\
  &=& [a_1^{v_1} a_2^{v_2}, a_1^{u_1} a_2^{u_2}] \\
  &=& [a_2, a_1]^{u_1 v_2 - v_2 u_1} \\
  &=& a_4^{u_1 v_2 - v_2 u_1},
\end{eqnarray*}
da alle anderen auftretenden Kommutatoren trivial sind und $[a_2, a_1]$ zentral ist.
\end{beweis}

\subsection{$p$"=Cover, Multiplikator und Nukleus}

\begin{lemma}
Für $G$ ergibt sich das $p$"=Cover, der Nukleus und der Multiplikator in
folgender Weise mit der Gewichtung $\omega$:
\begin{enumerate}
  \item[] $P(G)=\erz{a_1 \dd a_9 \mid [a_2, a_1] = a_4
  [a_3, a_1]=a_5, [a_3, a_2]=a_6, a_1^p= a_7, a_2^p=a_8, a_3^p=a_4 a_9}$.
  \item[] $M(G)=\erz{a_5 \dd a_9}$.
  \item[] $N(G)=\erz{1}$.
\end{enumerate}
Die Gruppe $P(G)$ hat die $p$"=Klasse 2.
\end{lemma}

\begin{beweis}
Die Präsentation von $P(G)$ ergibt sich aus \ref{cov} und \ref{newman} unter Verwendung des
Reduktionsverfahrens von Knuth"=Bendix, aus der sich $M(G)$ ablesen lässt. Die Gewichtung von
$P(G)$
\begin{enumerate}
  \item[] $1=\omega(a_1)=\omega(a_2)=\omega(a_3)$.
  \item[] $2=\omega(a_3)= \ldots =\omega(a_9)$.
\end{enumerate}
und damit $N(G)$ lässt sich nach \ref{gew} aus der Präsentation von $P(G)$ ablesen.
\end{beweis}

\begin{satz}
Die Gruppe $G$ ist abschließend.
\end{satz}

\begin{beweis}
Da $P(G)$ dieselbe $p$"=Klasse hat wie $G$, hat $G$ nach \ref{pfort} keine unmittelbaren
Nachfolger, d.\,h. $G$ ist abschließend.
\end{beweis}

\section{Nachfolger von $C_p^4$}

\begin{ver}
In diesem Abschnitt sei $$G =\erz{a_1, a_2, a_3, a_4 \mid [a_2,a_1], [a_3,a_1], [a_4, a_1],
[a_3,a_2], [a_4,a_2], [a_4,a_3], a_1^p, a_2^p, a_3^p, a_4^p}.$$ Damit ist $G$ isomorph zu
$C_{p}^4$.
\end{ver}

\subsection{Die Automorphismengruppe}

\begin{beme}
Die Gruppe $G$ hat die $p$"=Klasse 1 und ist zur additiven Gruppe von $\F_p^4$ isomorph. Die
Automorphismengruppe $Aut(G)$ kann daher mit $GL(4,p)=:A$ identifiziert werden.
\end{beme}

\subsection{$p$"=Cover, Multiplikator und Nukleus}

\begin{lemma}
Für $G$ ergibt sich das $p$"=Cover, der Nukleus und der Multiplikator in
folgender Weise mit der Gewichtung $\omega$:
\begin{enumerate}
  \item[] $P(G)=\erz{a_1 \dd a_{14} \mid [a_2,a_1]=a_{5}, [a_3,a_1]=a_{6}, [a_4, a_1]=a_{7},
  [a_3,a_2]=a_{8}, [a_4,a_2]=a_{9}, [a_4,a_3]=a_{10}, a_1^p=a_{11}, a_2^p=a_{12},
  a_3^p=a_{13}, a_4^p=a_{14}}$.
  \item[] $M(G)=\erz{a_5, \dd a_{14}}$.
  \item[] $N(G)=\erz{a_5, \dd a_{14}}$.
\end{enumerate}
Die Gruppe $P(G)$ hat die $p$"=Klasse 2, $G$ ist fortsetzbar und jede
Untergruppe von $M(G)$ ist zulässig.
\end{lemma}

\begin{beweis}
Die Präsentation von $P(G)$ ergibt sich aus \ref{cov} und \ref{newman} unter Verwendung des
Reduktionsverfahrens von Knuth"=Bendix, aus der sich $M(G)$ ablesen lässt. Die Gewichtung von
$P(G)$
\begin{enumerate}
  \item[] $1=\omega(a_1)=\omega(a_2)=\omega(a_3)=\omega(a_4)$.
  \item[] $2=\omega(a_5)= \ldots =\omega(a_{14})$.
\end{enumerate}
und damit $N(G)$ lassen sich nach \ref{gew} aus der Präsentation von $P(G)$ ablesen.
\end{beweis}

\subsection{Operation der Erweiterungsautomorphismen}

\begin{lemma}
Unter der Identifikation von $Aut(G)$ mit $GL(4,p)$ ist die Operation von
$A=GL(4,p)$ auf $\f^{10} \cong M(G)$ durch den Operationshomomorphismus $\vi$
gegeben: $$\vi: A \to Aut(M(G)): m=\begin{pmatrix}
  b_{1} & b_{2} & b_{3} & b_4 \\
  c_{1} & c_{2} & c_{3} & c_4 \\
  d_{1} & d_{2} & d_{3} & d_4 \\
  e_{1} & e_{2} & e_{3} & e_4
\end{pmatrix} \mt M = \begin{pmatrix}
  L   & 0 \\
  0   & m
\end{pmatrix},$$ wobei $$L=\begin{pmatrix}
  b_1 c_2 - c_1 b_2 & b_1 c_3 - c_1 b_3 & b_1 c_4 - c_1 b_4 & b_2 c_3 - c_2 b_3 & b_2 c_4 - c_2 b_4 & b_3 c_4 - c_3 b_4 \\
  b_1 d_2 - d_1 b_2 & b_1 d_3 - d_1 b_3 & b_1 d_4 - d_1 b_4 & b_2 d_3 - d_2 b_3 & b_2 d_4 - d_2 b_4 & b_3 d_4 - d_3 b_4 \\
  b_1 e_2 - e_1 b_2 & b_1 e_3 - e_1 b_3 & b_1 e_4 - e_1 b_4 & b_2 e_3 - e_2 b_3 & b_2 e_4 - e_2 b_4 & b_3 e_4 - e_3 b_4 \\
  c_1 d_2 - d_1 c_2 & c_1 d_3 - d_1 c_3 & c_1 d_4 - d_1 c_4 & c_2 d_3 - d_2 c_3 & c_2 d_4 - d_2 c_4 & c_3 d_4 - d_3 c_4 \\
  c_1 e_2 - e_1 c_2 & c_1 e_3 - e_1 c_3 & c_1 e_4 - e_1 c_4 & c_2 e_3 - e_2 c_3 & c_2 d_4 - e_2 c_4 & c_3 e_4 - e_3 c_4 \\
  d_1 e_2 - e_1 d_2 & d_1 e_3 - e_1 d_3 & d_1 e_4 - e_1 d_4 & d_2 e_3 - e_2 d_3 & d_2 e_4 - e_2 d_4 & d_3 e_4 - e_3 d_4
\end{pmatrix}$$
ist.
\end{lemma}

\begin{beweis}
Die Operation erhält man analog zu \ref{autcp3}.
\end{beweis}

\begin{beme}
Für die unmittelbaren Nachfolger der Ordnung $p^5$ werden statt der neundimensionalen
Unterräume von $M(G)$ ihre eindimensionalen Komplemente betrachtet und damit statt $\vi$ der
zu $\vi$ duale Operationshomomorphismus $\bar{\vi}$.
\end{beme}

\subsection{Bahnen zulässiger Untergruppen für Nachfolger der Ordnung $p^5$}

\begin{lemma}
Der Vektorraum $\f^{10} \cong M(G)$ zerfällt unter der Operation von $A=Aut(G)$ über
$\bar{\vi}$ in folgende Bahnen nichttrivialer Vektoren:
\begin{enumerate}
  \item[] $B_0=(0,0,0,0,0,0,0,0,0,0)^A$,
  \item[] $B_1=(1,0,0,0,0,0,0,0,0,0)^A$,
  \item[] $B_2=(1,0,0,0,0,1,0,0,0,0)^A$,
  \item[] $B_3=(0,0,0,0,0,0,1,0,0,0)^A$,
  \item[] $B_4=(1,0,0,0,0,0,1,0,0,0)^A$ mit $|B_4|=(p^2-1)(p^3-1)(p^2+1)$,
  \item[] $B_5=(1,0,0,0,0,0,0,0,1,0)^A$ mit $|B_5|=p^2(p^2-1)(p^3-1)(p^2+1)$,
  \item[] $B_6=(1,0,0,0,0,1,0,0,0,1)^A$ mit $|B_6|=p^2(p^2+1)(p^3-1)(p^2-1)(p-1)$.
\end{enumerate}
Zu den Vertretern der Bahnen $B_4$, $B_5$ und $B_6$ gehören die folgenden Stabilisatoren in
$Aut(G)$ mit ihrer jeweils angegebenen Ordnung:
\begin{eqnarray*}
  {S}_4 &=& \left\{
  \begin{pmatrix}
    1     & b_{2} & b_{3} & b_4 \\
    0     & 1     & c_{3} & c_4 \\
    0     & 0     & d_{3} & d_4 \\
    0     & 0     & e_{3} & e_4
  \end{pmatrix} \mid b_2, b_3, b_4, c_3, c_4 , d_3, d_4, e_3, e_4 \in \f, d_3 e_4 - e_3 d_4 \neq 0 \right\}  \\
    & & \text{mit }|{S}_4|=p^6(p^2-1)(p-1).\\
  {S}_5 &=& \left\{
  \begin{pmatrix}
    b_{1} & b_{2} & 0     & b_4 \\
    c_{1} & c_{2} & 0     & c_4 \\
    0     & 0     & 1     & d_4 \\
    0     & 0     & 0     & e_4
  \end{pmatrix} \mid b_1, b_2, b_4, c_1, c_2, c_4 , d_4, e_4 \in \f, b_1 c_2 - c_1 b_2 = 1, e_4 \neq 0 \right\}  \\
    & & \text{mit }|{S}_5|=p^4(p^2-1)(p-1).\\
    {S}_6 &=& \left\{
  \begin{pmatrix}
    b_{1} & b_{2} & b_1 e_2 - e_1 b_2     & 0 \\
    c_{1} & c_{2} & e_1 c_2 - c_1 e_2     & 0 \\
    0     & 0     & 1                     & 0 \\
    e_1   & e_2   & e_3                   & 1
  \end{pmatrix} \mid b_1, b_2, c_1, c_2, e_1, e_2, e_3 \in \f, b_1 c_2 - c_1 b_2 = 1 \right\}  \\
    & & \text{mit }|{S}_6|=p^4(p^2-1).\\
\end{eqnarray*}
\end{lemma}

\begin{beweis}
Es sei $\{v_1 \dd v_{10}\}$ die Standardbasis $V=\f^{10}$. Weiterhin seien $N=\{v_1 \dd v_6\}$
und $K=\{v_7 \dd v_{10}\}$. Damit ist $V$ die direkte Summe von $N$ und $K$ und $N$ und $K$
sind inavariante Unterräume unter der Operation von $A=Aut(G)$ über $\bar{\vi}$.

Zunächst wird gezeigt, daß $N$ unter der Operation von $A$ über $\bar{\vi}$ in die Bahnen
$B_0$, $B_1$ und $B_2$ zerfällt. Es sei $X=\{x \in M(4 \times 4,p) \mid -x=x^T\}$ die Menge
der schiefsymmetrischen $(4 \times 4)$"=Matrizen über $\f$. Man betrachte die die Abbildung
$$\theta: \f^6 \ra X: (x_1, x_2, x_3, x_4, x_5, x_6) \mt \begin{pmatrix}
  0      & x_1    & x_2   & x_3 \\
  -x_1   & 0      & x_4   & x_5 \\
  -x_2   & -x_4   & 0     & x_6 \\
  -x_3   & -x_5   & -x_6  & 0
\end{pmatrix}$$
Die Abbildung $\theta$ ist eine Bijektion zwischen $\f^6 \cong V$ und $X$. Durch
Ausmultiplikation stellt man fest, dass $A \cong GL(4,p)$ auf $L$ über $\bar{\vi}$ in
derselben Weise operiert wie $GL(4,p)$ auf $X$ vermittels $$\psi: (GL(4,p) \times X) \ra X:
(g,x) \mt g x g^T,$$ wobei $T$ die Transpositionsabbildung ist. Damit entspricht das Problem,
die Bahnen von $A$ in $L$ zu bestimmen, der Aufgabenstellung, die schiefsymmetrischen Matrizen
unter der Operation von $GL(4,p)$ über $\psi$ zu klassifizieren. Die Lösung dieser
Fragestellung ist beispielsweise in \cite{bri85}, S. 405 bis 414, dargestellt. Als Ergebnis
lässt sich festhalten, dass $X$ unter der Operation von $GL(4,p)$ über $\psi$ in drei Bahnen
zerfällt. Als Repräsentanten dieser drei Bahnen kann man die folgenden Matrizen angeben:
\begin{enumerate}
  \item $\begin{pmatrix}
    0      & 0      & 0     & 0 \\
    0      & 0      & 0     & 0 \\
    0      & 0      & 0     & 0 \\
    0      & 0      & 0     & 0
  \end{pmatrix}$
  \item $\begin{pmatrix}
    0      & 1      & 0     & 0 \\
    -1     & 0      & 0     & 0 \\
    0      & 0      & 0     & 0 \\
    0      & 0      & 0     & 0
  \end{pmatrix}$
  \item $\begin{pmatrix}
    0      & 1      & 0     & 0 \\
    -1     & 0      & 0     & 0 \\
    0      & 0      & 0     & 1 \\
    0      & 0      & -1    & 0
  \end{pmatrix}$
\end{enumerate}

In der Bahn der ersten Matrix liegt nur die Nullmatrix, in der Bahn der zweiten Matrix alle
ausgearteten Matrizen von $X$ ungleich der Nullmatrix und in der Bahn der dritten Matrix alle
nichtausgearteten Matrizen von $X$. Wendet man $\theta\1$ auf diese drei Matrizen an, dann
erhält man die Vertreter der Bahnen $B_0$, $B_1$ und $B_2$. Die Summe der Längen von $B_0$,
$B_1$ und $B_2$ ergibt die Mächtigkeit von $N$, also $p^6$.

Die Gruppe $GL(4,p)$ operiert transitiv auf $\f^4 \backslash\{0\}$ und $K$ ist unter der
Operation von $GL(4,p)$ über $\bar{\vi}$ invariant. Damit erhält man $B_3$ als weitere Bahn
und $|B_3|=p^4-1$.

Zu den übrigen drei Bahnen sind die Stabilisatoren angeben, die man über die Lösung einfacher
Gleichungssysteme berechnen kann. Mit den Bahn"=Stabilisator"=Sätzen und der Mächtigkeit von
$GL(4,p)$ nach \ref{glslz} erhält man die Längen dieser Bahnen. Alle drei Bahnenlängen sind
verschieden. Damit ist der Nachweis erbracht, dass $B_4$, $B_5$ und $B_6$ verschieden und
offensichtlich auch mit keiner der Bahnen $B_0$ bis $B_3$ identisch sind. Die Summe der Längen
aller sieben Bahnen liefert $p^{10}$. Damit ist gezeigt, dass $B_0$ bis $B_6$ sämtliche Bahnen
in $V$ unter der Operation von $A$ über $\bar{\vi}$ sind.
\end{beweis}

\begin{folg} \label{zu10}
Für die Nachfolger von $G$ der Ordnung $p^5$ ist die folgende Liste ein Vertretersystem
zulässiger Untergruppen.
\begin{enumerate}
  \item[] $M_1=\erz{a_{6},a_{7},a_{8},a_{9},a_{10},a_{11},a_{12},a_{13},a_{14}}$ entsprechend $\erz{(1,0,0,0,0,0,0,0,0,0)}^\bot$,
  \item[] $M_2=\erz{a_{5} a_{10}\1,a_{6},a_{7},a_{8},a_{9},a_{11},a_{12},a_{13},a_{14}}$ entsprechend $\erz{(1,0,0,0,0,1,0,0,0,0)}^\bot$,
  \item[] $M_3=\erz{a_{5},a_{6},a_{7},a_{8},a_{9},a_{10},a_{12},a_{13},a_{14}}$ entsprechend $\erz{(0,0,0,0,0,0,1,0,0,0)}^\bot$,
  \item[] $M_4=\erz{a_{5} a_{11}\1,a_{6},a_{7},a_{8},a_{9},a_{10},a_{12},a_{13},a_{14}}$ entsprechend $\erz{(1,0,0,0,0,0,1,0,0,0)}^\bot$,
  \item[] $M_5=\erz{a_{5} a_{13}\1,a_{6},a_{7},a_{8},a_{9},a_{10},a_{11},a_{12},a_{14}}$ entsprechend $\erz{(1,0,0,0,0,0,0,0,1,0)}^\bot$,
  \item[] $M_6=\erz{a_{5} a_{10}\1, a_5 a_{14}\1, a_{6},a_{7},a_{8},a_{9},a_{11},a_{12}, a_{13}}$ entsprechend $\erz{(1,0,0,0,0,1,0,0,0,1)}^\bot$.
\end{enumerate}
\end{folg}

\begin{beweis}
Nach \ref{suppl} ist jede Untergruppe von $M(G)$ zulässig. Daher liefert jede Bahn einen
Vertreter der zulässigen Untergruppen der Ordnung $p^9$.
\end{beweis}

\subsection{Nachfolger der Ordnung $p^5$}

\begin{satz} \label{ncp4p5}
In der folgenden Liste sind sämtliche Nachfolger von $G$ der Ordnung $p^5$ bis auf Isomorphie
angegeben:
\begin{enumerate}
  \item[] $G_1=\erz{g_1, g_2, g_3, g_4, g_5 \mid [g_2,g_1]=g_{5}}$,
  \item[] $G_2=\erz{g_1, g_2, g_3, g_4, g_5 \mid [g_2,g_1]=g_{5}, [g_4,g_3]=g_{5}}$,
  \item[] $G_3=\erz{g_1, g_2, g_3, g_4, g_5 \mid g_1^p=g_{5}}$,
  \item[] $G_4=\erz{g_1, g_2, g_3, g_4, g_5 \mid [g_2,g_1]=g_{5}, g_1^p=g_{5}}$,
  \item[] $G_5=\erz{g_1, g_2, g_3, g_4, g_5 \mid [g_2,g_1]=g_{5}, g_3^p=g_{5}}$,
  \item[] $G_6=\erz{g_1, g_2, g_3, g_4, g_5 \mid [g_2,g_1]=g_{5}, [g_4,g_3]=g_{5}, g_4^p=g_{5}}$.
\end{enumerate}
\end{satz}

\begin{beweis}
Dieser Satz ergibt sich nach \ref{bahnenzu}, indem man $P(G)$ nach jedem Repräsentanten des
Vertretersystems zulässiger Untergruppen faktorisiert, das in \ref{zu10} aufgelistet ist. Die
Präsentation der Faktorgruppe erhält man über das Verfahren aus \ref{faktor}.
\end{beweis}

\chapter{Zusammenfassung der Ergebnisse}

Von $C_p$, $C_p^2$, $C_p^3$ und $C_p^4$ aus sind im vorangegangen Kapitel endliche
Präsentationen von Gruppen berechnet worden, die zusammen mit $C_p^5$ ein Vertretersystem der
Isomorphieklassen der $p$"=Gruppen bis zur Ordnung $p^5$ ausmachen. Den elementarabelschen
Gruppen entsprechend, sind die Präsentationen in vier baumartigen Strukturen aufgetreten, die
gut die verwendete Methode widerspiegeln, aber schlecht zum Nachschlagen geeignet sind. Daher
werden die Präsenationen an dieser Stelle im Überblick aufgelistet. Sie werden nach zwei
Kriterien sortiert, nämlich erstens nach der minimalen Anzahl der Erzeuger und zweitens nach
der $p$"=Klasse. Die trivialen Relationen werden der Übersichtlichkeit halber weggelassen
(siehe Anmerkung \ref{kurz}).

\section{Die Isomorphieklassen der Gruppen der Ordnung $p$}

\begin{satz}
Ist $p$ eine Primzahl, so ist durch die folgende Liste ein Vertretersystem der
Isomorphieklassen der Gruppen der Ordnung $p$ angegeben:
\begin{enumerate}
  \item[] $G_1=C_p$.
\end{enumerate}
\end{satz}

\begin{beweis}
Dieser Satz ergibt sich unmittelbar aus dem Hauptsatz über endlich erzeugte abelsche Gruppen
(siehe etwa \cite{mey80}, S. 78).
\end{beweis}

\section{Die Isomorphieklassen der Gruppen der Ordnung $p^2$}

\begin{satz}
Ist $p$ eine Primzahl, so ist durch die folgende Liste ein Vertretersystem der
Isomorphieklassen der Gruppen der Ordnung $p^2$ angegeben:
\begin{enumerate}
  \item[] $G_1=\erz{g_1, g_2 \mid g_1^p=g_2}$,
  \item[] $G_2=C_p^2$.
\end{enumerate}
\end{satz}

\begin{beweis}
Dieser Satz ergibt sich unmittelbar aus dem Hauptsatz über endlich erzeugte abelsche Gruppen
(siehe etwa \cite{mey80}, S. 78). Alternativ ergibt sich $G_1$ nach \ref{ncp} als Nachfolger
von $C_p$.
\end{beweis}

\section{Die Isomorphieklassen der Gruppen der Ordnung $p^3$}

\begin{satz} \label{p3}
Ist $p$ eine Primzahl und $p>3$, so ist durch die folgende Liste ein Vertretersystem der
Isomorphieklassen der Gruppen der Ordnung $p^3$ angegeben:
\begin{enumerate}
  \item[] $G_1=\erz{g_1, g_2, g_3 \mid g_1^p=g_2, g_2^p=g_3} \cong C_{p^3}$, \ref{ncp},
  \item[] $G_2=\erz{g_1, g_2, g_3 \mid [g_2, g_1]=g_3} \cong D_p$, \ref{ncp2p3},
  \item[] $G_3=\erz{g_1, g_2, g_3 \mid g_1^p=g_3} \cong C_{p^2} \times C_p$, \ref{ncp2p3},
  \item[] $G_4=\erz{g_1, g_2, g_3 \mid [g_2, g_1]=g_3, g_1^p=g_3} \cong Q_p$, \ref{ncp2p3},
  \item[] $G_5=C_p^3$.
\end{enumerate}
\end{satz}

\begin{beweis}
Nach \ref{vollst} sind die Gruppen der Liste entweder elementarabelsch oder Nachfolger einer
elementarabelsche Gruppe $C$. Für den Fall, dass eine Gruppe $G$ unmittelbarer Nachfolger
einer anderen $p$"=Gruppe $H$ ist, ist anschließend an die endliche Präsentation von $G$ die
Nummer des Satzes angegeben, aus dem man entnehmen kann, von welcher Gruppe $H$ die Gruppe $G$
ein unmittelbarer Nachfolger ist und dass sie einen der Repräsentanten des Vertretersystems
unmittelbarer Nachfolger von $H$ darstellt. Damit lässt sich die endliche Folge unmittelbarer
Nachfolger von der elementarabelschen Gruppe $C$ bis zur Gruppe $G$ nachvollziehen und nach
\ref{bahnenzu} erkennen, dass $G$ zu keiner anderen Gruppe der Liste isomorph ist und die
Liste selbst ein Vertretersystem der Isomorphieklassen der Gruppen dieser Ordnung darstellt.
\end{beweis}

\section{Die Isomorphieklassen der Gruppen der Ordnung $p^4$}

\begin{satz} \label{p4}
Ist $p$ eine Primzahl und $p>3$, so ist durch die folgende Liste ein Vertretersystem der
Isomorphieklassen der Gruppen der Ordnung $p^4$ angegeben, wobei $w$ ein Erzeuger der
multiplikativen Gruppe von $\f$ ist (die Nummer hinter den endlichen Präsentationen gibt den
Satz an, aus dem man entnehmen kann, von welcher Gruppe der Listeneintrag ein unmittelbarer
Nachfolger ist):
\begin{enumerate}
  \item[] $G_1=\erz{g_1, g_2, g_3, g_4 \mid g_1^p=g_2, g_2^p=g_3, g_3^p=g_4} \cong C_{p^4}$, \ref{ncp}.
  \item[] $G_2=\erz{g_1, g_2, g_3, g_4 \mid g_1^p=g_3,
  g_2^p=g_4} \cong C_{p^2}^2$, \ref{ncp2p4}.
  \item[] $G_3=\erz{g_1, g_2, g_3, g_4 \mid [g_2, g_1]=g_3, g_1^p=g_4}$, \ref{ncp2p4} .
  \item[] $G_4=\erz{g_1, g_2, g_3, g_4 \mid [g_2, g_1]=g_4, g_1^p=g_3, g_2^p=g_4}$, \ref{ncp2p4}.
  \item[] $G_5=\erz{g_1,g_2,g_3,g_4 \mid g_1^p=g_3, g_3^p=g_4} \cong C_{p^3} \times C_p$,
  \ref{ncpp2p4}.
  \item[] $G_6=\erz{g_1,g_2,g_3,g_4, \mid g_1^p=g_3, g_3^p=g_4, [g_2,g_1]=g_4}$, \ref{ncpp2p4}.
  \item[] $G_7=\erz{g_1,g_2,g_3,g_4, \mid [g_2,g_1]=g_3, [g_3,g_1]=g_4}$, \ref{ndpp4}.
  \item[] $G_8= \erz{g_1,g_2,g_3,g_4, \mid [g_2,g_1]=g_3,
  [g_3,g_1]=g_4, g_1^p=g_4}$, \ref{ndpp4}.
  \item[] $G_9=\erz{g_1,g_2,g_3,g_4 \mid [g_2,g_1]=g_3,
  [g_3,g_1]=g_4, g_2^p=g_4}$, \ref{ndpp4}.
  \item[] $G_{10}=\erz{g_1,g_2,g_3,g_4 \mid [g_2,g_1]=g_3,
  [g_3,g_1]=g_4, g_2^p=g_4^w}$, \ref{ndpp4}.
  \item[] $G_{11}=\erz{g_1, g_2, g_3, g_4 \mid g_1^p=g_4}$, \ref{ncp3p4}.
  \item[] $G_{12}=\erz{g_1, g_2, g_3, g_4 \mid [g_2,g_1]=g_4}$, \ref{ncp3p4}.
  \item[] $G_{13}=\erz{g_1, g_2, g_3, g_4 \mid [g_2,g_1]=g_4, g_1^p=g_4}$, \ref{ncp3p4}.
  \item[] $G_{14}=\erz{g_1, g_2, g_3, g_4 \mid [g_2,g_1]=g_4, g_3^p=g_4}$, \ref{ncp3p4}.
  \item[] $G_{15}=C_p^4$.
\end{enumerate}
\end{satz}

\begin{beweis}
Nach \ref{vollst} sind die Gruppen der Liste entweder elementarabelsch oder Nachfolger einer
elementarabelsche Gruppe $C$. Für den Fall, dass eine Gruppe $G$ unmittelbarer Nachfolger
einer anderen $p$"=Gruppe $H$ ist, ist anschließend an die endliche Präsentation von $G$ die
Nummer des Satzes angegeben, aus dem man entnehmen kann, von welcher Gruppe $H$ die Gruppe $G$
ein unmittelbarer Nachfolger ist und dass sie einen der Repräsentanten des Vertretersystems
unmittelbarer Nachfolger von $H$ darstellt. Damit lässt sich die endliche Folge unmittelbarer
Nachfolger von der elementarabelschen Gruppe $C$ bis zur Gruppe $G$ nachvollziehen und nach
\ref{bahnenzu} erkennen, dass $G$ zu keiner anderen Gruppe der Liste isomorph ist und die
Liste selbst ein Vertretersystem der Isomorphieklassen der Gruppen dieser Ordnung darstellt.
\end{beweis}

\section{Die Isomorphieklassen der Gruppen der Ordnung $p^5$}

\begin{satz}
Ist $p$ eine Primzahl und $p>3$, so ist durch die folgende Liste aus $$61 + 2 \cdot p + ggT(4,
p-1) + 2 \cdot ggT(3, p-1)$$ endlichen Präsentationen ein Vertretersystem der
Isomorphieklassen der Gruppen der Ordnung $p^5$ angegeben. Dabei ist $w$ ein Erzeuger der
multiplikativen Gruppe von $\f$ sowie $W_3=\{x \in \f \mid x^3=1\}$ und $W_4=\{x \in \f \mid
x^4=1\}$. Weiterhin sei $a \in W_3$ sowie $b \in W_4$ und $k \in \{1 \dd \frac{p-1}{2}\}$. Die
Gruppen sind so angeordnet, dass einerseits die minimale Anzahl an Erzeugern zunimmt und dass
andererseits unter den Gruppen mit derselben Anzahl minimaler Erzeuger die $p$"=Klasse der
Gruppen ansteigt. (Die Nummer hinter den endlichen Präsentationen gibt den Satz an, aus dem
man entnehmen kann, von welcher Gruppe der Listeneintrag ein unmittelbarer Nachfolger ist.)
\begin{enumerate}
  \item[] $G_1=\erz{g_1, g_2, g_3, g_4, g_5 \mid g_1^p=g_2, g_2^p=g_3, g_3^p=g_4, g_4^p=g_5} \cong C_{p^5}$,  \ref{ncp},
  \item[] $G_2=\erz{g_1, g_2, g_3, g_4, g_5 \mid [g_2, g_1]=g_3, g_1^p=g_4, g_2^p=g_5}$,  \ref{ncp2p5},
  \item[] $G_3=\erz{g_1,g_2,g_3,g_4,g_5 \mid [g_2,g_1]=g_3, [g_3,g_1]=g_4, [g_3,g_2]=g_5}$, \ref{ndpp5},
  \item[] $G_4=\erz{g_1,g_2,g_3,g_4,g_5 \mid [g_2,g_1]=g_3, [g_3,g_1]=g_4, [g_3,g_2]=g_5, g_2^p=g_5}$, \ref{ndpp5},
  \item[] $G_5=\erz{g_1,g_2,g_3,g_4,g_5 \mid [g_2,g_1]=g_3, [g_3,g_1]=g_4, [g_3,g_2]=g_5, g_2^p=g_4}$, \ref{ndpp5},
  \item[] $G_6=\erz{g_1,g_2,g_3,g_4,g_5 \mid [g_2,g_1]=g_3, [g_3,g_1]=g_4, [g_3,g_2]=g_5, g_2^p=g_4^w}$, \ref{ndpp5},
  \item[] $G_7=\erz{g_1,g_2,g_3,g_4,g_5 \mid [g_2,g_1]=g_3, [g_3,g_1]=g_4, [g_3,g_2]=g_5, g_1^p=g_4, g_2^p=g_5}$, \ref{ndpp5},
  \item[] $G_8=\erz{g_1,g_2,g_3,g_4,g_5 \mid [g_2,g_1]=g_3, [g_3,g_1]=g_4 g_5, [g_3,g_2]=g_5, g_1^p=g_4, g_2^p=g_5}$, \ref{ndpp5},
  \item[] $G_9=\erz{g_1,g_2,g_3,g_4,g_5 \mid [g_2,g_1]=g_3, [g_3,g_1]=g_4 g_5^w, [g_3,g_2]=g_5, g_1^p=g_4, g_2^p=g_5}$, \ref{ndpp5},
  \item[] $G_{10}=\erz{g_1,g_2,g_3,g_4,g_5 \mid [g_2,g_1]=g_3, [g_3,g_1]=g_5^w, [g_3,g_2]=g_4, g_1^p=g_4, g_2^p=g_5}$, \ref{ndpp5},
  \item[] $G_{11}^k=\erz{g_1,g_2,g_3,g_4,g_5 \mid [g_2,g_1]=g_3, [g_3,g_1]=g_4, [g_3,g_2]=g_5^{w^k}, g_1^p=g_4, g_2^p=g_5}$, \ref{ndpp5},
  \item[] $G_{12}^k=\erz{g_1,g_2,g_3,g_4,g_5 \mid [g_2,g_1]=g_3, [g_3,g_1]=g_4 g_5^{w^{k}}, [g_3,g_2]=g_4^{w^{k-1}} g_5, g_1^p=g_4, g_2^p=g_5}$, \ref{ndpp5},
  \item[] $G_{13}=\erz{g_1,g_2,g_3,g_4,g_5 \mid g_1^p=g_3, g_2^p=g_4, g_3^p=g_5} \cong C_{p^3} \times
  C_{p^2}$, \ref{nt2},
  \item[] $G_{14}=\erz{g_1,g_2,g_3,g_4, g_5 \mid g_1^p=g_3,
  g_2^p=g_4, g_3^p=g_5, [g_2,g_1]=g_5}$, \ref{nt2},
  \item[] $G_{15}= \erz{g_1,g_2,g_3,g_4,g_5
  \mid [g_2,g_1]=g_3, g_1^p=g_4, g_4^p=g_5}$, \ref{nt3},
  \item[] $G_{16}= \erz{g_1,g_2,g_3,g_4,g_5 \mid [g_2,g_1]=g_3, g_1^p=g_4,
  [g_3,g_1]=g_5}$, \ref{nt3},
  \item[] $G_{17}= \erz{g_1,g_2,g_3,g_4, g_5
  \mid [g_2,g_1]=g_3, g_1^p=g_4, [g_3,g_1]=g_5, g_2^p=g_5}$, \ref{nt3},
  \item[] $G_{18}= \erz{g_1,g_2,g_3,g_4,g_5
  \mid [g_2,g_1]=g_3, g_1^p=g_4, [g_3,g_1]=g_5^w,g_2^p=g_5}$, \ref{nt3},
  \item[] $G_{19}= \erz{g_1,g_2,g_3,g_4,g_5 \mid [g_2,g_1]=g_3, g_1^p=g_4, [g_3,g_1]=g_5,
  g_4^p=g_5}$, \ref{nt3},
  \item[] $G_{20}= \erz{g_1,g_2,g_3,g_4,g_5
  \mid [g_2,g_1]=g_3, g_1^p=g_4, [g_3,g_1]=g_5, [g_3,g_2]=g_5}$, \ref{nt3},
  \item[] $G_{21}= \erz{g_1,g_2,g_3,g_4,g_5
  \mid [g_2,g_1]=g_3, g_1^p=g_4, [g_3,g_1]=g_5, [g_3,g_2]=g_5,g_2^p=g_5}$, \ref{nt3},
  \item[] $G_{22}= \erz{g_1,g_2,g_3,g_4,g_5
  \mid [g_2,g_1]=g_3, g_1^p=g_4, [g_3,g_1]=g_5, [g_3,g_2]=g_5,
  g_4^p=g_5}$, \ref{nt3},
  \item[] $G_{23}= \erz{g_1,g_2,g_3,g_4,g_5
  \mid [g_2,g_1]=g_3, g_1^p=g_4, [g_3,g_1]=g_5^w, [g_3,g_2]=g_5^w,
  g_4^p=g_5}$, \ref{nt3},
  \item[] $G_{24}=\erz{g_1,g_2,g_3,g_4,g_5 \mid
  [g_2,g_1]=g_3, g_2^p=g_3, g_1^p=g_4, [g_3,g_1]=g_5, [g_4,g_2]=g_5^{p-1}, g_3^p=g_5}$, \ref{nt4}
  \item[] $G_{25}=\erz{g_1,g_2,g_3,g_4,g_5 \mid
  [g_2,g_1]=g_3, g_2^p=g_3, g_1^p=g_4, g_4^p=g_5}$,  \ref{nt4}
  \item[] $G_{26}=\erz{g_1,g_2,g_3,g_4,g_5 \mid
  g_1^p=g_3, g_3^p=g_4, g_4^p=g_5} \cong C_{p^4} \times C_p$, \ref{nt5},
  \item[] $G_{27}=\erz{g_1,g_2,g_3,g_4,g_5\mid
  g_1^p=g_3, g_3^p=g_4, g_4^p=g_5, [g_2,g_1]=g_5}$, \ref{nt5},
  \item[] $G_{28}=\erz{g_1,g_2,g_3,g_4,g_5 \mid [g_2,g_1]=g_3, [g_3,g_1]=g_4,
  [g_4,g_1]=g_5}$, \ref{nt7},
  \item[] $G_{29}^a=\erz{g_1,g_2,g_3,g_4,g_5 \mid [g_2,g_1]=g_3, [g_3,g_1]=g_4,
  [g_4,g_1]=g_5^{w^a}, g_2^p=g_5}$, \ref{nt7},
  \item[] $G_{30}=\erz{g_1,g_2,g_3,g_4,g_5 \mid [g_2,g_1]=g_3, [g_3,g_1]=g_4,
  [g_4,g_1]=g_5, g_1^p=g_5}$, \ref{nt7},
  \item[] $G_{31}=\erz{g_1,g_2,g_3,g_4,g_5 \mid [g_2,g_1]=g_3, [g_3,g_1]=g_4,
  [g_4,g_1]=g_5, [g_3,g_2]=g_5}$, \ref{nt7},
  \item[] $G_{32}^a=\erz{g_1,g_2,g_3,g_4,g_5 \mid [g_2,g_1]=g_3, [g_3,g_1]=g_4,
  [g_4,g_1]=g_5^{w^a}, [g_3,g_2]=g_5^{w^a}, g_2^p=g_5}$, \ref{nt7},
  \item[] $G_{33}^b=\erz{g_1,g_2,g_3,g_4,g_5 \mid [g_2,g_1]=g_3, [g_3,g_1]=g_4,
  [g_4,g_1]=g_5^{w^b}, [g_3,g_2]=g_5^{w^b}, g_1^p=g_5}$, \ref{nt7},
  \item[] $G_{34}=\erz{g_1, g_2, g_3, g_4, g_5 \mid [g_2,g_1]=g_4, [g_3,g_1]=g_5}$, \ref{ncp3p5},
  \item[] $G_{35}=\erz{g_1, g_2, g_3, g_4, g_5 \mid [g_2,g_1]=g_4, g_3^p=g_5}$, \ref{ncp3p5},
  \item[] $G_{36}=\erz{g_1, g_2, g_3, g_4, g_5 \mid [g_2,g_1]=g_4,
  [g_3,g_2]=g_5, g_3^p=g_5}$, \ref{ncp3p5},
  \item[] $G_{37}=\erz{g_1, g_2, g_3, g_4, g_5 \mid  [g_3,g_2]=g_4, g_3^p=g_5}$, \ref{ncp3p5},
  \item[] $G_{38}=\erz{g_1, g_2, g_3, g_4, g_5 \mid  [g_3,g_1]=g_4,
  [g_3,g_2]=g_5, g_3^p=g_5}$, \ref{ncp3p5},
  \item[] $G_{39}=\erz{g_1, g_2, g_3, g_4, g_5 \mid [g_2,g_1]=g_4,
  [g_3,g_2]=g_5, g_3^p=g_4}$, \ref{ncp3p5},
  \item[] $G_{40}=\erz{g_1, g_2, g_3, g_4, g_5 \mid [g_3,g_2]=g_5, g_2^p=g_4, g_3^p=g_5}$, \ref{ncp3p5},
  \item[] $G_{41}=\erz{g_1, g_2, g_3, g_4, g_5 \mid [g_3,g_1]=g_4,
  [g_3,g_2]=g_5, g_2^p=g_4, g_3^p=g_5}$, \ref{ncp3p5},
  \item[] $G_{42}=\erz{g_1, g_2, g_3, g_4, g_5 \mid [g_2,g_1]=g_4,
  [g_3,g_2]=g_5, g_2^p=g_4, g_3^p=g_5}$, \ref{ncp3p5},
  \item[] $G_{43}=\erz{g_1, g_2, g_3, g_4, g_5 \mid g_2^p=g_4, g_3^p=g_5}$, \ref{ncp3p5},
  \item[] $G_{44}=\erz{g_1, g_2, g_3, g_4, g_5 \mid [g_2,g_1]=g_4, [g_3,g_1]=g_5,
  g_2^p=g_4, g_3^p=g_5}$, \ref{ncp3p5},
  \item[] $G_{45}=\erz{g_1, g_2, g_3, g_4, g_5 \mid [g_2,g_1]=g_5,
  g_2^p=g_4, g_3^p=g_5}$, \ref{ncp3p5},
  \item[] $G_{46}=\erz{g_1, g_2, g_3, g_4, g_5 \mid [g_2,g_1]=g_4 g_5, [g_3,g_1]=g_5,
  g_2^p=g_4, g_3^p=g_5}$, \ref{ncp3p5},
  \item[] $G_{47}=\erz{g_1, g_2, g_3, g_4, g_5 \mid [g_3,g_1]=g_5,
  g_2^p=g_4, g_3^p=g_5}$, \ref{ncp3p5},
  \item[] $G_{48}^k=\erz{g_1, g_2, g_3, g_4, g_5 \mid [g_2,g_1]=g_4, [g_3,g_1]=g_5^{w^k},
  g_2^p=g_4, g_3^p=g_5}$, \ref{ncp3p5},
  \item[] $G_{49}=\erz{g_1, g_2, g_3, g_4, g_5 \mid [g_2,g_1]=g_5^w, [g_3,g_1]=g_4,
  g_2^p=g_4, g_3^p=g_5}$, \ref{ncp3p5},
  \item[] $G_{50}^k=\erz{g_1, g_2, g_3, g_4, g_5 \mid [g_2,g_1]=g_4 g_5^{w^k},
  [g_3,g_1]=g_4^{w^{k-1}} g_5, g_2^p=g_4, g_3^p=g_5}$, \ref{ncp3p5},
  \item[] $G_{51}=\erz{g_1, g_2, g_3, g_4, g_5 \mid g_1^p=g_4, g_4^p=g_5}$, \ref{nt11},
  \item[] $G_{52}=\erz{g_1, g_2, g_3, g_4, g_5 \mid g_1^p=g_4, g_4^p=g_5,
  [g_3,g_2]=g_5}$, \ref{nt11},
  \item[] $G_{53}=\erz{g_1, g_2, g_3, g_4, g_5 \mid g_1^p=g_4, g_4^p=g_5, [g_3,g_1]=g_5}$, \ref{nt11},
  \item[] $G_{54}=\erz{g_1, g_2, g_3, g_4, g_5 \mid [g_2,g_1]=g_4, [g_4,g_2]=g_5}$, \ref{nt12},
  \item[] $G_{55}=\erz{g_1, g_2, g_3, g_4, g_5 \mid [g_2,g_1]=g_4, [g_4,g_2]=g_5, g_3^p=g_{5}}$, \ref{nt12},
  \item[] $G_{56}=\erz{g_1, g_2, g_3, g_4, g_5 \mid [g_2,g_1]=g_4, [g_4,g_2]=g_5, g_2^p=g_{5}}$, \ref{nt12},
  \item[] $G_{57}=\erz{g_1, g_2, g_3, g_4, g_5 \mid [g_2,g_1]=g_4, [g_4,g_2]=g_5, g_1^p=g_5}$, \ref{nt12},
  \item[] $G_{58}=\erz{g_1, g_2, g_3, g_4, g_5 \mid [g_2,g_1]=g_4, [g_4,g_2]=g_5^w, g_1^p=g_5}$, \ref{nt12},
  \item[] $G_{59}=\erz{g_1, g_2, g_3, g_4, g_5 \mid [g_2,g_1]=g_4, [g_4,g_2]=g_5, [g_3,g_1]=g_5}$, \ref{nt12},
  \item[] $G_{60}=\erz{g_1, g_2, g_3, g_4, g_5 \mid [g_2,g_1]=g_4, [g_4,g_2]=g_5, [g_3,g_1]=g_5, g_3^p=g_{5}}$, \ref{nt12},
  \item[] $G_{61}=\erz{g_1, g_2, g_3, g_4, g_5 \mid [g_2,g_1]=g_4, [g_4,g_2]=g_5, [g_3,g_1]=g_5, g_2^p=g_{5}}$, \ref{nt12},
  \item[] $G_{62}=\erz{g_1, g_2, g_3, g_4, g_5 \mid [g_2,g_1]=g_4, [g_4,g_2]=g_5, [g_3,g_1]=g_5, g_1^p=g_5}$, \ref{nt12},
  \item[] $G_{63}=\erz{g_1, g_2, g_3, g_4, g_5 \mid [g_2,g_1]=g_4, [g_4,g_2]=g_5^w, [g_3,g_1]=g_5^w, g_1^p=g_5}$, \ref{nt12},
  \item[] $G_{64}=\erz{g_1, g_2, g_3, g_4, g_5 \mid [g_2,g_1]=g_{5}}$, \ref{ncp4p5},
  \item[] $G_{65}=\erz{g_1, g_2, g_3, g_4, g_5 \mid [g_2,g_1]=g_{5}, [g_4,g_3]=g_{5}}$, \ref{ncp4p5},
  \item[] $G_{66}=\erz{g_1, g_2, g_3, g_4, g_5 \mid g_1^p=g_{5}}$, \ref{ncp4p5},
  \item[] $G_{67}=\erz{g_1, g_2, g_3, g_4, g_5 \mid [g_2,g_1]=g_{5}, g_1^p=g_{5}}$, \ref{ncp4p5},
  \item[] $G_{68}=\erz{g_1, g_2, g_3, g_4, g_5 \mid [g_2,g_1]=g_{5}, g_3^p=g_{5}}$, \ref{ncp4p5},
  \item[] $G_{69}=\erz{g_1, g_2, g_3, g_4, g_5 \mid [g_2,g_1]=g_{5}, [g_4,g_3]=g_{5},
  g_4^p=g_{5}}$, \ref{ncp4p5},
  \item[] $G_{70}=C_p^5$.
\end{enumerate}
In der folgenden Tabelle sind zu den oben aufgeführten Gruppen ihre $p$"=Klasse und ihre
minimale Anzahl an Erzeugern aufgelistet. In der vierten Spalte steht die Anzahl der Gruppen
mit den jeweiligen Eigenschaften.
\begin{table}[h] \centering
\begin{tabular}{l|l|l|l}
  Gruppen & Min. Erzsys. & $p$"=Klasse & Anzahl der Gruppen \\ \hline
   $G_1$                  & 1 & 5 & 1 \\
   $G_2$                  & 2 & 2 & 1 \\
   $G_3$ bis $G_{25}$     & 2 & 3 & $p+20$ \\
   $G_{26}$ bis $G_{33}$  & 2 & 4 & $ggT(4, p-1) + 2 \cdot ggT(3, p-1)+5$ \\
   $G_{34}$ bis $G_{50}$  & 3 & 2 & $p+14$ \\
   $G_{51}$ bis $G_{63}$  & 3 & 3 & 13  \\
   $G_{64}$ bis $G_{69}$  & 4 & 2 & 6  \\
   $G_{70}$               & 5 & 1 & 1  \\
\end{tabular}
\end{table}

\end{satz}

\begin{beweis}
Nach \ref{vollst} sind die Gruppen der Liste entweder elementarabelsch oder Nachfolger einer
elementarabelsche $p$"=Gruppe $C$. Für den Fall, dass eine Gruppe $G$ unmittelbarer Nachfolger
einer anderen $p$"=Gruppe $H$ ist, ist anschließend an die endliche Präsentation von $G$ die
Nummer des Satzes angegeben, aus dem man entnehmen kann, von welcher Gruppe $H$ die Gruppe $G$
ein unmittelbarer Nachfolger ist und dass sie einen der Repräsentanten des Vertretersystems
unmittelbarer Nachfolger von $H$ darstellt. Damit lässt sich die endliche Folge unmittelbarer
Nachfolger von der elementarabelschen Gruppe $C$ bis zur Gruppe $G$ nachvollziehen und nach
\ref{bahnenzu} erkennen, dass $G$ zu keiner anderen Gruppe der Liste isomorph ist und die
Liste selbst ein Vertretersystem der Isomorphieklassen der Gruppen dieser Ordnung darstellt.

Die minimale Anzahl der Erzeuger von $G$ ist nach der Definition des Nachfolgers und nach
\ref{des1} dadurch festgelegt, von welcher elementarabelschen $p$"=Gruppe $C$ die Gruppe $G$
ein Nachfolger ist. Ist nämlich $C=C_p^d$ und $G$ ein Nachfolger von $C$, so ist $d$ die
minimale Anzahl an Erzeugern von $G$. Aus der Folge unmittelbarer Nachfolger von der
elementarabelschen Gruppe $C$ bis zur Gruppe $G$, die im vorangegangenen Kapitel konstruiert
worden ist, lässt sich daher $d$ ablesen. Die $p$"=Klasse von $G$ lässt sich daran erkennen,
dass nach der Definition des unmittelbaren Nachfolgers mit jedem Glied dieser Folge die
$p$"=Klasse der jeweiligen Gruppe um Eins größer wird.
\end{beweis}

\section{Eine Bemerkung zu den Fällen $p=2$ und $p=3$}

\begin{beme}
In dieser Arbeit ist die Voraussetzung gemacht worden, dass $p$ eine Primzahl größer als 3
sei. Diese Einschränkung musste gemacht werden, da der theoretische Hintergrund des
algorithmischen Zuganges, auf den sich diese Arbeit stützt, nur unter dieser Einschränkung
gültig ist. Man vergleiche dazu \cite{hav77} und \cite{obr90}. Die Isomorphieklassen der
Gruppen der Ordnung $2^n$ und $3^n$ für $1 \leq n \leq 5$ sind beispielsweise in der
Gruppenbibliothek des Computer"=Algebra"=Systems \gap{} vorhanden (vgl. \cite{gap}).
\end{beme}

\newpage

\thispagestyle{empty}

\noindent Ich versichere, dass ich die Arbeit selbständig verfasst und keine anderen als die
angegebenen Hilfsmittel benutzt habe.


\addcontentsline{toc}{chapter}{\refname}

\begin{thebibliography}{999}

\bibitem{bag98} Bagnera, G.:
\emph{La composizione dei Gruppi finiti il cui grado é la quinta
potenza di un numero primo}, in \emph{Ann. mat. Pura Appl.}
\textbf{1}(3) (1898), 137 -- 228.

\bibitem{bes02}
Besche, Hans Ulrich, Eick, Bettina, und O'Brien, E.: \emph{A Millenium Project: Constructing
Small Groups}, in \emph{International Journal of Algebra and Computation} \textbf{12}(5)
(2002), 623 -- 644.

\bibitem{bri85} Brieskorn, Egbert:
\emph{Lineare Algebra und analytische Geometrie}, zweiter Band, Braunschweig
und Wiesbaden: Vieweg Verlag 1985.

\bibitem{cay54} Cayley, A.:
\emph{On the theory of groups, as depending on the symbolic equation
$\theta^n=1$}, in \emph{Philos. Mag.} \textbf{7}(4) (1854), 40 -- 47.

\bibitem{cay59} Cayley, A.:
\emph{On the theory of groups, as depending on the symbolic equation
$\theta^n=1$ -- Part III}, in \emph{Philos. Mag.} \textbf{18}(4) (1859), 34 --
37.

\bibitem{gap} The \gap{} Team: \emph{\gap{} --~Groups, algorithm, and programming, version 4},
Lehrstuhl D für Mathematik, RWTH Aachen, und School of Mathematical and Computational Science,
Universität von St Andrews, 1999.


\bibitem{eic99} Eick, Bettina, und O'Brien, E.: \emph{Enumerating p"=Groups} in \emph{J.
Austr. Math. Soc. (Series)} \textbf{67} (1999), 191 -- 205.

\bibitem{eic02} Eick, Bettina, Leedham"=Green, R., und O'Brien, E.: \emph{Computing
automorphism groups of $p$"=groups} in \emph{Comm. Alg.} \textbf{30} (2002), 2271 -- 2295.

\bibitem{hav77} Havas,
George, und Newman, Mike: \emph{Application of Computers to Question like those of Burnside}
in \emph{Burnside Groups} (Bielefelder Burnside"=Workshop 1977), in \emph{LNM} \textbf{806}
(1980), 211 -- 230 Berlin: Springer Verlag.

\bibitem{hal40}
Hall, P.: \emph{The classification of prime"=power groups}, in \emph{J. Reine
Angew. Math.} \textbf{182}(5) (1940), 613 -- 637.

\bibitem{hig60}
Higman, G.: \emph{Enumerating $p$"=groups. I: Inequalities}, in \emph{Proc. London Math. Soc.}
\textbf{10} (1960), 24 -- 30.

\bibitem{hoe93} Hölder, O.:
\emph{Die Gruppen der Ordnung $p^3, pq^2, pqr, p^4$}, in \emph{Math. Ann.}
\textbf{43} (1893), 301 -- 412.

\bibitem{hup67} Huppert, Bertram:
\emph{Endliche Gruppen}, erster Teilband, Berlin und Heidelberg: Springer
Verlag 1967.

\bibitem{mey80} Meyberg, Kurt:
\emph{Algebra}, erster Teilband, zweite Auflage, München und Wien: Carl Hanser
Verlag 1980.

\bibitem{obr90} O'Brien, E. A.:
\emph{The $p$"=Group Generation Algorithm} in \emph{J. Symb. Comp.} \textbf{9} (1990), 677 --
698.

\bibitem{obr94} O'Brien, E. A.:
\emph{Computing automorphism groups of $p$"=groups. Computational algebra and
number theory} in \emph{Math. Appl.} \textbf{325}, 83--90, Kluwer Acad. Publ.,
Dordrecht, 1995.

\bibitem{sim65}
Sims, Charles: \emph{Enumerating $p$"=groups}, in \emph{Proc. London Math. Soc.} \textbf{15}
(1965), 151 -- 166.


\bibitem{sim94} Sims, Charles:
\emph{Computing with finitely presented groups} aus der Reihe
\emph{Encyclopedia of mathematics and its applications}, Folge 48, Cambridge:
Cambridge University Press, 1994.

\bibitem{suz82} Suzuki, Michio:
\emph{Group Theory I} aus der Reihe \emph{Grundlehren der mathematischen
Wissenschaften --~A Series of Comprehensive Studies in Mathematics}, Folge 247,
Berlin, Heidelberg und New York: Springer, 1982.

\bibitem{syl72} Sylow, L.:
\emph{Théorèmes sur les groupes de substitutions}, in \emph{Math. Ann.}
\textbf{5} (1872), 584 -- 594.


\end{thebibliography}
\end{document}